\documentclass[11pt,letterpaper]{amsart}
\usepackage{amsmath,amsthm,xypic,graphics,amssymb}
\DeclareMathOperator{\conv}{conv}
\DeclareMathOperator{\Hom}{Hom}

\DeclareMathOperator{\sub}{sub}
\DeclareMathOperator{\svar}{svar}
\DeclareMathOperator{\csvar}{csvar}
\DeclareMathOperator{\var}{var}

\newcommand{\natu}{{\mathbb N}}
\newcommand{\real}{{\mathbb R}}
\newcommand{\qand}{\quad\text{and}\quad}

\newcommand{\updisjoint}{{\textstyle{.}}}

\newcommand{\ovline}{\underline{\phantom{aaa}}}
\newcommand{\ovtild}{\underline{\mathrm{form}}}
\renewcommand{\[}{\begin{equation}\notag}
\renewcommand{\]}{\end{equation}}

\newcommand{\xyomatrix}[1]{}
\newcommand{\leave}[1]{}

\swapnumbers
\theoremstyle{definition}
\newtheorem{point}{}[section]
\newtheorem{example}[point]{Example}

\newtheorem{defin}[point]{Definition}

\newtheorem{notat}[point]{Notation}
\newtheorem{termin}[point]{Terminology}
\newtheorem{remin}[point]{Reminder}
\newtheorem{remak}[point]{Remark}
\newtheorem{remark}[point]{Remark}

\newtheorem{conven}[point]{Convention}
\theoremstyle{plain}

\newtheorem{lemma}[point]{Lemma}

\newtheorem{cor}[point]{Corollary}

\newcommand{\marginextend}[1]{ \addtolength{\oddsidemargin}{-#1}  \addtolength{\evensidemargin}{-#1}
  \addtolength{\textwidth}{#1}\addtolength{\textwidth}{#1}}

\DeclareMathOperator*{\dire}{dir}
\newcommand{\fsubset}{\stackrel{\mathrm f}{\subset} }

\newcommand{\lrquote}{
Let $\mathfrak S\subset\mathfrak P(X)$ be an interval system,
$\mu:\mathfrak S\rightarrow\mathbb  R$ be a locally bounded measure. }
\newcommand{\quarterquote}{
Let $\mathcal  W$  be a commutative
topological group, $\mathfrak S\subset\mathfrak P(X)$ be an interval system,
$\mu:\mathfrak S\rightarrow\mathcal  W$ be a measure. }
\newcommand{\semiquote}{
Let $\mathfrak S\subset\mathfrak P(X)$ be an interval system, $\mathcal W$ be a locally convex space,
and $\mu:\mathfrak S\rightarrow\mathcal  W$ be a measure. }
\newcommand{\halfquote}{Let $\mathcal  V$, $\mathcal  W$, $\mathcal  Z$ be commutative
topological groups, $\mathfrak S\subset\mathfrak P(X)$ be an interval system,
$\mu:\mathfrak S\rightarrow\mathcal  W$ be a measure,
$L:\mathcal  V\times\mathcal  W\rightarrow\mathcal  Z$ be a biadditive pairing, which is
continuous in its second variable. }
\newcommand{\triquartquote}{Let $\mathcal  V$,  $\mathcal  Z$ be topological vector spaces,
$\mathcal W$ be a commutative topological group,
$\mathfrak S\subset\mathfrak P(X)$ be an interval system,
$\mu:\mathfrak S\rightarrow\mathcal  W$ be a measure,
$L:\mathcal  V\times\mathcal  W\rightarrow\mathcal  Z$ be a biadditive pairing, which is
continuous in its second variable and linear in the first variable. }
\newcommand{\intp}{\int^{\mathrm p}}
\newcommand{\intform}{\int^{\mathrm{form}}}
\marginextend{.6cm}
\begin{document}
\title{Notes on Lebesgue  integration}
\author{Gyula Lakos}
\date{September, 2004}
\address{Department of Geometry, E\"otv\"os University, P\'azm\'any P\'eter s.~1/C,  Budapest, H--1117, Hungary}
\email{lakos@cs.elte.hu}
\maketitle
\begin{abstract}
Alternative approaches to Lebesgue integration are considered.
\end{abstract}
\tableofcontents
\newpage
\section*{Introduction}
In what follows we discuss alternative interpretations of the Lebesgue integral
involving pairings of commutative topological group valued functions and
commutative topological group valued measures.
Several definitions will be considered.

In order to define integrals of Lebesgue type we will consider a measure
\[\mu:\mathfrak S\subset \mathfrak P(X)\rightarrow\mathcal  W,\]
and a biadditive pairing
\[L:\mathcal  V\times\mathcal  W\rightarrow\mathcal  Z.\]
Here $\mathcal  V$, $\mathcal  W$, $\mathcal  Z$ are commutative topological groups.
For the sake of simplicity, all commutative topological groups considered
are assumed to be complete, Hausdorff spaces.
In general, to define some kind of integral of a $\mathcal  V$-valued function,
it is sufficient to assume that the pairing $L$ is continuous in its second variable.
(Hence, we can choose the topology on $\mathcal  V$.)
This setting will essentially yield  a measure
\[\tilde\mu:\mathfrak S\rightarrow{\Hom}_{\mathrm{strong}}(\mathcal  V,\mathcal  Z),\]
ie. a strongly additive measure.
To keep the notation more suggestive, however, we keep using the pairing $L$.
Regarding the set system $\mathfrak S$, we will assume that it is an interval system.
These are structures slightly more general than semirings.

As a starting point for the integrals, one can consider the step functions, ie. functions of shape
\[s(x)=\sum_{j\in J}c_j\chi_{E_j}(x),\]
where $J$ is finite, $c_j\in\mathcal  V$, $E_j\in\mathfrak S$.
For such a step-function $s$ a pre-integral
\[\int^{\mathrm p}L(s,\mu)=\sum_{j\in J}L(c_j,\mu(E_j))\]
is defined.

Integration of step-functions is probably the very minimum we can  expect
from a Lebesgue type integration theory.
\\

 Now, let us consider the following definition:
\begin{defin}\label{def:01} (Lebesgue-McShane integral, or just the integral.)
A function $f:X\rightarrow\mathcal  V$ has integral
\[\int L(f,\mu)=a\]
if for each neighborhood $\mathcal  N$ of $a$ in  $\mathcal  Z$ there exists a set
\[\mathcal  A\subset\mathcal  V\times X\]
such that the following hold:
\begin{itemize}
\item[i.)] There are countable many sets $C_\lambda\subset\mathcal  V$, $D_\lambda\in\mathfrak S$
($\lambda\in\Lambda$) such that
\[\mathcal  A=\biggl(\bigcup_{\lambda\in\Lambda}C_\lambda\times D_\lambda  \biggr)\cup
\biggl(\{0\}\times\biggl(X\setminus \bigcup_{\lambda\in\Lambda}D_\lambda\biggr)\biggr).\]

\item[ii.)] The set $\mathcal  A$ contains the  graph
\[\{(f(x),x)\in\mathcal  V\times X\,: x\in X\}.\]
 of the function $f$.

\item[iii.)] The set $\mathcal  A$ contains the  graph of at least one step-function $s$.

\item[iv.)]  For all such step functions
\[\int^{\mathrm p}L(s,\mu)\in\mathcal  N.\]
\end{itemize}

[This definition requires a consistency statement for unicity, or we may just postulate that
the definition should be considered only in those cases when the integral is unique.]
\end{defin}

It turns out that this naive definition captures the Lebesgue integral as a special case
and it extends to the case of more general coefficients without difficulty.
In fact, I believe that the definition above gives \textit{the} correct definition of
the Lebesgue integral, which is the most conform to our expectations
(at least in the case when $X$ is a countable union of elements of $\mathfrak S$).

One can notice that the definition above is very geometric,
it can be considered as a generalized form of the Carath\'eodory extension process.
One the other hand, there is a functional-analytic interpretation:

In the definition above the function $f$ is approximated by step-functions $s$;
the role of the set $\mathcal  A$ is to specify an admissible set of step-functions.
For reasons which will be clear, later our notion of the integral can be interpreted
as the an extension of the pre-integral from
the set of step-functions with respect to the ``topological seminorm'' of the semivariation.

In particular, the good behavior of the semivariation of the measure $\mu$ is important for us.
The integral obtained according to the definition above has many of the special properties
of the classical Lebesgue integral.
Nevertheless, in order to obtain a realistic theory of integration, a stronger condition on $\mu$,
the property of being of locally finite semivariation is usually assumed.
This condition not only implies the continuity of $L(\cdot,\mu(\cdot))$ in the
first variable but establishes a proper notion of its boundedness.
One can do slightly better by asking for having a measure of $\sigma$-locally finite semivariation only.
\\

Nevertheless, Definition \ref{def:01} is by no means the most effective.
On one hand, we can be more generous in terms of approximations, ie. we can accommodate
a slightly larger class of functions. On the other hand, we can be more sophisticated
when choosing step-functions. Hence we obtain the extended Lebesgue-McShane,
Lebesgue-Kurzweil-Henstock and extended Lebesgue-Kurzweil-Henstock integrals.
(The names are picked up by similarity.)
These integrals are slightly more general, but they still require $\sigma$-additive measures.
Now, if the underlying space $X$ is a really nice topological space then with some cheating
one can assume that finitely additive set-functions are $\sigma$-additive.
In those cases we obtain the topological McShane and Kurzweil-Henstock integrals,
including the classical cases on finite intervals.

A common property of all of these integrals is that they can be obtained by
approximating the base measure $\mu$ by finite discrete measures.
(The convergence is non-sequential, of course, one has to use filters.)
This powerful property is also inherited by some other offsprings, e.~g.~
the Gaussian regularization of classical Lebesgue integrals.
The discrete-measure approximation principle hence results the single most
effective and general approach to integration according to my knowledge.

In fact, something more can be said:
On this very abstract level one can define Riemann/Lebesgue type
integration theories as filters $\mathcal  F$ on the set of finite formal sums
\[\sum_{j\in J} \delta_{x_j}\cdot X_j\,\]
where the $x_j$'s are from a set $X$, the $X_j$'s are from a set system
$\mathfrak S\subset\mathfrak P(X)$, and the integral is defined as the limit of
the filter $\mathcal  F$ under the image map
\[\sum_{j\in J} \delta_{x_j}\cdot X_j\quad\mapsto\quad\sum_{j\in J} L(f(x_j),\mu(X_j)).\]
(We might need to have some extra conditions on $f$, or; in the case of certain
regularized integrals, to pass to the vector space category.)
\\

A picture of different scope can, however, be considered by utilizing
other properties of the classical Lebesgue integral.
Instead of considering  sophisticated approximations one can go
``exactly'' after the function to be integrated.

Suppose that $\mathcal  V=\mathcal  W=\mathcal  Z=\mathbb R$, $L=\mathtt M$ is ordinary
multiplication, and $\mu\geq0$.
Now, assume that $f:X\rightarrow \mathbb R$ is a function,
\[f(x)='\sum_{\lambda\in\Lambda} c_\lambda\chi_{E_\lambda}(x)\qquad\text{for all }x\in X,\]
\[\sum_{\lambda\in\Lambda} c_\lambda\mu(E_\lambda)\quad\text{converges};\]
where $\Lambda$ is countable, $c_\lambda\in\mathbb R$, $E_\lambda\in\mathfrak S$;
and the prime sign ${}'$ in the equality indicates that we consider
equality only at points when the right side is convergent.

Then, one can prove that $f$ is integrable and
\[\int f\,\mu=\sum_{\lambda\in\Lambda} c_\lambda\mu(E_\lambda).\]
This is a simple consequence of Beppo Levi's theorem.
On the other hand, one can show that for each integrable function $f$ there is a
such an infinite step-function form. Hence, the observation above provides
a ``structure theorem'' for the classical integral, $\mu\geq 0$.

One can raise the question if an  alternative definition for the
integral using the infinite step-function form can be given.
The answer is obviously yes; the observations above prove the consistency of
such a definition. The real question is if the new definition
is convenient enough to build a theory of integration  from scratch.

We claim that a definition using infinite step-functions is convenient enough.
For the sake of simplicity we consider only measures of locally finite semivariation.
Moreover, using this approach, it is particularly convenient to handle those extended
cases when the integral is $+\infty$ or $-\infty$. The best way is to use the
extended set of real numbers $\mathbb R^*=\mathbb R\cup\{+\infty,-\infty,\pm\infty\}$,
where $\pm\infty$ is a short-hand notation for ``indefinite infinity''.
 That way every countable sum from $\mathbb R^*$ makes sense.
\begin{defin}\label{def:02} (Lebesgue-Riesz integral.)
Suppose that $\mu:X\rightarrow\mathbb R$ is a measure of locally finite semivariation,
$|\mu|^+$ and $|\mu|^-$ are its positive and negative parts.
We suppose that $f:X\rightarrow \mathbb R^*$ is a function,
\[f(x)='\sum_{\lambda\in\Lambda} c_\lambda\chi_{E_\lambda}(x)\qquad\text{for all }x\in X,\]
\[\sum_{\lambda\in\Lambda} c_\lambda|\mu|^+(E_\lambda)-
\sum_{\lambda\in\Lambda} c_\lambda|\mu|^-(E_\lambda)\neq\pm\infty;\]
where $\Lambda$ is countable, $c_\lambda\in\mathbb R$, $E_\lambda\in\mathfrak S$;
and the prime sign ${}'$ in the equality indicates that we consider
equality only at points when the right side is not $\pm\infty$.

Then one defines the  integral as
\[\int^{(LR)} f\,\mu=\sum_{\lambda\in\Lambda} c_\lambda|\mu|^+(E_\lambda)-
\sum_{\lambda\in\Lambda} c_\lambda|\mu|^-(E_\lambda).\]

[This definition requires a consistency statement for unicity, or we may just postulate that
this definition should be considered only in those cases when the integral is unique.]

[This definition differs from the example above in that (finite number)$='+\infty$, and
(finite number)$='-\infty$ are not allowed, but only (finite number)$='\pm\infty$.
That, however, makes no difference when the integral is supposed to be finite, in fact.
On the other hand, we gain that $+\infty$ and $-\infty$ will be acceptable values for the integral.]
\end{defin}

One can easily imagine several alternative definitions, but, for the sake of
simplicity, we will follow this one.
The definition above is pleasantly direct and it allows a parallel discussion with classical
extension theory.
\\
\\
About the structure of this paper:

In Sections \ref{sec:IS}--\ref{sec:INTMEAS} we define the  integral of
Definition \ref{def:01}, and discuss its basic properties. This part is more detailed.
In Sections \ref{sec:EXTLM}--\ref{sec:discrete} we consider some still fundamental
but more sophisticated aspects of the issue.

In Sections  \ref{sec:LKH}--\ref{sec:CLASSLEB} an admittedly erratic treatment of some
special cases follows. In those special cases the integral has nice additional properties.
In the last among those sections we discuss the case of the classical Lebesgue integrals,
in order to convince the reader that our definition for the integral is convenient enough.
Only the most elementary facts are demonstrated, after that one continue along the way
of classical treatments.

In Sections \ref{sec:INFORMAL}--\ref{sec:EXTENS} the integral of Definition \ref{def:02}
is considered. The discussion will be more detailed than before in the classical case.
In Section \ref{sec:EXTENS} the parallelism to extension theory is discussed.

The minimal amount of material which makes the point of these notes constitutes
Sections \ref{sec:IS}, \ref{sec:COUNT}, Subsections \ref{sec:MEASURE}.A,
 \ref{sec:envelopes}.A, and \ref{sec:LMcS}.A.
Sections \ref{sec:EXTLM}--\ref{sec:LKH} can be skipped when first reading.
Sections \ref{sec:INFORMAL}--\ref{sec:EXTENS} are essentially independent
from the other sections, expect some elementary tools mentioned earlier.

Otherwise, I would recommend reading the paper linearly, even if under the assumption
$\mathcal  V=\mathcal  W=\mathcal  Z=\mathbb R$.
\\

This paper is not a comprehensive survey  on integration.
There is a vast literature on general integration, cf.~\cite{dinc}, \cite{hens}, \cite{pang}
and related materials. There was no attempt made here to make a systematic comparison with
those other definitions.
My guiding objective was to show a sufficiently general, yet direct and practical
formulation with Definition \ref{def:01} for the Lebesgue integral.
As we will see, one can set this integral into a more general framework, but, in any case,
a certain concreteness cannot be avoided, hence spelling out specific approaches is
certainly useful, in my opinion.

The main reason behind preparing these notes is that I have not seen these specific approaches
expounded in other sources, or at least not on introductory level.
Especially, Sections \ref{sec:INFORMAL}--\ref{sec:EXTENS} were written for  fun,
their content should be easy exercise to anybody familiar with real analysis.
This latter part is a reformulation of Riesz's definition, cf. \cite{rszn},
for the Lebesgue integral.

I would like to thank Jared Wunsch for his good advices regarding this manuscript.

\newpage\section{Interval systems}\label{sec:IS}
\paragraph{\textbf{A. Fundamental properties}}
~\\

Interval systems are more or less the simplest set systems on
which measures can conveniently be considered.
\begin{defin}\label{def:IS}
A family of sets
$\mathfrak S$ is an interval system if for each $A,B\in\mathfrak S$
there exist countable families $\mathfrak C$ and $\mathfrak D$ of pairwise
disjoint sets from $\mathfrak S$ such that
\[A\setminus B=\bigcup_{C\in\mathfrak C}^\updisjoint C \qand A\cap B=
\bigcup_{D\in\mathfrak D}^\updisjoint D.\]
We can describe this phenomenon by saying that each element $A\in\mathfrak S$ can be decomposed
relative to any other element $B\in\mathfrak S$ using sets from $\mathfrak S$.
\end{defin}
\begin{termin}
More generally: Suppose that $\mathfrak A$ and $\mathfrak B$ are two families of sets.
We say that $\mathfrak B$ decomposes $\mathfrak A$ if
\begin{itemize}
\item[i.)] the family of sets $\mathfrak B$ contains
countably many pairwise disjoint sets;
\item[ii.)] every element of $\mathfrak A$ is a (necessarily countable) union of elements
of $\mathfrak B$.
\end{itemize}
We say that $\mathfrak B$ exactly decomposes $\mathfrak A$ if additionally
\begin{itemize}
\item[iii.)] $\bigcup\mathfrak A=\bigcup\mathfrak B$ holds.
\end{itemize}
The family $\mathfrak B$ [exactly] decomposes  the set $A$ if
it [exactly] decomposes  $\{A\}$.
\end{termin}
\begin{termin}
a.) The family of sets $\mathfrak A$ is finer than than the family of sets $\mathfrak B$
if each element of $\mathfrak B$ is a countable union of elements of $\mathfrak A$.

b.) The family of sets $\mathfrak E$ divides the family of sets $\mathfrak A$ if
for all $A\in\mathfrak A$ and $E\in\mathfrak E$ either $A\subset E$ or $A\cap E=\emptyset$.
\end{termin}
\begin{lemma}\label{lem:REFIN}
Let $N$ be an initial segment of $\natu$ and suppose that
$\mathfrak A=\{A_n\}_{n\in N}$ is an indexed family of sets from an interval
system $\mathfrak S$. Then, we  claim:

a). There exists a family of sets $\mathfrak A_{-1}\subset\mathfrak S$ which exactly
decomposes the family of sets
\[\mathfrak A'=\{A_0,\,\,A_1\setminus A_0,\,\,A_2\setminus( A_0\cup A_1),\,\,
A_3\setminus(A_0\cup A_1\cup A_2), \ldots\},\]
and so, in particular, $\bigcup\mathfrak A$.
If $\mathfrak E\subset\mathfrak S$ is finite then
we can assume that $\mathfrak A_{-1}$ is divided by $\mathfrak E$.

b). If $\mathfrak A_{-1}\subset\mathfrak S$ exactly decomposes $\bigcup\mathfrak A$
then there exist successive decompositions $\mathfrak A_n\subset\mathfrak S$
($n\in N$) such that $\mathfrak A_n$ exactly decomposes
\[\{\textstyle{\bigcup\mathfrak A}\},\,\mathfrak A_{-1},\,\ldots,\,\mathfrak A_{n-1}\]
and decomposes
\[\{A_0,\ldots, A_{n}\}.\]
If $\mathfrak A_{-1}$ was divided by $\mathfrak E$ then every $\mathfrak A_n$ is divided
by $\mathfrak E$.
\begin{proof}
a.) In an interval system a difference $A\setminus B$ decomposes, so, by induction,
(using decomposition set-wise) we can prove that the set
\[ (\ldots((A_n\setminus A_0)\setminus A_1)\setminus\ldots)\setminus A_{n-1}\]
exactly decomposes to a family $\mathfrak A'_{(n)}\subset\mathfrak S$.
Then let us do decompositions relatively to the elements of $\mathfrak E$.
That way the expression above exactly decomposes to $\mathfrak A_{(n)}$.
Now we can take
\[\mathfrak A_{-1}=\mathfrak A_{(0)}\cup\ldots\cup\mathfrak A_{(n)}\cup\ldots\,\,.\]

b.) We can obtain $\mathfrak A_n$ from $\mathfrak A_{n-1}$ by decomposing
each element $A\in\mathfrak A_{n-1}$ relative to $A_{n}$. If $\mathfrak E$ divided
$\mathfrak A_{n-1}$ then it will divide $\mathfrak A_{n}$, too.
\end{proof}
\end{lemma}
\begin{defin} \label{def:FOREST}
A family of sets $\mathfrak A$ is a forest if
\begin{itemize}
\item[i.)]
any two elements of the set $\mathfrak A$ are either disjoint or one
of them contains the other,
\item[ii.)]
 any subset of $\mathfrak A$ has a maximal element with respect
to containment;
\item[iii.)] any element of $\mathfrak A$ is contained in only finitely many other one.
\end{itemize}
A root of a forest means a maximal element of the forest.

The forest $\mathfrak A$ is fully branching if for each element $A\in\mathfrak A$ the maximal sets
contained in $A$ form an exact decomposition of $A$.
\end{defin}
\begin{lemma}\label{lem:FOREST} Suppose that $\mathfrak S$ is an interval system, and
$\mathfrak A\subset\mathfrak S$ is a countable family of sets.
Then, we claim, there exist a countable family of sets $\widetilde{\mathfrak A}\subset\mathfrak S$
such that
\begin{itemize}
\item[\textit{i.)}] The equality $\bigcup\mathfrak A=\bigcup \widetilde{\mathfrak A}$ holds.
\item[\textit{ii.)}] The countable family $\widetilde{\mathfrak A}$ is a fully branching forest.
\item[\textit{iii.)}]
Any set $C$ constructible from $\mathfrak A\cup \widetilde{\mathfrak A}$ is exactly decomposed by a
subset of $\widetilde{\mathfrak A}$. (Ie. any set which can be written as a finite expression of
sets from $\mathfrak A$, $\widetilde{\mathfrak A}$, and the operations $\cup,\,\cap,\,\setminus$
is  exactly decomposed by a subset of $\widetilde{\mathfrak A}$.)
\end{itemize}
Moreover, if $\mathfrak E\subset\mathfrak S$ is finite then we can assume that $\mathfrak E$
divides $\widetilde{\mathfrak A}$.
\begin{proof}
We can can number the elements of $\mathfrak A$ by a initial section $N$ of $\mathbb N$.
Let us apply Lemma \ref{lem:REFIN}.a and b. Then take
\[\widetilde{\mathfrak A}=\bigcup_{n\in N}\mathfrak A_{n}. \]
Then i.), ii.)  and the additional comment follow immediately.
Now, if a set $C$ is constructed from finitely many set
$C_1,\ldots,C_r\in\mathfrak A\cup\widetilde{\mathfrak A}$
then a set system $\mathfrak A_n$ of pairwise disjoint sets
will decompose all the sets $C_1,\ldots,C_r$, and hence any set constructible from them.
\end{proof}
\end{lemma}

\paragraph{\textbf{B. Some constructible classes}}
\begin{notat}\label{def:SETUN}
If $\mathfrak S$ is a family of sets then

let $\boldsymbol\Sigma\mathfrak S$ contain
the sets which occur as a countable union of elements of $\mathfrak S$;

let $\boldsymbol\Sigma_0\mathfrak S$ contain the sets which occur as a finite union
of elements of $\mathfrak S$; and

let $\boldsymbol\Sigma_{00}\mathfrak S$ contain the sets which
occur as a finite disjoint union of elements of $\mathfrak S$.
Furthermore,

let $\boldsymbol\Sigma_{\mathrm c}\mathfrak S$ contain
the sets which occur as expressions of finitely many  elements of $\mathfrak S$ and
$\cup,\,\cap$ and $\setminus$ (ie. the constructible $(\cup,\cap,\setminus)$-closure.)

\end{notat}
\begin{lemma}\label{lem:SETUN}
Assume that $\mathfrak S$ is an interval system. Then we, claim:

a.) Any element $A\in\boldsymbol\Sigma\mathfrak S$ can obtained as a countable disjoint union
of elements of $\mathfrak S$ (ie. it is exactly decomposed by a countable subset $\mathfrak D$ of
$\mathfrak S$).

b.) If $\mathfrak D_0$ contains finitely many pairwise disjoint elements from $\mathfrak S$ and
$\bigcup \mathfrak D_0\subset A\in\boldsymbol\Sigma\mathfrak S$
then $\mathfrak D_0$ can be extended to en exact decomposition $\mathfrak D$ of $A$.

c.) If $A\in\boldsymbol\Sigma\mathfrak S$ is exactly decomposed by $\mathfrak D_1$ and $\mathfrak D_2$
then there is an exact decomposition $\mathfrak D$ of $A$, which also exactly decomposes $\mathfrak D_1$
and $\mathfrak D_2$.

d.) Each element of $\boldsymbol\Sigma_{\mathrm c}\mathfrak S$ can be countably decomposed in $\mathfrak S$,
ie. belongs to $\boldsymbol\Sigma\mathfrak S$.
\begin{proof}
a.) Suppose that $A=\bigcup\mathfrak A$, where $\mathfrak A$ is a countable subset of $\mathfrak A$.
List the elements of $\mathfrak A$ is some order and apply Lemma \ref{lem:REFIN}.a.

b.) Assume things for $A$ similarly. Then list the elements of $\mathfrak D_0\cup\mathfrak A$
in some order but take the elements of $\mathfrak D_0$ ahead. Then apply Lemma \ref{lem:REFIN}.a.
In the final decomposition, the elements not contained in $\bigcup\mathfrak D_0$ will provide
the required extension.

c.) Decompose each set $D_1\cap D_2$ ($D_1\in\mathfrak D_1,\,D_2\in\mathfrak D_2$) separately
and take the union of the decompositions.

d.) That follows from Lemma \ref{lem:FOREST} applied to the generating sets of a
$(\cup,\cap,\setminus)$-constructible set.
\end{proof}
\end{lemma}
\newpage
\section{Countable sums}\label{sec:COUNT}
We collect certain elementary definitions and observations here.
See \cite{enge} for notions of general topology.\\

\paragraph{\textbf{A. Algebraic sums.}}
\begin{defin}
Let $\mathcal  V$ be a commutative group.

a.) Suppose that $\{a_\lambda\}_{\lambda\in\Lambda}$ is a countable indexed family
of elements of $\mathcal  V$.
Assume that $a_\lambda=0$ for all $\lambda\in\Lambda$ except for finitely many.
Only in that case we define the algebraic sum
\[\sum_{\lambda\in\Lambda}a_\lambda\]
as any finite sum consisting all the nonzero elements.

b.) More generally, suppose $\{A_\lambda\}_{\lambda\in\Lambda}$ is a countable
indexed family of subsets of $\mathcal  V$.
Assume that $0\in A_\lambda$ for all $\lambda\in\Lambda$ except for finitely many.
Only in that case we define the algebraic sum
\[\sum_{\lambda\in\Lambda}A_\lambda\]
as
\[\biggl\{\sum_{\lambda\in\Lambda}a_\lambda\,:\,a_\lambda\in A_\lambda,\, a_\lambda=0
\text{ for all $\lambda\in\Lambda$ except for finitely many} \biggr\}. \]

c.) When we add elements $a_\lambda\in\mathcal  V$ and subsets $A_\lambda\subset\mathcal  V$
then we, in fact, consider the sets $\{a_\lambda\}$ instead of $a_\lambda$.
\end{defin}
These sums have the contracting property:
\begin{lemma}
Let $\mathcal  V$ be a commutative group.
Let $\{A_\lambda\}_{\lambda\in\Lambda}$ be a countable indexed family of subsets of $\mathcal  V$.
Suppose that
\[\Lambda=\bigcup^\updisjoint_{\gamma\in\Gamma}\Lambda_\gamma\]
is a countable decomposition.
Then, we claim,
\[\sum_{\lambda\in\Lambda}A_\lambda=\sum_{\gamma\in\Gamma}\biggl(\sum_{\lambda\in\Lambda_\gamma}A_\lambda\biggr),\]
if left side makes sense or if the right side makes sense and $\Gamma$ is finite. \qed
\end{lemma}

\paragraph{\textbf{B. Uniform topology.}}
\begin{conven}
In what follows if $\mathcal T$ is supposed to be a (not necessarily open) neighborhood of $x$
in $\mathcal V$ then will just write that $0\in\mathcal T\subset\mathcal V$ is a neighborhood.
\end{conven}
The following lemma will be used to get rid of closure signs:
\begin{lemma}
Let $\mathcal  V$ be a commutative continuous group. Suppose that $H\subset\mathcal  V$ is arbitrary.
Assume that $0\in\mathcal  T\subset \mathcal  V$ is a neighborhood. Then, we claim,
\[\overline H\subset H+\mathcal  T.\]

More generally,
\[\overline H=\bigcap_{\substack{ 0\in\mathcal  T\subset \mathcal  V\\\text{neighborhood}}}H+\mathcal  T\]
holds. \qed
\end{lemma}
The following lemma will help us to replace the $\varepsilon/\delta$ formalism:
\begin{lemma}
Let $\mathcal  V$ be a commutative topological group.
Assume that $0\in\mathcal  T\subset \mathcal  V$ is a neighborhood.

a.) We claim that for any natural number $n$ there exist a
neighborhood $0\in\mathcal  U\subset \mathcal  V$ such that
\[\mathcal  U-\mathcal  U+\mathcal  U-\ldots\pm\mathcal  U_{\quad(n\text{ terms})\quad}\subset\mathcal  T.\]

b.) If $\Lambda$ is any countable set then there exists a family
neighborhoods $0\in\mathcal  T_{\lambda}\subset \mathcal  V$ (``division by $\Lambda$'') such that
\[\sum_{\lambda\in\Lambda}\mathcal  T_\lambda\subset\mathcal  T.\]
\begin{proof}
a.) That follows from the continuity of the expression
\[x_1-x_2+x_3-\ldots\pm x_n\]
around $0$.

b.) We can assume that $\Lambda$ is an initial segment of $\mathbb N$.
Let us set $\mathcal  T_{-1}=\mathcal  T$ and let us choose $\mathcal  T_{n+1}$ recursively such that
\[\mathcal  T_{n+1}+\mathcal  T_{n+1}\subset\mathcal  T_n.\]
Then one can see that for each $m\in\mathbb N$
\[\sum_{0\leq n\leq m}\mathcal  T_n\subset\mathcal  T.\]
Also, taking union we obtain that
\[\sum_{n\in\mathbb N}\mathcal  T_n\subset\mathcal  T.\]
\end{proof}
\end{lemma}
A consequence is:
\begin{lemma}\label{lem:SUMCLO}
Let $\mathcal  V$ be a commutative topological group.
Let $\{A_\lambda\}_{\lambda\in\Lambda}$ be a countable
indexed family of subsets of $\mathcal  V$.
Assume that $0\in A_\lambda$ for all $\lambda\in\Lambda$ except for finitely many.
Then, we claim,
\[\sum_{\lambda\in\Lambda}\overline{A_\lambda}\subset \overline{\sum_{\lambda\in\Lambda}A_\lambda}\]
Consequently,
\[\overline{\sum_{\lambda\in\Lambda}\overline{A_\lambda}}=\overline{\sum_{\lambda\in\Lambda}A_\lambda}\]
holds. \qed
\end{lemma}

\paragraph{\textbf{C. Convergence.}}
\begin{remin}
In the case of general topological spaces convergence is discussed in terms of
filters, filter bases, and filter subbases (these latter ones are set systems which
satisfy the finite nonempty intersection property).

Each filter base $\mathfrak F$ on a set $X$
(in the sense that $\subset\mathfrak P(X)$)
generates a filter $\mathfrak F^{\mathrm{f}\,\,\mathrm{in}\,\,X}$.
If no ambiguity occurs we just write $\mathfrak F^{\mathrm{f}}$.
Each filter subbase $\mathfrak F$ generates a filter base $\mathfrak F^{\mathrm{b}}$.

A filter base $\mathfrak F_1$ is finer than the filter base $\mathfrak F_2$ if
each element $F_1\in \mathfrak F_1$ is contained in an element $F_2\in \mathfrak F_2$.
It is denoted by $\mathfrak F_1\succ\mathfrak F_2$.
In terms of generated filters this means $\mathfrak F_1^{\mathrm f}\supset \mathfrak F_2^{\mathrm f}$.
Two filter bases are equivalent if they are both finer than each other.
\end{remin}
However, sometimes it is useful to consider set systems other than above.
\begin{defin}\label{def:CONV}
Let $\mathcal  V$ be a commutative topological group.
A family of sets $\mathfrak H\subset\mathfrak P(\mathcal  V)$
is convergent to $a\in\mathcal  V$ if the following conditions hold:
\begin{itemize}
\item[\texttt{(V1)}]
For each neighborhood
$0\in\mathcal  T\subset\mathcal  V$ there exists a set $H\in\mathfrak H$
such that
\[H\subset a-\mathcal  T.\]
\item[\texttt{(V2)}] For each set $H\in\mathfrak H$ and each neighborhood
$0\in\mathcal  T\subset\mathcal  V$
\[H\cap(a-\mathcal  T)\neq 0.\]
\end{itemize}
\qquad Or, expressing the same thing in different ways:
\begin{itemize}
\item For each set $H\in\mathfrak H$ and each neighborhood
$0\in\mathcal  T\subset\mathcal  V$
\[a\in H+\mathcal  T.\]
\end{itemize}
\begin{itemize}
\item For each set $H\in\mathfrak H$
\[a\in\overline H.\]
\end{itemize}
\end{defin}
\begin{defin}\label{def:CAUCHY}
Let $\mathcal  V$ be a commutative topological group.
A family of sets $\mathfrak H\subset\mathfrak P(\mathcal  V)$ is
a Cauchy system if the following conditions holds:
\begin{itemize}
\item[\texttt{(C1)}] For each neighborhood $0\in\mathcal  T\subset\mathcal  V$ there exists a set
$H\in\mathfrak H$ such that
\[H-H\subset \mathcal  T.\]
\item[\texttt{(C2)}] For any $H_1,H_2\in\mathfrak H$ and neighborhoods
$0\in\mathcal  T_1,\mathcal  T_2\subset\mathcal  V$
\[(H_1+\mathcal  T_1)\cap( H_2+\mathcal  T_2)\neq\emptyset.\]
\end{itemize}
\end{defin}
\begin{lemma}\label{lem:CONV}
Let $\mathcal  V$ be a commutative topological group.
Then every convergent set system is a Cauchy set system. \qed
\end{lemma}
\begin{lemma}\label{lem:CONVEQUI}
Let $\mathcal  V$ be a commutative topological group.
Suppose that $\mathfrak H\subset\mathfrak P(\mathcal V)$ is a set system.
We can consider the ``smeared'' set system
\[\mathfrak H^{\mathrm t}=\{H+\mathcal  T\,:\,H\in\mathfrak H,\,0\in\mathcal  T\subset\mathcal  V
\text{ is a neighborhood}\}.\]
Then, we claim, the set systems $\mathfrak H$, $\mathfrak H^{\mathrm t}$,
$\mathfrak H^{\mathrm t\mathrm f}$ are equiconvergent. \qed
\end{lemma}
\begin{lemma}\label{lem:CAUCHY}
Let $\mathcal  V$ be a commutative topological group.
Suppose that $\mathfrak H\subset\mathfrak P(\mathcal  V)$ is a Cauchy system.

Then, we claim that for any finitely many $H_1,\ldots,H_n\in\mathfrak H$, neighborhoods
$0\in\mathcal  T_0,\mathcal  T_1,\ldots,\mathcal  T_n\subset\mathcal  V$, and
natural numbers $r,s\in\mathbb N$ there exists
a set $H\in\mathfrak H$ and a neighborhood $0\in\mathcal  T\subset\mathcal  V$ such that
\[(H+\mathcal  T)-(H+\mathcal  T)+\ldots+(H+\mathcal  T)_{\quad(2r+1\text{ \rm{terms}})\quad}
\subset(H_1+\mathcal  T_1)\cap\ldots\cap( H_r+\mathcal  T_n)\]
and
\[(H+\mathcal  T)-(H+\mathcal  T)+(H+\mathcal  T)-\ldots-(H+\mathcal  T)_{\quad(2s\text{ \rm{terms}})\quad}
\subset\mathcal  T_0.\]

In particular, the ``smeared'' set system
\[\mathfrak H^{\mathrm t}=\{H+\mathcal  T\,:\,H\in\mathfrak H,\,0\in\mathcal  T\subset\mathcal  V
\text{ is a neighborhood}\}\]
forms a filter base, generating a Cauchy filter $\mathfrak H^{\mathrm t\mathrm f}$.
\qed
\end{lemma}
In order to make our life simpler we adopt the following:
\begin{conven}
In what follows ``commutative topological group'' means
complete, Hausdorff commutative topological group. Similarly, all
topological vector spaces are assumed to be complete, Hausdorff.
\end{conven}
Completeness makes Cauchy filters convergent, the Hausdorff property implies
uniqueness for limits.
\\

\paragraph{\textbf{D. Absolute convergent sums.}}
\begin{defin}\label{def:SUMdef}
Let $\mathcal  V$ be a commutative topological group.
Suppose that $\{a_\lambda\}_{\lambda\in\Lambda}$ is a countable indexed family of elements of
$\mathcal  V$. We say that $a$ is the absolute sum of the family $\{a_\lambda\}_{\lambda\in\Lambda}$
if for each neighborhood $0\in\mathcal  T\subset\mathcal  V$ there
 exist a finite set $\Xi\subset\Lambda$
such that
\[\sum_{\lambda\in\Xi}a_\lambda+\sum_{\lambda\in\Lambda\setminus\Xi}\{0,a_\lambda\}
\equiv
\biggl\{\sum_{\omega\in\Omega}a_\omega\,:\,
\Omega\text{ is finite},\,\Xi\subset\Omega\subset\Lambda\biggr\}\subset a-\mathcal  T.\]
\end{defin}
\begin{lemma}\label{lem:SUMdef}
Let $\mathcal  V$ be a commutative topological group.
Suppose that $\{a_\lambda\}_{\lambda\in\Lambda}$ is a countable indexed family of elements of
$\mathcal  V$. If the sum $a$ exists then it is unique. The sum $a$ exists
if and only if for each neighborhood $0\in\mathcal  T'\subset\mathcal  V$ there
exist a finite set $\Xi\subset\Lambda$ such that
\[\sum_{\lambda\in\Lambda\setminus\Xi}\{0,a_\lambda\}
\equiv\biggl\{\sum_{\omega\in\Omega'}a_\omega\,:\,
\Omega'\text{ is finite},\,\Omega'\subset\Lambda\setminus\Xi\biggr\}\subset \mathcal  T'.\]
In particular, partial sums exist.
\begin{proof}
These statements follow from the Hausdorff and completeness properties, respectively.
\end{proof}
\end{lemma}
\begin{notat}
For such an absolute sum as in Definition \ref{def:SUMdef}  the  notation
\[a=\sum_{\lambda\in\Lambda}^{\ovline} a_\lambda\]
will be used.
\end{notat}
\begin{remark}\label{rem:SUMdef}
If $\Lambda$ is finite then the notion of the absolute sum is the same as the algebraic one.
If $\Lambda$ is infinite then one can simply prove that
\[\sum_{\lambda\in\Lambda}^{\ovline} a_\lambda\]
exists if and only if for any ordering $p:\mathbb N\xrightarrow{\simeq}\Lambda$ the sequence
\[\lim_{n\rightarrow\infty}\sum_{m=0}^n a_{p(m)}\]
converges. In the latter case there will be a common limit, the sum.
\end{remark}
An immediate consequence of the definition is:
\begin{lemma}
Let $\mathcal  V$ be a commutative topological group.
Suppose that $\{a_\lambda\}_{\lambda\in\Lambda}$ is a countable indexed family of elements of
$\mathcal  V$, such that the corresponding sum is convergent.
Then, we claim, for each $\Xi\subset\Lambda$ finite
\[\sum_{\lambda\in\Lambda}^{\ovline} a_\lambda\in
\overline{\sum_{\lambda\in\Xi}a_\lambda+\sum_{\lambda\in\Lambda\setminus\Xi}\{0,a_\lambda\}}.\]
On the other hand, for each $0\in\mathcal  T\subset\mathcal  V$ neighborhood there exists
a finite $\Xi\subset\Lambda$
such that
\[\overline{\sum_{\lambda\in\Xi}a_\lambda+\sum_{\lambda\in\Lambda\setminus\Xi}\{0,a_\lambda\}}
\subset \sum_{\lambda\in\Lambda}^{\ovline} a_\lambda-\mathcal  T .\]
\qed
\end{lemma}
Quite the most important property of sums is the contraction property:
\begin{lemma}\label{lem:CONTR}
Let $\mathcal  V$ be a commutative topological group.
Suppose that $\{a_\lambda\}_{\lambda\in\Lambda}$ is a countable indexed family of elements of
$\mathcal  V$. Also assume that we have a disjoint decomposition
\[\Lambda=\bigcup_{\gamma\in\Gamma}^\updisjoint\Lambda_\gamma.\]
Then, we have
\[\sum_{\lambda\in\Lambda}^{\ovline} a_\lambda=\sum_{\gamma\in\Gamma}^{\ovline}
\biggl(\sum_{\lambda\in\Lambda_\gamma}^{\ovline} a_\lambda\biggr)\]
if the left side exists or if the right side exists and $\Gamma$ is finite.
\begin{proof}
Suppose that the left side exists. According to the Lemma \ref{lem:SUMdef} the partial sums
exist. Let
\[a=\sum_{\lambda\in\Lambda}^{\ovline}a_\lambda,\qquad
a_{(\gamma)}=\sum_{\lambda\in\Lambda_\gamma}^{\ovline} a_\lambda. \]
Consider an arbitrary neighborhood $0\in\mathcal  T'\subset\mathcal  V$.

Let $0\in\mathcal  T\subset V$ be a neighborhood such that $\mathcal  T-\mathcal  T\subset \mathcal  T'$.
Let $\Xi$ be a set as in Definition \ref{def:SUMdef} and
\[\Xi'=\{\gamma\,:\,\gamma\in \Gamma,\,\Lambda_\gamma\cap\Xi\neq\emptyset\}.\]
Now, let us divide $\mathcal  T$ by $\Gamma$.
For each $\gamma\in\Gamma$ let us consider a set $\Xi_\gamma$ as in Definition
\ref{def:SUMdef} but with $\Lambda_\gamma$ and $\mathcal  T_\gamma$.
Then
\begin{multline}
\sum_{\gamma\in\Xi'}a_{(\gamma)}+\sum_{\gamma\in\Gamma\setminus\Xi'}\{0,a_{(\gamma)}\}
\subset \sum_{\gamma\in\Xi'}\biggl(\bigg(\sum_{\lambda\in\Xi_{\gamma}\cup(\Xi\cap \Lambda_\gamma)}
a_\lambda\biggr)+\mathcal  T_\gamma\biggr)+\\
+\sum_{\gamma\in\Gamma\setminus\Xi'}\biggl(\bigg(\sum_{\lambda\in\Xi_{\gamma}}\{0,a_\lambda\}\biggr)+
\mathcal  T_\gamma\biggr)\subset\sum_{\lambda\in\Xi\cup{\bigcup_{\gamma\in\Xi'}\Xi_\gamma}}a_\lambda
+\sum_{\lambda\in{\bigcup_{\gamma\in\Gamma\setminus\Xi'}\Xi_\gamma}}\{0,a_\lambda\}
+\sum_{\gamma\in\Gamma}\mathcal  T_\gamma\\
\subset\biggl(\sum_{\lambda\in\Xi}a_\lambda+\sum_{\lambda\in\Lambda\setminus\Xi}\{0,a_\lambda\}\biggr)
+\sum_{\gamma\in\Gamma}\mathcal  T_\gamma\subset(a-\mathcal  T)+\mathcal  T\subset a-\mathcal  T'.
\notag
\end{multline}
That proves that $\Xi'$ satisfies the requirements of Definition \ref{def:SUMdef}
with respect to the family $\{a_{(\gamma)}\}_{\gamma\in\Gamma}$, $a$ and $\mathcal  T'$.
Being $\mathcal  T'$ arbitrary that establishes the our case.

On the other hand, suppose that  the right side exists and $\Gamma$ is finite.
We can still use the notation
\[a_{(\gamma)}=\sum_{\lambda\in\Lambda_\gamma}^{\ovline} a_\lambda. \]

Consider an arbitrary neighborhood $0\in\mathcal  T\subset\mathcal  V$.
We can divide $\mathcal  T$ by $\Gamma$.
Now, let $\Xi_{\gamma}$ be as the one in the Definition
\ref{def:SUMdef} but with $\Lambda_\gamma$ and $\mathcal  T_\gamma$.
Then we claim,
\[\Xi=\bigcup_{\gamma\in\Gamma}\Xi_\gamma\]
satisfies the requirements in Definition \ref{def:SUMdef}.
Indeed,
\[\sum_{\lambda\in\Xi}a_\lambda+\sum_{\lambda\in\Lambda\setminus\Xi}\{0,a_\lambda\}
=\sum_{\gamma\in\Gamma}\biggl(
\sum_{\lambda\in\Xi_\gamma}a_\lambda+\sum_{\lambda\in\Lambda_\gamma\setminus\Xi_\gamma}
\{0,a_\lambda\}\biggr)\subset\]\[\subset\sum_{\gamma\in\Gamma} (a_{(\gamma)}-\mathcal  T_\gamma)=
\sum_{\gamma\in\Gamma} a_{(\gamma)}-\sum_{\gamma\in\Gamma}\mathcal  T_\gamma\subset
\sum_{\gamma\in\Gamma} a_{(\gamma)}-\mathcal  T.\]
\end{proof}
\end{lemma}
There are other notable properties, including behavior with respect to continuous
linear functionals, which are quite straightforward, so we finish the discussion here.
\newpage
\section{Measures}\label{sec:MEASURE}
\paragraph{\textbf{A. Fundamentals.}}
\begin{defin}\label{def:MEASURE}
Suppose that $\mathfrak S$ is an interval system and
$\mathcal  W$ is a commutative topological group. Then the function
\[\mu:\mathfrak S\rightarrow\mathcal  W\]
is called a $\mathcal  W$-valued measure on $\mathfrak S$ if it
is $\sigma$-additive, ie. for any countable disjoint decomposition
\[A=\bigcup_{\lambda\in\Lambda}^\updisjoint A_\lambda \]
in $\mathfrak S$ the equality
\[\mu(A)=\sum_{\lambda\in\Lambda}^{\ovline}\mu(A_\lambda) \]
holds.
\end{defin}
\begin{example}
Let $\mathfrak I$ be the set  of possibly degenerate finite intervals on $\real$.
One can simply check that $\mathfrak I$ is an interval system.
Furthermore, we claim that the interval length function
\[l:\mathfrak I\rightarrow\mathbb R.\]
is a measure. Indeed, the $\sigma$-additive property
follows from the usual compact-open argument.
\end{example}
\begin{remark}
As Remark \ref{rem:SUMdef} shows, it would be enough to ask $\sigma$-additivity
in the sense of ordered sums, then absolute convergence follows automatically.
\end{remark}
\begin{conven}
If $E\subset X$ then then we can consider its characteristic function $\chi_E$, which is
supposed to be a $\mathbb Z$-valued function. The natural pairing $\mathtt m$ between $\mathbb Z$
and any commutative topological group $\mathcal  V$ will be used all the time.
\end{conven}
\begin{defin}\label{def:STEP}
Assume that $\mathfrak S\subset \mathfrak P(X)$ is an interval system and $\mathcal  V$ is
commutative topological group.

a.) A function
\[s:X\rightarrow\mathcal  V\]
is a step-function if has form
\[s=\sum_{j\in J}c_j\chi_{E_j},\]
where $J$ is finite, $c_j\in\mathcal  V$, $E_j\in\mathfrak S$ ($j\in J$).
A set function $s$ is simple if there is a form, where the sets $E_j$ are pairwise disjoint.

Ie. step functions come from ``finite non-disjoint sums'', while simple step functions
come from more special ``finite disjoint sums''.

b.) If $s$ is a step function and it takes a constant  value on the set $D$ then this
constant value is denoted by $s^D$.
\end{defin}
\begin{lemma}\label{lem:SUM}
\halfquote
Suppose that $J$ is finite, $c_j\in\mathcal  V$, $E_j\in\mathfrak S$ $(j\in J)$.

a.) If $\mathfrak D\subset\mathfrak S$ decomposes $\{E_j\,:\,j\in J\}$
then, we claim,
\[\sum_{j\in J}L(c_j,\mu(E_j))=
\sum_{D\in\mathfrak D}^{\ovline}L\biggl(s^D,\mu(D)\biggr)\]
(including the existence of the right side).

b.) The sum \[\sum_{j\in J}L(c_j,\mu(E_j))\]
can be recovered from the step-function
\[s=\sum_{j\in J} c_j\chi_{E_j}.\]
\begin{proof} a.) We can obtain the right side from the left side through the
equalities
\[=\sum_{j\in J}\biggl(
\sum^{\ovline}_{\substack{D\in\mathfrak D\\D\subset E_j}}L(c_j,\mu(D))\biggr)=
\sum^{\ovline}_{\substack{j\in J,\,D\in\mathfrak D\\D\subset E_j}}L(c_j,\mu(D))=
\sum^{\ovline}_{D\in\mathfrak D}\biggl( \sum_{\substack{j\in J\\D\subset E_j}}L(c_j,\mu(D)) \biggr)=.\]
The first equality follows from the $\sigma$-additivity of $\mu$ and the continuity
of $L$ in the second variable.
The second one is legal because a finite sum of infinite sums still makes up a
valid sum. The third one is legal because one can contract terms in an infinite sum.
The fourth one follows from the additivity of $L$ in the first variable.

b.) Suppose that we have an other sum form such that
$J'$ is finite, $c_j'\in\mathcal  V$, $E_j'\in\mathfrak S$ $(j\in J')$.
Apply Lemma \ref{lem:REFIN}.a and b. in order to find an exact decomposition $\mathfrak D$ of
\[\{E_j\,:\,j\in J\}\cup\{E_j'\,:\,j\in J'\}.\]
Now apply part a.~of this lemma in these two cases.
On the left sides we have the two original sums, while on
the right we have the very same expression in those two cases,
 because according to our assumptions
the coefficients paired with $\mu(D)$ are the same.
Hence the equality of the original sums is implied.
\end{proof}
\end{lemma}
~

\paragraph{\textbf{B. Semivariation.}}~\\

The following statements of regularity are not necessary to define the integral itself,
but they are important when we study its properties.
\begin{defin}\label{def:XMEAVAR}
Let $\mathfrak S$ be an interval system, $\mathcal  Z$ be a commutative topological group
and $\tilde\mu:\mathfrak S\rightarrow\mathcal  Z$ be a measure.
Suppose that $A$ is a union of countable many elements of $\mathfrak S$
(ie. $A\in\boldsymbol\Sigma\mathfrak S$).
Then we define the semivariation of $A$ with respect to $\tilde\mu$ as
\begin{multline}
\svar(A,\tilde\mu)=\\\overline{\biggl\{\sum_{j\in J}\tilde\mu(A_j)\,:\,
J\text{ is finite},\,A_j\in\mathfrak S,\,
A_j\subset A,\text{ the $A_j$ are pairwise disjoint}\biggr\}}.\notag
\end{multline}
\end{defin}
\begin{lemma}\label{lem:XMEAVAR}
Let $\mathfrak S$ be an interval system, $\mathcal  Z$ be a commutative topological group
and $\tilde\mu:\mathfrak S\rightarrow\mathcal  Z$ be a measure. Suppose that $A$ is a countable  disjoint union
\[A=\bigcup_{\lambda\in\Lambda}^\updisjoint A_\lambda\]
of elements $A_\lambda\in\mathfrak S$. Then, we claim, for a fixed set $T\subset \mathcal  Z$ the
following two statements are equivalent:

i.) For each neighborhood  $0\in\mathcal  U\subset\mathcal  Z$ and for all countable decomposition
\[A=\bigcup_{\omega\in\Psi}^\updisjoint C_\psi,\]
($C_\psi\in\mathfrak S$),
there exist a finite set $\Omega\subset\Psi$ such that
\[\sum_{\psi\in\Psi\setminus\Omega}\{0,\tilde\mu(C_\psi) \}\subset T+\mathcal  U. \]

ii.) For each neighborhood
$0\in\mathcal  T\subset\mathcal  Z$
there exists a finite set $\Xi\subset\Lambda$ such that
\[\svar\biggl(\bigcup_{\lambda\in\Lambda\setminus\Xi} A_\lambda,\tilde\mu\biggr)\subset T+\mathcal  T .\]
\begin{proof}
i.$\Rightarrow$ii.)
Suppose that it is otherwise. We can assume that $\Lambda=\mathbb N$.
Let $\mathcal  U$ be such a neighborhood of $0$ such that
$\mathcal  U+\mathcal  U+\mathcal  U\subset\mathcal  T$.
By induction we will construct a monotone increasing sequence of natural numbers
$s_k$ $(k\in\mathbb N)$, $s_{-1}={-1}$ such that
\[\bigcup_{s_{k-1}<n\leq s_k} A_n\]
exactly decomposes into a family  $\{C_\psi\,:\,\psi\in\Psi_k \}\subset\mathfrak S$
and for a finite subset $\Omega_k\subset\Psi_k$
\[\sum_{\omega\in\Omega_k }\tilde\mu(C_\omega)\notin T+\mathcal  U.\]
That will contradict to our assumption with respect to the decomposition
\[\{C_\psi\,:\,\psi\in\Psi\}\subset\mathfrak S,\qquad
\Psi=\bigcup_{k\in\mathbb N}\Psi_k\]
of $A$.

As for the inductive construction:
Suppose that $s_k$ is defined. Then according to our assumption,
there exists an element
\[t_k\in\svar\biggl(\bigcup_{s_{k-1}<n}A_n,\tilde\mu\biggr)\setminus (T+\mathcal  T).\]
Now, from the definition of semivariation follows that
there exist finitely many pairwise disjoint set $B_{\xi}$ $(\xi\in\Xi_k)$,
such that $B_{\xi}\subset \bigcup_{n>s_{k-1}}A_n$ and
\[\sum_{\xi\in\Xi_k}\tilde\mu(B_\xi)\in t_k-\mathcal  U.\]

Now, each $B_\xi$ exactly decomposes into
\[\{B_\xi\cap A_{s_k+1},\,B_\xi\cap A_{s_k+2},\,B_\xi\cap A_{s_k+3},\ldots, \},\]
whose elements  exactly decompose in $\mathfrak S$ further.
Thus, we obtain that $B_\xi$ exactly decomposes into a family
$\mathfrak B_\xi\subset\mathfrak S$, such
that each element is contained in one of the $A_n$'s $(n>s_{k-1})$.
Then, being $\Xi_k$ finite,
\[\sum_{\xi\in\Xi_k} \tilde\mu(B_\xi)=\sum_{\xi\in\Xi_k} \biggl(\sum_{B\in\mathfrak  B_\xi}^{\ovline}
\tilde\mu(B)\biggr)=\sum_{\xi\in\Xi_k,\,B\in\mathfrak B_\xi}^{\ovline}\tilde\mu(B).\]
Hence, by the definition of sum there exist finitely many non-empty elements $C_\omega$
$(\omega\in\Omega_k)$ from
\[\bigcup_{\xi\in\Xi_k}\mathfrak B_\xi\]
such that
\[\sum_{\omega\in\Omega_k} \tilde\mu(C_\omega)\in
\sum_{\xi\in\Xi_k,\,B\in\mathfrak  B_\xi}^{\ovline}\tilde\mu(B)-\mathcal  U
\subset t_k-\mathcal  U-\mathcal  U.\]
So, in particular,
\[\sum_{\omega\in\Omega_k} \tilde\mu(C_\omega)\notin T+\mathcal  U, \]
because otherwise
\[t_k\in \sum_{\omega\in\Omega_k} \tilde\mu(C_\omega)+\mathcal  U+\mathcal  U
\subset (T+\mathcal  U)+\mathcal  U+\mathcal  U
\subset T+\mathcal  T\]
would hold.
Let $s_k$ be the highest index such that $C_\omega\subset A_{s_k}$ for an
$\omega\in\Omega_k$.

According Lemma \ref{lem:SETUN}.b we can find countable many elements which complement the $C_\omega$'s to
form an exact decomposition $\{C_\psi\,:\,\psi\in\Psi_k\}$ of
$\{A_n\,:\,s_{k-1}<n<s_k\}$ in $\mathfrak S$.

That completes the construction step.

ii$\Rightarrow$i.) Let us choose a neighborhood $0\in\mathcal  T\subset\mathcal  Z$ such that
$\mathcal  T+\mathcal  T+\mathcal  T\subset\mathcal  U$. Now, for each set $C_\psi$ the set
\[\{C_\psi\cap A_0,\,C_\psi\cap A_1,\,\ldots \}\]
decomposes $C_\psi$; moreover, each element decomposes further in $\mathfrak S$.

Hence we obtain a decomposition $\mathfrak B_\psi\subset\mathfrak S$ of $C_\psi$ such that each
element of $\mathfrak B_\psi$ is contained in one of the $A_n$'s.
Let us choose $\Xi$ such that
\[\svar\biggl(\bigcup_{\lambda\in\Lambda\setminus\Xi} A_\lambda,\tilde\mu\biggr)\subset T+\mathcal  T .\]
We define
\[\mathfrak B=\bigcup_{\psi\in\Psi}\mathfrak B_\psi.\]
Let
\[\mathfrak B_{(1)}=\{B\in\mathfrak B\,:\,B\subset A_\lambda,\,\lambda\in\Xi\},\]
\[\mathfrak B_{(2)}=\{B\in\mathfrak B\,:\,B\subset A_\lambda,\,\lambda\in\Lambda\setminus\Xi\}.\]
Now, being $\Xi$ finite
\[\sum_{\xi\in\Xi}\tilde\mu(A_\xi)=\sum_{\xi\in\Xi}\,\,\sum_{B\in\mathfrak B,\,B\subset A_\xi}^{\ovline}
\tilde\mu(B)=\sum_{B\in\mathfrak B_{(1)}}^{\ovline}\tilde\mu(B).\]
According to the definition of sums
let $\mathfrak B_{(0)}\subset\mathfrak B_{(1)} $ a finite subset  such that
\[ \sum_{B\in\mathfrak B_{(1)}\setminus \mathfrak B_{(0)}}\{0,\tilde\mu(B)\}\subset \mathcal  T\]
Let $\Omega$ be the set of such indices such that $C_\omega$ contains an element of $\mathfrak B_{(0)}$.
We claim that this $\Omega$ satisfies our requirements.

Indeed, if $J\subset \Psi\setminus\Omega$ is finite then
\[\sum_{j\in J}\tilde\mu(C_j)=\sum_{j\in J}\quad\sum_{B\in\mathfrak B,\,B\subset C_j }^{\ovline}\tilde\mu(B)=
\sum_{B\in\mathfrak B,\,B\subset C_j,\,j\in J }^{\ovline}\tilde\mu(B)\]\[\subset
\sum_{B\in\mathfrak B,\,B\subset C_j,\,j\in J }\{0,\tilde\mu(B)\}+\mathcal  T
\subset\sum_{B\in\mathfrak B_{(1)}\setminus \mathfrak B_{(0)} }  \{0,\tilde\mu(B)\}+
\sum_{B\in\mathfrak B_{(2)}}\{0,\tilde\mu(B)\}+\mathcal  T\]
\[\subset \mathcal  T+(T+\mathcal  T)+\mathcal  T\subset T+\mathcal  U.\]
\end{proof}
\end{lemma}
\begin{defin}\label{def:MEAVAR}
\halfquote

Suppose that $A$ is a union of countable many elements of $\mathfrak S$
and $c\in\mathcal  V$.
Then we define the semivariation of $A$ with respect $L,c,\mu$ as
\begin{multline}
\svar(L,c,A,\mu)=\\\overline{\biggl\{\sum_{j\in J}L(c,\mu(A_j))\,:\,
J\text{ is finite},\,A_j\in\mathfrak S,\,
A_j\subset A,\text{ the $A_j$ are pairwise disjoint}\biggr\}}.\notag
\end{multline}
\end{defin}
\begin{lemma}\label{lem:MEAVAR}
\halfquote
Suppose that
\[A=\bigcup^{\updisjoint}_{\lambda\in\Lambda}A_\lambda\]
is a disjoint decomposition in $\mathfrak S$. Also suppose that
$\mathcal  T$ is a neighborhood of $0$ in $\mathcal  Z$.
Then we claim that there exists a finite set $\Xi$ such that
\[\svar\biggl(L,c,\bigcup_{\lambda\in\Lambda\setminus\Xi}A_\lambda,\mu\biggr)\subset
\mathcal  T.\]
\begin{proof}
Let us consider the function
$\tilde\mu:\mathfrak S\rightarrow\mathcal  Z$
defined by
\[\tilde\mu(A)=L(c,\mu(A)).\]
The fact that $\tilde\mu$ is a measure follows from the continuity of $L$ in the
second variable. Let us apply Lemma \ref{lem:XMEAVAR} with $T=0$. Here Lemma \ref{lem:XMEAVAR}.i
holds trivially because $\tilde\mu$ is a measure and $A\in\mathfrak S$.
Then \ref{lem:XMEAVAR}.ii  yields our statement.
\end{proof}
\end{lemma}
\begin{remark}
One can notice that
\[\svar(A,\mu)=\svar(\mathtt m,1,A,\mu)\]
where $\mathtt m$ is the natural pairing between $\mathbb Z$ and $\mathcal  W$.
\end{remark}
~

\paragraph{{\textbf{C. Constructive extension of measures}}}
~\\

The following statement will immediately follow from the properties of
integrable sets later (cf. Lemma \ref{lem:INTRING}).
However it also has an elementary proof.

\begin{lemma}
Let $\mathfrak S$ be an interval system, $\mathcal  W$ be a commutative topological group
and $\tilde\mu:\mathfrak S\rightarrow\mathcal W$ be a measure.
Then there is a unique measure extension
\[\mu_{\mathrm c}:\boldsymbol\Sigma_{\mathrm c}\mathfrak S\rightarrow\mathcal W\]
to the $(\cup,\cap,\setminus)$-constructible sets.
\begin{proof} We prove this in two steps. First we extend it to the $(\cap,\setminus)$-constructible
closure $\boldsymbol\Sigma_{\mathrm d}\mathfrak S$.

1. Suppose that $A$ is $(\cap,\setminus)$-constructible. Then according to
Lemma \ref{lem:SETUN}.d the set $A$ decomposes in $\mathfrak S$:
\[A=\bigcup_{\lambda\in\Lambda}^\updisjoint A_\lambda,\]
where $\Lambda$ is countable, $A_\lambda\in\mathfrak S$. We claim that we can set
\[\mu_{\mathrm d}(A)=\sum_{\lambda\in\Lambda}^{\ovline} \mu(A_\lambda),\]
and this will yield a measure extension
\[ \mu_{\mathrm d}:\boldsymbol\Sigma_{\mathrm d}\mathfrak S\rightarrow\mathcal W.\]
Indeed, there is an element $B\in\mathfrak S$ such that $B\subset A$. Then $C=B\setminus A$ is also
$(\cap,\setminus)$-constructible, which may be decomposed similarly to A:
\[C=\bigcup_{\omega\in\Omega}^\updisjoint C_\omega.\tag{\S}\]
Then
\[B=\bigcup_{\lambda\in\Lambda}^\updisjoint A_\lambda
\,\dot\cup\,\bigcup_{\omega\in\Omega}^\updisjoint C_\omega. \]
Applying the $\sigma$-additive property of $\mu$ to this decomposition we can notice that the corresponding
partial sum as in (\S) exists. Also,
\[ \sum_{\lambda\in\Lambda}^{\ovline} \mu(A_\lambda)= \mu(B)-\sum_{\omega\in\Omega}^{\ovline} \mu(C_\omega)\]
shows that $\mu'(A)$ does not depend on the actual decomposition by $A_\lambda$'s.

Now, if
\[A=\bigcup_{\lambda\in\Lambda}^\updisjoint A_\lambda,\]
is a countable decomposition by the $(\cap,\setminus)$-constructible sets $A_\lambda$, and
$\mathfrak A_\lambda $ decomposes $A_\lambda$ then
\[\mu_{\mathrm d}(A)=\sum_{\lambda\in\Lambda,\,S\in\mathfrak A_\lambda}^{\ovline} \mu(S)=
\sum_{\lambda\in\Lambda}^{\ovline}\sum_{S\in\mathfrak A_\lambda}^{\ovline}\mu(S)=
\sum_{\lambda\in\Lambda}^{\ovline}\mu_{\mathrm d}(A_\lambda). \]
That proves $\sigma$-additivity.

2. Suppose that $A$ is $(\cup,\cap,\setminus)$-constructible. Then there is a finite decomposition
\[A=\bigcup_{j\in J}^\updisjoint A_j,\]
where $J$ is finite and the $A_j$ are $(\cap,\setminus)$-constructible. We claim that
\[\mu_{\mathrm c}(A)=\sum_{j\in J}\mu_{\mathrm d}(A_j)\]
will yield a measure extension. Indeed, if $A'_{j'}$ ($j'\in J'$) is an other decomposition then
\[ \sum_{j\in J}\mu_{\mathrm d}(A_j)=\sum_{j\in J,\,j'\in J'}\mu_{\mathrm d}(A_j\cap A'_{j'})=
\sum_{j'\in J'}\mu_{\mathrm d}(A'_{j'}),\]
so the extension is well-defined. The $\sigma$-additivity can be proven as follows:
If
\[A=\bigcup_{\lambda\in\Lambda}^\updisjoint B_\lambda,\]
is a countable decomposition then
\[\mu_{\mathrm c}(A)=\sum_{j\in J}\mu_{\mathrm d}(A_j)=
\sum_{j\in J}\sum_{\lambda\in\Lambda}^{\ovline}\mu_{\mathrm d}(A_j\cap B_\lambda)=
\sum_{j\in J,\,\lambda\in\Lambda}^{\ovline}\mu_{\mathrm d}(A_j\cap B_\lambda)=\]
\[=\sum_{\lambda\in\Lambda}^{\ovline}\sum_{j\in J}\mu_{\mathrm d}(A_j\cap B_\lambda)=
\sum_{\lambda\in\Lambda}^{\ovline}\mu_{\mathrm c}(B_\lambda)\]
proves $\sigma$-additivity.

The unicity statements follow from the fact that every $(\cup,\cap,\setminus)$-constructible set can
be decomposed in $\mathfrak S$.
\end{proof}
\end{lemma}
\begin{point}
This is a very reasonable extension, where the underlying set
system $\boldsymbol\Sigma_{\mathrm c}\mathfrak S$ is a ring, and for the most part
it is harmless.

However, in what follows we will \textit{not} use this extension,
because it emphasizes set-theoretical concerns.
\end{point}

\newpage
\section{Envelopes}\label{sec:envelopes}
\paragraph{{\textbf{A. Fundamentals}}}
\begin{notat} Suppose that $Z\subset Y\times X$.

i.)
If $x\in X$ then we use the notation
\[Z^x=\{y\in Y\,:\,(y,x)\in Z\}.\]
Similarly, for $R\subset X$,
\[Z^R=\bigcap_{x\in R}Z^x=\{y\in Y\,:\,(y,x)\in Z\,\text{ for all }\,x\in R\}.\]

Sometimes we also use the notation
\[{}^yZ=\{x\in X\,:\,(y,x)\in Z\}.\]

ii.)
We say that the function $f:X\rightarrow Y$ is approximated by $Z$ if its graph
\[\{(y,x)\in Y\times X\,:\,y=f(x)\}\]
is contained in the set $Z$.
\end{notat}
\begin{defin}
Suppose that $\mathcal  A\subset \mathcal  V\times X$,  and $\mathcal  V$ is a commutative topological group.
Then, for a countable set system $\mathfrak D\subset \mathfrak P(X)$ we define
\[\langle\mathfrak D\rangle_{\mathcal  A}=\biggl(\bigcup_{D\in\mathfrak D}\mathcal  A^D\times D  \biggr)\cup
\biggl(\{0\}\times\biggl(X\setminus \bigcup_{D\in\mathfrak D}D\biggr)\biggr).\]
\end{defin}
\begin{conven}
In what follows a function $f:X\rightarrow\mathcal  V$ and its graph
\[\{(v,x)\in\mathcal  V\times X\,:\,v=f(x)\}\]
will be identified.
\end{conven}
\begin{defin}\label{def:APPNEIGH}
Suppose that $\mathcal  V$ is a topological vector space and
$\mathfrak S\subset\mathfrak P(X)$ is an interval system. Then we say that
\[\mathcal  A\subset \mathcal  V\times X\]
is a sectioned envelope ( or just, simply, envelope)
of sets in $\mathcal  V$ with respect to $\mathfrak S$ if the following conditions hold:
\begin{itemize}
\item[\texttt{(S$_{\mathtt s}$)}] There is a step-function $s:X\rightarrow\mathcal  V$ (with respect to $\mathfrak S$ and
$\mathcal  V$) such that its graph is contained in the set $\mathcal  A$; ie.
\[s\subset\mathcal  A.\]
\item[\texttt{(G)}] There exists a countable generating set $\mathfrak D\subset\mathfrak S$
for $\mathcal  A$; ie. a countable family of sets $\mathfrak D\subset\mathfrak S$ such that
\[\mathcal  A=\langle\mathfrak D\rangle_{\mathcal  A}.\]
\end{itemize}

We say that the envelope $\mathcal  C$ is pointed if it contains the graph of the zero step-function $0$.
(Or, equivalently,  $0\in\mathcal  C^x$ for all $x\in X$.)
\end{defin}
\begin{remark}
It is easy to see then that the conditions for a pointed case can be summed by:
\begin{itemize}
\item[\texttt{(S$_{\mathtt 0}$G)}]
There exists a countable generating set $\mathfrak D\subset\mathfrak S$
for $\mathcal  C$ such that
\[\mathcal  C=\biggl(\bigcup_{D\in\mathfrak D}\mathcal  C^D\times D  \biggr)\cup
\biggl(\{0\}\times X\biggr).\]
\end{itemize}
\end{remark}
\begin{lemma}\label{lem:APPNEIGH}
Let $\mathfrak S\subset \mathfrak P(X)$ be a interval system
and $\mathcal  V$ be a commutative topological group.
Suppose that $\mathcal  A$ is  an envelope of sets in $\mathcal  V$ with respect to $\mathfrak S$,
which contains the graph of the step-function
\[s=\sum_{j\in J}c_j\chi_{E_j},\]
where $J$ is finite, $c_j\in\mathcal  V$, $E_j\in\mathfrak S$; and $\mathfrak D$
is a generator system for $\mathcal  A$. Assume that $\mathfrak D'\subset\mathfrak S$
is an other countable family of sets. Then we claim:

a.) If $\mathfrak D'$ is finer than
$\mathfrak D$ and $\mathfrak D'$ is divided by $\mathfrak E=\{E_J\,:\,j\in J\}$; or

b.) if $\mathfrak D'$ is finer than
$\mathfrak D$ and $\bigcup\mathfrak D=\bigcup\mathfrak D'$; or

c.) if $\mathfrak D'$ is finer than
$\mathfrak D$ and $\mathcal  A$ is a pointed envelope;
\\
then $\mathfrak D'$ is also a generator system for $\mathcal  A$.
\begin{proof}
The statement is nontrivial only in the case a.) when we have to prove that
for any
\[x\in\left({\textstyle\bigcup}\mathfrak D'\right)\setminus\left({\textstyle\bigcup}\mathfrak D\right)\]
there is an element $D\in\mathfrak D'$ such that $0\in\mathcal  A^D$. However, by our assumptions
there is an element $D\in\mathfrak D'$, such that $x\in D$, and being divided by $\mathfrak E$ the value
of $s$ is constant on $D$. That value must be $0$. Hence $0\in\mathcal  A^D$.
\end{proof}
\end{lemma}

\begin{defin}\label{def:BLIN}
\halfquote

Suppose that $\mathcal  A$ is  an envelope of sets in $\mathcal  V$ with respect to $\mathfrak S$.
Then we define the bilinear sum
\[\intp L(\mathcal A,\mu)=\overline{\biggl\{\sum_{j\in J}L(c_j,\mu(E_j))\,:\,J\text{ is finite},
\,c_j\in\mathcal  V,\,E_j\in\mathfrak S,\sum_{j\in J} c_j\chi_{E_j}\subset\mathcal  A\biggr\}}.\notag\]
\end{defin}
The following observation will be used so often
that we do not make separate references to it.
\begin{lemma}
\halfquote

Suppose that $\mathcal  A_1$, $\mathcal A_2$ are envelopes
of sets in $\mathcal  V$ with respect to $\mathfrak S$ such that
\[\mathcal A_1\subset\mathcal A_2.\]
Then, we claim,
\[\intp L(\mathcal A_1,\mu)\subset \intp L(\mathcal A_2,\mu).\]
\begin{proof}
That is immediate from Definition \ref{def:BLIN}.
\end{proof}
\end{lemma}

The following statement is crucial.
\begin{lemma} \label{lem:NEIGH}
Let $\mathfrak S\subset \mathfrak P(X)$ be a interval system
and $\mathcal  V$ be a commutative topological group.
Suppose that $\{\mathcal  A_\lambda\}_{\lambda\in\Lambda}$ is an indexed family
of  envelopes of sets in $\mathcal  V$ with respect to $\mathfrak S$, such that
$\mathcal  A_\lambda$ is pointed for all $\lambda\in\Lambda$ except for finitely many.
Consider the pointwise sum
\[\sum_{\lambda\in\Lambda}\mathcal  A_\lambda=\biggl\{(v,x)\,:\,x\in X,\,v\in
\sum_{\lambda\in\Lambda}\mathcal  A_\lambda{}^x \biggr\}.\]
Then, we claim, the (pointwise) sum above is  an envelope of sets in $\mathcal  V$ with respect to $\mathfrak S$.
The sum is pointed if all the summands are pointed.

Moreover, if $\mathcal  W$, $\mathcal  Z$ are commutative topological groups,
$\mu:\mathfrak S\rightarrow\mathcal  W$ is a measure,
and  $L:\mathcal  V\times\mathcal  W\rightarrow\mathcal  Z$ is a biadditive pairing
which is continuous in its second variable then
\[\intp L\biggl(\sum_{\lambda\in\Lambda}\mathcal  A_\lambda,\mu\biggr)=
\overline{\sum_{\lambda\in\Lambda}\intp L(\mathcal A_\lambda,\mu)}.\]
\begin{proof}
Let
\[\mathcal A=\sum_{\lambda\in\Lambda}\mathcal  A_\lambda.\]

First we have to prove the envelope properties.
Let $\Xi\subset\Lambda$ be a finite set
such that for all $\lambda\in\Lambda\setminus\Xi$ the envelope
$\mathcal  A_\lambda$ is pointed.
For each $\xi\in\Xi$ let us choose a step-function
\[s_\xi=\sum_{j\in J_\xi} c_{\xi,j}\chi_{E_{\xi,j}}\]
approximated by $\mathcal  A_\xi$.
Then the step-function
\[s=\sum_{\xi\in\Xi}s_\xi=\sum_{\xi\in\Xi,\,j\in J_\xi} c_{\xi,j}\chi_{E_{\xi,j}} \]
will certainly be approximated by $\mathcal A$,
hence the first condition for envelopes is satisfied.

Next, we have to find a generator system for $\mathcal A$.
Let
\[\mathfrak E=\{E_{\xi,j}\,:\,\xi\in\Xi,\,j\in J_\xi\},\]
and  $\mathfrak D_\lambda$ be a generator system for $\mathcal  A_\lambda$.

Let us apply Lemma \ref{lem:FOREST} to the countable family of sets
\[\bigcup_{\lambda\in\Lambda}\mathfrak D_\lambda.\]
such that the result $\widetilde{\mathfrak D}$ should be divided by $\mathfrak E$.
We claim that $\widetilde{\mathfrak D}$ will be a generator system for $\mathcal  A$.

Notice that according to Lemma \ref{lem:APPNEIGH} $\widetilde{\mathfrak D}$ will be a generator
system for all $\mathcal  A_\lambda$.

If, $x\notin\bigcup\widetilde{\mathfrak D}$ implies $\{0\}=\mathcal  A_\lambda{}^x$ for
all $\lambda$, hence $\{0\}=\mathcal  A^x$.
Now assume that $x\in\bigcup\widetilde{\mathfrak D}$ and $v\in\mathcal  A^x$.
What we have to show that there is a set $D\in \widetilde{\mathfrak D}$ such that $v\in\mathcal  A^D$.

By definition, there is a finite set $\Xi\subset\Omega_x \subset\Lambda$ and elements
\[v_\lambda\in\mathcal  A_\lambda{}^x \quad (\lambda\in \Omega_x)\]
such that
\[v=\sum_{\lambda\in \Omega_x}v_\lambda.\]
According to the fact that $\widetilde{\mathfrak D}$ is a generator for all $\mathcal  A_\lambda$
we see that for $\lambda\in\Omega_x$ there is an element $D_{x,\lambda}\in\widetilde{\mathfrak D}$
such that $v_\lambda\in \mathcal  A^{D_{x,\lambda}}$. Let $D_x$ be the smallest of these $D_{x,\lambda}$'s,
ie. the intersection. Then $v_\lambda\in \mathcal  A^{D_{x}}$ for $\lambda\in\Omega_x$.
Hence,
\[v=\sum_{\lambda\in \Omega_x}v_\lambda\in\sum_{\lambda\in \Omega_x}\mathcal  A_\lambda{}^{D_x}.\]

Being $\mathcal  A_\lambda$ pointed for $\lambda\in\Lambda\setminus\Omega_x$ we find that
\[v\in\sum_{\lambda\in\Lambda}\mathcal  A_\lambda{}^{D_x}.\tag{\dag} \]
But then
\[v\in\sum_{\lambda\in\Lambda}\mathcal  A_\lambda{}^{D_x}\subset
\biggl(\sum_{\lambda\in\Lambda}\mathcal  A_\lambda\biggr)^{D_x},\]
proving our statement about the generator set $\widetilde{D}$.
The comment about the sum of pointed envelopes should be obvious.

As for the second part: The $\supset$ direction. It is enough to prove
\[\sum_{\lambda\in\Lambda} \intp L(\mathcal  A_\lambda,\mu)\subset \intp L(\mathcal  A,\lambda).\]
But that immediately follows from Lemma \ref{lem:SUMCLO} and the definition of the the  bilinear sums.

The $\subset$ direction.
It is enough to prove that for each step-function
\[s'=\sum_{j\in J'}c_j'\chi_{E_j'},\]
whose graph is contained in the sum $\mathcal  A$ it yields
\[\sum_{j\in J'}L(c_j',\mu(E_j'))\in\overline{\sum_{\lambda\in\Lambda}\intp L(\mathcal  A_\lambda,\mu)}.\]
Let
\[\mathfrak E'=\{E'_j\,;\,j\in J'\}.\]
Let us apply Lemma \ref{lem:FOREST} to
\[\mathfrak E\cup\mathfrak E'\cup\bigcup_{\lambda\in\Lambda}\mathfrak D_\lambda,\]
such that the resulted set $\widetilde{\mathfrak D}$ is divided by $\mathfrak E\cup\mathfrak E'$.
Notice that both $s$ and $s'$ are constant on elements $D\in\widetilde{\mathfrak D}$.
(And those values are $s^D$ and $s'{}^D$.)

According to our previous arguments, for each $x\in\bigcup\widetilde{\mathfrak D}$ there exist
$D\in \widetilde{\mathfrak D}$ such that
\[s(x),s'(x)\in\sum_{D\in\widetilde{\mathfrak D}}\mathcal  A^{D}.\]
In particular, those sets $D\in \widetilde{\mathfrak D}$ for which
\[s^D,s'{}^D\in\sum_{D\in\widetilde{\mathfrak D}}\mathcal  A^{D}\]
cover $\bigcup\widetilde{\mathfrak D}$. Let $\mathfrak D'$ be the set of the maximal such
elements (remember: $\widetilde{\mathfrak D}$ is forest). Then $\mathfrak D'$ decomposes
$\bigcup\widetilde{\mathfrak D}$, and, by our assumptions, all elements or $\mathfrak E\cup\mathfrak E'$.

According to Lemma \ref{lem:SUM}
\[\sum_{j\in J'}L(c_j',\mu(E_j'))=
\sum_{j\in J}L(c_j,\mu(E_j))+\sum_{D\in\mathfrak D'}^{\ovline} L(s'{}^D,\mu(D))-
\sum_{D\in\mathfrak D'}^{\ovline} L(s^D,\mu(D)),\]
and making contractions in the right side we obtain
\[=\sum_{j\in J}L(c_j,\mu(E_j))+\sum_{D\in\mathfrak D'}^{\ovline} L(s'{}^D-s^D,\mu(D)).\]
Hence, it is enough to prove that for any finite set $K\subset \mathfrak D'$
\[\sum_{j\in J}L(c_j,\mu(E_j))+\sum_{D\in K} L(s'{}^D-s^D,\mu(D))\in
\sum_{\lambda\in\Lambda}\intp L(\mathcal  A_\lambda,\mu),\]
because then the big sum will  certainly be in the closure.

Now, we can finish our proof as follows:

Assume that for $D\in K$
\[s'{}^D=\sum_{\lambda\in\Omega_D}s_{D,\lambda}' \]
where $\Omega_D\supset\Xi$ is finite and $s_{D,\lambda}\in \mathcal  A_\lambda{}^D$.
Set $s_\lambda=0$ for $\lambda\in\Lambda\setminus\Xi$.
Then let
\[h_\lambda=s_\lambda+\sum_{D\in K}-s_{\lambda}|_D+\sum_{D\in K,\,\lambda\in\Omega_D} s'_{D,\lambda}\chi_D. \]
for $\lambda\in\Omega=\bigcup_{K\in D}\Omega_D$ and let $h_\lambda=0$ otherwise.
Now, notice that
\[s+\sum_{D\in K} (s'-s)|_D=\sum_{\lambda\in\Omega}h_\lambda,\]
$h_\lambda$ is a step-function,
and the graph of $h_\lambda$ is contained in  $\mathcal  A_\lambda$.

In particular the evaluated bilinear sum corresponding to $h_\lambda$ (cf. Lemma \ref{lem:SUM})
belongs to $\intp L(\mathcal  A_\lambda,\mu)$,
and, consequently,
\[\sum_{j\in J}L(c_j,\mu(E_j))+\sum_{D\in K} L(s'{}^D-s^D,\mu(D))\in
\sum_{\lambda\in\Lambda}\intp L(\mathcal  A_\lambda,\mu). \]
\end{proof}
\end{lemma}
\begin{lemma}\label{lem:NEGAT}
Let $\mathfrak S\subset \mathfrak P(X)$ be a interval system
and $\mathcal  V$ be a topological vector space.
Suppose that $\mathcal  C$ is  an envelope of pointed  sets in
$\mathcal  V$ with respect to $\mathfrak S$.

Consider
\[-\mathcal  C=\{(x,v)\,:\,x\in X,\,v\in-\mathcal  C^x \}.\]
Then, we claim that the set above is  an envelope of pointed  sets in $\mathcal  V$ with respect to $\mathfrak S$.

Moreover, if $\mathcal  W$, $\mathcal  Z$ are topological spaces,
$\mu:\mathfrak S\rightarrow\mathcal  W$ is a measure,
and  $L:\mathcal  V\times\mathcal  W\rightarrow\mathcal  Z$ is a bilinear pairing
continuous in its second variable then
\[\intp L(-\mathcal  C,\mu)=-\intp L(\mathcal C,\mu).\]
\begin{proof}
Trivial.
\end{proof}
\end{lemma}

\begin{lemma}\label{lem:APPROX}
\halfquote
Let $f:X\rightarrow\mathcal  V$ be a function.

Suppose that $\mathcal  A_1$ and $\mathcal  A_2$ are two envelopes approximating $f$.
Then, we claim,
\[0\in \intp L(\mathcal  A_1-\mathcal  A_2,\mu)= \overline{\intp L(\mathcal  A_1,\mu)-\intp L(\mathcal  A_2,\mu)}.\]
In particular, if $0\in\mathcal  T_1,\mathcal  T_2\subset\mathcal  Z$ are neighborhoods then
\[\biggl(\intp L(\mathcal  A_1,\mu)+\mathcal  T_1 \biggr)\cap
\biggl(\intp L(\mathcal  A_2,\mu)+\mathcal  T_2\biggr)\neq\emptyset. \]

\begin{proof}
The envelope $\mathcal  A_1-\mathcal  A_2$ approximates $f-f=0$. That is a step-function, hence
$0\in \intp L(\mathcal  A_1-\mathcal  A_2,\mu)$. The equality follows from Lemma \ref{lem:NEIGH} and
\ref{lem:NEGAT}.
\end{proof}
\end{lemma}
The following lemma yields a more transparent picture about envelopes:
\begin{lemma}\label{lem:DECO}
Let $\mathfrak S\subset \mathfrak P(X)$ be a interval system
and $\mathcal  V$ be a commutative topological group.
\begin{itemize}
\item[\textit{a.)}] Suppose that
\[s=\sum_{j\in J}c_j\chi_{E_j}\]
is a step-function. Then its graph $s$ is  an envelope.
\item[\textit{b.)}] If the envelope $\mathcal  A$ contains the graph of the step-function $s$, then
\[\mathcal  C=-s+\mathcal  A\]
is a pointed envelope. In this case
\[\mathcal  A=s+\mathcal  C\]
holds.
\end{itemize}
Hence, envelopes containing the graph of $s$ can be specified by pointed envelopes
to be added to $s$.

Moreover, if $\mathcal  W$, $\mathcal  Z$ are commutative topological groups,
$\mu:\mathfrak S\rightarrow\mathcal  W$ is a measure,
and  $L:\mathcal  V\times\mathcal  W\rightarrow\mathcal  Z$ is a biadditive pairing
which is continuous in its second variable then in the cases above
\begin{itemize}
\item[\textit{a'.)}]
\[\intp L(s,\mu)=\sum_{j\in J}L(c_j,\mu(E_j));\]
\item[\textit{b'.)}]
\[\intp L(\mathcal  A,\mu)=\intp L(s,\mu)+\intp L(\mathcal  C,\mu)=
\sum_{j\in J}L(c_j,\mu(E_j))+\intp L(\mathcal  C,\mu)\]
\end{itemize}
hold.
\begin{proof} a.) We can use Lemma \ref{lem:REFIN} for the finite set system
\[\{E_j\,:\,j\in J\}\]
in order to get an exact decomposition. That will be certainly be generator system.

b.) That follows from Lemma \ref{lem:NEIGH} and \ref{lem:NEGAT}.

a'.) According to Lemma \ref{lem:SUM} every sum in the definition of $L(s,\mu)$ gives this single value.

b'.) That follows from Lemma \ref{lem:NEIGH} and \ref{lem:NEGAT}.
\end{proof}
\end{lemma}
~

\paragraph{\textbf{B. More about structure}}~\\

In general, pointed envelopes are much nicer than envelopes . Not only because we
can  add them up freely, but, also, the corresponding bilinear sums are simpler to evaluate:
\begin{lemma}\label{lem:CEVAL}
\halfquote

Suppose that $\mathcal  C$ is a pointed envelope of sets in $\mathcal  V$ with respect to $\mathfrak S$.
Then, we claim,
\begin{multline}
\intp L(\mathcal  C,\mu)=\overline{\biggl\{\sum_{j\in J}L(c_j,\mu(E_j))\,:\,J\text{ is finite},
\,c_j\in\mathcal  V,\,E_j\in\mathfrak S,\ldots}\\
\overline{\ldots\text{ the }E_j\text{ are pairwise disjoint, }
\sum_{j\in J} c_j\chi_{E_j}\subset\mathcal  C\biggr\}}\notag
\end{multline}
One can see that this is a closed pointed set.
\begin{proof}
The $\supset$ part is trivial. On the other hand, consider any finite ``non-disjoint'' sum.
According to Lemma \ref{lem:SUM}.a, such a finite sum can be written as a convergent sum
\[\sum_{\omega\in\Omega} L(c_\omega,E_\omega),\]
where $\Omega$ is countable, $E_\omega\in\mathfrak S$ are pairwise disjoint,
$c_\omega\in\mathcal  C^{E_\omega}$. Being $\mathcal  C$ a pointed envelope,
the finite partial sums belong to the right side, hence, by closure, the full sum.
Taking closure again we obtain $\subset$.

The fact that the zero step-function is approximated guarantees $0\in \intp L(\mathcal  C,\mu)$.
\end{proof}
\end{lemma}
\begin{remark}
Hence, for a pointed envelope $\mathcal  C$ we bilinear sum $\intp L(\mathcal  C,\mu)$ can be
defined by ``finite disjoint'' sums. That can be taken as an alternative definition.
\end{remark}
\begin{defin}
Let $\mathfrak S\subset\mathfrak P(X)$ be an interval system, $\mathcal  V$ be a
commutative topological group. Suppose that $c\in\mathcal  V$, and $A$ is a countable
union of elements of $\mathfrak S$ (ie. $A\in\boldsymbol\Sigma\mathfrak S$).
Then we define the pointed envelope
\[\var(c,A)=(\{c\}\times A)\cup(\{0\}\times X).\]
\end{defin}
\begin{cor}\label{cor:SVAR}
\halfquote

Suppose that $c\in\mathcal  V$, and $A$ is a countable union of elements of $\mathfrak S$. Then
\[\svar(L,c,A,\mu)=\intp L(\var(c,A),\mu).\]
\begin{proof}
We apply Lemma \ref{lem:CEVAL} to the right side.
That immediately yields the definition of the semivariation on the left, however.
\end{proof}
\end{cor}

The following lemma is not needed in order to define the integral but it is a useful tool otherwise.
It says the we can essentially assume that an envelope contains a simple  step-function.
\begin{lemma}\label{lem:APPREFIN}
\halfquote

Suppose that
\[s=\sum_{j\in J}c_j\chi_{E_j}\]
is a step-function, $J$ is finite, $c_j\in\mathcal  V$, $E_j\in\mathfrak S$ $(j\in J)$;
and $\mathcal  C$ is a pointed envelope, so that
\[s+\mathcal  C\]
is  an envelope.
Also suppose that $\mathfrak D\subset\mathfrak S$ decomposes
\[\{E_j\,:\,j\in J\}.\]
(Notice that $s$ takes a constant value $s^D$ on each $D\in\mathfrak D$.)
Let $0\in\mathcal  T'\subset \mathcal  Z$ be a neighborhood.

Then, we claim, there exists a finite set $J'\subset\mathfrak D$,
such that with the choice of the step-function
\[s'= \sum_{D\in J'}s_D\chi_{D},\]
and the pointes envelopes
\[\mathcal  D=\sum_{D\in \mathfrak D\setminus J'}\var(s^D,D),\qquad
\mathcal  C'=\mathcal  C+\mathcal  D,\]
it yields
\[s+\mathcal  C\subset s-\mathcal  D+\mathcal  C=s'+\mathcal  D+\mathcal  C=s'+\mathcal  C'\]
and
\[L(s+\mathcal  C,\mu)\subset L(s'+\mathcal  C',\mu)\subset L(s+\mathcal  C,\mu)+\mathcal  T.\]

A property of the envelope $s'+\mathcal  C'$ above is that $s'$ is a simple step-function.
\begin{proof}
Everything holds independently from the choice of $J'$ except
\[L(s'+\mathcal  C',\mu)\subset L(s+\mathcal  C,\mu)+\mathcal  T.\]
For that reason we have to choose $J'$ well.

Let $0\in\mathcal  T\subset Z$ be a neighborhood such that
\[-\mathcal  T+\mathcal  T+\mathcal  T\subset \mathcal  T'.\]
Let us divide $\mathcal  T$ by $J$. Let
\[\mathfrak E_j=\{D\,:\,D\in\mathfrak D,\,D\subset E_j\}.\]
As Lemma \ref{lem:MEAVAR} shows,
for each $E_j$ there exists a finite set $\Omega_j\subset\mathfrak E_j$ such that
\[\svar\biggl(L,c_j,\mu,\bigcup_{D\in\mathfrak E_j\setminus\Omega_i}D\biggr)
\subset\mathcal  T_j.\]
Let
\[J'=\bigcup_{j\in J}\Omega_j.\]
Define
\[\mathcal  D'=\sum_{j\in J}\var\biggl(c_j,\bigcup_{D\in\mathfrak E_j\setminus\Omega_i}D\biggr).\]
Then
\[\mathcal  D\subset\mathcal  D'.\]
Thus,
\[s'+\mathcal  C'=s+\mathcal  C-\mathcal  D\subset s+\mathcal  C-\mathcal  D'\]
and
\[\intp  L(s'+\mathcal  C',\mu)\subset\intp L(s+\mathcal  C-\mathcal  D',\mu)\]\[=
\overline{\intp  L(s+\mathcal  C,\mu)-
\sum_{j\in J}\intp L\biggl(\var\biggl(c_j,\bigcup_{D\in\mathfrak E_j\setminus\Omega_i}D\biggr) ,\mu\biggr) }\]
\[=\overline{\intp  L(s+\mathcal  C,\mu)-
\sum_{j\in J}\svar\biggl(L,c_j,\bigcup_{D\in\mathfrak E_j\setminus\Omega_i,\mu}D\biggr) }
=\overline{\intp  L(s+\mathcal  C,\mu)-
\sum_{j\in J}\mathcal  T_j }\]
 \[=\overline{ \intp L(s+\mathcal  C,\mu)-
\mathcal  T }\subset\intp  L(s+\mathcal  C,\mu)-\mathcal  T+\mathcal  T
\subset \intp L(s+\mathcal  C,\mu)+\mathcal  T'.\]
\end{proof}
\end{lemma}
\begin{lemma}\label{lem:APPFILTER}
\halfquote

Suppose that $f:X\rightarrow\mathcal  V$ is a function and $\mathcal  A_r$ ($r\in R$)
are finitely many envelopes approximating $f$.
Then, we claim, for any neighborhood $0\in\mathcal  T\subset\mathcal  Z$
there exits  an envelope $\mathcal  A$ approximating $f$ such that
\[\intp L(\mathcal  A,\mu)\subset \bigcap_{r\in R}\biggl(\intp L(\mathcal  A_j,\mu)+\mathcal  T\biggr).\]
\begin{proof}
Let us write the envelopes as
\[\mathcal  A_r=s_r+\mathcal  C_r=\sum_{j\in J_r}c_{r,j}\chi_{E_{r,j}}+\mathcal  C_r.\]
Consider
\[\mathfrak E=\{E_{r,j}\,:\,r\in R,\, j\in J_r\},\]
and a generator system $\mathfrak D_r$ for each $\mathcal  A_r$.
Take a forest refinement $\widetilde{\mathfrak D}$ for
\[\mathfrak D=\bigcup_{r\in R}\mathfrak D_r\]
such that it is divided by $\mathfrak E$.
Notice the for each $D\in\widetilde{\mathfrak D}$ any step-function $s_r$
takes a constant value $s_r{}^D$.

Also notice that $\widetilde{\mathfrak D}$ will be a common generator system for all $\mathcal  A_r$.
Then, for each $x\in\bigcup\widetilde{\mathfrak D}$ and for all $r\in R$
there exists an element $D_{x,r}$ such that
\[f(x)\in\mathcal  A_r{}^{D_{x,r}}.\]
Taking the minimal of those (ie. the intersection) as $D_x$ we find
\[f(x)\in\bigcap_{r\in R}\mathcal  A_r{}^{D_{x}}.\]

In particular the elements $D\in\widetilde{\mathfrak D}$ such that
\[\bigcap_{r\in R}\mathcal  A_r{}^D\neq\emptyset\]
cover $\bigcup \widetilde{\mathfrak D}$. Let $\mathfrak D'$ be the set of the maximal such $D$'s.
Then $\mathfrak D'$ decomposes $\mathfrak E$.

As the previous lemma shows for each $r\in R$ there is a set $\Xi_r\subset \widetilde{\mathfrak D}$ such that for
the envelope
\[ \mathcal  A_r'=s+\mathcal  C_r-\sum_{ D\in\widetilde{\mathfrak D}\setminus\Xi_r}\var(s_r{}^D,D ) \]
we have
\[\intp L(\mathcal  A_r',\mu)\subset \intp L(\mathcal  A_r,\mu)+\mathcal  T.\]
Let
\[\Xi=\bigcup_{r\in R} \Xi_r.\]

We claim that
\[\mathcal  A=\bigcap_{r\in R} \mathcal  A_r'\]
is  an envelope.
Indeed, by definition, for each $D\in\mathfrak D'$ we can choose an element
\[s_D\in\bigcap_{r\in R}\mathcal  A_r{}^D.\]
Then, one can see that $\mathcal  D$ contains the graph of
\[\sum_{\xi\in\Xi} s_D\chi_D.\]
Furthermore, $\widetilde{\mathfrak D}$ will be a generator system
for $\mathcal  A$. That proves that $\mathcal  A$ is  an envelope.

Obviously, $\mathcal  A$ contains the graph of $f$, hence it is an approximation.
 \end{proof}
\end{lemma}

\begin{lemma}\label{lem:UNIOLIM}
Let $\mathfrak S\subset \mathfrak P(X)$ be a interval system
and $\mathcal  V$ be a commutative topological group.
Suppose that $\{\mathcal  A_n\}_{n\in\mathbb N}$ is an indexed family
of  envelopes of sets in $\mathcal  V$ with respect to $\mathfrak S$, such that
\[\mathcal A_0\subset\mathcal A_1\subset \mathcal A_2\subset \mathcal A_3\subset\ldots\,.  \]
Then, we claim, the union
\[\bigcup_{n\in\mathbb N}\mathcal A_n\]
is  an envelope of sets in $\mathcal  V$ with respect to $\mathfrak S$.

Moreover, if $\mathcal  W$, $\mathcal  Z$ are commutative topological groups,
$\mu:\mathfrak S\rightarrow\mathcal  W$ is a measure,
and  $L:\mathcal  V\times\mathcal  W\rightarrow\mathcal  Z$ is a biadditive pairing
which is continuous in its second variable then
\[\intp L\biggl(\bigcup_{n\in\mathbb N}\mathcal  A_n,\mu\biggr)=
\overline{\bigcup_{n\in\mathbb N}\intp L(\mathcal A_n,\mu)}.\]
\begin{proof}
Let us denote the union of $\mathcal A_n$'s by $\mathcal A$.

a.) If $s$ is a step-function contained in $\mathcal A_0$, then $s$ is contained in the union of the
$\mathcal A$. Similarly, if $\mathfrak D_n$ is a generator system for $\mathcal A_n$, then
one can immediately see that the
union of the $\mathfrak D_n$'s will be a generator system for the union of the $\mathcal A$.
That proves our first statement.

b.) The second statement. The ``$\supset$'' part follows from the closedness of the left side.
The ``$\subset$'' can be proven as follows. If the step function $s$ contained in $\mathcal A_0$
then
\[-s+\mathcal A_0\subset-s+\mathcal A_1\subset -s+\mathcal A_2\subset -s+\mathcal A_3\subset\ldots\,, \]
\[\bigcup_{n\in\mathbb N}(-s+\mathcal A_n)=-s+\bigcup_{n\in\mathbb N}\mathcal A_n,\]
and Lemma \ref{lem:DECO}
shows that it is enough to prove the statement for pointed envelopes.

So, assume that $\mathcal A_n=\mathcal C_n$ are pointed envelopes.
It is enough to show that for each step-function
\[s=\sum_{j\in J}c_j\chi_{E_j}\subset\mathcal A\]
we have
\[\sum_{j\in J}L(c_j,\mu({E_j}))\in\intp L(\mathcal A,\mu).\]

Let $\mathfrak E=\{E_j\,:\,j\in J\}$.
Let $\mathfrak D_n$ be a generator system for $\mathcal A$.
Let us use Lemma \ref{lem:FOREST} in order to obtain a forest refinement
 $\widetilde{\mathfrak D}$ of
\[\mathfrak D=\bigcup_{n\in\mathbb N}\mathfrak D_n\cup\mathfrak E, \]
such that it is divided by $\mathfrak E$.

Then, for each $x\in\bigcup\mathfrak E$ there is a maximal set $D\in \widetilde{\mathfrak D} $
such that there exists an element $n\in\mathbb N$ such that
\[s^D\in (\mathcal A_n)^D.\]
These sets form an exact decomposition $\mathfrak D'$ of $\mathfrak E$.
According to Lemma \ref{lem:SUM} we see that
\[ \sum_{j\in J}L(c_j,\mu({E_j}))=\sum^{\ovline}_{D\in\mathfrak D'}L(s^D,\mu(D)).\tag{\S}\]
Now, any finite partial sum of the right side is contained in a set
\[\intp L(\mathcal A_n,\mu),\]
henceforth in
\[\bigcup_{n\in\mathbb N}\intp L(\mathcal A_n,\mu).\]
Consequently,  (\S) is contained in
\[\overline{\bigcup_{n\in\mathbb N}\intp L(\mathcal A_n,\mu)}.\]
And that is what we wanted to prove.
\end{proof}
\end{lemma}

\newpage
\section{The Lebesgue-McShane integral}\label{sec:LMcS}
\paragraph{{\textbf{A. Fundamentals}}}
\begin{defin}\label{def:LMINT}
\halfquote
Suppose that $f: X\rightarrow\mathcal  V$ is function.

We say is that the value $a\in\mathcal  Z$ is the integral of $f$ with respect to $L,\mu$ if
for each neighborhood $0\in\mathcal  T\subset \mathcal  Z$ there exists  an envelope $\mathcal  A$ approximating  $f$ such that
\[\intp L(\mathcal  A,\mu)\subset a+\mathcal  T.  \]
\end{defin}
\begin{lemma}\label{lem:LMINT}
\halfquote
Suppose that $f: X\rightarrow\mathcal  V$ is function.

The integral of $f$, if exists, is unique. Sufficient and necessary condition for the
existence of the integral is that
for each neighborhood $0\in\mathcal  T\subset \mathcal  Z$ there exists
 an envelope $\mathcal  A$ approximating $f$ such that
\[\intp L(\mathcal  A,\mu)-\intp L(\mathcal  A,\mu)\subset \mathcal  T.\]
\begin{proof}

The necessity part of the existence statement immediately follows from the definition.
The sufficiency can be proven as follows:
From our assumption and Lemma \ref{lem:APPROX} it follows that
\[\mathfrak H=\biggl\{\intp L(\mathcal  A,\mu)\,:\,\text{ $\mathcal  A$ approximates $f$ }\biggr\},\]
forms a Cauchy system.
Then Lemma \ref{lem:CAUCHY} and completeness imply convergence of $\mathfrak H$.
The limit will necessarily be  an integral.
The convergence of $\mathfrak H$ and the Hausdorff property implies uniqueness.
(Remark: Here we did not use the more advanced Lemma \ref{lem:APPFILTER}.)
\end{proof}
\end{lemma}
\begin{notat}
We use the notation
\[\int^{(LM)} L(f,\mu),\]
or rather just
\[\int L(f,\mu)\]
for the integral.
\end{notat}
\begin{remark}
One can see that Definition \ref{def:01} and \ref{def:LMINT} are equivalent.
In the first case  certain sets $B$ are supposed to form a finer system then
neighborhood filter of $a$, while in the second case that should happen with sets $\overline B$.
But this distinction does not matter, because if $0\in\mathcal U\subset\mathcal  Z$ is a neighborhood
and $\mathcal U-\mathcal U\subset\mathcal T$
then $\overline B\subset a-\mathcal U$ implies $\overline B\subset B+\mathcal U-\mathcal U\subset a-\mathcal T.$

(Also, more generally, the corresponding equiconvergent systems $\mathfrak H^{\mathrm{tf}}$
are the same.)
\end{remark}
\begin{lemma}\label{lem:INTCONT}
\halfquote

Suppose that $f: X\rightarrow\mathcal  V$ is an integrable function.
Suppose that the envelope  $\mathcal  A$ is approximates of $f$
Then, we claim,
\[\int L(f,\mu)\in \intp L(\mathcal  A,\mu). \]
\begin{proof}
Follows from the definition of convergence of set systems (applied to the system $\mathfrak H$ as in the proof
of Lemma \ref{lem:LMINT}) and the fact that the right side is closed.
\end{proof}
\end{lemma}
\begin{lemma}\label{lem:LMADD}
\halfquote
Suppose that $f,g:X\rightarrow\mathcal  V$ are integrable functions.
Then, we claim, $f+g$, $-f$ are integrable, and
\[\int L(f+g,\mu)=\int L(f,\mu)+\int L(g,\mu),\qquad \int L(-f,\mu)=-\int L(f,\mu).\]
\begin{proof}
Let $0\in\mathcal  U\subset\mathcal  Z$ be a neighborhood
such  that $\mathcal  U+\mathcal  U+\mathcal  U\subset\mathcal  T$.
Let $\mathcal  A$ and $\mathcal  B$ envelopes approximating $f$ and $g$ respectively such that
\[\intp L(\mathcal  A,\mu )\subset \int L(f,\mu)+\mathcal  U\qand
 \intp L(\mathcal B,\mu )\subset \int L(g,\mu)+\mathcal  U.\]
Then $\mathcal  A+\mathcal  B$ approximates $f+g$ and
\begin{multline}
\intp L(\mathcal  A+\mathcal  B,\mu)=\overline{\intp L(\mathcal  A,\mu )+ \intp L(\mathcal  B,\mu )}
\subset \intp L(\mathcal  A,\mu )+ \intp L(\mathcal  B,\mu )+\mathcal  U\\\subset
\int L(f,\mu)+\int L(g,\mu)+\mathcal  U+\mathcal  U+\mathcal  U\subset
\int L(f,\mu)+\int L(g,\mu)+\mathcal  T.\notag
\end{multline}
Hence,
\[\int L(f,\mu)+\int L(g,\mu)\]
satisfies the definition of integral for $f+g$.
For $-f$, invert everything.
\end{proof}
\end{lemma}
The integral is monotone:
\begin{notat}
Suppose that $\mathcal  V$ is a commutative group and $\mathsf V$ is an additive
cone there.

a.) We write $a\leq_{\mathsf V}b$ if $b-a\in\mathsf V$.

b.) The interval $[a,b]_{\mathsf V}$ contains those elements $c$
for which $a\leq_{\mathsf V}c\leq_{\mathsf V}b$.
\end{notat}
\begin{lemma}\label{lem:MON}
\halfquote

Suppose that $\mathsf V$ and $\mathsf Z$
are additive cones in $\mathcal  V$ and $\mathcal  Z$ respectively, $\mathsf Z$
is closed, and for each  $v\in\mathsf V$, $E\in\mathfrak S$,
\[L(v,\mu(E))\in \mathsf Z.\]
Suppose that $f:X\rightarrow\mathcal  V$ is integrable and $f(x)\in\mathsf V$ (pointwise).
Then, we claim,
\[\int L(f,\mu )\in\mathsf Z.\]
\begin{proof}
Being $f$ integrable, there exist a set $A$, which is a countable union of element
of $\mathfrak S$, and $f$ vanishes outside $A$. Consider the
pointed envelope
\[\mathcal  C=(\mathsf V\times A)\cup(\{0\}\times X).\]
Then $\mathcal  C$  approximates on $f$, hence
\[\int L(f,\mu )\in \intp L(\mathcal  C,\mu).\]
On the other hand,
\[\intp L(\mathcal  C,\mu)\subset\mathsf Z;\]
because the finite ``disjoint'' sums for $\intp L(\mathcal  C,\mu)$
yield elements of $\mathsf Z$, and $\mathsf Z$ is closed.
\end{proof}
\end{lemma}
\begin{lemma}\label{lem:ALGBP}
\halfquote

Suppose that $\Lambda$ is countable,  $f_\lambda:X\rightarrow\mathcal  V$ ($\lambda\in\Lambda$)
are integrable functions, and
\[f=\sum_{\lambda\in\Lambda}f_\lambda\]
pointwise, algebraically, ie. for each $x\in X$ at most finitely many $f_\lambda(x)$ are nonzero.

Also assume that
for each neighborhood $0\in\mathcal  T\subset\mathcal  Z$
there is exist a a finite set $\Xi\subset\Lambda$ and envelopes $\mathcal  C_\lambda$
($\lambda\in\Lambda\setminus\Xi$) such that $\mathcal  C_\lambda$ is an approximation of $f_\lambda$ and
\[\sum_{\lambda\in\Lambda\setminus\Xi}\intp L(\mathcal  C_\lambda,\mu)\subset\mathcal  T.\]

Then, we claim that $f$ is integrable and
\[\int L(f,\mu)=\sum_{\lambda\in\Lambda}^{\ovline}\int L(f_\lambda,\mu).\]
\begin{proof}
Our fundamental assumption and
\[\int L(f_\lambda,\mu)\in \intp L(\mathcal  C_n,\mu)\]
immediately implies that
\[\sum_{\lambda\in\Lambda}^{\ovline}\int L(f_\lambda,\mu).\]
is convergent according to the Cauchy criterium.

Let $0\in\mathcal  U\subset\mathcal  Z$ be a neighborhood
such that $\mathcal  U+\mathcal  U+\mathcal  U\subset\mathcal  T$.

We can find a $\Xi$ and pointed envelopes $\mathcal  C_\lambda$ ($\lambda\in\Lambda\setminus\Xi$)
approximating $f_\lambda$ such that
\[\sum_{\lambda\in\Lambda\setminus\Xi}\intp L(\mathcal  C_\lambda,\mu)\subset\mathcal  U.\]
We can assume that $\Xi$ is so large that
\[\sum_{\lambda\in\Xi}\int L(f_\lambda,\mu)\in
\sum_{\lambda\in\Lambda}^{\ovline}\int L(f_\lambda,\mu)+\mathcal  U. \]

Let us divide $\mathcal  U$ by $\Xi$ and choose
envelopes $\mathcal  A_\xi$ approximating $f_\xi$ such that
\[\intp L(\mathcal  A_\xi,\mu)\subset \int L(f_\xi,\mu)+\mathcal  U_\xi\]
for $\xi$. Consider the envelope
\[\mathcal  A=\sum_{\xi\in\Xi}\mathcal  A_\xi+\sum_{\lambda\in\Lambda\setminus\Xi }\mathcal  C_\lambda.\]
Clearly, $\mathcal  A$ approximates $f$.
On the other hand,
\[\intp L(\mathcal  A,\mu)\subset\overline{\sum_{\xi\in\Xi}\intp L(\mathcal  A_\xi,\mu)+
\sum_{\lambda\in\Lambda\setminus\Xi}\intp L(\mathcal  C_\lambda,\mu)}
\]\[\subset \sum_{\xi\in\Xi}\intp L(\mathcal  A_\xi,\mu)+
\sum_{\lambda\in\Lambda\setminus\Xi}\intp L(\mathcal  C_\lambda,\mu)+\mathcal  U\]
\[\subset\sum_{\xi\in\Xi}\int L(f_\xi,\mu)+\sum_{\xi\in\Xi}\mathcal  U_\xi+
\sum_{\lambda\in\Lambda\setminus\Xi}\mathcal  U_\lambda+\mathcal  U
\subset \sum_{\xi\in\Xi}\int L(f_\xi,\mu)+\mathcal  U+\mathcal  U\]
\[\subset\sum_{\lambda\in\Lambda}^{\ovline}\int L(f_\lambda,\mu)+\mathcal  U+\mathcal  U+\mathcal  U
\subset\sum_{\lambda\in\Lambda}^{\ovline}\int L(f_\lambda,\mu)+\mathcal  T. \]
\end{proof}
\end{lemma}

\paragraph{{\textbf{B. More about structure}}}
\begin{lemma}\label{lem:ADDSUPER}
\halfquote

Suppose that $\Lambda$ is countable,  $f,f_\lambda:X\rightarrow\mathcal  V$ ($\lambda\in\Lambda$)
are integrable functions and
\[\sum_{\lambda\in\Lambda}f_\lambda=f,\]
such that at each point $x\in X$ at most one of the $f_\lambda(x)$'s is nonzero.

Then we claim, for each neighborhood $0\in\mathcal  T\subset\mathcal  Z$ there exist an finite set
$\Xi$, and a family of pointed envelopes $\mathcal  C_\lambda$ ($\lambda\in\Lambda\setminus\Xi$) such that
$\mathcal  C_\lambda$ approximates $f_\lambda$ and
\[\intp L\biggl(\sum_{\lambda\in\Lambda\setminus\Xi}\mathcal  C_\lambda,\mu\biggr)\subset\mathcal  T.\]
In particular,
\[\int L(f,\mu)=\sum^{\ovline}_{\lambda\in\Lambda}\int L(f_\lambda,\mu). \]
\begin{proof}
Consider a neighborhood $0\in\mathcal T\subset\mathcal Z$.
Let $0\in\mathcal U\subset\mathcal Z$ be a neighborhood such that
\[\mathcal U+\mathcal U-\mathcal U+\mathcal U-\mathcal U+\mathcal U\subset\mathcal T.\]
Let us divide $\mathcal U$ by $\Lambda$.

Let us consider  an envelope
\[\mathcal  A_\lambda=\sum_{j\in J_\lambda}c_{\lambda,j}\chi_{E_{\lambda,j}}+\mathcal  C_\lambda\]
for each $f_\lambda$ such that
\[\intp L(\mathcal  C_\lambda,\mu)\subset\mathcal  U_\lambda.\]
And similarly,  an envelope
\[\mathcal  A=\sum_{j\in J}c_j\chi_{E_{j}}+\mathcal  C\]
for $f$ such that
\[\intp L(\mathcal  C,\mu)\subset\mathcal  U.\]
According to Lemma \ref{lem:APPREFIN} we can assume that for a fixed $\lambda$ the sets
$E_{\lambda,j}$ are pairwise disjoint; and similarly that the sets $E_j$ are pairwise disjoint.
Let $\mathfrak E=\{E_j\,:\,j\in J\}$

Let $\mathfrak D$ be a generator system for $\mathcal  A$ and $\mathfrak D_\lambda$ be
a generator system for $\mathcal  A_\lambda$.  Apply Lemma \ref{lem:FOREST} for
\[\mathfrak D\cup\bigcup_{\lambda\in\Lambda}\mathfrak D_\lambda\cup\mathfrak E\]
such that the resulted forest $\widetilde{\mathfrak D}$ is divided by the set $\mathfrak E$.

Now, for each point $x\in E_j$ there exists an index $\lambda\in\Lambda$ such that $f(x)=f_\lambda(x)$.
Then, there exists a set $D\in\widetilde{\mathfrak D}$ such that
$x\in D\subset E_j$ and  $f(x)\in \mathcal  A^D\cap\mathcal  A_\lambda{}^D$.

In particular, the sets $D\in\widetilde{\mathfrak D}$ such that
there exists $\lambda\in\Lambda$ such that
\[\bigcup_{\lambda\in\Lambda}\mathcal  A^D\cap\mathcal  A_\lambda{}^D\neq\emptyset\]
cover all sets $E_j$.
For $J\in J$ let $\mathfrak E_j$ be the set of maximal such $D$'s contained in $E_j$.
Then we know that $\mathfrak E_j$ exactly decomposes $E_j$.

Let us divide $\mathcal  U$ by $J$.
Now, applying Lemma \ref{lem:MEAVAR} we can obtain finite families
$\mathfrak B_j\subset\mathfrak E_j$ for all $j\in J$ such that
\[\svar\left(L,c_j,{\textstyle\bigcup}(\mathfrak E_j\setminus\mathfrak B_j) ,\mu\right)
\subset\mathcal U_j.\]
Let
\[\mathfrak B=\bigcup_{j\in J}\mathfrak B_j.\]
We know that for each $D\in\mathfrak B$ there exists an index $\lambda_D\in\Lambda$
such that
\[\mathcal  A^D\cap\mathcal  A_{\lambda_D}{}^D\neq\emptyset.\]
Let
\[\Xi=\{\lambda_D\,:\,D\in\mathfrak B\}.\]
Now, let
\[\mathcal  C_\lambda'=(\{0\}\times X)\cup\mathcal  A_\lambda.\]
We claim that these $\mathcal  C'_\lambda$'s satisfy our requirements.
First of all, it is clear that $\mathcal  C'_\lambda$ is a pointed envelope
approximating $f_\lambda$.

Now let us examine the behavior of these $\mathcal  C'_\lambda$'s over each point $x\in X$.
If $f_\lambda(x)=0$ then $0\in \mathcal  A_\lambda$, hence
\[\mathcal  C_\lambda'{}^x\subset \mathcal  A_\lambda{}^x-\mathcal  A_\lambda{}^x=
\mathcal  C_\lambda{}^x-\mathcal  C_\lambda{}^x.\tag{x}\]
In the single exceptional case $f_\lambda(x)\neq 0$, $\lambda\in\Lambda\setminus\Xi$
there are various cases:

a.) If $x$ is not in any of the the $E_j$'s:
Then $f(x)\in\mathcal  C^x$ and $f(x)\in\mathcal  A_\lambda{}^x$.
Now,
\[\mathcal  A_\lambda{}^x\subset f(x)-\mathcal  C_\lambda{}^x+\mathcal  C_\lambda{}^x
\subset \mathcal  C^x-\mathcal  C_\lambda{}^x+\mathcal  C_\lambda{}^x,\]
implies
\[\mathcal  C_\lambda'{}^x\subset\mathcal  C^x-\mathcal  C_\lambda{}^x+\mathcal  C_\lambda{}^x.
\tag{a}\]

b.) If $x\in E_j$, $x\in D\in\mathfrak E_j\setminus\mathfrak B_j$ :
Then $f(x)\in c_j+\mathcal  C^x$ and $f(x)\in\mathcal  A_\lambda{}^x$.
Now
\[\mathcal  A_\lambda{}^x\subset f(x)-\mathcal  C_\lambda{}^x+\mathcal  C_\lambda{}^x
\subset c_j+\mathcal  C^x-\mathcal  C_\lambda{}^x+\mathcal  C_\lambda{}^x,\]
implies
\[\mathcal  C_\lambda'{}^x\subset
\svar\left(c_j,{\textstyle\bigcup}(\mathfrak E_j\setminus\mathfrak B_j)\right) {}^x+
\mathcal  C^x-\mathcal  C_\lambda{}^x+\mathcal  C_\lambda{}^x.\tag{b}\]

c.) If $x\in E_j$, $x\in D\in\mathfrak B_j$ :
Then $f(x)\in c_j+\mathcal  C^x$, $f(x)\in\mathcal  A_\lambda{}^x$ and
$\mathcal  A^D\cap\mathcal  A_{\lambda_D}{}^D\neq\emptyset$.
The latter fact implies that there
is a common element $v\in\mathcal  A^x\cap\mathcal  A_{\lambda_D}{}^x$.
Then
\[\mathcal  A_\lambda{}^x\subset f(x)-\mathcal  C_\lambda{}^x+\mathcal  C_\lambda{}^x
\subset  v+\mathcal  C^x-\mathcal  C^x-\mathcal  C_\lambda{}^x+\mathcal  C_\lambda{}^x.\]

On the other hand, $\Xi\ni\lambda_D\neq\lambda\in\Lambda\setminus\Xi$, which implies
 $0\in\mathcal  A_{\lambda_D}$, hence
 \[v\in \mathcal  A_{\lambda_D}{}^x\subset \mathcal  A_{\lambda_D}{}^x-\mathcal  A_{\lambda_D}{}^x
 =\mathcal  C_{\lambda_D}{}^x-\mathcal  C_{\lambda_D}{}^x\subset
\sum_{\xi\in\Xi}(\mathcal  C_{\xi}{}^x-\mathcal  C_{\xi}{}^x) .\]
Consequently,
\[\mathcal  A_\lambda{}^x\subset  \sum_{\xi\in\Xi}(\mathcal  C_{\xi}{}^x-\mathcal  C_{\xi}{}^x)
+\mathcal  C^x-\mathcal  C^x-\mathcal  C_\lambda{}^x+\mathcal  C_\lambda{}^x.\]
or, more generally
\[\mathcal  C_\lambda'{}^x\subset  \sum_{\xi\in\Xi}(\mathcal  C_{\xi}{}^x-\mathcal  C_{\xi}{}^x)
+\mathcal  C^x-\mathcal  C^x-\mathcal  C_\lambda{}^x+\mathcal  C_\lambda{}^x.\tag{c}\]
Overall, from (a), (b), (c), in the exceptional case we can make the crude statement,
\[\mathcal  C_\lambda'{}^x\subset
\sum_{j\in J}\svar\left(c_j,{\textstyle\bigcup}(\mathfrak E_j\setminus\mathfrak B_j)\right) {}^x
+\sum_{\xi\in\Xi}(\mathcal  C_{\xi}{}^x-\mathcal  C_{\xi}{}^x)
+\mathcal  C^x-\mathcal  C^x-\mathcal  C_\lambda{}^x+\mathcal  C_\lambda{}^x.\tag{y}\]
Summing for $\lambda\in\Lambda\setminus\Xi$ the non-exceptional cases (x) and the exceptional case (y) we obtain
\[\sum_{\lambda\in\Lambda\setminus\Xi}
\mathcal  C_\lambda'{}^x\subset
\sum_{j\in J}\svar\left(c_j,{\textstyle\bigcup}(\mathfrak E_j\setminus\mathfrak B_j)\right) {}^x
+\sum_{\lambda\in\Lambda}(\mathcal  C_{\lambda}{}^x-\mathcal  C_{\lambda}{}^x)
+\mathcal  C^x-\mathcal  C^x.\]
This is true for all $x\in X$, so it yields
\[\sum_{\lambda\in\Lambda\setminus\Xi}
\mathcal  C_\lambda'\subset
\sum_{j\in J}\svar\left(c_j,{\textstyle\bigcup}(\mathfrak E_j\setminus\mathfrak B_j)\right)
+\sum_{\lambda\in\Lambda}(\mathcal  C_{\lambda}-\mathcal  C_{\lambda})
+\mathcal  C-\mathcal  C.\]
Evaluating,
\[\intp L\biggl(\sum_{\lambda\in\Lambda\setminus\Xi}\mathcal  C_\lambda',\mu\biggr)=
 \overline{\sum_{j\in J}L(\svar\left(c_j,{\textstyle\bigcup}(\mathfrak E_j\setminus\mathfrak B_j)\right),\mu)
+\ldots}\]\[\overline{\ldots+
\sum_{\lambda\in\Lambda}\biggl(\intp L(\mathcal  C_{\lambda},\mu)-\intp L(\mathcal  C_{\lambda},\mu)\biggr)
+\intp L(\mathcal  C,\mu)-\intp L(\mathcal  C,\mu)}  \]
\[\subset\sum_{j\in J}\mathcal  U_j+\sum_{\lambda\in\Lambda}\mathcal  U_\lambda-
\sum_{\lambda\in\Lambda}\mathcal  U_\lambda+\mathcal  U-\mathcal  U+\mathcal  U
\subset\mathcal  U+\mathcal  U-\mathcal  U+\mathcal  U-\mathcal  U+\mathcal  U\subset\mathcal  T.\]
That proves our main statement.

In particular, we see that the assumptions of the previous lemma hold, so the statement
about the integral follows.
\end{proof}
\end{lemma}
\begin{lemma}\label{lem:LOCI}
\halfquote

Suppose that $f:X\rightarrow\mathcal  V$ is an integrable function.
Then, we claim that for each $E\in\mathfrak S$ the function $f|_E$, given by
\[f|_E(x)=\begin{cases} f(x)& \qquad\text{ if }x\in E \\
0& \qquad\text{ if }x\notin E\end{cases}\]
is integrable. Moreover,
\[\mu_{L,f}(E)=\int L(f|_E,\mu)\]
is a measure on $\mathfrak S$.
\begin{proof}
Consider  an envelope
\[\mathcal  A=\sum_{j\in J}c_j\chi_{E_j}+\mathcal  C\]
approximating $f$
such that
\[\intp L(\mathcal  A,\mu)-\intp L(\mathcal  A,\mu)\subset\mathcal  T.\]
According to Lemma \ref {lem:APPREFIN} we can assume that
the sets $E_j$ are pairwise disjoint and each of them is contained in $E$ or it is
disjoint from it.
(We take a forest refinement of the original $E_j$'s such that it is divided by $E$ and that yields a set
system $\mathfrak D$ to ``support'' the simple step function.)
Then
\[\mathcal  A_E=\sum_{j\in J,\,E_j\subset E}c_j\chi_{E_j}+\mathcal  C\]
approximates $f|_E$ and
\[\intp L(\mathcal  A_E,\mu)-\intp L(\mathcal  A_E,\mu)=
\intp L(\mathcal  A,\mu)-\intp L(\mathcal  A,\mu)\subset\mathcal  T.\]
Hence the integrability of $f|_E$ follows from the Cauchy criterium.

Now, the $\sigma$-additive property of the measure immediately follows from
the previous lemma. (But one can easily give a direct proof, too.)
\end{proof}
\end{lemma}

\begin{defin}
Suppose that $\mathcal  Z$ is a commutative topological group.
A neighborhood subbasis $\mathfrak B$ is a family of neighborhoods
of $0$ in $\mathcal  Z$ such that the following holds:
\begin{itemize}
\item[i.)] For any neighborhood $0\in\mathcal  T\subset\mathcal  Z$ there are
finitely many elements $B_1,\ldots,B_n$ in $\mathfrak B$ such that
\[B_1\cap\ldots\cap B_n\subset \mathcal  T.\]
\end{itemize}
The neighborhood subbasis $\mathfrak B$ is divisible if
\begin{itemize}
\item[ii.)] For any $B\in\mathfrak B$ there exist $B_1,B_2\in\mathfrak B$ such that
\[B_1+B_2\subset B.\]
\end{itemize}
\end{defin}
\begin{lemma}\label{lem:APPDIV}
\halfquote

Suppose that $f:X\rightarrow\mathcal  V$ is a function.
Assume that $\mathfrak B$ is a neighborhood subbasis.
Then we claim:

a.) $f$ has integral $a$ if and only if for each $B\in\mathfrak B$ there is  an envelope
$\mathcal  A$ approximating $f$ such that
\[\intp L(\mathcal  A,\mu)\subset a- B.\]

b.) If $\mathfrak B$ is divisible then
$f$ is integrable if and only if for each $B\in\mathfrak B$ there is  an envelope
$\mathcal  A$ approximating $f$ such that
\[\intp L(\mathcal  A,\mu)-\intp L(\mathcal  A,\mu)\subset B.\]
\begin{proof}
a.) For an arbitrary neighborhood $0\in\mathcal  T\subset\mathcal  Z$ we can choose elements
from $\mathfrak B$ such that
\[B_1\cap\ldots\cap B_n\subset\mathcal  T.\]
Lemma \ref{lem:APPFILTER} implies that there is  an envelope $\mathcal  A$ such that
\[\intp L(\mathcal  A,\mu)\subset (a-B_1)\cap\ldots\cap(a-B_n).\]
Then
\[a-\intp L(\mathcal  A,\mu)\subset B_1\cap\ldots\cap B_n\subset\mathcal  T,\]
hence
\[\intp L(\mathcal  A,\mu)\subset a-\mathcal  T.\]

b.) Similarly, for an arbitrary neighborhood $0\in\mathcal  T\subset\mathcal  Z$ we can choose elements
from $\mathfrak B$ such that
\[B_1\cap\ldots\cap B_n\subset\mathcal  T.\]
Then, we can divide $B_r$ such that
\[B_{r,1}+B_{r,2}\subset B_r,\]
and we can find a neighborhood $0\in\mathcal  U\subset\mathcal  Z$ such that
\[\mathcal  U-\mathcal  U\subset B_{r,2}\]
for all $r$.
By our assumptions, we can choose envelopes $\mathcal  A_r$ such that
\[\intp L(\mathcal  A_r,\mu)-\intp L(\mathcal  A_r,\mu)\subset B_{r,1}.\]
Then, by Lemma \ref{lem:APPFILTER}, we can find  an envelope $\mathcal  A$ such that
\[\intp L(\mathcal  A,\mu)\subset \intp L(\mathcal  A_r,\mu)+\mathcal  U\]
for all $r$. That implies
\[\intp L(\mathcal  A,\mu)-\intp L(\mathcal  A,\mu)\subset
\intp L(\mathcal  A_r,\mu)-\intp L(\mathcal  A_r,\mu)+B_{r,2}\subset B_{r,1}+B_{r,2}\subset B_r\]
for all $r$, hence
\[L(A,\mu)\subset B_1\cap\ldots\cap B_n\subset\mathcal  T.\]
According to the Cauchy criterium this implies integrability.
\end{proof}
\end{lemma}
\begin{remark}
The divisibility of $\mathfrak B$ can be weakened to the  condition that
for each $B\in\mathfrak B$ there is an element $B_0\in\mathfrak B$ and a neighborhood
$0\in\mathcal  U\subset\mathcal  Z$ such that
\[B_0+\mathcal U\subset B.\]
\end{remark}
\begin{remark}
Lemma \ref{lem:APPDIV}  is useful in the case of induced topologies, including the weak and
strong topologies. It says the integration  componentwise is the same as
integration in the induced topology.
\end{remark}

\newpage\section{Semivariation and negligible sets}
\paragraph{\textbf{A. Semivariation of integrable functions}}
\begin{defin}
\halfquote

A function $f:X\rightarrow \mathcal  V$ is locally $L,\mu$-integrable if for all
$A\in\mathfrak S$ the function $f|_A$ is $L,\mu$-integrable.
(We see that according to Lemma \ref{lem:LOCI} all integrable functions are locally integrable.)
\end{defin}
\begin{defin}
\halfquote
Suppose that $f:X\rightarrow \mathcal  V$ is locally integrable.
Then the semivariation of $f$ is
\begin{multline}
\svar(L,f,\mu)=\\
\overline{\biggl\{\sum_{j\in J}\int L(f|_{A_j},\mu)\,:
\,J\text{ is finite},\,A_j\in\mathfrak S,\,
\text{ the $A_j$ are pairwise disjoint}\biggr\}}.\notag
\end{multline}
Ie., in other terms, this is just the semivariation of the measure $\mu_{L,f}$.
\end{defin}
\begin{lemma}\label{lem:SVCHAR}
\halfquote
Then we claim:

a.) If $f,g$ are locally integrable functions then
\[\svar(L,f+g,\mu)\subset\overline{\svar(L,f,\mu)+\svar(L,g,\mu)}\]
and
\[\svar(L,-f,\mu)=-\svar(L,f,\mu).\]

b.) If $f,g$ are locally integrable functions and $f-g\subset\mathcal  C$ for the pointed envelope $\mathcal  C$
then
\[\svar(L,f,\mu)\subset \overline{\svar(L,g,\mu)+\intp L(\mathcal  C,\mu)}\]

c.) If $s:\rightarrow\mathcal  V$ is a step-function then
\[\svar(L,s,\mu)=L(s\cup(\{0\}\times X),\mu).\]

d.) If $f$ is integrable and $\mathcal  A$ is  an envelope
approximating $f$ then
\[\svar(L,f,\mu)\subset \intp L(\mathcal  A\cup (\{0\}\times X),\mu)\subset\]
\[\subset\overline{\svar(L,f,\mu)+\intp L(\mathcal  A,\mu)-\intp L(\mathcal  A,\mu)}.\]

e.) If $f$ is integrable and $\mathcal  C$ is a pointed envelope
approximating $f$ then
\[\svar(L,f,\mu)\subset \intp L(\mathcal  C,\mu). \]

f.) On the other hand, if $f$ is integrable and
$0\in\mathcal  T\subset\mathcal  Z$ is an arbitrary neighborhood
then, we claim, there exists a pointed envelope $\mathcal  C$  approximating $f$, such that
\[\intp L(\mathcal  C,\mu)\subset \svar(L,f,\mu)+\mathcal  T.\]

g.) If $f$ is integrable then
\[\int L(f,\mu)\in\svar(L,f,\mu).\]
\begin{proof}
a.) Integrating
\[\sum_{j\in J}(f+g)|_{A_j}=\sum_{j\in J}f|_{A_j}+\sum_{j\in J}g|_{A_j}\]
we obtain that
\[\sum_{j\in J}\int L((f+g)|_{A_j},\mu)=\sum_{j\in J}\int L(f|_{A_j},\mu)
+\sum_{j\in J}\int L(g|_{A_j},\mu)\]\[\in\svar(L,f,\mu)+\svar(L,g,\mu).\]
Taking closure yields our first statement. The statement about $-f$ is straightforward.

b.) The function
\[\sum_{j\in J}(f-g)|_{A_j}\]
is approximated by $\mathcal  C$, hence, by Lemma \ref{lem:INTCONT}
\[\sum_{j\in J}\int L(f|_{A_j},\mu)=\sum_{j\in J}\int L(g|_{A_j},\mu)+
\sum_{j\in J}\int L((f-g)|_{A_j},\mu)\]
\[\in \svar(L,g,\mu)+\intp L(\mathcal  C,\mu).\]
Taking closure yields our statement.

e.) This is the special case of point b. with $g=0$.

c.) On simply check that $s\cup(\{0\}\times X)$ is a pointed envelope.
Then $\subset$ follows from point e. On the other hand, every finite disjoint sum
as in Lemma \ref{lem:CEVAL} is of form
\[\sum_{j\in J} L(s|_{A_j},\mu).\]
That implies $\supset$.

d.) The first containment follows from point e. We have to prove the second one.
We can assume that
\[\mathcal  A=s+\mathcal  C,\]
such that $s$ is step-function and $\mathcal  C$ is  an envelope. Now,
\[\intp L(\mathcal  A\cup (\{0\}\times X),\mu)\subset\intp L((s\cup(\{0\}\times X))+\mathcal  C,\mu  )=
\overline{\svar(s,\mu)+\intp L(\mathcal  C,\mu)}.\]
Then, by b.,
\[\subset\overline{\overline{\svar(f,\mu)+L(-\mathcal  C,\mu)}+\intp L(\mathcal  C,\mu)  }=
\overline{{\svar(f,\mu)-\intp L(\mathcal  C,\mu)}+\intp L(\mathcal  C,\mu)  }\]
\[=\overline{{\svar(f,\mu)-\intp L(\mathcal  A,\mu)}+\intp L(\mathcal  A,\mu)  }\]
f.) For a neighborhood $0\in\mathcal  T\subset Z$ let $\mathcal  A$ be  an envelope
approximating $f$ such that
\[\intp L(\mathcal  A,\mu)\subset \int L(f,\mu)+\mathcal  T.\]
According to the previous point
\[\intp L(\mathcal  A\cup (\{0\}\times X),\mu)\subset \svar(f,\mu)-\mathcal  T+\mathcal  T+\mathcal  T.\]
Being $\mathcal  T$ arbitrary, our statement follows with $\mathcal  C=\mathcal  A\cup (\{0\}\times X)$.

g.) This follows from f. and
\[ \int L(f,\mu)\in \intp L(\mathcal  C,\mu),\]
ie. Lemma \ref{lem:CEVAL}.
\end{proof}
\end{lemma}
\begin{lemma}\label{lem:SVTEST}
\halfquote

Suppose that $f:X\rightarrow\mathcal  V$ is locally integrable.
We say that a function $g:X\rightarrow\mathcal  V$ is a restriction of $f$ if for each $x\in X$
we have $g(x)=f(x)$ or $g(x)=0$.

We claim:
\[\svar(L,f,\mu)=\overline{\biggl\{\int L(g,\mu)\,:
\,g \text{ is integrable, }g \text{ is a restriction of $f$ }\biggr\}}.\]
In particular, if the locally integrable function $h:X\rightarrow\mathcal  V$
is a restriction of $f$ then
\[\svar(L,h,\mu)\subset\svar(L,f,\mu).\]

\begin{proof}
0. The part $\subset$ of the equality is trivially true because $\svar(L,f,\mu)$ is,
by definition, generated by some special restrictions.

1. First, we prove the part $\supset$ of the equality under the restriction that $f$ is integrable.
Under this assumption we we can find a pointed envelope $C$ to $f$ such that
\[\intp L(\mathcal  C,\mu)\subset \svar(L,f,\mu)+\mathcal  T,\]
where $0\in\mathcal  T\subset \mathcal  Z$ is an arbitrarily small neighborhood.
But then $\mathcal  C$ also approximates $g$, hence
\[\int L(g,\mu)\in \intp L(\mathcal  C,\mu)\subset \svar(L,f,\mu)+\mathcal  T.\]
Because of the arbitrary choice of $\mathcal  T$ we conclude
\[\int L(g,\mu)\in\svar(L,f,\mu),\]
and this implies our statement.

At this point we know that under the additional assumption that $f,h$ are integrable
\[\svar(L,h,\mu)\subset\svar(L,f,\mu)\]
holds.

2. Here, we prove the part $\supset$ of the equality if  $f$ is locally integrable.
Lemma \ref {lem:SVCHAR}.g implies that it is enough to show that
\[ \sum_{j\in J}\int L(g|_{A_j},\mu)\in \svar(L,f,\mu)\]
whenever $A_j\in\mathfrak S$ ($j\in J$) are pairwise disjoint.
On the other hand, Lemma \ref {lem:SVCHAR}.g and the end of the previous point implies that
\[\sum_{j\in J}\int L(g|_{A_j},\mu)=\int L\biggl(\sum_{j\in J}g|_{A_j},\mu\biggr)\in
\svar\biggl(L,\sum_{j\in J}f|_{A_j},\mu\biggr).\]
Then, it is enough to show that
\[\svar\biggl(L,\sum_{j\in J}f|_{A_j},\mu\biggr)\subset \svar(L,f,\mu).\]
For that reason it is enough to show that
\[\int L\biggl(\sum_{j'\in J'}\biggl(\sum_{j\in J}f|_{A_j}\biggr)\biggr|_{A_{j'}'},\mu\biggr)
\in\svar(L,f,\mu)\tag{*}\]
whenever $A_{j'}'\in\mathfrak S$ ($j'\in J'$) are pairwise disjoint.
Let us decompose each set $A_j\cap A'_{j'}$ ($j\in J,\,j'\in J'$ ) in $\mathfrak S$, and take the union $\mathfrak B$
of these decompositions. Then
\[ \sum_{j'\in J'}\biggl(\sum_{j\in J}f|_{A_j}\biggr)\biggr|_{A_{j'}'}=\sum_{B\in\mathfrak B}f|_B.\]
This and the additivity of the integral (e. g. Lemma \ref{lem:ADDSUPER}) implies that
\[\int L\biggl(\sum_{j'\in J'}\biggl(\sum_{j\in J}f|_{A_j}\biggr)\biggr|_{A_{j'}'},\mu\biggr)=
\sum_{B\in\mathfrak B}^{\ovline}\int L(f|_B,\mu). \]
The right side is in the closure of its finite partial sums, which are in $\svar(L,f,\mu)$
by definition. That implies (*).
\end{proof}
\end{lemma}
\begin{lemma}\label{lem:SVSUM}
\halfquote

Suppose that $f,f_\lambda:X\rightarrow\mathcal  V$ ($\lambda\in\Lambda$) are locally integrable and
\[f(x)=\sum_{\lambda\in\Lambda}f_\lambda(x)\]
pointwise, in algebraical sense.
Then, we claim,

a.)
\[\svar(L,f,\mu)\subset \overline{\sum_{\lambda\in\Lambda}\svar(L,f_\lambda,\mu) }.\]

b.) If  for each point $x\in X$ at most one $\lambda\in\Lambda$ such that $f_\lambda(x)\neq 0$ then
\[\svar(L,f,\mu)= \overline{\sum_{\lambda\in\Lambda}\svar(L,f_\lambda,\mu) }. \]
\begin{proof}

a.) We prove this in two steps:

1. First we prove our statement under the assumption that $f$ and $f_\lambda$ are integrable.
Then, it  follows from Lemma \ref{lem:SVCHAR}:
We can slightly enlarge the right side (cf. Lemma \ref{lem:SVCHAR}.f) to get approximating envelopes,
which we can add up; then we know that this slightly enlarged ($+\mathcal  T$)
sum contains the left side (cf. Lemma \ref{lem:SVCHAR}.e).

2. The general statement follows from
\[\sum_{j\in J}\int L(f|_{A_j},\mu)\in\svar\biggl(L,\sum_{j\in J}f|_{A_j},\mu\biggr)
\]\[\subset
\overline{\sum_{\lambda\in\Lambda}\svar\biggl(L,\sum_{j\in J}f_\lambda|_{A_j},\mu\biggr) }
\subset \overline{\sum_{\lambda\in\Lambda}\svar(L,f_\lambda,\mu) }\]

b.) The direction $\supset$ follows from Lemma \ref{lem:SVTEST}.
\end{proof}
\end{lemma}
\begin{remak}
For the purposes of Lemma \ref{lem:SVSUM}.a the condition can be relaxed to
\[f(x)\in\sum_{\lambda\in\Lambda}\{0,f_\lambda(x)\}\]
(pointwise).
\end{remak}
\begin{lemma}\label{lem:SVSIG}
\halfquote

If $f:X\rightarrow\mathcal V$ is locally integrable,$A\in\boldsymbol\Sigma\mathfrak S$ and
$0\in\mathcal  T\subset\mathcal  Z$ is an arbitrary neighborhood
then, we claim, there exists a pointed envelope $\mathcal  C$  approximating $f|_A$, such that
\[\intp L(\mathcal  C,\mu)\subset \svar(L,f,\mu)+\mathcal  T.\]
\begin{proof}
Assume (cf. Lemma \ref{lem:SETUN}.a) that $A$ decomposes as
\[A=\bigcup_{\lambda\in\Lambda}^\updisjoint A_\lambda.\]
Let $0\in\mathcal U\subset\mathcal Z$ be a neighborhood such that $\mathcal U+\mathcal U\subset \mathcal T$.
Let us divide $\mathcal U$ by $\lambda$.
According to Lemma \ref{lem:SVCHAR}.f we can find pointed envelopes $\mathcal C_\lambda$
approximating $f|_{A_\lambda}$ such that
\[\intp L(\mathcal C_\lambda,\mu)\subset\svar(f|_{A_\lambda},\mu)+\mathcal U_\lambda.\]
Then the pointed envelope
\[\mathcal C=\sum_{\lambda\in\Lambda}\mathcal C_\lambda\]
approximates $f|_A$ and
\[\intp L(\mathcal C,\mu)=\overline{\sum_{\lambda\in\Lambda}\intp L(\mathcal C_\lambda,\mu)}\subset
\overline{\sum_{\lambda\in\Lambda}\svar(f|_{A_\lambda},\mu)+\mathcal U_\lambda}\subset
\overline{\svar(f|_A,\mu)+\mathcal U}\subset\]
\[\subset\svar(f|_A,\mu)+\mathcal U+\mathcal U\subset \svar(f|_A,\mu)+\mathcal T\subset \svar(f,\mu)+\mathcal T.\]

The second containment relation is a consequence of Lemma \ref{lem:SVSUM}.b, just like the last one.
\end{proof}
\end{lemma}
\begin{lemma}\label{lem:SVMON}
\halfquote

Suppose that $\mathsf V$ and $\mathsf Z$ are additive cones in $\mathcal  V$ and $\mathcal  Z$
respectively, $\mathsf Z$ is closed, and for each  $v\in\mathsf V$, $A\in\mathfrak S$,
\[L(v,\mu(A))\in\mathsf Z.\]
Suppose that $f:X\rightarrow\mathcal  V$ is integrable and $0\leq_{\mathsf V}f$ (pointwise).

Then,
\[\svar(L,f,\mu)\subset \left[0,\int L(f,\mu)\right]_{\mathsf Z},\]
endpoints contained.
\begin{proof}
This follows from the monotonicity of the integral, Lemma \ref{lem:MON}.
\end{proof}
\end{lemma}
\begin{lemma}\label{lem:SVTRANS}
\halfquote

Assume that $f_\lambda:X\rightarrow\mathcal  V$ ($\lambda\in\Lambda$) is a countable
family of integrable functions.
Then the following two statements are equivalent:
\begin{itemize}
\item[i.)] For each neighborhood $0\in\mathcal  T\subset\mathcal  Z$ there is a finite set
$\Xi\subset\Lambda$ such that
\[\sum_{\lambda\in\Lambda\setminus\Xi}\svar(L,f_\lambda,\mu)\subset \mathcal  T.\]
\item[ii.)] For each neighborhood $0\in\mathcal  T\subset\mathcal  Z$ there are envelopes
$\mathcal  C_\lambda$ approximating $f_\lambda$, and a finite set
$\Xi\subset\Lambda$ such that
\[\sum_{\lambda\in\Lambda\setminus\Xi}\intp L(\mathcal  C_\lambda,\mu)\subset \mathcal  T.\]
\end{itemize}
\begin{proof}
That immediately follows from Lemma \ref{lem:SVCHAR}.e and f.
\end{proof}
\end{lemma}

\paragraph{\textbf{B. Negligible sets}}
\begin{defin}
\halfquote

i.) A set $A\subset X$ is a weakly $L,\mu$-negligible set if  each function
$f:X\rightarrow\mathcal  V$ supported on $A$ (ie. $f=f|_A$) is integrable and
\[\int L(f,\mu)=0.\]
Ie. in other terms: for each function $f:X\rightarrow\mathcal  V$ vanishing outside of $A$
and for each neighborhood
$0\in\mathcal  T\subset Z$ there exists  an envelope $\mathcal  C$ of pointed sets in $\mathcal  V$ with respect to $\mathfrak S$
 such that
\[ f\subset \mathcal  C\quad\qand\quad \intp L(\mathcal  C,\mu)\subset \mathcal  T.\]

ii.) A set $A\subset X$ is a strongly $L,\mu$-negligible set if for each neighborhood
$0\in\mathcal  T\subset Z$ there exists  an envelope $\mathcal  C$ of pointed sets in $\mathcal  V$ with respect to $\mathfrak S$
 such that
\[\mathcal  V\times A\subset\mathcal  C\quad\qand\quad \intp L(\mathcal  C,\mu)\subset \mathcal  T.\]
\end{defin}
\begin{lemma}
\halfquote Then we claim:

a.) Strongly negligible sets are weakly negligible.

b.) Subsets and countable unions of weakly (strongly) $L,\mu$-negligible sets are
weakly (strongly)  $L,\mu$-negligible sets.
\begin{proof}
a.) That immediately follows from the definition.

b.) For subsets the statement is trivial.
In the weak case we can prove the countable union property as follows:

Suppose that $A_\lambda$ ($\lambda\in\Lambda$) are weakly negligible and $f$ is supported on
\[A=\bigcup_{\lambda\in\Lambda}A_\lambda.\]
Let $0\in\mathcal  T\subset\mathcal  Z$ be an arbitrary neighborhood.
Let us choose a neighborhood  $0\in\mathcal  U\subset\mathcal  Z$ such that
\[\mathcal  U+\mathcal  U\subset \mathcal  T\]
Let us divide $\mathcal  U$ by $\Lambda$. Let us choose pointed envelopes $\mathcal  C_\lambda$
\[f|_{A_\lambda}\subset\mathcal  C_\lambda
\quad\qand\quad \intp L(\mathcal  C_\lambda,\mu)\subset \mathcal  U_\lambda.\]
Take
\[\mathcal  C=\sum_{\lambda\in\Lambda}\mathcal  C_\lambda.\]
Then
\[f\subset \mathcal  C\]
but
\[\intp L(\mathcal  C,\mu)=\overline{ \sum_{\lambda\in\Lambda}\intp L(\mathcal  C_\lambda,\mu)}=
\overline{ \sum_{\lambda\in\Lambda}\mathcal  U_\lambda}=\overline{\mathcal  U}\subset \mathcal  T.\]
That proves the weak negligibility statement. The strong negligibility statement
can be proven similarly, except instead of
\[f|_{A_\lambda}\subset \mathcal  C_\lambda\]
rather
\[(\mathcal  V\times A_\lambda)\cup(\{0\}\times X)\subset \mathcal  C_\lambda\]
should be written, and similarly for $f$.
\end{proof}
\end{lemma}
\begin{lemma}
\halfquote

Suppose that $f:X\rightarrow \mathcal  V$ is a function and the (pointed) envelope $\mathcal  A$
approximates $f$.
Assume that $f'$ differs from $f$ on a weakly negligible set only.
Then we claim,
for each neighborhood $0\in\mathcal  T\subset\mathcal  Z$
there exists a (pointed) envelope $\mathcal  A'\supset\mathcal  A$ such that
$\mathcal  A'$ approximates $f'$ but
\[\intp L(\mathcal  A',\mu)\subset \intp L(\mathcal  A,\mu)+\mathcal  T.\]

\begin{proof}
Let us choose a neighborhood  $0\in\mathcal  U\subset\mathcal  Z$ such that
\[\mathcal  U+\mathcal  U\subset \mathcal  T.\]
Let $\mathcal  C$ be a pointed envelope such that
$f'-f$ is approximated by $\mathcal  C$ but
\[\intp L(\mathcal  C,\mu)\subset\mathcal  U.\]
Now take $\mathcal  A'=\mathcal  A+\mathcal  C$. Then
\[\intp L(\mathcal  A',\mu)=\overline{\intp L(\mathcal  A,\mu)+\intp L(\mathcal  C,\mu)}
\subset\intp L(\mathcal  A,\mu)+\mathcal U+\mathcal U\subset  \intp L(\mathcal  A,\mu)+\mathcal T.\]
\end{proof}
\end{lemma}
\begin{lemma}\label{lem:NEGCOSH}
\halfquote

Suppose that $C$ is negligible. Then, we claim for any envelope $\mathcal  A$
\[\intp L(\mathcal  A,\mu)=\biggl\{\sum_{j\in J} L(c_j,\mu(E_j)\,:\,J\text{ is finite,\, etc.,}
 \sum_{j\in J} c_j\chi_{E_j}(x)\in \mathcal  A^x
\text{ for }x\notin C ) \biggr\}.\]
\begin{proof}
The $\subset$ part of the statement is trivial. On the other hand if
\[s=\sum_{j\in J} c_j\chi_{E_j}\]
is such a step-function the we can find a function $f$ supported on $C$ such that
\[s+f\subset \mathcal A.\]
But then, because of the weak negligibility,
\[\int L(s,\mu)=\int L(s+f,\mu)\in\intp L(\mathcal A,\mu).\]
Taking closure yields our statement.

\end{proof}
\end{lemma}
\begin{point}
Now, all the semivariation and summation statements above can be rewritten
``up to negligible sets'', because we can just assume that all the function
involved are so that they are all $0$ on a large negligible set.
We refrain from spelling out the exact statements.
\end{point}

\newpage
\section{Integrable and extended measurable sets}\label{sec:INTMEAS}
\paragraph{\textbf{A. Integrable sets}}~\\

Previously, we have discussed integrable functions.
Using characteristic functions we may look to extend the notion of integrability to sets.

The following lemma helps to characterize the integrable sets.
\begin{lemma}\label{lem:INTCHAR}
\halfquote

Suppose that $c\in\mathcal  V$, $S\subset X$ and $c\chi_S$ is integrable.
Then we, claim,  that for each neighborhood $0\in\mathcal  T'\subset\mathcal  Z$ there exists a
set  $F\in\boldsymbol\Sigma_{00}\mathfrak S$, and sets
$A,B\in\boldsymbol\Sigma\mathfrak S$,
\[A\subset F,\qquad B\subset X\setminus F \]
such that
\[\tilde{\mathcal A}=c\chi_F-\var(c,A)+\var(c,B)\]
approximates $c\chi_E$ and
\[\intp L(\tilde{\mathcal A},\mu)\subset\int L(c\chi_F,\mu)+\mathcal  T'.\]

A property of $\tilde{\mathcal A}$ is a that is takes its values only from the set
$\{0,c\}\subset \mathcal V$.
\begin{proof}
Let $0\in\mathcal  T\subset\mathcal  Z$ be a  neighborhood such that
$\mathcal   T-\mathcal   T\subset\mathcal  T'$. Let
\[\mathcal  A=\sum_{j\in J}c_j\chi_{E_j}+\mathcal  C\]
be  an envelope approximating $c\chi_S$ but
\[\intp L(\mathcal A,\mu)\subset \int L(c\chi_F,\mu)+\mathcal  T.\]
Let $\mathfrak D$ be a a generator system for $\mathcal  A$.
Let $\mathfrak E=\{E_j\,:\,j\in J\}$.
Apply Lemma \ref{lem:FOREST} to $\mathfrak D\cup \mathfrak E$, such that the result
$\widetilde{\mathfrak D}$ is divided by $\mathfrak E$.
Then $\widetilde{\mathfrak D}$ is a generator system for $\mathcal  A$.

The approximating property implies for each element $x\in E_j$
there exist $D\in \widetilde{\mathfrak D}$ such that $x\in \widetilde{\mathfrak D}\subset E_j$
and $c\chi_S(x)\in\mathcal  A^D$.
In particular,
\[0\in\mathcal  A^D\quad\text{or}\quad c\in \mathcal  A^D.\]
Let $\mathfrak D'$ be the family of the maximal sets $D\in \widetilde{\mathfrak D}$
such that
\[0\in\mathcal  A^D\quad\text{or}\quad c\in \mathcal  A^D.\]
Then $\mathfrak D'$ decomposes $\{E_j\,:\,j\in J\}$. Hence, by Lemma \ref{lem:APPREFIN}
there is envelope
\[ \mathcal  A'=\sum_{j\in J'}c_j'\chi_{E_j'}+\mathcal  C' \]
such that the $E_j'$ are pairwise disjoint, they are from $\mathfrak D'$,
$\mathcal  A'\supset \mathcal  A$ approximates $c\chi_S$ and
\[\intp L(\mathcal  A',\mu)\subset \intp L(\mathcal  A,\mu)+\mathcal  T.\]
Let
\[F=\bigcup_{j\in J',\,c\in\mathcal  A'{}^{E'_j}} E_j', \]
\[A=\bigcup_{D\in \widetilde{\mathfrak D},\,D\subset F,\, 0\in\mathcal  A'{}^D}D,\qquad
B=\bigcup_{D\in \widetilde{\mathfrak D},\,D\not\subset F,\, c\in\mathcal  A'{}^D}D\]
Then
\[\tilde{\mathcal A}=\chi_F-\var(c,A)+\var(c,B)\]
will be  an envelope because it contains the simple step-function $\chi_F$, and we see that
\[\chi_F-\var(c,A)+\var(c,B)\subset \mathcal  A',\]
proving the containment in the statement of the lemma.
On the other hand, $\tilde{\mathcal A}$ approximates $c\chi_S$ by  construction.
 \end{proof}
\end{lemma}

\begin{defin}\label{def:MEASTRONG}
\halfquote

Then we define the measure
\[\tilde\mu:\mathfrak S\rightarrow{\Hom}_{\mathrm{strong}}(\mathcal  V,\mathcal  Z),\]
in the strong topology, by
\[[\tilde\mu(A)]v=L(v,\mu(A)). \]
It is easy to check that $\tilde\mu$ is indeed a measure.
\end{defin}
\begin{remark}
Notice that ${\Hom}_{\mathrm{strong}}(\mathcal  V,\mathcal  Z)$ remains complete and
Hausdorff according to our conventions.
In theory, for the purpose of integration it enough the keep the measure $\tilde\mu$ and we
can forget about $L$, $\mu$ and $\mathcal  W$.
\end{remark}
\begin{lemma}\label{lem:INTBICHAR}
\halfquote
We also consider the natural pairing $\mathtt m:\mathbb Z\times\mathcal  Z\rightarrow\mathcal  Z$.

Suppose that $S\subset X$. The we claim
\[\text{$c\chi_S$ is integrable for any $c\in\mathcal  V$ with respect to $L$ and $\mu$}\]
if and only if
\[\text{$\chi_S$ is integrable  with respect to $\mathtt m$ and $\tilde\mu$}.\]
\begin{proof}
This follows from the fact that special envelopes as in Lemma \ref {lem:INTCHAR} can be used
to decide integrability (cf. also Lemma \ref{lem:APPDIV}).
\end{proof}
\end{lemma}
Hence, in what follows we take the liberty of discussing integrability of sets
in the case when the measure is paired with $\mathcal  V=\mathbb Z$, $\mathcal  W=\mathcal  Z$.

\begin{defin}\label{def:INTGBLE}
\quarterquote

A set $S\subset X$ is called $\mu$-integrable if $\chi_S$ is integrable with respect to natural
pairing $\mathtt m$ to the integers and $\mu$.
We denote the family all integrable sets  by $\mathfrak S^\mu$.
\end{defin}
\begin{lemma}\label{lem:INTRING}
\quarterquote

The family $\mathfrak S^\mu$ of the $\mu$-integrable sets  forms a ring, ie. it is closed
for $\cap,\cup,\setminus$.
\begin{proof}
Suppose that $S_1$ and $S_2$ are integrable sets. Let $L=\mathtt m$.
Let $0\in\mathcal  T'\subset\mathcal Z$ be an arbitrary neighborhood.
Then chose a neighborhood $0\in\mathcal  T\subset\mathcal Z$  such that
\[\mathcal  T+\mathcal  T+\mathcal  T+\mathcal  T+\mathcal  T\subset\mathcal T'.\]
According to Lemma \ref{lem:INTCHAR} there exist envelopes
\[\mathcal  A_k=\sum_{j\in J_k}\chi_{E_{k,j}}+\mathcal  C_k\]
approximating $\chi_{S_k}$, $k=1,2$,
such that the sets $E_{k,j}$ are pairwise disjoint for a fixed $k$, and
\[\intp L(\mathcal  C_k,\mu)\subset\mathcal  T.\]

Applying Lemma \ref{lem:REFIN} we can find a
a countable family of sets $\mathfrak D$ decomposing the sets $E_{k,j}$.
Then according to Lemma \ref{lem:APPREFIN} we can find
there exist envelopes
\[\mathcal  A_k'=\sum_{E\in\mathfrak C_k}\chi_{E}+\mathcal  C_k'\]
approximating $\chi_{S_k}$, $k=1,2$,
such that $\mathfrak C_k\subset \mathfrak D$ are finite sets, and
\[\intp L(\mathcal  C_k',\mu)\subset\mathcal  T+\mathcal  T.\]
Then one can see that
\[\mathcal  A_\cap=\sum_{E\in\mathfrak C_1\cap \mathfrak C_2}\chi_{E}+\mathcal  C_1'+\mathcal  C_2'\]
approximates $\chi_{S_1\cap S_2}$ and
\[\intp L(\mathcal  C_1'+\mathcal  C_2',\mu)\subset \overline{\intp L(\mathcal  C_1',\mu)
+\intp L(\mathcal  C_2',\mu)}\subset
\mathcal  T+\mathcal  T+\mathcal  T+\mathcal  T+\mathcal  T\subset\mathcal T'.\]
According to the Cauchy criterium the integral exists. The other cases follow similarly
or just by linearity.
\end{proof}
\end{lemma}
\begin{lemma}\label{lem:INTCLOS}
\halfquote

Suppose that the set $S\subset X$ is integrable with respect to $\mu$.
Then, we claim, it is integrable with respect to
\[\tilde\mu:\mathfrak S\rightarrow{\Hom}_{\mathrm{strong}}(\mathcal  V,\mathcal  Z).\]
\begin{proof}
Let $\mathcal A$ be  an envelope over $\mathfrak S$ with values in $\mathbb Z$.
Then we define $c\mathcal A$ to be the envelope obtained from $\mathcal A$ such that
each fiber $\mathcal A^x$ is multiplies by $c$. Then one can prove that
It follows that for  an envelope $\mathcal  A$ over $\mathbb Z$
\[\intp L(c\mathcal  A,\mu)=\overline{L\biggl(c,\intp \mathtt m\,(\mathcal  A,\mu)\biggr)}.\]
From this and the continuity of $L$ in the second variable the statement follows.
\end{proof}
\end{lemma}

\paragraph{\textbf{B. Extended measureable sets}}
\begin{defin}\label{def:EMDEF}
\quarterquote

Let $\widetilde{\mathfrak S^\mu}$ denote the set of those $M\subset X$ such that
\[E\in \mathfrak S \quad \Rightarrow \quad E\cap M\in \mathfrak S^\mu,\]
ie. which intersects elements of $\mathfrak S$  in integrable sets.
The elements of $\widetilde{\mathfrak S^\mu}$ are called extended measurable sets (with respect to $\mu$).
\end{defin}
\begin{lemma}\label{lem:EMRING}
\quarterquote

a.) Integrable sets are extended measurable, ie. $\mathfrak S^\mu\subset \widetilde{\mathfrak S^\mu}$.

b.) The extended measurable sets, $\widetilde{\mathfrak S^\mu}$, form an algebra, ie. it is closed for
$\cap,\cup,\setminus$ and $X\in\widetilde{\mathfrak S^\mu}$.
\begin{proof}
That immediately follows from Lemma \ref{lem:INTRING}.
\end{proof}
\end{lemma}

\begin{lemma}\label{lem:EMMEAS}
\halfquote

Suppose that $f:X\rightarrow\mathcal  V$ is an integrable function.
Then, we claim that for each $M\in\widetilde{\mathfrak S^\mu}$
(or $M\in\widetilde{\mathfrak S^{\tilde \mu}}$) the function $f|_M$ is integrable. Moreover,
\[\hat \mu_{L,f}(M)=\int L(f|_M,\mu)\]
is a measure on $\widetilde{\mathfrak S^\mu}$.
\begin{proof}
Let $\mathcal  T+\mathcal  T+\mathcal  T\subset \mathcal  T'$. Let
\[\mathcal  A=\sum_{j\in J} c_J\chi_{E_j}+\mathcal  C\]
be  an envelope approximating $f$ such that for the pointed envelope $C$
\[\intp L(\mathcal  C,\mu)\subset\mathcal  T.\]
Then, according to our assumptions the sets $c_j\chi_{E_j\cap M}$ are integrable.
So, there is  an envelope
\[\mathcal  A'=s'+\mathcal  C'\]
approximating
\[\sum_{j\in J} c_j\chi_{E_j\cap M},\]
and such that
\[\intp L(\mathcal  C',\mu)\subset \mathcal T.\]
Then
\[s'+\mathcal  C'+\mathcal  C\]
will approximate $f|_M$ while
\[\intp L(\mathcal  C'+\mathcal  C,\mu)\subset \mathcal  T+\mathcal  T+\mathcal  T\subset\mathcal  T'.\]
Here $\mathcal  T'$ was arbitrary, hence the integrability follows from the Cauchy criterium.

The $\sigma$-additive property of the measure immediately follows from Lemma \ref{lem:ADDSUPER}.
\end{proof}
\end{lemma}
\begin{remak}
In the statement the special case with extended $\tilde\mu$-measurability in is not really necessary, because
this is indeed a special case of the statement with a plain $\mu$.
\end{remak}
\begin{lemma}\label{lem:EMCEAS}
\quarterquote

A set $M$ is extended measurable if and only if
\[E\in \mathfrak S^\mu \quad \Rightarrow \quad E\cap M\in \mathfrak S^\mu,\]
ie. it intersects integrable sets in integrable sets.
In particular,
 $\widetilde{\mathfrak S}^\mu$ is the measurable closure of ${\mathfrak S}^\mu$ in
 set-arithmetical sense.
\begin{proof} Sufficiency is obvious, while necessity follows from
the statement of Lemma \ref {lem:EMMEAS} about $f|_{M}$, such that
$f$ is characteristic function and $L=\mathtt m$.
\end{proof}
\end{lemma}

\paragraph{\textbf{C. Controlled measures}}
~\\

While a $\delta$-ring structure is not necessary for
a measure, it might certainly be useful.
\begin{defin}\label{def:CMEA}
A measure $\mu:\mathfrak S\rightarrow\mathcal  V$ is called
controlled if for any countable pairwise disjoint family of sets of $\mathfrak S$
contained in an other element of $\mathfrak S$,
\[\bigcup_{\lambda\in\Lambda}^\updisjoint A_\lambda\subset A,\]
implies that
\[\sum_{\lambda\in\Lambda}^{\ovline} \mu(A_\lambda)\]
exits.
\end{defin}

In ordinary circumstances this condition is quite unnoticeable:
\begin{lemma}\label{lem:CMEAEX}
The measure $\mu:\mathfrak S\rightarrow\mathcal  W$ is controlled, if
\begin{itemize}
\item[\textit{a.})] $\mu$ is $[0,+\infty)$-valued; or
\item[\textit{b.)}]  $\mathcal  V$ is a locally convex vector space and $\mu$ is of locally finite variation,
     ie. for each continuous seminorm $p$ on $\mathcal  V$ the set
     \[\mathrm{var}_p(\mu,A)=\hspace{10cm}\]\[
     \overline{\biggl\{\sum_{j\in J}p(\mu(A_j))\,:\, J\text{ is finite  }\,A_j\in\mathfrak S,\,
     A_j\subset A,\text{ the $A_j$ are pairwise disjoint}\biggr\}}\]
     is bounded; or
\item[\textit{c.)}] $\mathfrak S$ is a $\delta$-ring, ie.
\[\bigcup^{\updisjoint}_{\lambda\in\Lambda}A_\lambda\subset A\qquad\qquad
\text{implies}\qquad\qquad
\bigcup^{\updisjoint}_{\lambda\in\Lambda}A_\lambda\in\mathfrak S.\]
\end{itemize}
\begin{proof}
These statements are straightforward.
\end{proof}
\end{lemma}
\begin{lemma}\label{lem:CMEAVAR}
\halfquote
Suppose that
\[\bigcup^{\updisjoint}_{\lambda\in\Lambda}A_\lambda\subset A\]
in $\mathfrak S$. Also suppose that
$\mathcal  T$ is a neighborhood of $0$ in $\mathcal  Z$, and $c\in\mathcal V$ is arbitrary
Then we claim that there exists a finite set $\Xi$ such that
\[\svar\biggl(L,c,\bigcup_{\lambda\in\Lambda\setminus\Xi}A_\lambda,\mu\biggr)\subset
\mathcal  T.\]
\begin{proof}
Apply Lemma \ref{lem:XMEAVAR} with the choice of $T=0$, $\tilde \mu(A)=L(c,\mu(A))$.
\end{proof}
\end{lemma}
\begin{lemma}\label{lem:CADDSUPER}
\halfquote

Assume that $\mu$ is controlled.
Suppose that $\Lambda$ is countable,  $f,f_\lambda:X\rightarrow\mathcal  V$ ($\lambda\in\Lambda$)
are integrable functions and for each $x\in X$
\[\sum_{\lambda\in\Lambda}f_\lambda(x)=f(x) \text{ or }=0,\]
such that at each point $x\in X$ at most one of the $f_\lambda(x)$'s is nonzero.

Then we claim, for each neighborhood $0\in\mathcal  T\subset\mathcal  Z$ there exist an finite set
$\Xi$, and a family of pointed envelopes $\mathcal  C_\lambda$ ($\lambda\in\Lambda\setminus\Xi$) such that
$\mathcal  C_\lambda$ approximates $f_\lambda$ and
\[\intp L\biggl(\sum_{\lambda\in\Lambda\setminus\Xi}\mathcal  C_\lambda,\mu\biggr)\subset\mathcal  T.\]
In particular,
\[\sum_{\lambda\in\Lambda}f_\lambda. \]
is integrable.
\begin{proof}
We can repeat the proof of Lemma \ref{lem:ADDSUPER}, with minor changes.
The difference is that that $\mathfrak D'$
will not necessarily cover $\bigcup\widetilde{\mathfrak D}$ where the sum of the functions is zero.
So, the sets $\mathfrak E_i$ will not cover $E_i$ and we have to apply
to Lemma \ref{lem:CMEAVAR} instead of Lemma \ref{lem:MEAVAR}.
\end{proof}
\end{lemma}
\begin{lemma}
\quarterquote

Suppose that $\mu$ is controlled. Then, we claim, $\mathfrak S^\mu$ is a $\delta$-ring
and $\widetilde{\mathfrak S^\mu}$ is a $\sigma$-algebra.
\begin{proof}
It is enough to prove the statement for integrable case.
That follows from Lemma \ref{lem:CADDSUPER} applied to characteristic functions of sets.
\end{proof}
\end{lemma}

\paragraph{\textbf{D. Measure extensions by integrable sets}}
\begin{conven}\label{def:EXTINT}
If $\mu:\mathfrak S\rightarrow \mathcal  W$ is a measure then let
\[\hat\mu:\mathfrak S^\mu\rightarrow\mathcal  W\]
be the extension to integrable set.
\end{conven}
\begin{lemma}\label{lem:MECHAR}
\quarterquote

Then, we claim, $\hat\mu$ is the unique extension $\nu$ of $\mu$ to $\mathfrak S^\mu$ such that
for any $A\in\mathfrak S$
\[\svar(A,\nu)=\svar(A,\mu).\]

Similar statement holds for any sub- interval system of $\mathfrak S^\mu$.
\begin{proof}
From the additive property of the semivariation it follows that
\[\svar(A,\nu)=\svar(A,\mu)\]
holds for any $A\in\mathfrak S$ if and only if it holds for any $A\in\boldsymbol\Sigma\mathfrak S$.

Then the statement follows from Lemma \ref{lem:INTCHAR}.
(One can shorten the proof by using Lemma \ref{lem:SVTEST}.)
\end{proof}
\end{lemma}

\begin{lemma}\label{lem:MECOS}
\halfquote

Suppose that $\mathcal  A$ is a (pointed) envelope over $\mathfrak S$.
Then, we claim, $\mathcal  A$ is a (pointed) envelope over $\mathfrak S^\mu$ and
\[\intp L(\mathcal  A,\mu)=\intp L(\mathcal  A,\hat\mu).\]
\begin{proof}
The (pointed) envelopeness over $\mathfrak S^\mu$ is trivial, also the $\subset$ part of the
equality. On the other hand, if $\hat s$ is a step-function over $\mathfrak S^\mu$ then it is
$L,\mu$-integrable and by Lemma \ref {lem:INTCONT} the corresponding sum is contained in
\[\intp L(\mathcal  A,\mu).\]
Taking closure, this proves $\supset$, hence our statement.
\end{proof}
\end{lemma}
\begin{lemma}\label{lem:MEINT}
\halfquote

Suppose that $f$ is $L,\mu$-integrable. Then $f$ is $L,\hat\mu$-integrable, the integrals are equal,
and
\[\svar(L,f,\mu)=\svar(L,f,\hat\mu).\]
\begin{proof}
The first two statement follows from Lemma \ref{lem:MECOS}, while the last one follows from
Lemma \ref{lem:EMMEAS} and Lemma \ref{lem:SVTEST}.
\end{proof}
\end{lemma}
\begin{lemma}\label{lem:MECONTINT}
\halfquote

If $\mathcal  A$ is  an envelope over $\mathfrak S^\mu$ and there is countable set $R\subset \mathcal  V$
such that $\mathcal  A\subset R\times X$ (ie. $\mathcal  A$ takes only countable many values) then
for each neighborhood $0\in\mathcal  T\subset \mathcal  Z$ there is  an envelope $\tilde{\mathcal  A}$ over $\mathfrak S$
such that
\[\mathcal  A\subset\tilde{ \mathcal  A}\]
but
\[\intp L(\tilde{\mathcal  A},\mu)
\subset
\intp L(\mathcal  A,\hat\mu)+\mathcal  T.\]
\begin{proof}
Assume that
\[\mathcal  A=\hat s+\mathcal  C\]
such that $\hat s$ is a step-function over $\mathfrak S^\mu$ contained in $\mathcal  A$.
One can notice that $\mathcal  C$ has a similar properties to $\mathcal  A$ except with
an extended set $R'\supset R$.

Let $0\in\mathcal  T'\subset\mathcal Z$ be a neighborhood such that
\[\mathcal  T'+\mathcal  T'+\mathcal  T'-\mathcal T'\subset \mathcal  T.\]

Now, $\hat s$ is $L,\mu$-integrable, hence there is  an envelope
\[\mathcal A_0=s+\mathcal C_0  \] containing the step-function $s$ over $\mathfrak S$, and
 approximating $\hat s$, such that
\[\intp L(\mathcal  C_0,\mu)\subset \mathcal  T'.\]

Let $\mathfrak D\subset\mathfrak S^\mu$ be a generator system for $\mathcal  C$.
Let us divide $\mathcal  T'$ by $R'\times \mathfrak D$.
Then Lemma \ref{lem:INTCHAR} implies that
for each $r\in R'$, $D\in\mathfrak D,\,r\in\mathcal  A^D$
there is a set $F_{r,D}\in\boldsymbol\Sigma_{00}\mathfrak S$
and a set $A_{r,D}\in \boldsymbol\Sigma\mathfrak S$ such that
$D$ differs from $F_{r,D}$ only on $A_{r,D}$ but
\[\svar(L,r,A_{r,D},\mu)\subset \mathcal  T'_{r,D}.\]
Let
\[\mathcal  C'=\sum_{r\in R',\,D\in\mathfrak D,\,r\in\mathcal  A^D} \var(r,A_{r,D}).\]
Then
\[\intp L(\mathcal  C',\mu)\subset
\overline{\sum_{r\in R',\,D\in\mathfrak D,\,r\in\mathcal  A^D}\mathcal  T'_{r,D}}
\subset \overline{\mathcal  T'}.\]

Let us define
\[\tilde{\mathcal  C}=\bigcup_{D\in\mathfrak D,\,r\in\mathcal  A^D}\{r\}\times(F_{r,D} \cup A_{r,D})
\cup \{0\}\times X.\]
The set $\mathcal C'$ is a pointed  envelope over $\mathfrak S$ because if we decompose
the sets $F_{r,D}$ and $A_{r,D}$ in $\mathfrak S$
and we take the union of these decompositions then it will yield a generator system for $\mathcal C'$.

Let us define
\[\mathcal  A=\mathcal  A_0+\tilde{\mathcal  C}.\]
Then $\mathcal  A$ is a clearly  an envelope over $\mathfrak S$ because it is a sum of envelope
and a pointed envelope.
On the other hand,
\[\hat s\subset\mathcal  A_0\quad\text{and}\quad \mathcal  C\subset \tilde{\mathcal  C}\]
implies that
\[\mathcal  A=\hat s+\mathcal  C\subset \mathcal  A_0 +\tilde{\mathcal  C}=\tilde{\mathcal  A}.\]

On the other hand, $\tilde{\mathcal C}\subset \mathcal C+\mathcal C'$ implies that
\[\tilde{\mathcal  A}=\mathcal  A_0 +\tilde{\mathcal  C}\subset\mathcal  A_0+\mathcal  C+\mathcal  C'.\]
Then, continuing,
\[\subset \hat s -\mathcal C_0+\mathcal C_0+\mathcal  C+\mathcal  C'
\subset \mathcal A-\mathcal C_0+\mathcal C_0+\mathcal C'\]
implies that
\[\intp L(\tilde{\mathcal  A},\mu)=\intp L(\tilde{\mathcal  A},\hat\mu)
\subset\intp L({\mathcal  A}-\mathcal C_0+\mathcal C_0+\mathcal C',\hat\mu)\subset \]\[\subset
\overline{\intp L(\mathcal A,\hat\mu) +\mathcal T'-\mathcal T'+\overline{\mathcal T'}}\subset
\intp L(\mathcal A,\hat\mu) +\mathcal T'-\mathcal T'+\mathcal T'+\mathcal T'
\subset \intp L(\mathcal A,\hat\mu) +\mathcal T.\]
\end{proof}
\end{lemma}
\begin{remak}
This statement will be used to prove the conversion of Lemma \ref{lem:MEINT} in cases when
countable-valued envelopes can be used to ``approximate'' general envelopes.
\end{remak}

\paragraph{\textbf{E. Measure extensions by constructible sets}}
\begin{lemma}\label{lem:CNINTAH}
\quarterquote

Consider the constructible measure extension
\[\mu_{\mathrm c}:\boldsymbol\Sigma_{\mathrm c}\mathfrak S\rightarrow\mathcal W.\]
Then we claim that $\mu_{\mathrm c}$ is a restriction of
\[\hat\mu:\mathfrak S^\mu\rightarrow\mathcal  W.\]
\begin{proof}
That follows from Lemma \ref{lem:INTRING} and the unicity statement for $\mu_{\mathrm c}$,
or from the corresponding version of Lemma \ref{lem:MECHAR}.
\end{proof}
\end{lemma}

\begin{lemma}\label{lem:CNCO}
\halfquote
Consider the constructible measure extension
\[\mu_{\mathrm c}:\boldsymbol\Sigma_{\mathrm c}\mathfrak S\rightarrow\mathcal W.\]
Then we claim:

a.) A set $\mathcal C$ is pointed envelope over $\boldsymbol\Sigma_{\mathrm c}\mathfrak S$ if and only
if it is a pointed envelope over $\mathfrak S$. In this case
\[\intp L(\mathcal C,\mu_{\mathrm c})=\intp L(\mathcal C,\mu).\]

b.) Any envelope $\mathcal A$ over $\mathfrak S$ is  an envelope over
$\boldsymbol\Sigma_{\mathrm c}\mathfrak S$, and
\[\intp L(\mathcal A,\mu_{\mathrm c})=\intp L(\mathcal A,\mu).\]

c.) For any envelope $\mathcal A$ over $\boldsymbol\Sigma_{\mathrm c}\mathfrak S$
and any neighborhood $0\in\mathcal T\subset\mathcal Z$ there is  an envelope
$\mathcal A'$ over $\mathfrak S$ such that
\[\mathcal A\subset\mathcal A'\]
and
\[\intp L(\mathcal A',\mu)\subset\intp L(\mathcal A,\mu_{\mathrm c})+\mathcal T.\]
\begin{proof}
a.) The first statement follows from Lemma \ref{lem:SETUN}.d and the other part follows from
Lemma \ref{lem:INTCONT}.

b.) The first statement is trivial, the second part follows from Lemma \ref{lem:INTCONT}.

c.) We can consider a usual decomposition
\[\mathcal A=s+\mathcal C\]
over $\boldsymbol\Sigma_{\mathrm c}\mathfrak S$. From part a. we see that it is enough to
prove the statement in the case $\mathcal A=s$. That follows from the integrability of $s$
over $\mathfrak S$.
\end{proof}
\end{lemma}
\begin{lemma}\label{lem:CNINT}
\halfquote

Then we claim a function $f:X\rightarrow\mathcal V$ is (locally)
$L,\mu$-integrable if and only if it is (locally) $L,\mu_{\mathrm c}$-integrable.
If $f$ is integrable then
\[\int L(f,\mu)=\int L(f,\mu_{\mathrm c}).\]
If $f$ is locally integrable then
\[\svar(L,f,\mu)=\svar(L,f,\mu_{\mathrm c}).\]
\begin{proof}
That follows from Lemma \ref{lem:CNCO}
\end{proof}
\end{lemma}

\newpage
\section{Extended integration}\label{sec:EXTLM}

\paragraph{\textbf{A. Local properties of the ordinary integral}}
\begin{defin}\label{def:EOSUP}
Let $\mathfrak S\subset\mathfrak P(X)$ be an interval system and $\mathcal V$ be a commutative topological
group. Let $\mathcal A$ be  an envelope of sets in $\mathcal V$ over $\mathfrak S$.

We say that $\mathcal A$ is supported on $A\in\boldsymbol\Sigma\mathfrak S$
if there is a generator system $\mathfrak D$ of $\mathcal A$ such that $\bigcup\mathfrak D=A$.
\end{defin}
\begin{lemma}\label{lem:EOSUP}
Let $\mathfrak S\subset\mathfrak P(\mathfrak S)$ be an interval system, and $\mathcal  V$
be a commutative topological groups. Assume that $\mathcal  A$ is  an envelope of sets in $\mathcal V$ over $\mathfrak S$.
Then we claim:

a.) If $\mathcal A$ is supported on $A$ then it is also supported on any set
$A\in\boldsymbol\Sigma\mathfrak S$ such that $B\supset A$.

b.) If  $X\in\boldsymbol\Sigma\mathfrak S$ then $\mathcal A$ is supported on $X$.
\begin{proof}
a.) This follows from Lemma \ref{lem:APPNEIGH}, we can add up any set to $\bigcup\mathfrak D$
but we have to decompose everything by $\mathfrak E$.
b.) This is a special case of a.
\end{proof}
\end{lemma}
\begin{remark}
In particular, in the case   $X\in\boldsymbol\Sigma\mathfrak S$
we can simplify things by assuming $\bigcup\mathfrak D=X$.
\end{remark}
\begin{lemma}\label{lem:EOINT}
\halfquote
Assume that $f:X\rightarrow\mathcal V$ is integrable.

Then we claim that
there is set $A\in\boldsymbol\Sigma\mathfrak S$ such that $f|_A=f$, ie. $f$ vanishes outside of $A$.
\begin{proof}
That immediately follows from the existence of an approximating envelope.
\end{proof}
\end{lemma}

On the other hand, if $X\notin \boldsymbol\Sigma\mathfrak S$,
then this is a serious restriction excluding some functions which are intuitively integrable.
The objective of this section is to get around this technical difficulty.

\begin{defin}\label{def:EORES}
 a.) If $f:X\rightarrow Y$ is a function and $S\subset X$
then let
\[f\|_S:S\rightarrow Y\]
be the set-theoretic restriction of $f$ to $S$. (Ie.
we not just set the set the function value to $0$ outside of $S$ but
we restrict the actual domain).

b.) If $\mathfrak S\subset\mathfrak P(\mathfrak S)$ is a set system and
and $A\subset X$ then we define
\[\mathfrak S\|_A=\{S\in\mathfrak S\,:\,S\subset A\}.\]
\end{defin}
\begin{lemma}\label{def:EOCON}
Suppose that $\mathfrak S\subset\mathfrak P(\mathfrak S)$ is an interval system.

We claim that if $A\in\boldsymbol\Sigma\mathfrak S$ or
$X\setminus A\in \boldsymbol\Sigma_{\mathrm c}\mathfrak S$
then $\mathfrak S\|_{A}$ is an interval system.
\begin{proof}
That immediately follows from the properties of interval systems.
\end{proof}
\end{lemma}
\begin{lemma}\label{lem:EORES}
\halfquote

Assume that $A\in\boldsymbol\Sigma\mathfrak S$ or
$X\setminus A\in \boldsymbol\Sigma_{\mathrm c}\mathfrak S$.
Then we claim:

a.) If $\hat{\mathcal A}$ is  an envelope over $\mathfrak S\|_A$ , then
\[\mathcal A=\hat{\mathcal A}\cup((X\setminus A)\times\{ 0\})\]
will be  an envelope of sets in $\mathcal V$ over $\mathfrak S$ such that
\[\intp L(\hat{\mathcal A},\mu\|_{\mathfrak S\|_A})=\intp L(\mathcal A,\mu).\]

In particular, if $f:X\rightarrow\mathcal V$ vanishes outside of $A$  and  $\hat{\mathcal A}$
approximates $f\|_A$ then $\mathcal A$ approximates $f$.

b.) If $\mathcal A$ is  an envelope of sets in $\mathcal V$ over $\mathfrak S$
containing $\{0\}\times(X\setminus A)$ then
for every neighborhood $0\in\mathcal T\subset\mathcal Z$
there is  an envelope $\hat{\mathcal A}$  over $\mathfrak S\|_A$
such that
\[\mathcal A\cap (\mathcal V\times A)\subset\hat{\mathcal A}\]
but
\[\intp L(\hat{\mathcal A},\mu\|_{\mathfrak S\|_A})\subset\intp L(\mathcal A,\mu)+\mathcal T.\]

In particular, if $f:X\rightarrow\mathcal V$ vanishes outside of $A$  and  $\mathcal A$
approximates $f$ then $\hat{\mathcal A}$ approximates $f\|_A$.
\begin{proof}
a.) The envelope properties are trivial to check. Now, assume that
\[\hat{\mathcal A}=\hat s+\hat{\mathcal C}\]
is sum of a step-function and a pointed envelope over ${\mathfrak S\|_A}$.
We can apply the construction to $\hat{\mathcal C}$, too.

Then
\[\intp L(\hat s,\mu\|_{\mathfrak S\|_A}))=\intp L(s,\mu)\]
because the decomposition of $\hat s$ to characteristic functions can be used not only over
$\mu\|_{\mathfrak S\|_A}$ but also over $\mu$.
Furthermore,
\[\intp L(\hat{\mathcal C},\mu\|_{\mathfrak S\|_A})=\intp L(\mathcal C,\mu)\]
because the simple step-functions as in Lemma \ref{lem:CEVAL}
(omitting the $0$ coefficients) will be the same for $\hat{\mathcal C}$ and $\mathcal C$.
That implies
\[\intp L(\hat{\mathcal A},\mu\|_{\mathfrak S\|_A})=
\intp L(\hat s,\mu\|_{\mathfrak S\|_A}))+\intp L(\hat{\mathcal C},\mu\|_{\mathfrak S\|_A})=\]\[
=\intp L(s,\mu)+\intp L(\mathcal C,\mu)=\intp L(\mathcal A,\mu).\]

b.)
Suppose that
\[\mathcal A=\sum_{j\in J}c_j\chi_{E_j}+\mathcal C\]
is a decomposition of $\mathcal A$ to a step-function and a pointed envelope.
Let $\mathfrak E=\{E_j\,:\,j\in J\}$.

1. First we prove the statement under the assumption
$A\in\boldsymbol\Sigma\mathfrak S$.

Suppose that $A$ exactly decomposes into $\mathfrak A$ over $\mathfrak S$.
Let $\mathfrak D$ be a generator system for $\mathcal A$.
Then  let $\widetilde{\mathfrak D}$ be a forest refinement (cf. Lemma \ref{lem:FOREST}) of
$\mathfrak D_0\cup\mathfrak E\cup\mathfrak A$ divided by $\mathfrak E$.

Then $\{0\}\times (X\setminus A)\subset \mathcal A$ implies that
  for each $x\in\bigcup\mathfrak D$ there is a set $D\in\widetilde {\mathfrak D}$ such that
$x\in D$ and
\[0\in\mathcal A^D \qquad \text{or}\qquad D\subset A\]
The set system $\widetilde{\mathfrak D}$ is divided by $\mathfrak E$, so $x\in E_j$ implies $D\subset E_j$.
In particular the family $\mathfrak D'$ of the maximal sets of $\widetilde{\mathfrak D} $
such that
\[0\in\mathcal A^D \qquad \text{or}\qquad D\subset A\]
will decompose  $\mathfrak E$. Then applying Lemma \ref{lem:APPREFIN}
we can replace $\mathcal A$ by  an envelope $\mathcal A'\supset \mathcal A$ such that
it contains a simple step-function
\[s'=\sum_{j\in J'}c_j'\chi_{E_j'},\]
the sets $E'_j$ are pairwise disjoint, the sets $E'_j$ are from $\mathfrak D'$,
but
\[\intp L(\mathcal A',\mu)\subset \intp L(\mathcal A,\mu)+\mathcal T.\]
Omitting  the coefficients $c'_j$ which can be replaced by $0$ we can assume that
$E'_j\subset A$ for all $j\in J'$. Then from Lemma \ref{lem:SETUN}.c one can see that
\[\hat{\mathcal A}=\mathcal A'\cap(\mathcal V\times A)\]
is  an envelope. According to point a. we see that
\[\intp L(\hat{\mathcal A},\mu\|_{\mathfrak S\|_A})=
\intp L(\hat{\mathcal A}\cup ((X\setminus A)\times\{0\}),\mu)\subset
\intp L(\mathcal A',\mu)\subset\intp L(\mathcal A,\mu)+\mathcal T.\]

2. We prove that statement under the assumption $X\setminus A\in \boldsymbol\Sigma_{\mathrm c}\mathfrak S$.
Assume that $X\setminus A$ is a $(\cup,\cap,\setminus)$-expression of $B_1,\ldots,\,B_n$.
Then can exactly decompose $\mathfrak E$ by $\mathfrak D$ a such that it
is divided $\{B_1,\ldots,\,B_n\}$.

 Then applying Lemma \ref{lem:APPREFIN} and omitting the coefficients which can be set to $0$ we obtain
 a similar $\mathcal A'$ with a similar simple step-function, such that $E_j'\cap A=0$ for all $j\in J'$.

 The proof ends in the same way.
\end{proof}
\end{lemma}

\begin{lemma}\label{lem:EOBER}
\halfquote

Assume that $A\in\boldsymbol\Sigma\mathfrak S$ or
$X\setminus A\in \boldsymbol\Sigma_{\mathrm c}\mathfrak S$.
Suppose that $f:X\rightarrow\mathcal V$  vanishes outside of $A$.
Then we claim:

a.) The function $f$ is (locally) $L,\mu$-integrable if and only if $f\|_A$ is (locally)
$L,\mu\|_{\mathfrak S\|_A}$-integrable. The integrals, if exist, are equal.

b.) If $f$ is locally integrable then
\[\svar(L,f,\mu)=\svar(L,f\|_A,\mu\|_{\mathfrak S\|_A})=\]
\[\overline{\biggl\{\sum_{j\in J}\int L(f|_{A_j},\mu)\,:
\,J\text{ is finite},\,A_j\in\mathfrak S\|_A,\,
\text{ the $A_j$ are pairwise disjoint}\biggr\}}.\]
\begin{proof}
That immediately follows from Lemma \ref{lem:SVTEST} and \ref{lem:EORES}.
\end{proof}
\end{lemma}
\begin{lemma}\label{lem:EOCUT}
\halfquote

Suppose that $\mathcal A$ is  an envelope of sets in $\mathcal V$ over $\mathfrak S$.
 Assume $A\in\boldsymbol\Sigma_{\mathrm c}\mathfrak S$.
Then, we claim, for each neighborhood $0\in\mathcal T\subset\mathcal Z$ there are envelopes
$\mathcal A_1$ and $\mathcal A_2$ such that
\[(\mathcal A\cap (\mathcal V\times A))\cup( \{0\}\times (X\setminus A))\subset\mathcal A_1,\]
\[(\mathcal A\cap (\mathcal V\times (X\setminus A)))\cup( \{0\}\times A)\subset\mathcal A_2,\]
but
\[\overline{\intp L(\mathcal A_1,\mu)+\intp L(\mathcal A_2,\mu)
}\subset\intp L(\mathcal A,\mu)+\mathcal T.\]

If $\mathcal A$ approximates $f:X\rightarrow\mathcal V$ then $\mathcal A_1$ approximates $f|_A$
and $\mathcal A_2$ approximates $f|_{X\setminus A}$.
\begin{proof}
This follows from Lemma \ref{lem:EORES}.
\end{proof}
\end{lemma}
In what follows we use these locality properties and their consequences without much reference
to the original statements.
~\\

\paragraph{\textbf{B. The extended integral}}
\begin{lemma}\label{lem:EQUIGEN}
\halfquote

Suppose that $f:X\rightarrow \mathcal  V$ is a function.
We claim that the following conditions for an element  $a\in\mathcal  Z$ are equivalent:

$\mathtt{(\overline E^s)}$ For each neighborhood $0\in\mathcal  T\subset\mathcal  Z$ and
set $A_0\in \boldsymbol\Sigma_{00}\mathfrak S$
there is a set $A'\in\boldsymbol\Sigma_{00}\mathfrak S$ such that
\[A_0\cap A'=\emptyset\]
and
for all $A_0\cup A'\subset B\in\boldsymbol\Sigma\mathfrak S$ there is envelope $\mathcal  B$ such that
\[f|_B\subset \mathcal  B\qquad\text{and}\qquad \intp L(\mathcal  B,\mu)\subset a+\mathcal  T.\]

$\mathtt{(\overline E^m)}$
 For each neighborhood $0\in\mathcal  T\subset\mathcal  Z$
there is a set $A\in\boldsymbol\Sigma_{00}\mathfrak S$ such that
for all $A\subset B\in\boldsymbol\Sigma\mathfrak S$ there is envelope $\mathcal  B$ such that
\[f|_B\subset \mathcal  B\qquad\text{and}\qquad \intp L(\mathcal  B,\mu)\subset a+\mathcal  T.\]

$\mathtt{(\overline E_\Sigma)}$
 For each neighborhood $0\in\mathcal  T\subset\mathcal  Z$
there is a set $A\in\boldsymbol\Sigma\mathfrak S$ such that
for all $A\subset B\in\boldsymbol\Sigma\mathfrak S$ there is envelope $\mathcal  B$ such that
\[f|_B\subset \mathcal  B\qquad\text{and}\qquad \intp L(\mathcal  B,\mu)\subset a+\mathcal  T.\]

$\mathtt{(\overline E_0)}$
 For each neighborhood $0\in\mathcal  T\subset\mathcal  Z$
there is a set $A\in\boldsymbol\Sigma_{0}\mathfrak S$ such that
for all $A\subset B\in\boldsymbol\Sigma_{0}\mathfrak S$ there is envelope $\mathcal  B$ such that
\[f|_B\subset \mathcal  B\qquad\text{and}\qquad \intp L(\mathcal  B,\mu)\subset a+\mathcal  T.\]

$\mathtt{(\overline E_{00})}$
  For each neighborhood $0\in\mathcal  T\subset\mathcal  Z$ and
set $A_0\in \boldsymbol\Sigma_{00}\mathfrak S$
there is a set $A\in\boldsymbol\Sigma_{00}\mathfrak S$ such that
for all $A'\in\boldsymbol\Sigma_{00}\mathfrak S$ such that $A\cap A'=\emptyset$
there is envelope $\mathcal  B$ such that
\[f|_{(A\cup A')}\subset \mathcal  B\qquad\text{and}\qquad \intp L(\mathcal  B,\mu)\subset a+\mathcal  T.\]

$\mathtt{(\overline L^m)}$
 The function $f$ is locally integrable and  for each neighborhood $0\in\mathcal  T\subset\mathcal  Z$
there is a set $A\in\boldsymbol\Sigma_{00}\mathfrak S$ such that
for all $A\subset B\in\boldsymbol\Sigma_{\mathrm c}\mathfrak S$
\[\int L(f|_B,\mu)\subset a+\mathcal  T.\]

$\mathtt{(\overline L^w)}$
 The function $f$ is locally integrable and  for each neighborhood $0\in\mathcal  T\subset\mathcal  Z$
there is a set $A\in\boldsymbol\Sigma_{\mathrm c}\mathfrak S$ such that
for all $A'\in\boldsymbol\Sigma_{00}\mathfrak S$ such that $A\cap A'=\emptyset$
\[\int L(f|_{A},\mu)+\int L(f|_{A'},\mu)\subset a+\mathcal  T.\]

$\mathtt{(\overline S_1)}$
 The function $f$ is locally integrable and for each neighborhood $0\in\mathcal  T\subset\mathcal  Z$
there is a set $A\in\boldsymbol\Sigma_{00}\mathfrak S$ such that
\[\int L(f|_A,\mu)+\svar(L,f-f|_A,\mu)\subset a-\mathcal  T.\]

$\mathtt{(\overline S_2)}$
The function $f$ is locally integrable and for each neighborhood $0\in\mathcal  T\subset\mathcal  Z$
there is a simple step-function $s$ such that
\[\int L(s,\mu)+\svar(L,f-s,\mu)\subset a-\mathcal  T.\]

$\mathtt{(\overline S_3)}$
The function $f$ is locally integrable and for each neighborhood $0\in\mathcal  T\subset\mathcal  Z$ and
there is an integrable function $g$ such that
\[\int L(g,\mu)+\svar(L,f-g,\mu)\subset a-\mathcal  T.\]

Moreover, such an $a\in\mathcal Z$, if exists, is unique.
Integrable functions functions satisfy the conditions
above and the corresponding value is the integral.

\begin{proof}
The equivalence statement:

$\mathtt{(\overline E^s)}$$\Rightarrow$$\mathtt{(\overline E^m)}$ is obvious, while
$\mathtt{(\overline E^m)}$$\Rightarrow$$\mathtt{(\overline E^s)}$ follows from
 Lemma \ref{lem:EORES} and \ref{lem:EOCUT}.

$\mathtt{(\overline E^m)}$$\Rightarrow$$\mathtt{(\overline E_\Sigma)}$ is obvious, while
$\mathtt{(\overline E_\Sigma)}$$\Rightarrow$$\mathtt{(\overline E^m)}$ follows from
Lemma \ref{lem:APPREFIN}.

$\mathtt{(\overline E^m)}$$\Rightarrow$$\mathtt{(\overline E_0)}$ is obvious, just like
 $\mathtt{(\overline E^s)}$$\Rightarrow$$\mathtt{(\overline E_{00})}$.

$\mathtt{(\overline E_0)}$ implies local integrability because $E\subset B$
can be assumed.
Similarly, $\mathtt{(\overline E_{00})}$
implies the local integrability of $f$ because if $E\in\mathfrak S$  then $A_0=E$
can be chosen. In all of those cases we can apply Lemma \ref{lem:EOCUT} to
separate $E$ and apply the Cauchy criterium to the envelope over $E$.
Having the local integrability established,
$\mathtt{(\overline E_0)}$ immediately implies $\mathtt{(\overline L^w)}$,
$\mathtt{(\overline E_{00})}$ immediately implies $\mathtt{(\overline L^w)}$.

$\mathtt{(\overline L^m)}$$\Rightarrow$$\mathtt{(\overline L^w)}$ is obvious, while
$\mathtt{(\overline L^w)}$$\Rightarrow$$\mathtt{(\overline L^m)}$ follows
from Lemma \ref{lem:SVCHAR}.g.

$\mathtt{(\overline L^w)}$$\Rightarrow$$\mathtt{(\overline S_1)}$ is obvious, while
$\mathtt{(\overline S_1)}$$\Rightarrow$$\mathtt{(\overline L^w)}$ follows
from Lemma \ref{lem:APPREFIN} and \ref{lem:SVCHAR}.b.

$\mathtt{(\overline S_1)}$$\Leftrightarrow$$\mathtt{(\overline S_2)}$$\Leftrightarrow$$\mathtt{(\overline S_3)}$
 follows from Lemma \ref{lem:SVCHAR}.b.

$\mathtt{(\overline S_1)}$$\Rightarrow$$\mathtt{(\overline E_\Sigma)}$
is a consequence of Lemma \ref{lem:SVSIG}.
\\
By that we have proved the equivalence statements.

The statement about unicity and integrable functions should be straightforward.
\end{proof}
\end{lemma}

\begin{defin}\label{def:EXTENDINT}
\halfquote
Suppose that $f:X\rightarrow\mathcal V$ is a function.

If there is an $a$ as in Lemma \ref{lem:EQUIGEN} above then it is called the
extended Lebesgue-McShane integral.
We can use the notation
\[\int^{\mathrm{LM}} L(f,\mu)=a. \]
\end{defin}
\begin{point}
In one form or other (weaker or stronger) all the previous statements about integrable functions
can be transferred into statements about extended integrable functions.
(Cf. especially characterization Lemma \ref{lem:EQUIGEN}.$\mathtt{(\overline S_2)}$.)
This is mainly a technical matter so will not spend time with it.
\end{point}
\begin{lemma}\label{lem:OLDCHAR}
\halfquote
Suppose that $f:X\rightarrow\mathcal V$.

Then, we claim, $f$ is integrable if and only if $f$ is extended integrable and there is a set
$A\in \boldsymbol\Sigma\mathfrak S$ such that $f=f|_A$ (ie. $f$ vanishes outside of $A$).
\begin{proof}
Necessity is obvious. Sufficiency follows from
characterization version Lemma  \ref{lem:EQUIGEN}.$\mathtt{(\overline S_2)}$
and Lemma \ref{lem:SVSIG}.
\end{proof}
\end{lemma}

\paragraph{\textbf{C. Elementary constructions on filters}}
\begin{remin}\label{rem:FDE}
a.) If $\mathfrak F$ is a filter base on a set $X$ and $f:X\rightarrow Y$ then the induced
set system
\[f_*\mathfrak F=\{\{f(x)\,:\,x\in F\}\,:\,F\in\mathfrak F\}\]
will be a filter base on $Y$. This construction is
compatible with taking the generated filters (if we extend the image).

b.)
If $P$ is a set and $\leq$ is a preordering on $S$ (ie. it is reflexive and transitive) then
a subset $S\subset P$ is said to be  $\leq$-top if
\begin{itemize}
\item[\texttt{(T1)}] for each $p\in P$ there is an element $s\in S$ such that $p\leq s$, and
\item[\texttt{(T2)}] if $s\in S$ and $s\leq t$ then $t\in S$.
\end{itemize}
The set of $\leq$-tops form a filter base
\[\dire^{\leq}_{p\in P}.\]

Special case: If  $(P,\leq)$ is upward directed, ie. for any two $p_1,p_2\in P$ there exists an element
$p\in P$ such that $p_1,p_2\leq p$ the filter base above is equivalent to the filter base
\[\{\{s\in P\,:\,p\leq s\}\,:\,p\in P\}.\]

c.) If
\[a:P\rightarrow Y\]
is a map to a topological space $Y$ then we can define
\[\lim_{p\in P}^{\leq}a_p\]
as the limit of the filter base
\[\dire^{\leq}_{p\in P} a_p=a_*\dire^{\leq}_{p\in P}.\]
\end{remin}
\begin{point}
The limit construction ``$\dire$'' can be generalized to filter bases.
The guiding principle is that an element $a$ can be substituted by the single set-element filter base
\[\{\{a\}\}.\]
\end{point}
\begin{remin}\label{rem:FDS}
a.)
If $\{\mathfrak F_\lambda\}_{\lambda\in\Lambda}$ is  a family of filter bases then we can consider
their finest common coarsening
\[\prod_{\lambda\in\Lambda}\mathfrak F_\lambda=\left\{\bigcup_{\lambda\in\Lambda} F(\lambda)\,:\, F
\text{ is a choice function for } \{\mathfrak F_\lambda\}_{\lambda\in\Lambda}\right\}.\]
(In terms of the generated filters this is just their intersection.)

This filter base will converge to a point if and only if every component filter base
converges to that point.

b.) If $\{\mathfrak F_\lambda\}_{\lambda\in\Lambda}$ is  a family of filter bases such that
for any $\alpha,\beta\in\Lambda$ there exist $\gamma\in\Lambda$ such that
$\mathfrak F_\gamma \succ\mathfrak F_\alpha,\, \mathfrak F_\beta$
 then we can consider their their coarsest common refinement
\[\coprod_{\lambda\in\Lambda}\mathfrak F_\lambda= \bigcup_{\lambda\in\Lambda}\mathfrak F_\lambda.\]
(It also yields union in terms of the generated filters.)

This filter base will converge to a point if one of the component filter bases
converges to that point, but there is no implication in the other direction.

c.) If $(P,\leq)$ is a preordered set and $\{\mathfrak F_{p}\}_{p\in P}$ is an indexed set
of filter bases then we can consider the filter base
\[\dire^{\leq}_{p\in P} \mathfrak F_p=
\coprod_{\Lambda\in\scalebox{.8}{$\displaystyle{\dire^{\leq}_{p\in P} }$}}\biggl(
\prod_{p\in\Lambda}\mathfrak F_p\biggr).\]


d.)
These  constructions above are compatible with filter base equivalence,
ie. if we substitute the filter bases $\mathfrak F_p$ by equivalent filter bases
$\mathfrak F_p'$ then the constructions above will yield equivalent filter bases
as results. Also, the constructions above are compatible with direct image
as in Reminder \ref{rem:FDE}.a. More precisely, the constructions above and taking direct image commute.
\end{remin}
\begin{remak}
We use ``$\dire$'' for the limit construction of filter bases because ``$\lim$''
is reserved for the (topological)  limit of a concrete filter base.
\end{remak}
What follows will be used only in the next section.
\begin{remin}
\label{rem:FILADD}
Assume that $\mathcal V$ is a commutative group

a.) If $\mathfrak F_1$ and $\mathfrak F_2$ are filter bases on a commutative group $\mathcal V$
then we define their sum
\[\mathfrak F_1+\mathfrak F_2=\{F_1+F_2\,:\,F_1\in\mathfrak F_1,\,F_2\in\mathfrak F_2\}.\]
Similarly, for a filter base $\mathfrak F$ on $\mathcal V$ we can define
\[-\mathfrak F=\{-F\,:\,F\in\mathfrak F\}.\]

b.) More generally, if $\{\mathfrak F_\lambda\}_{\lambda\in\Lambda}$ is a countable indexed family
of filter bases on a commutative group $\mathcal V$ then we can define the formal sum
\[\sum_{\lambda\in\Lambda}^{\mathrm{form}}\mathfrak F_\lambda=
\biggl\{\sum_{\xi\in\Xi} F_\xi+\sum_{\lambda\in\Lambda\setminus\Xi} (F_\lambda\cup\{0\})
 \,:\,F_\lambda\in\mathfrak F_\lambda,\,\Xi\subset\Lambda\text{ is finite}\biggr\}. \]

c.) Each element $v\in\mathcal V$ defines a canonical filter base
\[\{\{v\}\}\]
on $\mathcal V$.

In the light of this, for a countable indexed family
$\{v_\lambda\}_{\lambda\in\Lambda}$ of elements of $\mathcal V$
we can define
\[\sum_{\lambda\in\Lambda}^{\mathrm{form}} v_\lambda=
\sum_{\lambda\in\Lambda}^{\mathrm{form}}\{\{ v_\lambda\}\}=
\biggl\{\sum_{\xi\in\Xi} v_\xi+\sum_{\lambda\in\Lambda\setminus\Xi} \{0,v_\lambda\}
 \,:\,\Xi\subset\Lambda\text{ is finite}\biggr\}. \]

d.)
Suppose that $\mathfrak F$ is a filter base on the commutative group $\mathcal V$.
Then we define the Cauchy test  filter base $C(\mathfrak F)$ of $\mathfrak F$ as the
finest common coarsening of the filter bases
\[\mathfrak F-\mathfrak F+\mathfrak F-\mathfrak F+\ldots+\mathfrak F-\mathfrak F_{\quad(2n \text{\,terms})}.\]
\end{remin}
\begin{lemma}\label{lem:FORMALCRIT}
Suppose that $\mathfrak F$ is a filter base on the commutative group $\mathcal V$.

a.) Then, we claim,
$C(\mathfrak F)$ is a $0$-neighborhood base of
a (not necessarily Hausdorff or complete) topology $\mathtt T_{\mathfrak F}$
on $\mathcal V$ compatible with the group structure.

b.)  A filter base  $\mathfrak G$ on $\mathcal V$ is a Cauchy filter base  with respect to the topology
 $\mathtt T_{\mathfrak F}$ above
if and only
\[C(\mathfrak G)\succ C(\mathfrak F).\]

c.) For each such $\mathtt T_{\mathfrak F}$-Cauchy filter base $\mathfrak G$ there is a
one coarsest filter  $\mathfrak G'$, which is $\mathtt T_{\mathfrak F}$-Cauchy
and it is coarser than $\mathfrak G$.
This is the the filter
\[\mathfrak G'=\mathfrak G^{\mathrm{t\,\,with\,\,\,\mathtt T_{\mathfrak F}\,,f}}
=(\mathfrak G+C(\mathfrak F))^{\mathrm f}.\]

d.) The space
\[{}_{\mathfrak F}\mathcal V\]
of the coarsest $\mathtt T_{\mathfrak F}$-Cauchy filters forms a commutative group under
the operations ``$+$'' and  ``$-$''.

The space ${}_{\mathfrak F}\mathcal V$ has a topology
${}_{\mathfrak F}\mathtt T$, with the prescription that for each $\mathtt T_{\mathfrak F}$-open
set $U$ the set
\[\mathcal U=\{x\in{}_{\mathfrak F}\mathcal V\,:\,U\in x\},\]
is ${}_{\mathfrak F}\mathtt T$-open
and these are exactly the open sets of ${}_{\mathfrak F}\mathcal V$.

The topology ${}_{\mathfrak F}\mathtt T$ is complete, Hausdorff;
$({}_{\mathfrak F}\mathcal V,{}_{\mathfrak F}\mathtt T)$ it is the (Hausdorff-Cauchy)
reduction-completion
of $(\mathcal V,\mathtt T_{\mathfrak F})$.

e.) There is a natural continuous map
\[\iota_{\mathfrak F}:\mathcal V\rightarrow{}_{\mathfrak F}\mathcal V\]
given by
\[x\mapsto (\{\{x\}\}+C(\mathfrak F))^{\mathrm f}.\]

f.) If $\mathcal Z$ is a commutative topological group (complete, Hausdorff), and
\[\phi:\mathcal V\rightarrow \mathcal Z\]
is a homomorphism such that
\[\phi_*\mathfrak F\]
is convergent then there exists a unique continuous map
$\psi: {}_{\mathfrak F}\mathcal V\rightarrow \mathcal Z$ such that the diagram
\[\xymatrix{\mathcal V\ar[r]^{\iota_{\mathfrak F}}\ar[dr]_\phi & {}_{\mathfrak F}\mathcal V\ar[d]^\psi \\ &
\mathcal Z}\]
commutes. Ie. ${}_{\mathfrak F}\mathcal V$ is the crudest, universal reduction-completion of $\mathcal V$
such that the image of $\mathfrak F$ is convergent. \qed
\end{lemma}

\paragraph{\textbf{D. The extended integral as a limit}}
\begin{defin}\label{def:SPECLIM}
Suppose that
$\mathfrak S\subset \mathfrak P(X)$ is an interval system. For
$A,B\in\boldsymbol\Sigma_{00}\mathfrak S$ we define the relation
\[A\,\dot\subset\,B\]
to hold if and only if there is a set $C\in\boldsymbol\Sigma_{00}\mathfrak S$
such that
\[A\,\dot\cup\,C=B\]
(ie. $B$ is the disjoint union of $A$ and $C$.)
\end{defin}
\begin{defin}\label{def:EQUILIM}
\halfquote

Suppose that $f:X\rightarrow \mathcal V$ is a function, and
$A\in\boldsymbol\Sigma\mathfrak S.$

Then we define
\[{}^{\mathrm{LM}}\widetilde{\mathfrak F}_{L,f,\mu}^A =
\biggl\{\intp L(\mathcal A,\mu)\,:\,\mathcal A\text{ is  an envelope approximating }f|_A
\biggr\}^{\mathrm{tf}}.\]
Notice that there are envelopes like that: for example $\mathcal A=(\mathcal V\times A)\cup(\{0\}\times X)$
is a possible envelope. Notice the ``smearing'' and the filter extension.
Lemma \ref{lem:APPFILTER} shows after smearing we have a filter base, so the filter extension
can be taken.
\end{defin}
\begin{lemma}\label{lem:EQUILIM}
\halfquote
Suppose that $f:X\rightarrow\mathcal V$ is a function.

Then, we claim, the following statements are equivalent:

$\mathtt{(\overline E)}$ (cf. Lemma \ref{lem:EQUIGEN} and Definition \ref{def:EXTENDINT})
\[\int^{\mathrm{LM}} L(f,\mu)=a. \]

$\mathtt{(\widetilde E_{00})}$
\[\lim\dire_{A\in\boldsymbol\Sigma_{00}\mathfrak S}^{\dot\subset}
{}^{\mathrm{LM}}\widetilde{\mathfrak F}_{L,f,\mu}^A=a.\]

$\mathtt{(\widetilde E_{0})}$
\[\lim\dire_{A\in\boldsymbol\Sigma_{0}\mathfrak S}^{\subset}
{}^{\mathrm{LM}}\widetilde{\mathfrak F}_{L,f,\mu}^A=a.\]

$\mathtt{(\widetilde E_{\Sigma})}$
\[\lim\dire_{A\in\boldsymbol\Sigma\mathfrak S}^{\subset}
{}^{\mathrm{LM}}\widetilde{\mathfrak F}_{L,f,\mu}^A=a.\]

$\mathtt{(\widetilde L_{00})}$ The function $f$ is locally integrable and
\[\lim_{A\in\boldsymbol\Sigma_{00}\mathfrak S}^{\dot\subset}
\int L(f|_A,\mu)=a.\]

$\mathtt{(\widetilde L_{0})}$ The function $f$ is locally integrable and
\[\lim_{A\in\boldsymbol\Sigma_{0}\mathfrak S}^{\subset}
\int L(f|_A,\mu)=a.\]
\begin{proof}
$\mathtt{(\widetilde E_{00})}$ is equivalent to $\mathtt{(\overline E_{00})}$.
$\mathtt{(\widetilde E_{0})}$ is equivalent to $\mathtt{(\overline E_{0})}$.
$\mathtt{(\widetilde E_{\Sigma})}$ is equivalent to $\mathtt{(\overline E_{\Sigma})}$.
$\mathtt{(\overline L^m)}$$\Rightarrow$$\mathtt{(\widetilde L_{00})}$$\Rightarrow$$\mathtt{(\overline L^w)}$.
$\mathtt{(\overline L^m)}$$\Rightarrow$$\mathtt{(\widetilde L_{0})}$$\Rightarrow$$\mathtt{(\overline L^w)}$.
\end{proof}
\end{lemma}

\newpage
\section{Discrete approximations}\label{sec:discrete}

In this section we show that the conditions on envelopes can be relaxed.
The price we have to pay is that $\intp L(\mathcal  A,\mu)$ cannot be defined as
a single set but we have to think in terms of filter bases.
While this requires a more obscure notion of convergence, it yields way to a more powerful approach.
\\

\paragraph{\textbf{A. Generalized sums}}
\begin{defin}
Let $\mathcal  V$ be a commutative topological group, and let
$\{A_\lambda\}_{\lambda\in\Lambda}$ be a countable family of subsets of $\mathcal  V$.
Then we say that the sum of $\{A_\lambda\}_{\lambda\in\Lambda}$ is formally contained in the closed set $T$
\[\sum^{\ovtild}_{\lambda\in\Lambda} A_\lambda\fsubset T \]
if for each neighborhood $0\in\mathcal  T\subset\mathcal  V$ there exist a finite set $\Xi\subset\Lambda$ such that
\[\sum_{\lambda\in\Xi} A_\lambda+\sum_{\lambda\in\Lambda\setminus\Xi}(A_\lambda\cup\{0\})\subset T+\mathcal  T. \]

If $Z\subset\mathcal  Z$ is a general set then we write
\[\sum^{\ovtild}_{\lambda\in\Lambda} A_\lambda\fsubset Z \]
if $Z$ contains a closed set $T$ as above.

If $a_\lambda\in\mathcal  V$ then we write
\[\sum^{\ovtild}_{\lambda\in\Lambda} a_\lambda\fsubset Z\]
if the statements holds with $A_\lambda=\{a_\lambda\}$.
\end{defin}
\begin{remark}\label{rem:URREM}
For a closed set $T$ the situation can be described as
\[\sum^{\mathrm{form}}_{\lambda\in\Lambda} A_\lambda\succ\{T\}^{\mathrm t}, \]
but we rather use the more suggestive notation above.
\end{remark}
\begin{lemma}\label{lem:URCONT}
Let $\mathcal  V$ be a commutative topological group.
Suppose that $\{a_\lambda\}_{\lambda\in\Lambda}$ is a countable indexed family of elements of $\mathcal  V$.
Assume that
\[\Lambda=\bigcup_{\gamma\in\Gamma}^{\updisjoint}\Lambda_\gamma\]
is a decomposition.

a.) Assume that
\[a_{(\gamma)}=\sum_{\lambda\in\Lambda_\gamma}^{\ovline}a_\gamma\]
is convergent for all $\gamma\in\Gamma$. Then, we claim,
\[\sum^{\ovtild}_{\lambda\in\Lambda} a_\lambda\fsubset Z\]
implies
\[\sum^{\ovtild}_{\gamma\in\Gamma} a_{(\gamma)}\fsubset Z.\]

b.) If $\Gamma$ is finite and
\[\sum^{\ovtild}_{\lambda\in\Lambda_\gamma} a_\lambda\fsubset Z_\gamma \]
then
\[\sum^{\ovtild}_{\lambda\in\Lambda} a_\lambda\fsubset
\overline{\sum_{\gamma\in\Gamma} Z_\gamma }.\]
\begin{proof} This follows in the usual manner.
\end{proof}
\end{lemma}

\paragraph{\textbf{B. Unrestricted envelopes}}
\begin{defin}\label{def:URAPPNEIGH}
Suppose that $\mathcal  V$ is a commutative topological group and
$\mathfrak S\subset\mathfrak P(X)$ is an interval system. Then we say that
\[\mathcal  H\subset \mathcal  V\times X\]
is an unrestricted envelope of sets in $\mathcal  V$ with respect to $\mathfrak S$ if
\begin{itemize}
\item[\texttt{(S$_{\mathtt f}$)}] For all $x\in X$
\[\emptyset\neq\mathcal  H^x.\]
\item[\texttt{(G)}]
There exists a countable generating set $\mathfrak D\subset\mathfrak S$
for $\mathcal  A$; ie. a countable family of sets $\mathfrak D\subset\mathfrak S$ such that
\[\mathcal  H=\langle\mathfrak D\rangle_{\mathcal  H}.\]
\end{itemize}

We say that $\mathcal  H$ is supported on a set $A\in\boldsymbol\Sigma\mathfrak S$
is there is a generator system $\mathfrak D$ for $\mathcal  H$ such that $A=\bigcup\mathfrak D$.
\end{defin}
\begin{remark}
So it is not necessarily true that $\mathcal H$ is supported on $A\in\boldsymbol\Sigma\mathfrak S$
implies that $\mathcal  H$ is supported on $B\in\boldsymbol\Sigma\mathfrak S$ if $B\supset A$.
\end{remark}
There are plenty of examples:
\begin{lemma}\label{lem:URCO}
Suppose that $\mathcal  V$ is a commutative topological group and
$\mathfrak S\subset\mathfrak P(X)$ is an interval system.

Assume that $\mathcal  A$ is  an envelope with a generator system $\mathfrak D$.
Then, we claim, $\mathcal  A$ is an unrestricted envelope on any set $A\in\boldsymbol\Sigma\mathfrak S$
such  that $A\supset \bigcup\mathfrak D$.
\begin{proof} Suppose that $\mathfrak B\subset\mathfrak S$ is countable,
$B=\bigcup\mathfrak B$,  and
$s=\sum_{j\in J}c_j\chi_{E_j}$ is a step-function such that its graph is contained in $\mathcal  A$,
\[\mathfrak E=\{E_j\,:\,j\in J\}.\]
Then apply Lemma \ref {lem:FOREST} to $\mathfrak D\cup\mathfrak B$ such that  the result should be
divided by $\mathfrak E$. One can see that it yields a generator system for $\mathcal  A$.
\end{proof}
\end{lemma}
\begin{lemma}\label{lem:URAPPNEIGH}
Let $\mathfrak S\subset \mathfrak P(X)$ be a interval system
and $\mathcal  V$ be a commutative topological group.
Assume that $\mathcal  H$ is an unrestricted envelope of sets in $\mathcal  V$ with respect to $\mathfrak S$.

a.) Suppose that $\mathfrak D$ is a generator system for $\mathcal  H$.
Assume that $\mathfrak D'\subset\mathfrak S$
is an other countable family of sets which is finer than $\mathfrak D$ but
\[\bigcup\mathfrak D'=\bigcup\mathfrak D.\]
Then, we claim, $\mathfrak D'$ is also a generator system  for $\mathcal  H$.

b.) The family supporting sets of $\mathcal  H$ is closed for
finite intersections.
In particular, if $\mathcal  H$ is supported on $A_1$ and $A_2$ then it is also supported on
$A_1\cap A_2$.

c.) The family supporting sets of $\mathcal  H$ is closed for countable unions.
In particular, if $\mathcal  H$ is supported on $A_1$ and $A_2$ then it is also supported on
$A_1\cup A_2$.

d.) If $\mathcal  H$ is supported on $A$ and $\{A_\omega\}_{\omega\in\Omega}$ is a
countable family of
sets such that $\mathcal  H^{A_\omega}\neq\emptyset$ for all $\omega\in\Omega$ then
$\mathcal  H$ is also supported on
\[A\cup\bigcup_{\omega\in\Omega}A_\omega.\]
\begin{proof}
These statements are all straightforward, except b.). In that case let $\mathfrak D_1$ and
$\mathfrak D_2$ be the corresponding generator systems. Then for each pair of sets
$D_1\in\mathfrak D_1$,
$D_2\in\mathfrak D_2$, the set $D_1\cap D_2$ decomposes in $\mathfrak S$. Let $\mathfrak D$
be the union of all these decompositions. Then one can see that $\mathfrak D$ is a generator
system for $\mathcal  H$ and $\bigcup\mathfrak D=A_1\cap A_2$.
\end{proof}
\end{lemma}
\begin{defin}\label{def:URBLIN}
\halfquote

Suppose that $\mathcal  H$ is an unrestricted envelope of sets in $\mathcal  V$ with respect to $\mathfrak S$,
supported on $A\in\boldsymbol\Sigma\mathfrak S$.
Then for a closed set $T\subset\mathcal  Z$ we write
\[\intform_A L(\mathcal  H,\mu)\fsubset T\]
if
\[\sum_{\omega\in\Omega}^{\ovtild}L(c_\omega,\mu(E_\omega))\fsubset T\]
whenever $\Omega$ is countable, the sets $E_\omega\in\mathfrak S$ ($\omega\in\Omega$) are pairwise disjoint,
$A=\displaystyle \bigcup^{\updisjoint}_{\omega\in\Omega}E_{\omega}$,
and $c_\omega\in \mathcal  H^{E_\omega}$.
We extend this notation for all sets $Z\subset \mathcal Z$ which contain a closed set $T$ as above.
\end{defin}
\begin{remak}
For a closed set $T$,
again, we can also make sense out this notation by taking the finest common coarsening
(cf. Reminder \ref {rem:FDS}.a) of a bunch of filter bases and proceed as in Remark
\ref{rem:URREM}.
\end{remak}
\begin{lemma} \label{lem:URNEIGH}
Let $\mathfrak S\subset \mathfrak P(X)$ be a interval system and $\mathcal  V$ be a
commutative topological group. Suppose that $\mathcal  H_1$ and $\mathcal  H_2$ are unrestricted envelopes
supported on $A\in\boldsymbol\Sigma\mathfrak S$. Consider the pointwise sum
\[\mathcal  H_1+\mathcal  H_2=\biggl\{(v,x)\,:\,x\in X,\,v\in
\mathcal  H_1{}^x+ \mathcal  H_2{}^x\biggr\}.\]
Then, we claim, the (pointwise) sum above is an unrestricted envelope of sets in $\mathcal  V$
with respect to $\mathfrak S$ supported on $A$.

Moreover, if $\mathcal  W$, $\mathcal  Z$ are commutative topological groups,
$\mu:\mathfrak S\rightarrow\mathcal  W$ is a measure,
and  $L:\mathcal  V\times\mathcal  W\rightarrow\mathcal  Z$ is a biadditive pairing
which is continuous in its second variable
and
\[\intform_A L(\mathcal  H_1,\mu)\fsubset Z_1\quad\qand\quad \intform_A L(\mathcal  H_2,\mu)\fsubset Z_2\]
then
\[\intform_A L(\mathcal  H_1+\mathcal  H_2,\mu)\fsubset\overline{ Z_1+Z_2}.\]
\begin{proof}
Regarding  the first part: If $\mathfrak D_1$ and $\mathfrak D_2$ are corresponding generator systems
the let us apply Lemma \ref{lem:FOREST} to  $\mathfrak D_1\cup\mathfrak D_2$. One can see that the
resulted forest yields a generator system for $\mathcal  H_1+\mathcal  H_2$.

Regarding the second part:
Consider any sum as in the Definition \ref{def:URBLIN}.  Let
\[\mathfrak E=\{E_\omega\,:\,\omega\in\Omega\}.\]
Let us apply Lemma \ref{lem:FOREST} to  $\mathfrak D_1\cup\mathfrak D_2\cup\mathfrak E$,
let us denote the result by $\mathfrak R$.
Then for all such $x\in A\cap\bigcup\mathfrak E$ there exist sets
$D_1\in\mathfrak D_1,\,D_2\in\mathfrak D_2,E_\omega\in\mathfrak E$ containing $x$
such that $c_\omega\in\mathcal  H_1{}^{D_1}+\mathcal  H_2{}^{D_2}$.
Hence, there exists a set $D\in\mathfrak R$ such that $x\in D\subset E_\omega\cap D_1\cap D_2$,
and so $c_\omega\in
\mathcal  H_1{}^{D}+\mathcal  H_2{}^D$. Let $\mathfrak  R'$ be the set of maximal such elements
$D\in\mathfrak R$
such that $D$ is contained in some
 sets $E_\omega\in\mathfrak E$, $D_1\in\mathfrak D_1$, $D_2\in\mathfrak D_2$ and
$c_\omega\in\mathcal  H_1{}^{D}+\mathcal  H_2{}^D$.

Now, $\mathfrak R'$ exactly decomposes  $\mathfrak E$.
 Let us use the notation
$c_D'=c_\omega$ if $D\subset E_{\omega}$.

Hence, according to Lemma \ref{lem:URCONT}.a it is enough to prove that
\[\sum_{D\in\mathfrak R'}^{\ovtild} L(c_D',\mu(D))\fsubset \overline{Z_1+Z_2}. \]

Let $c_D'=c_{D,1}+c_{D_2}$ where $c_{D,1}\in \mathcal H_1{}^D$ and
$c_{D,2}\in\mathcal H_2{}^D$.
Then
\[\sum_{D\in\mathfrak R'}^{\ovtild} L(c_{D,1},\mu(D))\fsubset Z_1\quad\qand\quad
\sum_{D\in\mathfrak R'}^{\ovtild} L(c_{D,2},\mu(D))\fsubset Z_2. \]
implies our statement according to Lemma \ref{lem:URCONT}.b and a.
\end{proof}
\end{lemma}
\begin{lemma} \label{lem:UREXTEND}
\halfquote

a.) If $\mathcal  A$ is  an envelope, which is an unrestricted envelope
supported on $A\in\boldsymbol\Sigma\mathfrak S$ then
\[ \intform_A L(\mathcal  A,\mu)\fsubset T\quad\Leftrightarrow\quad \intp L(\mathcal  A,\mu)\subset T.\]

b.) If $\mathcal  H$ is an unrestricted envelope supported on $A$ and
\[\intform_A L(\mathcal  H,\mu)\fsubset T\]
then for each neighborhood $0\in\mathcal  U\subset\mathcal  Z$
there exists  an envelope $\mathcal  A\supset\mathcal  H$, which is
an unrestricted envelope supported on $A$,  such that
\[\intp L(\mathcal  A,\mu)\subset T-T+T+\mathcal  U.\]
\begin{proof}
a.) It is easy to see this if $\mathcal  A=s$ is a (graph of a) step-function.
Also, the statement is clear if $\mathcal  A=\mathcal  C$ is a pointed envelope.
For a sum
\[\mathcal  A=s+\mathcal  C\]
the statement follows from the previous lemma.

b.) Let $\mathfrak D$ be a generator system for $\mathcal H$.
Then, by the usual methods, we can find a exact decomposition of $\mathfrak D'$ of
$\bigcup\mathfrak D$ such that for each element $D\in\mathfrak D'$ there is an element
$c_D\in\mathcal  A^D$. Let $\mathfrak S'$ contain the elements $S\in\mathfrak S$ such that
$S\subset D$ for an $D\in\mathfrak D'$. We see that
$\mathfrak S'$ is an interval system on $\bigcup\mathfrak D$.

Let us define a measure $\tilde\mu:\mathfrak S'\rightarrow\mathcal  Z$ by
\[\tilde\mu(S)=L(c_D,\mu(S))\]
if $S\subset D\in\mathfrak D'$.
Notice that $\bigcup\mathfrak D$ countably decomposes in $\mathfrak S'$

According to our assumptions for each neighborhood $0\in\mathcal  T\subset\mathcal  Z $
 for each countable decomposition $\mathfrak E$
 of $\bigcup\mathfrak D$ in $\mathfrak S'$ we can find a finite subset
 $\Omega\subset \mathfrak E $ such that
 \[\sum_{E\in\Omega}\tilde\mu(E)+\sum_{E\in\mathfrak E\setminus\Omega}\{0,\tilde\mu(E)\}
 \subset T+\mathcal  T.\]
 In particular,
 \[\sum_{E\in\mathfrak E\setminus\Omega}\{0,\tilde\mu(E)\}
 \subset T-T+\mathcal  T-\mathcal  T.\]
Applying Lemma \ref{lem:XMEAVAR} we find that there are finitely many elements
$\Xi\subset \mathfrak D'$ such that
\[\svar\biggl(\bigcup_{D\in\mathfrak D'\setminus\Xi} D,\tilde\mu\biggr)
\subset T-T+\mathcal  T-\mathcal  T+\mathcal  T.\]
We can assume that $\mathcal T$ was chosen such that
\[\mathcal T-\mathcal T+\mathcal T+\mathcal T\subset\mathcal U,\]
and $\Xi$ is chosen accordingly.

Notice that
\[ \svar\biggl(\bigcup_{D\in\mathfrak D'\setminus\Xi}D ,\tilde\mu\biggr)=
\overline{\sum_{D\in\mathfrak D'\setminus\Xi}\svar(D ,\tilde\mu)}\]
\[=\overline{\sum_{D\in\mathfrak D'\setminus\Xi}\svar(L,c_D,D ,\mu)}=
\intp L\biggl(\sum_{D\in\mathfrak D'\setminus\Xi}\svar(c_D,D),\mu\biggr).\]
Hence, defining
\[\mathcal  C=\sum_{D\in\mathfrak D'\setminus\Xi}\svar(c_D,D)\]
as a pointed envelope we find
that
\[\intp L(\mathcal  C,\mu)\subset T-T+\mathcal  T-\mathcal  T+\mathcal  T.\]
Now, we see that
$\mathcal  H-\mathcal  C$
is  an envelope, containing the graph of
\[\sum_{D\in\Xi}c_D\chi_D.\]
Meanwhile, from Lemma \ref {lem:URNEIGH},
\[\intform L(\mathcal  H-\mathcal  C,\mu)\fsubset\overline{T-(T-T+\mathcal  T-\mathcal  T+\mathcal  T)}\]
\[\subset T-T+T-\mathcal  T+\mathcal  T-\mathcal  T+\mathcal  T\subset T-T+T+\mathcal  U. \]
According to point a. that implies
\[\intp L(\mathcal  H-\mathcal  C,\mu)\subset T-T+T+\mathcal  U.\]
\end{proof}
\end{lemma}

\begin{lemma}\label{lem:UEQUILIM}
\halfquote
Suppose that $f:X\rightarrow\mathcal V$ is a function.

Then, we claim, the following statements are equivalent:

$\mathtt{(\overline E)}$ (cf. Lemma \ref{lem:EQUIGEN} and Definition \ref{def:EXTENDINT})
\[\int^{\mathrm{LM}} L(f,\mu)=a. \]

$\mathtt{(\overline E_{\Sigma}')}$  For each neighborhood $0\in\mathcal  T\subset\mathcal  Z$
there is a set $A\in\boldsymbol\Sigma\mathfrak S$ such that
for all $A\subset B\in\boldsymbol\Sigma\mathfrak S$
there is an unrestricted envelope $\mathcal  B$ supported on $B$ such that
\[f|_B\subset \mathcal  B\qquad\text{and}\qquad \intform_B L(\mathcal  B,\mu)\fsubset a+\mathcal  T.\]

$\mathtt{(\overline E_{0}')}$  For each neighborhood $0\in\mathcal  T\subset\mathcal  Z$
there is a set $A\in\boldsymbol\Sigma_{0}\mathfrak S$ such that
for all $A\subset B\in\boldsymbol\Sigma_{0}\mathfrak S$ there is an unrestricted envelope $\mathcal  B$
 supported on $B$  such that
\[f|_B\subset \mathcal  B\qquad\text{and}\qquad \intform_B L(\mathcal  B,\mu)\fsubset a+\mathcal  T.\]

$\mathtt{(\overline E_{00}')}$   For each neighborhood $0\in\mathcal  T\subset\mathcal  Z$ and
set $A_0\in \boldsymbol\Sigma_{00}\mathfrak S$
there is a set $A\in\boldsymbol\Sigma_{00}\mathfrak S$ such that
for all $A'\in\boldsymbol\Sigma_{00}\mathfrak S$ such that $A\cap A'=\emptyset$
there is an unrestricted envelope $\mathcal  B$  supported on $B$ such that
\[f|_{(A\cup A')}\subset \mathcal  B\qquad\text{and}\qquad \intform_{(A\cup A')} L(\mathcal  B,\mu)
\fsubset a+\mathcal  T.\]
\begin{proof}
The implications
$\mathtt{(\overline E_{\Sigma})}$$\Rightarrow$$\mathtt{(\overline E_{\Sigma}')}$,
$\mathtt{(\overline E_{0})}$$\Rightarrow$$\mathtt{(\overline E_{0}')}$,
$\mathtt{(\overline E_{00})}$$\Rightarrow$$\mathtt{(\overline E_{00}')}$, follow from
Lemma \ref{lem:EORES}, Lemma \ref{lem:URCO} and Lemma  \ref{lem:UREXTEND}.a.

The statements
$\mathtt{(\overline E_{\Sigma}')}$$\Rightarrow$$\mathtt{(\overline E_{\Sigma})}$,
$\mathtt{(\overline E_{0}')}$$\Rightarrow$$\mathtt{(\overline E_{0})}$,
$\mathtt{(\overline E_{00}')}$$\Rightarrow$$\mathtt{(\overline E_{00})}$, follow from
Lemma \ref{lem:UREXTEND}.b.
\end{proof}
\end{lemma}

\paragraph{\textbf{C. Tag systems}}
\begin{defin}\label{def:TAGDEF}
Suppose that $\mathfrak S\subset\mathfrak P(X)$ is an interval system. A tag system on $\mathfrak S$
is a pair $(\mathfrak D,h)$ such that
\begin{itemize}
\item[\texttt{(H1)}] $\mathfrak D$ is a countable subset of $\mathfrak S$;
\item[\texttt{(H2)}]
$h:\bigcup\mathfrak D\rightarrow\mathfrak D$ is a function such that
 \[x\in h(x)\] for all
$x\in \bigcup\mathfrak D$.
\end{itemize}
The tag system $(\mathfrak D,h)$ is supported on $A$ if $A=\bigcup\mathfrak D$.
\end{defin}
\begin{defin}\label{def:TAGFIN}
Suppose that $\mathfrak S\subset\mathfrak P(X)$ is an interval system.
We say that the tag system $(\mathfrak D',h')$ is finer than
$(\mathfrak D,h)$, if
\[\bigcup\mathfrak D=\bigcup\mathfrak  D'\]
and for all $x\in \bigcup\mathfrak D=\bigcup\mathfrak D'$
\[h'(x)\subset h(x).\]
\end{defin}
\begin{lemma}\label{lem:TAGFIN}
Suppose that $\mathfrak S\subset\mathfrak P(X)$ is an interval system.
Then we claim:

For any finitely many tag system supported on $A\in\boldsymbol\Sigma\mathfrak S$
there is a tag system supported on $A$ which is finer than any of them.
\begin{proof}
Suppose that $(\mathfrak D_1,h_1),\ldots, (\mathfrak D_n,h_n )$ are those tag systems.
Let us apply Lemma \ref{lem:FOREST} to $\mathfrak D_1\cup\ldots\cup\mathfrak D_n$, let us denote the result
by $\widetilde{\mathfrak D}$. Then for each $x\in A$ there is a set $D\in\widetilde{\mathfrak D}$
such that
\[x\in D\subset h_1(x)\cap\ldots\cap h_n(x).\]
Let $h(x)$ be such the maximal such $D$ in $\widetilde{\mathfrak D}$.
Then $(\widetilde{\mathfrak D},h)$ will be a tag system finer than any $(\mathfrak D_j,h_j)$.
\end{proof}
\end{lemma}
\begin{defin}\label{def:TAGDIV}
Suppose that $\mathfrak S\subset\mathfrak P(X)$ is an interval
system, and $(\mathfrak D,h)$ is a tag system on $\mathfrak S$.

A tagged Lebesgue-McShane division (or just tagged division) associated to $(\mathfrak D,h)$ is a pair
$(\mathfrak E, c)$ such that
\begin{itemize}
\item[\texttt{(D1)}]
 The set $\mathfrak E$ is countable set of pairwise disjoint nonempty elements of $\mathfrak S$;
\item[\texttt{(D2)}] $c:\mathfrak E\rightarrow \bigcup\mathfrak D$ is function such that
\[E\subset h(c(E))\]
 for each $E\in\mathfrak E$.
\end{itemize}
\end{defin}
\begin{lemma}\label{lem:TAGREFIN}
Suppose that $\mathfrak S\subset\mathfrak P(X)$ is an interval
system.

a.) Any tag system $(\mathfrak D,h)$ allows a tagged  Lebesgue-McShane  division $(\mathfrak E,c)$.

b.) If the tag system $(\mathfrak D',h')$ is finer than $(\mathfrak D,h)$ and
$(\mathfrak E,c)$ is a tagged division associated to $(\mathfrak D',h')$, then
$(\mathfrak E,c)$ is also associated to $(\mathfrak D,h)$.
\begin{proof}
a.) Let us apply the Lemma \ref{lem:FOREST} to $\mathfrak D$.
Then for each $x\in\bigcup\mathfrak D$ there is a set $D\in\widetilde{\mathfrak D}$ such that
\[x\in D\subset h(x).\]
Hence the sets $D\in \widetilde{\mathfrak D}$ for which there is an element $x\in D$ so
that $D\subset h(x)$ cover $\bigcup\mathfrak D$.
Let $\mathfrak E$ be the set of all the maximal sets as above.
For $E\in\mathfrak E$ let  $c(E)\in E$ be a such point that $E\subset h(c(E))$.

b.) That is immediate from the definitions.
\end{proof}
\end{lemma}

\begin{defin}\label{def:TAGFILTER}
Suppose that $\mathfrak S\subset\mathfrak P(X)$ is an interval system.

a.) A (finite) formal sum over  $\mathfrak S\subset\mathfrak P(X) $ is a finite expression
\[\sum_{j\in J}\delta_{x_j}\cdot E_j,\]
where $J$ is finite, $x_j\in X$, $E_j\in\mathfrak S$.
They form a commutative group $\mathsf F(X,\mathfrak S)$.

b.) If $(\mathfrak E,c)$ is a tagged division then its associated filter base is
\[\mathfrak F_{(\mathfrak E,c)}=\sum^{\mathrm{form}}_{E\in\mathfrak E}\delta_{c(E)}\cdot E\]

c.) If $(\mathfrak D,h)$ is a tag system then its associated Lebesgue-McShane filter base is
\[{}^{\mathrm{LM}}\mathfrak F^{(\mathfrak D,h)}=\prod_{
\substack{(\mathfrak E,c)\text{ Lebesgue-McShane tagged }
\\\text{ division associated to }(\mathfrak D,h)}}\mathfrak F_{(\mathfrak E,c)}\]
ie. the coarsest filter base generated by the tagged divisions associated to $(\mathfrak D,h)$.

d.) If $A\in\boldsymbol\Sigma\mathfrak S$ then its associated  Lebesgue-McShane filter base is
\[{}^{\mathrm{LM}}\mathfrak F^A=\coprod_{ (\mathfrak D,h)
\text{ supported on  }A}{}^{\mathrm{LM}}\mathfrak F^{(\mathfrak D,h)}, \]
ie. the finest filter base generated by tag systems on $A$.
(Cf. Lemma \ref{lem:TAGREFIN}.b for the applicability of the definition.)

e.) The associated Lebesgue-McShane filter base of countable kind is defined as
\[{}^{\mathrm{LM}}\mathfrak F^{\mathfrak S}=\dire_{A\in \boldsymbol\Sigma\mathfrak S}^{\subset}
{}^{\mathrm{LM}}\mathfrak F^A,\]
ie. the limit of the filter bases ${}^{\mathrm{LM}}\mathfrak F^A$ as $A\in\boldsymbol\Sigma\mathfrak S$ increases.

 The associated Lebesgue-McShane filter base of finite kind is defined as
\[{}^{\mathrm{LM}}\mathfrak F_0^{\mathfrak S}=\dire_{A\in \boldsymbol\Sigma_0\mathfrak S}^{\subset}
{}^{\mathrm{LM}}\mathfrak F^A,\]
ie. the limit of the filter bases ${}^{\mathrm{LM}}\mathfrak F^A$ as $A\in\boldsymbol\Sigma_0\mathfrak S$ increases.

 The associated Lebesgue-McShane filter base of disjoint finite kind is defined as
\[{}^{\mathrm{LM}}\mathfrak F_{00}^{\mathfrak S}=\dire_{A\in \boldsymbol\Sigma_0\mathfrak S}^{\dot\subset}
{}^{\mathrm{LM}}\mathfrak F^A,\]
ie. the limit of the filter bases ${}^{\mathrm{LM}}\mathfrak F^A$ as $A\in\boldsymbol\Sigma_{00}\mathfrak S$ increases.

\end{defin}
\begin{remak}
Despite of the multiple steps, the Lebesgue-McShane filter bases are just filter bases on finite formal sums.
\end{remak}
\begin{lemma}
Suppose that $\mathfrak S\subset\mathfrak P(X)$ is an interval system. We claim:

a.) If $X\in\boldsymbol\Sigma\mathfrak S$ then ${}^{\mathrm{LM}}\mathfrak F^{\mathfrak S} $
is equivalent to ${}^{\mathrm{LM}}\mathfrak F^{X}$.

b.) If $X\in\boldsymbol\Sigma_0\mathfrak S$ then ${}^{\mathrm{LM}}\mathfrak F^{\mathfrak S} $
and ${}^{\mathrm{LM}}\mathfrak F_0^{\mathfrak S} $
are both equivalent to ${}^{\mathrm{LM}}\mathfrak F^{X}$.

b.) If $X\in\boldsymbol\Sigma_{00}\mathfrak S$ then
${}^{\mathrm{LM}}\mathfrak F^{\mathfrak S} $,
${}^{\mathrm{LM}}\mathfrak F_0^{\mathfrak S} $, and ${}^{\mathrm{LM}}\mathfrak F_{00}^{\mathfrak S} $
are all equivalent to ${}^{\mathrm{LM}}\mathfrak F^{X}$.

\begin{proof}
This immediately follows from the construction.
\end{proof}
\end{lemma}

\begin{defin}
\halfquote
Suppose that $f:X\rightarrow \mathcal  V$ is a function.

We define the action $L(f,\mu)$ on finite formal sums by
\[L(f,\mu)\biggl(\sum_{j\in J}\delta_{x_j}\cdot E_j\biggr)=\sum_{j\in J}L(f(x_j),\mu( E_j)).\]
\end{defin}

\paragraph{\textbf{D. Tag systems vs unrestricted envelopes}}
\begin{lemma}\label{lem:TAGRES}
Let $\mathfrak S\subset\mathfrak P(X)$ be an interval system and $\mathcal V$ be a
commutative topological group.
Suppose  $f:X\rightarrow\mathcal  V$ is a function.

a.) Assume that $(\mathfrak D,h)$ is a tag system
such that
 that $f|_{\bigcup\mathfrak D}=f$. (Ie. ${\bigcup\mathfrak D}$
 covers the non-zero places of $f$.)
Then we claim,
\[\mathcal  H_{f,(\mathfrak D,h)}=\bigcup_{x\in \bigcup\mathfrak D}\{f(x)\}\times h(x)\cup
\left(\{0\}\times \left(X\setminus\bigcup\mathfrak D\right)\right)\]
is an unrestricted envelope approximating $f$.

b.) If $\mathcal  H$ is an unrestricted envelope approximating $f$
with generator system $\mathfrak D$
then we can find
a tag system $(\mathfrak D,h)$ as above such that
\[f\subset \mathcal  H_{f,(\mathfrak D,h)}\subset \mathcal  H.\]
\begin{proof}
a.) This is straightforward; $\mathfrak D$ will be a generator system.

 b.)
Assume that $\mathfrak D$ is a generator system for $\mathcal  H$.
For each $x\in\bigcup\mathfrak D$ there is a set $D\in\mathfrak D$ such that
$x\in D$ and $f(x)\in\mathcal  H^D$. For each $x\in\bigcup\mathfrak D$
we can choose such a $D$ as $h(x)$. Then the statement is obvious.
\end{proof}
\end{lemma}

\begin{lemma}\label{lem:TAGVAL}
\halfquote

Suppose that $f:X\rightarrow\mathcal  V$ is a function.
Then, for a tag system $(\mathfrak D,h)$ supported on $B$, and
 a closed set $T\subset \mathcal  V$ the following two statements are equivalent:

i.)
\[\intform_B L(\mathcal  H_{f,(\mathfrak D,h)},\mu)\fsubset T\]

ii.) For each $0\in\mathcal  T\subset\mathcal  Z$ there is set $G$
\[L(f,\mu)_*\,\, {}^{\mathrm{LM}}\mathfrak F^{(\mathfrak D,h)}\ni G\subset T+\mathcal  T.\]
\begin{proof}
That follows from the construction of the associated Lebesgue-McShane filter bases.
\end{proof}
\end{lemma}

\begin{lemma}\label{lem:TAGINT}
\halfquote

Suppose that $f:X\rightarrow\mathcal  V$ is a function.
Then the following two statements are equivalent:

$\mathtt{(\overline E)}$ The extended Lebesgue-McShane integral of $f$ is $a$.

$\mathtt{(E_{\Sigma})}$  The filter base
\[L(f,\mu)_*\,\, {}^{\mathrm{LM}}\mathfrak F^{\mathfrak S}\]
converges to $a$.

$\mathtt{(E_{0})}$ The filter base
\[L(f,\mu)_*\,\, {}^{\mathrm{LM}}\mathfrak F_0^{\mathfrak S}\]
converges to $a$.

$\mathtt{(E_{00})}$  The filter base
\[L(f,\mu)_*\,\, {}^{\mathrm{LM}}\mathfrak F_{00}^{\mathfrak S}\]
converges to $a$.
\begin{proof}
Consider the statement of Lemma \ref{lem:EQUILIM}. First,
according to Lemma \ref{lem:TAGRES} we can formulate
the corresponding statements by only unrestricted envelopes of special kind.
Then according to Lemma \ref {lem:TAGVAL} we can write the statements in terms of
direct images of Lebesgue-McShane filter bases associated to tag systems.
Then we can reformulate the statements in terms of direct images of
the big Lebesgue-McShane filter bases.
\end{proof}
\end{lemma}
\begin{remak}
One can immediately see that $L(f,\mu)_*\,\, {}^{\mathrm{LM}}\mathfrak F^{A}$ is finer than
$ {}^{\mathrm{LM}}\widetilde{\mathfrak F}^{A}_{L,f,\mu}$. On the other hand,
this is not necessary the true in the other directioneven if after taking a smearing ``$\mathrm t$''.
The reason for that is that Lemma \ref{lem:UREXTEND}.b can be used only if we can guarantee that
there are enough ``Cauchy small'' sets around.
\end{remak}
\begin{lemma}\label{lem:TAGSIG}
\halfquote

Suppose that $f:X\rightarrow\mathcal V$ is a function and $X\in\boldsymbol\Sigma\mathfrak S$.
Then, we claim, the following statements are equivalent:

i.) The Lebesgue-McShane integral of $f$ is $a$.

ii.) The extended Lebesgue-McShane integral of $f$ is $a$.

iii.) The filter base
\[L(f,\mu)_*\,\, {}^{\mathrm{LM}}\mathfrak F^{X}\]
converges to $a$.
\begin{proof}
The equivalence of i. and ii. immediately follows from Lemma \ref{lem:OLDCHAR}.
The equivalence of ii. and iii. follows from characterization Lemma
\ref{lem:TAGINT}.$\mathtt{(E_{\Sigma})}$.
\end{proof}
\end{lemma}
\begin{remak}
One aspect of the Lebesgue-McShane filter bases that we can consider integration
as a limit of integration with discrete measures. A formal sum is an actual prescription
how to ``concentrate'' $\mu$ to a measure supported on a finite set.
The effect of $L(f,\mu)$ on a formal sum is nothing else but the integral of $f$
on the the concentrated measure. Here is the title of this section from.
\end{remak}
\begin{remark}
If one wants to define the extended Lebesgue-Shane integral quickly then
after reviewing Subsection \ref{sec:EXTLM}.C
he should start from Subsection \ref{sec:discrete}.C and
then define the extended Lebesgue-McShane integral as in Lemma \ref{lem:TAGINT}.$\mathtt{(E_{0})}$.
\end{remark}

One has to notice that the filter bases ${}^{\mathrm{LM}}\mathfrak F^{\mathfrak S}$,
${}^{\mathrm{LM}}\mathfrak F^{\mathfrak S}_{0}$,
${}^{\mathrm{LM}}\mathfrak F^{\mathfrak S}_{00}$ are not necessarily equivalent,
yet under favorable conditions they induce equiconvergent filters.
One may wonder if we can state this on a formal level.
The answer is affirmative.
~\\

\paragraph{\textbf{E. The integral as a formal construction}}
\begin{defin}\label{def:FRMSIGADD}
Suppose that $\mathfrak S\subset\mathfrak P(X)$ is an interval system.
a.) For each element $x\in X$ and each countable decomposition
\[A=\bigcup^{\updisjoint}_{\lambda\in\Lambda}A_\lambda\]
$(A,A_\lambda\in\mathfrak S)$ we can consider the filter base
\[\{\{\delta_x\cdot A\}\}-\sum^{\mathrm{form}}_{\lambda\in\Lambda}\delta_x\cdot A_\lambda\]
on $\mathsf F(X,\mathfrak S)$.

Let us define $\mathfrak F^{X,\mathfrak S}_{\sigma}$, the associated $\sigma$-additivity
filter base, as the finest common coarsening of the filter bases above \textit{and} the
filter base $\{\{0\}\}$.

Then we can define the $\sigma$-additive completion of $\mathsf F(X,\mathfrak S)$ as
\[\mathsf F_\sigma(X,\mathfrak S)={}_{(\mathfrak F^{X,\mathfrak S}_{\sigma})}\mathsf F(X,\mathfrak S).\]

b.) One can also consider the free abelian groups $\mathsf F(X)$ and $F(\mathfrak S)$.
While $\mathsf F(X)$ does very well with the discrete topology, $F(\mathfrak S)$
can be completed to a $\sigma$-additive group $F_{\sigma}(\mathfrak S)$.
The construction is similar, one should complete according to
a $\sigma$-additivity filter base $\mathfrak F^{\mathfrak S}_{\sigma}$, which is a finest common
coarsening of certain filter bases on $\mathsf F(\mathfrak S)$ plus the filter base $\{\{0\}\}$.
\end{defin}
\begin{lemma}\label{lem:FRMSIGADD}
 Suppose that $\mathfrak S\subset \mathfrak P(X)$ is an interval system.
There is a  unique pairing
\[L_{ \mathfrak F^{X,\mathfrak S}_{\sigma} }:\mathsf F(X)\times \mathsf F_{\sigma}({\mathfrak S})\rightarrow
\mathsf F_\sigma(X,\mathfrak S)\]
such that
\[L_{ \mathfrak F^{X,\mathfrak S}_{\sigma} }(x,\iota_{ \mathfrak F^{\mathfrak S}_{\sigma} }(S))
=\iota_{ \mathfrak F^{X,\mathfrak S}_{\sigma} }(\delta_x\cdot S).\]
This pairing is the crudest, universal pairing to complete, Hausdorff commutative topological groups.

If $\mathfrak J$ is a filter base coarser than $\mathfrak F_\sigma^{X,\mathfrak S}$ then this
factorizes further to a pairing
\[L_{ \mathfrak J }:\mathsf F(X)\times
\mathsf F_{\sigma}({\mathfrak S})\rightarrow
{}_{\mathfrak J}\mathsf F(X,\mathfrak S)\]

\begin{proof}
That follows from the universality of the constructions.
\end{proof}
\end{lemma}
\begin{lemma}\label{def:FRMINTEG}
 Suppose that $\mathfrak S\subset \mathfrak P(X)$ is an interval system.

a.) Then, we claim,
\[C(\mathfrak F^{X,\mathfrak S}_{\sigma}+{}^{\mathrm{LM}}\mathfrak F^{\mathfrak S})^{\mathrm f}=
C(\mathfrak F^{X,\mathfrak S}_{\sigma}+{}^{\mathrm{LM}}\mathfrak F_{0}^{\mathfrak S})^{\mathrm f}=
C(\mathfrak F^{X,\mathfrak S}_{\sigma}+{}^{\mathrm{LM}}\mathfrak F_{00}^{\mathfrak S})^{\mathrm f}.\]
This common filter
\[{}^{\mathrm{LM}}\mathfrak F_{\mathrm{Cauchy}}^{\mathfrak S}\]
is the finest $0$-neighborhood inducing  (cf. Lemma \ref{lem:FORMALCRIT}) filter $\mathfrak J$
on $\mathsf F(X,\mathfrak S)$ which is coarser than $C(\mathfrak F_{\sigma}^{X,\mathfrak S})$
and the tautological function
\[\mathsf x:X\rightarrow \mathsf F(X)\]
\[x\mapsto x\]
is extended Lebesgue-McShane integrable with respect to
the measure
\[\mu_{\mathfrak S}:\mathfrak S\rightarrow\mathsf F_{\sigma}(\mathfrak S)\]
\[S\mapsto \iota_{ \mathfrak F^{\mathfrak S}_{\sigma} }(S)\]
and the formal pairing
\[L_{ \mathfrak J }:\mathsf F(X)\times
\mathsf F_{\sigma}({\mathfrak S})\rightarrow
{}_{\mathfrak J}\mathsf F(X,\mathfrak S).\]

b.) Moreover, we claim,
\[({}^{\mathrm{LM}}\mathfrak F^{\mathfrak S}_{\mathrm{Cauchy}}+
{}^{\mathrm{LM}}\mathfrak F^{\mathfrak S})^{\mathrm f}=
({}^{\mathrm{LM}}\mathfrak F^{\mathfrak S}_{\mathrm{Cauchy}}+
{}^{\mathrm{LM}}\mathfrak F_{0}^{\mathfrak S})^{\mathrm f}=
({}^{\mathrm{LM}}\mathfrak F^{\mathfrak S}_{\mathrm{Cauchy}}+
{}^{\mathrm{LM}}\mathfrak F_{00}^{\mathfrak S})^{\mathrm f},
\]
and this common filter
\[{}^{\mathrm{LM}}\mathfrak F_{\mathrm{full}}^{\mathfrak S}\in
 {}_{{}^{\mathrm{LM}}\mathfrak F_{\mathrm{Cauchy}}^{\mathfrak S}}
 \mathsf F(X,\mathfrak S)\]
is the filter representing the integral of the tautological function $\mathsf x$.
Ie.
\[{}^{\mathrm{LM}}\mathfrak F_{\mathrm{full}}^{\mathfrak S}=\int^{\mathrm{LM}}
 L_{ {}^{\mathrm{LM}}\mathfrak F_{\mathrm{Cauchy}}^{\mathfrak S}  }(\mathsf x,\mu_{\mathfrak S}) .\]

c.) Let $\mathcal  V$, $\mathcal  W$, $\mathcal  Z$ be commutative
topological groups,
$\mu:\mathfrak S\rightarrow\mathcal  W$ be a measure,
$L:\mathcal  V\times\mathcal  W\rightarrow\mathcal  Z$ be a biadditive pairing, which is
continuous in its second variable. Assume that $f$ is extended Lebesgue-McShane integrable
with respect to $L,\mu$. Then, we, claim there exists a unique
continuous homomorphism
\[H_{L,f,\mu}: {}_{{}^{\mathrm{LM}}\mathfrak F_{\mathrm{Cauchy}}^{\mathfrak S}}
 \mathsf F(X,\mathfrak S)\rightarrow \mathcal Z\]
 such that
 \[H_{L,f,\mu}(\iota_{ \mathfrak F^{X,\mathfrak S}_{\sigma} }(\delta_x\cdot A))=L(f(x),\mu(A)).\]
 In this case
 \[H_{L,f,\mu}({}^{\mathrm{LM}}\mathfrak F_{\mathrm{full}}^{\mathfrak S})=
 \int^{\mathrm{LM}} L(f,\mu).\]
\begin{proof}
That follows from the universality of the construction and Lemma \ref{lem:TAGINT}.
\end{proof}
\end{lemma}
\begin{remak}
Humanely speaking: over an interval system $\mathfrak S\subset\mathfrak P(X)$
the most general extended Lebesgue-McShane integrable function is
the tautological function $\mathsf x$; the integral takes values in the commutative topological group
${}_{{}^{\mathrm{LM}}\mathfrak F_{\mathrm{Cauchy}}^{\mathfrak S}}
 \mathsf F(X,\mathfrak S)$; the actual value of the integral is the filter
 ${}^{\mathrm{LM}}\mathfrak F_{\mathrm{full}}^{\mathfrak S}$.
\end{remak}

\begin{remak}
In particular. Lemma \ref{lem:OLDCHAR} shows that  in case of
$X\in\boldsymbol\Sigma\mathfrak S$ there is a universal class of envelopes recovering the integral.
\end{remak}

\begin{remak}
One can continue  discussing the (extended) Lebesgue-McShane integral on the formal level.
\end{remak}
~

\paragraph{\textbf{F. The Lebesgue-Kurzweil-Henstock integral}}
~\\

We can change Definition \ref{def:TAGDIV} slightly:
\begin{defin}\label{def:TAGDIVKH}
Suppose that $\mathfrak S\subset\mathfrak P(X)$ is an interval
system, and $(\mathfrak D,h)$ is a tag system on $\mathfrak S$.

A tagged Lebesgue-Kurzweil-Henstock division
associated to $(\mathfrak D,h)$ is a pair
$(\mathfrak E, c)$ such that
\begin{itemize}
\item[\texttt{(D1)}]
 The set $\mathfrak E$ is countable set of pairwise disjoint nonempty elements of $\mathfrak S$
\item[\texttt{(D2')}] $c:\mathfrak E\rightarrow \bigcup\mathfrak D$ such that
 for each $E\in\mathfrak E$
\[c(E)\in E\subset h(c(E))\]
 holds.
\end{itemize}
\end{defin}
\begin{point}
Then the corresponding version Lemma \ref{lem:TAGREFIN} holds with the same proof.
Then, in Definition \ref{def:TAGFILTER} the terms ``Lebesgue-McShane'' can be
substituted by ``Lebesgue-Kurzweil-Henstock''.

That way one obtains the filter bases
\[{}^{\mathrm{LKH}}\mathfrak F^A\]
and the corresponding derivatives.
\end{point}
\begin{defin}
\halfquote

Suppose that $f:X\rightarrow\mathcal V$ is a function. Then we say that the extended
Lebesgue-Kurzweil-Henstock integral of $f$ is $a$ if
\[L(f,\mu)_*\,\, {}^{\mathrm{LKH}}\mathfrak F_0^{\mathfrak S}\]
converges to $a$.

In the case when $f$ vanishes outside of a set $A\in\boldsymbol\Sigma\mathfrak S$
we also say that it is the ordinary Lebesgue-Kurzweil-Henstock integral.
\end{defin}
\begin{remark}
One can immediately see that the Lebesgue-Kurzweil-Henstock filter bases are finer than the
corresponding Lebesgue-McShane filter bases, so the Lebesgue-Kurzweil-Henstock integral
is, in general,  more effective than the Lebesgue-McShane integral.

That, however, very much depends on the situation. Over $\{0,1\}^{\mathbb N}$, with respect to the
dyadic measure, the Lebesgue-Kurzweil-Henstock integral presents a dramatic improvement
(it is very much as the ordinary Kurzweil-Henstock integral). On the other hand, over $[0,1]$,
with respect to the
interval measure, there is no significant improvement.
In particular, the Lebesgue-Kurzweil-Henstock will \textit{not} yield
the ordinary Kurzweil-Henstock integral.
This is at least surprising, if
we consider the fact that these measure spaces above are essentially the same.

In Section \ref{sec:LKH}, however, we will see  how to extend the benefits of
Lebesgue-Kurzweil-Henstock divisions to the interval measure.
\end{remark}

Nevertheless, the filter base ${}^{\mathrm{LKH}}\mathfrak F_0^{\mathfrak S}$ is already sufficiently
fine to make the following comment here:
\begin{remark}
We can use formal sums like in ${}^{\mathrm{LKH}}\mathfrak F_0^{\mathfrak S}$
to define the integral of noncommutative semigroup valued functions. In this case
it is reasonable to assume that $X$ has a preordering $\precsim$ and during
taking the direct image by $L(f,\mu)$ we consider only evaluations
\[L(f(x_1),\mu(E_1))L(f(x_2),\mu(E_2))\ldots L(f(x_n),\mu(E_n))\]
when $i\leq j$ implies $x_i \not\succsim x_j$. If it is possible then we might also ask
 for the sets $E_j$ to be compatible with the preordering.
Here the linearity of the measure pairing, and a lot of other properties,
will be lost but the direct image makes sense.

This setting has a limited value but it is realistic to ask about integrable, say,  $\mathbb Z$-valued
functions  on $\mathbb R$
with respect to a dynamical system, which can be thought as a noncommutative measure on intervals.
\end{remark}
\begin{remark}
If $\mathfrak S$ is cofinal in $\boldsymbol\Sigma_0\mathfrak S$ with respect to $\subset$
then one can give an extra boost to convergence by defining the ``improper integral'' filter bases
\[{}^{\mathrm{LM}}\mathfrak F_{\mathrm{imp}}^{\mathfrak S}=\dire_{A\in \mathfrak S}^{\subset}
{}^{\mathrm{LM}}\mathfrak F^A\quad\text{and}\quad
{}^{\mathrm{LKH}}\mathfrak F_{\mathrm{imp}}^{\mathfrak S}=\dire_{A\in \mathfrak S}^{\subset}
{}^{\mathrm{LKH}}\mathfrak F^A.\]
\end{remark}
~

\paragraph{\textbf{G. Measurable functions}}~\\

Unrestricted envelopes allow us to define an important function class:
\begin{defin}\label{def:URMEAS}
\halfquote
Suppose that $f:X\rightarrow\mathcal  V$ is a function.

i.) We say that $f$ is $L,\mu$-measurable
if $f$ vanishes outside of a set $A\in\boldsymbol\Sigma\mathfrak S$
and for all set $A\subset B\in\boldsymbol\Sigma\mathfrak S$
and for each neighborhood $0\in\mathcal  T\subset \mathcal Z$
there is an unrestricted envelope $\mathcal  H$ supported on $B$ which
approximates $f$ and
\[\intp L(\mathcal  H-\mathcal  H,\mu)\subset \mathcal  T.\]

ii.) The function $f$ is extended measurable if for each $E\in\mathfrak S$ the function
$f|_E$ is measurable.
\end{defin}

\begin{remak}
So the difference is that for integrable functions we demand the ``Cauchy property'' with envelopes,
while in for measure functions we allow the ``Cauchy property'' with unrestricted envelopes
\end{remak}
\begin{lemma}\label{lem:UREMEAS}
\quarterquote

Then , we claim, a set S is extended $\mu$-measurable if and only if
its characteristic function $\chi_S$ is an extended $\mathtt m, \mu$-measurable function.
\begin{proof}
Suppose that $E\in\mathfrak S$ and
Let $\mathcal H$ be an unrestricted envelope approximating $\chi_S|_E=\chi_{E\cap S}$ such that
the support of $\mathcal H$ contains $E$.
Then we can pass to the the unrestricted envelope
\[\mathcal H'=((\{0,1\}\times E)\cap\mathcal H)\cup(\{0\}\times (X\setminus E)).\]
By the usual methods there is an exact decomposition $\mathfrak D$ of $E$ such that
for every $D\in\mathfrak D$
\[(\mathcal H')^D\neq0.\]
Using Lemma \ref{lem:MEAVAR} we can find a finite subset $K\subset\mathfrak D$.
such that
\[\svar(\mathtt m,1,\bigcup(\mathfrak D\setminus K),\mu)\subset \mathcal T.\]
Then
\[\tilde{\mathcal H}=\mathcal H'-\svar(1,\bigcup(\mathfrak D\setminus K) )
+\svar(1,\bigcup(\mathfrak D\setminus K) )\]
will be  an envelope approximating $\chi_S|_E=\chi_{E\cap S}$
while with $L=\mathtt m$
\[\intp L(\tilde{\mathcal H}-\tilde{\mathcal H} ,\mu)\subset
\overline{\intp L(\mathcal  H-\mathcal  H,\mu)+\mathcal T-\mathcal T+\mathcal T-\mathcal T }
\subset \mathcal T-\mathcal T+\mathcal T-\mathcal T+\mathcal T-\mathcal T.\]
Then, according to the Cauchy criterium, integrability of $S\cap E$ follows.
Ie. extended measurability in new sense implies extended measurabiliy in old sense.
The other direction is obvious.
\end{proof}
\end{lemma}
\newpage
\section{Locally compact integrals (on base)} \label{sec:LKH}

\paragraph{\textbf{A. Locally compact interval systems}}
\begin{defin}
An interval system $\mathfrak S$ is locally compact if any countable disjoint decomposition
\[A=\bigcup^{\updisjoint}_{\lambda\in\Lambda}A_\lambda\]
in $\mathfrak S$ is finite, ie. with the exception of finitely many terms all components must
be the empty set $\emptyset$.
\end{defin}
\begin{example} Coin-tossing, ie. the cylinder sets of  $\{0,1\}^{\mathbb N}$ form a
locally compact interval system. More generally, any set system, such that
\begin{itemize}
\item[i.)] it is a finitely rooted, finitely branching forest,
\item[ii.)] descending chains have non-empty intersection,
\end{itemize}
forms a locally compact interval system. That also applies to the
$(\cup,\cap,\setminus)$-constructible closure of such a system.
\end{example}
\begin{lemma} Suppose that $\mathcal  V$ is a commutative topological group,
and $\mathfrak S\subset\mathfrak P(X)$ is a locally compact interval system.

Then, we claim, any finitely additive function $\theta:\mathfrak S\rightarrow\mathcal  V$
will be a measure, regardless the topology on $\mathcal  V$.
\begin{proof}
All sum reduces to finite sums if we disregard the $0$'s.
\end{proof}
\end{lemma}
\begin{remak}
Having a locally compact interval system has nice consequences with respect to integration.

a.) In Definition \ref{def:TAGFILTER}.b, if we have a tag system $(\mathfrak E,c)$ such that
$\bigcup \mathfrak E \in\boldsymbol\Sigma_0\mathfrak S$ then the associated filter base
$\mathfrak F_{(\mathfrak E,c)}$ can be substituted by the equivalent filter base
\[\dot{\mathfrak F}_{(\mathfrak E,c)}=
\biggl\{\biggl\{\sum_{E\in \mathfrak E}\delta_{c(E)}\cdot E\biggr\}\biggr\}.\]

As a consequence, for $A\in\boldsymbol\Sigma_0\mathfrak S$ the filter base
${}^{\mathrm{LM}}\mathfrak F^A$ can be substituted by an equivalent filter
\[{}^{\mathrm{LM}}\dot{\mathfrak F}^A,\]
which is a filter base on formal Riemann sums on $A$. In particular,
\[ L(f,\mu)_*{}^{\mathrm{LM}}\dot{\mathfrak F}^A\]
will be a filter on Riemann sums.

b.)   It yields
\[C(\mathfrak F^{X,\mathfrak S}_\sigma)=\{G\}^{\mathrm f},\]
where $G$ is the subgroup of $F(X,\mathfrak S)$ generated by elements
\[\delta_x\cdot A-\sum_{\lambda\in\Lambda}\delta_x\cdot A_\lambda,\]
where $x\in X$ and $A=\bigcup_{\lambda\in\Lambda}A_\lambda$ is a finite decomposition.
That way $\mathsf F_\sigma(X,\mathfrak S)$ is just a quotient group of $\mathsf F(X,\mathfrak S)$.
Similar comments apply to $\mathsf F_\sigma(\mathfrak S)$.
\end{remak}
\begin{remak}
The examples above seem to be of limited use. Notably, even the interval system $\mathfrak I$
of the finite intervals on $\mathbb R$ fails to be locally compact.
\end{remak}
~

\paragraph{\textbf{B. Stone completions}}
~\\

However, one can produce locally compact interval systems in great numbers:
\begin{conven}
In this section, for a set system $\mathfrak W$ we will use the term
``filter base/ filter $\mathfrak F$ on $\mathfrak W$''.
This will not be understood in the sense that $\mathfrak F\subset \mathfrak P(\mathfrak W)$
plus $\mathfrak F$ is a filter base / filter base upward
closed to the containment relation in $\mathfrak P(\mathfrak W)$; but it will be understood in the sense that
$\mathfrak F\subset \mathfrak W$ plus $\mathfrak F$ is a filter base / filter base upward
closed to the containment relation in $\mathfrak W$.
(Cf. \cite{grat} for a quite general setting.)
\end{conven}
\begin{remin} If $\mathfrak S\subset\mathfrak P(X)$ is a set system then we may
consider its the associated ring $\boldsymbol\Sigma_{\mathrm c}\mathfrak S$.
Then, we may consider its Stone completion. More specifically:

i.) We define $ \mathfrak X_{\mathfrak S}$ to be the set of maximal filters on $\boldsymbol\Sigma_{\mathrm c}\mathfrak S$.

ii.) For each $A\in\boldsymbol\Sigma_{\mathrm c}\mathfrak S$ we define
\[\beta A=\{\mathbf x\in \mathfrak X_{\mathfrak S}\,:\,A\in\mathbf x\}.\]

iii.) We define $\beta\mathfrak S$ ($\beta\boldsymbol\Sigma_{\mathrm c}\mathfrak S$)
as the image of the set $\mathfrak S$ ($\boldsymbol\Sigma_{\mathrm c}\mathfrak S$) under
the $\beta$.
\end{remin}
\begin{lemma}\label{lem:CS1}
 Suppose that $\mathfrak S\subset\mathfrak P(X)$ is a set system.
The we claim:

a.) For $\mathbf x\in \mathfrak X_{\mathfrak S}$ and $A\in\mathfrak S$
\[A\notin\mathbf x\quad\Leftrightarrow\quad \exists B\in\mathbf  x\text{ s. t. }B\cap A=\emptyset;\]
and if $A=\bigcup_{j\in J}A_j,\, J \text{ is finite},\, A_j\in \boldsymbol\Sigma_{\mathrm c}\mathfrak S$ then
\[A\in\mathbf x \quad\Leftrightarrow\quad \exists j\in J\text{ s. t. }A_j\in\mathbf x .\]

b.) The set system $\beta\boldsymbol\Sigma_{\mathrm c}\mathfrak S$ is algebraically isomorphic
to $\boldsymbol\Sigma_{\mathrm c}\mathfrak S$ through $\beta$
(ie. isomorphic with respect to $\subset$ and finite $\cup,\cap,\setminus$).

c.) For $A,A_\lambda\in\boldsymbol\Sigma_{\mathrm c}\mathfrak S$, $(\lambda\in\Lambda)$,
$\lambda$ is countable, we have
\[A=\bigcup_{\lambda\in\lambda} \beta A_\lambda\quad\Leftrightarrow\quad\exists\Xi\subset
\Lambda\text{ finite, s. t. }
A=\bigcup_{\lambda\in\Xi} \beta A_\lambda=  \bigcup_{\lambda\in\lambda} \beta A_\lambda.\]
In particular, $\beta\boldsymbol\Sigma_{\mathrm c}\mathfrak S$ is a locally compact interval system.
\qed
\end{lemma}
\begin{defin}
An interval system $\mathfrak S$ is locally finite if for each $A,B\in\mathfrak S$, $A\subset B$
the set $B\setminus A$ will be a finite disjoint union of elements of $\mathfrak S$.
\end{defin}
\begin{example}
The interval system $\mathfrak I$ of the finite intervals on $\mathbb R$.
\end{example}
\begin{lemma}\label{lem:CSLF}
If the interval system $\mathfrak S$ is locally finite then every element of its
constructible $(\cap,\cup,\setminus)$-closure is a finite disjoint union of elements of $\mathfrak S$.

In particular, $\boldsymbol\Sigma_{00}\mathfrak S=\boldsymbol\Sigma_{0}\mathfrak S=
\boldsymbol\Sigma_{\mathrm c}\mathfrak S$.
\begin{proof}
One can prove this by induction.
\end{proof}
\end{lemma}
\begin{remark}
In particular, step-functions and simple step-functions over locally finite intervals are the same.
That makes some statements like \ref{lem:APPREFIN} less important.
\end{remark}
\begin{lemma}\label{lem:CS2}
Let $\mathfrak S\subset \mathfrak P(X)$ be a locally finite interval system.
Then we claim

a.) By cofinality, the set of maximal filters on $\mathfrak S$ corresponds to the
set of maximal filters on $\boldsymbol\Sigma_{\mathrm c}\mathfrak S$.

Ie., the Stone completion $\mathfrak X_{\mathfrak S}$
can be defined as the set of maximal filters on $\mathfrak S$.

b.) For $\mathbf x\in \mathfrak X_{\mathfrak S}$ and $A\in\mathfrak S$
\[A\notin\mathbf  x\quad\Leftrightarrow\quad \exists B\in\mathbf  x\text{ s. t. }
B\cap A=\emptyset;\]
and if $A=\bigcup_{j\in J}A_j,\, J \text{ is finite},\, A_j\in\mathfrak S$ then
\[A\in\mathbf x \quad\Leftrightarrow\quad \exists j\in J\text{ s. t. } A_j\in\mathbf x .\]

c.) $\beta\mathfrak S$ is  a locally finite interval system which is algebraically isomorphic to $\mathfrak S$.
\begin{proof}
a.) It is sufficient to prove that the elements of $\mathfrak S$ are cofinal in the
maximal filters on $\boldsymbol\Sigma_{\mathrm c}\mathfrak S$.
This follows from Lemma \ref{lem:CSLF} and Lemma \ref{lem:CS1}.a.

The points b., c.) follow immediately.
\end{proof}
\end{lemma}
\begin{remak}
Lemma \ref{lem:CS2}.a is quite useful when one actually computes the completion.
\end{remak}
One idea to extend integration to finitely additive set functions is that
instead of $\mathfrak S$ we use the Stone completion $\beta\mathfrak S$,
where $\sigma$-additivity is trivially guaranteed.
An obstacle to that approach that we have to extend a function $f:X\rightarrow\mathcal  V$
to a function $f':\mathfrak X_{\mathfrak S}\rightarrow \mathcal  V$.
For that reason we need a map $\mathfrak X_{\mathfrak S}\rightarrow X$.
\begin{defin}\label{def:CONSIS}
Let $\mathfrak S\subset \mathfrak P(X)$ be a set system.
Assume that $X$ is a topological space.

We say that the topology of $X$ is compatible to $\mathfrak S$ if
every element of $\mathfrak X_{\mathfrak S}$ is convergent to one unique element of $X$ (as a filter base).
Ie. there is a well-defined limit map
\[\lim: \mathfrak X_{\mathfrak S}\rightarrow X.\]
\end{defin}
\begin{lemma}\label{lem:CONSIS1}
Let $\mathfrak S\subset \mathfrak P(X)$ be a set system.
Suppose that $X$ is a topological space, which is compatible to $\mathfrak S$.
Assume that $S\in\mathfrak S$ and $\mathbf x\in \mathfrak X_{\mathfrak S}$.

Then, we claim, $\lim \mathbf x\in S^\circ$ (interior) implies $\mathbf x\in\beta S$.
\begin{proof}
The convergence of $\mathbf x$ implies that there is an element $A\in \mathbf x$ such that
 $A\subset S^\circ$.
That, however implies $S\in \mathbf x$.
\end{proof}
\end{lemma}
\begin{lemma}\label{lem:CONSIS2}
Let $\mathfrak S\subset \mathfrak P(X)$ be a set system.
Assume that $X$ is a topological space, which is compatible to $\mathfrak S$.
Then for $A\in\boldsymbol\Sigma_{\mathrm c}\mathfrak S$
\[\lim (\beta A)\subset\overline A.\]
(Here we mean the image set under $\lim$.)
\begin{proof}
This follows from that $\mathbf x\in\beta A$ implies $A\in\mathbf x$, which implies
that the limit must be in $\overline A$.
\end{proof}
\end{lemma}
\begin{lemma}\label{lem:CONSIS3}
Let $\mathfrak S\subset \mathfrak P(X)$ be an locally finite interval system.
Assume that the following conditions hold:

i.) For all $A\in\mathfrak S$ the space $\overline A$ is a compact Hausdorff space.

ii.) For each open  cover $\mathcal  Q$ of $\overline A$ there is a finite decomposition $\mathfrak B$ of $A$
such that each element of $B\in\mathfrak B$ is contained in an element $Q $ of  $ \mathcal  Q$.
(Remark: in fact $\overline B\subset Q$ can be assumed.)

Then, we claim, the topology of $X$ is compatible to $\mathfrak S$.

Moreover,
\[\lim(\beta A)=\overline A\]
holds.
\begin{proof}
Assume that $A\in\mathbf x$. Condition i. implies that $x$ has a limit point $y$ in $\overline A$,
 the intersection of
all sets $\overline B$, where $B\in\mathbf x$. The Hausdorff property and ii.) implies that
$\mathbf x$ converges to $y$.

The additional property can be proven as follows:
For each element $y\in\overline A$ and each finite decomposition $\mathfrak B$ of $A$ we can consider
the set
\[S_{y,\mathfrak B}=\bigcup_{B\in\mathfrak B,\, y\in\overline B}B. \]
Taking common refinements of finite decompositions we see that for a fixed $y$ these sets
$S_{y,\mathfrak B}$ form a
filter base. This filter base a has $y$ as a limit point and according to ii.) the filter base
converges to it. We can extend the convergent filter base to a maximal one.
\end{proof}
\end{lemma}
\begin{example}\label{exa:SCOM1}
Consider the locally finite interval system $\mathfrak I$, which is
the  system of finite intervals $\mathfrak I$ on $\mathbb R$.

Let us define the following filters
 $[a]$, $[a+]$, $[a-]$ ($a\in\mathbb R$), ``contained element'', ``right limit'',
``left limit'' on $\mathfrak I$: Let
\[I\in[a]\quad \Leftrightarrow\quad a\in I,\qquad
I\in[a+]\quad \Leftrightarrow\quad a\in\overline{ I\cap(a,+\infty)},\]
\[I\in[a-]\quad \Leftrightarrow\quad a\in\overline{ I\cap(-\infty,a)}.\]
One can see that these are the maximal filters on $\mathfrak I$, hence they form
the Stone completion of $\mathfrak I$.

Furthermore, $\mathfrak I$ is compatible with the ordinary topology of $\mathbb R$, and
 the limit map is given by
\[\lim\, [a]=a,\quad\lim\, [a+]=a,\quad\lim\, [a-]=a.\]
\end{example}
\begin{example}\label{exa:SCOM2}
Let $V$ be a finite dimensional real vector space.
Let a $\mathbb SV$ be its radial space, ie. the space which contains the
rays $\lceil v\rceil$ belonging to vectors $0\neq v\in V$.
This is a sphere. If $f\neq 0$ is a functional on $V$ the we can consider the
corresponding open and closed half-spaces
\[\{\lceil v\rceil\,:\,v\neq 0,\,f(v)>0\}
\quad\text{and}\quad \{\lceil v\rceil\,:\,v\neq 0,\,f(v)\geq 0\}\]

We can consider the system of simplices $\mathfrak O_{\Delta}$, system of convex polyhedrons
 $\mathfrak O_{\mathrm c\mathrm o}$,
and the system of polyhedrons $\mathfrak O_{\mathrm c}$.
Here $\mathfrak O_{\mathrm c}$ is the $(\cup,\cap,\setminus)$-constructible closure
of the open and closed half-spaces; $\mathfrak O_{\mathrm c\mathrm o}$ contains the
finite intersections of open and closed half-spaces; while $\mathfrak O_{\Delta}$
contains those convex polyhedrons which are strictly contained in an open half-space,
simplices, and open in the non-degenerate directions.
(But one can come with alternative versions.)

One can see that each of them is an locally finite interval system, in fact their
common constructible closure is $\mathfrak O_{\mathrm c}$, so their Stone completion is
the same.

For the sake of convenience, let us fix a positive scalar product $g$ on $V$.
That way we represent rays as unit vectors.
One can notice that for two orthonormal vectors $v,w$ and angle $\alpha$
we can associate the unit vector
\[\Lambda(v,\alpha,w)=v\cos\alpha+w\sin\alpha.\]
More generally, for orthonormal vectors $v_0,\ldots,v_n$ and angles $\alpha_1,\ldots,\alpha_n$
we can consider the vector
\[\Lambda_{v_0,\ldots,v_n}(\alpha_1,\ldots,\alpha_n)=\Lambda(v_0,\alpha_1,\Lambda(v_1,\alpha_2,\ldots
\Lambda(v_{n-1},\alpha_n,v_n))\ldots ).\]
(For $n=0$ it is just $\Lambda_{v_0}=v_0$.)

Then, for orthonormal vectors $v_0,\ldots,v_n$ (where $n=0,\ldots,\dim V$) we can define the filters
\[[v_0,\ldots,v_n]\]
on either of the set systems above by the prescription
\[I\in[v_0,\ldots,v_n]\quad \Leftrightarrow\quad\exists\, 0<\alpha_1,\ldots,\alpha_n
\text{ s. t. } \Lambda_{v_0,\ldots,v_n}((0,\alpha_1)\times \ldots\times(0,\alpha_n))\subset I.\]
One can see that these filters are all the possible maximal filters.

Furthermore, the natural topology on $\mathbb SV$ is compatible and
\[\lim\,[v_0,\ldots,v_n]=\lceil v_0\rceil.\]
(The Stone completion does not depend on the scalar product of course, but
one has to look for a  convenient way to describe it.)
\end{example}
\begin{example}\label{exa:SCOM3}
Consider the interval $[0,1]$ and the interval system $\mathfrak I\|_{[0,1]}$ on it.
One can see that the Stone completion $\mathfrak X_{\mathfrak I\|_{[0,1]}}$of $\mathfrak I\|_{[0,1]}$ contains the
filters $[a+]$ $(a\in[0,1))$, $[a]$ $(a\in[0,1])$, $[a-]$ $(a\in(0,1])$, similarly to
Example \ref{exa:SCOM1}.

Let $\Omega$ be an arbitrary set and consider the set $[0,1]^\Omega$.
Let $(\mathfrak I\|_{[0,1]})^{\Omega,[0,1]}$ be the interval system of the
cylinder sets on $[0,1]^\Omega$, ie. the set of sets
\[C{}^{I_1}_{\omega_1}\ldots{}^{I_n}_{\omega_n}=\{x\in[0,1]^\Omega\,:\, x(\omega_k)\in I_k\},\]
where $\omega_k\in\Omega$ are different from each other, $I_k\in \mathfrak I\|_{[0,1]}$.

Then, to each element $y\in (\mathfrak X_{\mathfrak I\|_{[0,1]}})^\Omega$ we can assign a
filter
\[[y]\]
on $(\mathfrak I\|_{[0,1]})^{\Omega,[0,1]}$ by the prescription
\[C\in [y]\quad\Leftrightarrow\quad \exists \{\omega_1,\ldots,\omega_n\}\subset \Omega, I_k\in y(\omega_k)
\text{ s. t. }
C{}^{I_1}_{\omega_1}\ldots{}^{I_n}_{\omega_n}\subset \Omega.\]
One can see that these are all the maximal filters on $(\mathfrak I\|_{[0,1]})^{\Omega,[0,1]}$.
\end{example}

\paragraph{\textbf{C. The McShane integral}}
~\\

One idea to generalize  the notion of the (extended) Lebesgue-McShane integrals
is to get rid of the condition of $\sigma$-additivity of measures.

\begin{defin}
Let $\mathcal  V,\mathcal  W,\mathcal  Z$ be commutative topological groups,
$X$ be a topological space, $\mathfrak S\subset\mathfrak P(X)$ be a
locally finite interval system compatible
to $X$, $\kappa:\mathfrak S\rightarrow\mathcal  V$ be a finitely additive
function, $L:\mathcal  V\times\mathcal  W\rightarrow\mathcal  Z$ be a biadditive pairing,
which is continuous in its second variable.

Suppose that $f:X\rightarrow\mathcal  V$ is a function.
Then we define its McShane integral
\[\int^{\mathrm M} L(f,\kappa)\]
as the extended Lebesgue-McShane integral
\[\int^{\mathrm{LM}} L(\lim{}^*f,\beta_*\kappa),\]
if it exists. Here $\lim{}^*f=f\circ\lim:\mathfrak X_{\mathfrak S}\rightarrow\mathcal V$
and
$\beta_*\kappa=\kappa\circ\beta^{-1}:\beta\mathfrak S\rightarrow\mathcal W$.
\end{defin}
\begin{remark}
Ie. instead of the pairing \[L(f,\kappa):X\times \mathfrak S\rightarrow\mathcal  Z\] we use
\[L(f,\kappa)\circ(\lim,\beta^{-1}):\mathfrak X_{\mathfrak S}\times \beta\mathfrak S\rightarrow\mathcal  Z.\]
\end{remark}

Then one can write down the associated McShane filter
bases obtained from the Lebesgue-MacShane filter bases.
For the sake of simplicity we consider only the one of finite kind.
\begin{defin}\label{def:MTAGFIN}
Suppose that $\mathfrak S\subset\mathfrak P(X)$ is a locally finite interval system,
$X$ is a topological space compatible to $\mathfrak S$.
Assume that $A\in\boldsymbol\Sigma_0\mathfrak S$.
A topological tag system supported on  $A$
is a pair $(\mathfrak D,h)$ such that
\begin{itemize}
\item[\texttt{(HT1)}] $\mathfrak D$ is a countable subset of $\mathfrak S$ such that
$\bigcup\mathfrak D=A$, moreover, $A\in\boldsymbol\Sigma_0\mathfrak D$;
\item[\texttt{(HT2)}]
$h:\beta A\rightarrow\mathfrak D$ is a function such that \[\mathbf x\in \beta h(\mathbf x)\]
for all $\mathbf x\in \beta A$.
\end{itemize}
\end{defin}

\begin{defin}\label{def:MTAGDIV}
Suppose that $\mathfrak S\subset\mathfrak P(X)$ is a locally finite interval
system, $X$ is a topological space compatible to $\mathfrak S$.
Assume that  $(\mathfrak D,h)$ is a tag system on $\mathfrak S$.

A tagged McShane division  associated to $(\mathfrak D,h)$ is a pair
$(\mathfrak E, c)$ such that
\begin{itemize}
\item[\texttt{(DT1)}]
 $\mathfrak E$ is a finite set of pairwise disjoint nonempty elements of $\mathfrak S$;
\item[\texttt{(DT2)}] $c:\mathfrak E\rightarrow \beta A$ is a function such that
\[E\subset h(c(E))\]
for each $E\in\mathfrak E$.
\end{itemize}
\end{defin}
\begin{lemma}\label{lem:MTAGREFIN}
Suppose that $\mathfrak S\subset\mathfrak P(X)$ is a locally finite interval system,
$X$ is a topological space compatible to $\mathfrak S$.
Then, we claim:

a.) Any tag system $(\mathfrak D,h)$ allows a tagged  McShane  division $(\mathfrak E,c)$.

b.) If the tag system $(\mathfrak D',h')$ is finer than $(\mathfrak D,h)$ and
$(\mathfrak E,c)$ is a tagged McShane division associated to $(\mathfrak D',h')$, then
$(\mathfrak E,c)$ is also associated to $(\mathfrak D,h)$.
\begin{proof}
This is just Lemma \ref{lem:TAGREFIN} under special circumstances.
\end{proof}
\end{lemma}
\begin{defin}\label{def:MTAGFILTER}
Suppose that $\mathfrak S\subset\mathfrak P(X)$ is a locally finite interval system,
$X$ is a topological space compatible to $\mathfrak S$.

a.) If $(\mathfrak E,c)$ is a tagged McShane division then its associated filter base is
as the  set
\[{}^{\mathrm t}\mathfrak F_{(\mathfrak E,c)}=
\biggl\{\biggl\{\sum_{E\in\mathfrak E}\delta_{\lim c(E)}\cdot E\biggr\}\biggr\}.\]

b.) If $(\mathfrak D,h)$ is a topological tag system the its associated McShane filter base is
\[{}^{\mathrm{M}}\mathfrak F^{(\mathfrak D,h)}=\prod_{
\substack{(\mathfrak E,h)\text{ McShane tagged }
\\\text{ division associated to }(\mathfrak D,h)}}{}^{\mathrm t}\mathfrak F_{(\mathfrak E,c)}\]
ie. the coarsest filter base generated by the tagged divisions associated to $(\mathfrak D,h)$.

c.) If $A\in\boldsymbol\Sigma_0\mathfrak S$ then its associated  McShane filter base is
\[{}^{\mathrm{M}}\mathfrak F^A=\coprod_{\substack{ (\mathfrak D,h) \text{ topological tag}\\
\text{ system supported on  }A }}{}^{\mathrm{M}}\mathfrak F^{(\mathfrak D,h)}, \]
ie. the finest filter base generated by tag systems on $A$.
(Cf. Lemma \ref{lem:MTAGREFIN}.b for the applicability of the definition.)

d.) The associated McShane filter base of finite kind is defined as
\[{}^{\mathrm{M}}\mathfrak F_0^{\mathfrak S}=\dire_{A\in \boldsymbol\Sigma_0\mathfrak S}^{\subset}
{}^{\mathrm{M}}\mathfrak F^A,\]
ie. the limit of the filter bases ${}^{\mathrm{M}}\mathfrak F^A$ as $A\in\boldsymbol\Sigma_0\mathfrak S$ increases.
\end{defin}
\begin{remak}
Despite of the multiple steps, the Lebesgue-McShane filter bases are just filter bases on finite formal sums.
\end{remak}
The following statements make connection between the Lebesgue-McShane and McShane integrability.
\begin{lemma}\label{lem:DEAT1}
\halfquote
Let $\mathfrak S$ be locally finite and
$X$ be a topological space compatible to $\mathfrak S$.

Assume that $A\in\boldsymbol\Sigma_0\mathfrak S$ and $f:X\rightarrow\mathcal V$ is a function.

a.) Suppose that
 for each $S\in\mathfrak S$, $S\subset A$ the set
\[(A\cap\lim(\beta S)\setminus S)\cup(S\setminus\lim(\beta S) )\]
is weakly $L,\mu$-negligible.
 Also assume that $\mathcal  A'$ is  an envelope over $\beta\mathfrak S$ approximating
 $\lim{}^* f|_{\beta A}$.

Then we claim, for each neighborhood $0\in\mathcal  T\subset\mathcal  Z$
there exist  an envelope $\tilde{\mathcal  A}$ over $\mathfrak S$ approximating $f|_A$ such that
\[\intp L(\tilde{\mathcal  A},\mu)\subset \intp L(\mathcal  A',\beta_* \mu)+\mathcal  T.\]

b.)
Suppose that for each $S\in\mathfrak S$, $S\subset A$ the set
\[((\beta S)\setminus \lim{}^{-1}(S))\cup
(\beta A\cap\lim{}^{-1}(S)\setminus (\beta S))\]
is weakly $L,\beta_*\mu$-negligible.

Then, we claim, for each envelope $\mathcal  A$ over $\mathfrak S$
approximating $f|_A$, and neighborhood
$0\in\mathcal  T\subset \mathcal Z$
there exists  an envelope $\tilde {\mathcal  A}$ approximating $\lim^* f|_{\beta A}$ such that
\[\intp L(\hat{\mathcal  A},\beta_*\mu)\subset \intp L(\mathcal  A,\mu)+\mathcal  T.\]

\begin{proof}
a.) 1. Let
\[\mathcal  A'=\sum_{j\in J} c_j\chi_{\beta E_j}+\mathcal  C'\]
be  an envelope over $\beta\mathfrak S$ approximating $\lim{}^* f|_{\beta A}$, such that
$\mathcal C'$ is a pointed envelope.

For $\mathcal A'$ and $\mathcal C'$ we take intersection by
$\mathcal V\times\beta A$, and then we can extend them back to $X$.
(Being the the interval system locally finite the step-functions will not get too fragmented.)
This and the local finiteness of $\mathfrak S$ allow us the assumptions
\[A=\bigcup_{j\in J}^{\updisjoint}E_j\]
and
\[\mathcal C'\subset (\mathcal V\times \beta A)\cup \{0\}\times \mathfrak X_{\mathfrak S}.\]
Then there is a generator system   $\beta\mathfrak D=\{\beta D\,:\,D\in\mathfrak D\}$
for the envelope $\mathcal  C'$, and we can assume that
$\bigcup\mathfrak D=A$, $\mathfrak D$ is a forest,
and $\mathfrak D$ is divided by $\mathfrak E=\{E_j\,:\,j\in J\}$.
This latter one also implies that $\beta\mathfrak D$ is a generator system for $\mathcal A'$.

Then let us define the pointed envelope
\[\mathcal  C=\bigcup_{D\in\mathfrak D}\mathcal  C^{\beta D}\times  D
\,\,\cup\,\,(\{0\}\times X),\]
and the envelope
\[\mathcal  A=\sum_{j\in J} c_j\chi_{E_j}+\mathcal  C.\]

One can see that
\[\mathcal  A=\bigcup_{D\in\mathfrak D}\mathcal  A^{\beta D}\times  D
\,\,\cup\,\,(\{0\}\times (X\setminus A)).\]

2. Here we prove  that
\[\intp L(\mathcal  A,\mu)\subset \intp L(\mathcal  A',\beta_* \mu).\]
Now,
\[\intp L\biggl(\sum_{j\in J} c_j\chi_{\beta E_j},\beta^*\mu\biggr)=\sum_{j\in J} L(c_j,\mu( E_j))
=\intp L\biggl(\sum_{j\in J} c_j\chi_{E_j},\mu\biggr)\]
shows that it is enough to prove that
\[\intp L(\mathcal  C,\mu)\subset \intp L(\mathcal  C',\beta_*\mu) \tag{v1}.\]
Now, assume that
\[\sum_{j\in\tilde J} \tilde c_j\chi_{\beta \tilde E_j}\subset \mathcal  C\]
is a finite disjoint sum. Then, according to the definition of $\mathcal  C$ we find that
\[\tilde E_j\subset \bigcup_{D\in\mathfrak D,\, c_j\in (\mathcal A')^{\beta D}} D. \]

From this it follows that $E_j'$ is decomposed by a countable set subset $\mathfrak D_j$ of
$\mathfrak S$ such that for each $D'\in\mathfrak D_j$ there is set $D\in\mathfrak D$ such that
\[D'\subset D\quad\text{ and }\quad c_j\subset (\mathcal  A')^{\beta D}.\]
Then
\[\sum_{j\in J'} L(c_j',\mu( E_j'))=
\sum_{j\in J'}\sum^{\ovline}_{D\in\mathfrak D_j} L(c_j,\mu(D))=
\sum_{j\in J'\, D\in\mathfrak D_j}^{\ovline} L(c_j,\mu(D)).\]
Now, every finite partial sum of right side belongs to $ \intp L(\mathcal  C',\beta_*\mu)$. Hence, by closure,
\[\sum_{j\in J'} L(c_j',\mu( E_j'))\in \intp L(\mathcal  C',\beta_*\mu).\]
the disjoint sum was arbitrary hence this implies (v1).

3. Let
\[C=\bigcup_{D\in\mathfrak D}
(A\cap\lim(\beta D)\setminus D)\cup(D\setminus\lim(\beta D)).\]
This is a weakly $L,\mu$-negligible set.

Assume that $x\in A\setminus C$. Being $A=\bigcup\mathfrak D$ there is an element
$D\in\mathfrak D$ such that $x\in D$. Being $D\setminus\lim(\beta D)\subset C$
there is an element $\mathbf x\in\beta D$ such that $\lim \mathbf x=x$.
Then, being $\beta D$ a generator system  and a forest we find that there is a
possibly smaller $D'\in\mathfrak D$ such that $\mathbf x\in \beta D'$ and
$\lim^*f(\mathbf x)\in(\mathcal A')^{\beta D'}$. That implies $f(x)\in\mathcal A^{D'}$.
Being $A\cap\lim(\beta D')\setminus D'\subset C$ implies that $x\in D'$.
This and the previous fact implies that $f(x)\in\mathcal A^x $.

If $x\notin A$ then $f|_A(x)=0\in\mathcal A^x$ is obvious.
Ultimately, we see that $f|_A(x)\in \mathcal A^{x}$ except if $x\in C$.

Then, the weak negligibility of $C$ implies that we can a add a
small pointed envelope $\tilde{\mathcal  C}$ to $\mathcal  A$ such that
\[\tilde{\mathcal  A}=\mathcal  A+\tilde{\mathcal  C}\]
would approximate $f|_A$ and satisfy our statement.

b.) 1. Consider  an envelope
\[\mathcal  A=\sum_{j\in J} c_j\chi_{E_j}+\mathcal  C\]
 approximating $f|_A$. Again we can restrict to $A$ and extend again.
This and the local finiteness of $\mathfrak S$ allow us the assumptions
\[A=\bigcup_{j\in J}^{\updisjoint}E_j\]
and
\[\mathcal C\subset (\mathcal V\times  A)\cup \{0\}\times X.\]
Then there is a generator system   $\mathfrak D\subset\mathfrak S$
for the envelope $\mathcal  C$, and we can assume that
$\bigcup\mathfrak D=A$, $\mathfrak D$ is a forest,
and $\mathfrak D$ is divided by $\mathfrak E=\{E_j\,:\,j\in J\}$.
This latter one also implies that $\mathfrak D$ is a generator system for $\mathcal A$.

Then let us define the pointed envelope
\[\mathcal  C'=\bigcup_{D\in\mathfrak D}\mathcal  C^D\times \beta D
\,\,\cup\,\,(\{0\}\times \mathfrak X_{\mathfrak S})\]
and the envelope
\[\mathcal  A'=\sum_{j\in J} c_j\chi_{\beta E_j}+\mathcal  C'.\]
(Here the finiteness of $\mathfrak E$ is important for $\mathcal C'$ in order to be a pointed envelope,
cf. Lemma
\ref{lem:CS1}.c).

One can see that
\[\mathcal  A'=\bigcup_{D\in\mathfrak D}\mathcal  A^D\times \beta D
\,\,\cup\,\,(\{0\}\times (\mathfrak X_{\mathfrak S}\setminus\beta A)).\]
2. Then, we claim,
\[\intp L(\mathcal  A',\beta^*\mu)\subset \intp L(\mathcal  A,\mu).\]
Now,
\[\intp L\biggr(\sum_{j\in J} c_j\chi_{\beta E_j},\beta_*\mu\biggl)=\sum_{j\in J} L(c_j,\mu( E_j))
=\intp L\biggr(\sum_{j\in J} c_j\chi_{E_j},\mu\biggl)\]
shows that it is enough to prove that
\[\intp L(\mathcal  C',\beta^*\mu)\subset \intp L(\mathcal  C,\mu).\tag{v2}\]
Let
\[\sum_{j\in J'} c_j'\chi_{\beta E_j'}\subset \mathcal  C'\]
be a finite disjoint sum. Then Lemma \ref{lem:CS1}.c and the algebraic equivalence implies
\[\sum_{j\in J'} c_j'\chi_{ E_j'}\subset \mathcal  C.\]
Hence,
\[\sum_{j\in J'} L(c_j',\mu( E_j'))\in \intp L(\mathcal  C,\mu).\]
The finite disjoint sum was arbitrary, hence we have proved (v2).

3.
Let
\[C=\bigcup_{S\in\mathfrak D\cup\mathfrak E}
((\beta S)\setminus \lim{}^{-1}(S))\cup
(\beta A\cap\lim{}^{-1}(S)\setminus (\beta S)).\]
That is a weakly $L,\beta_*$-negligible set.

Let $\mathbf x\in\beta A\setminus C$ be arbitrary. Assume $\lim \mathbf x=x$.
Then $(\lim^{*} f)(\mathbf x)=f(x)$.
We know that there is an element $E\in\mathfrak E$ such that $\mathbf x\in\beta E$.
Then
\[(\beta E)\setminus \lim{}^{-1}(E)\subset C\]
implies that $x\in E$.
That implies that  is an element $D\in\mathfrak D$ such that
$x\in D$ and $f(x)\in\mathcal A^D$. Now, $x\in D$ and
\[\beta A\cap \lim{}^{-1}(D)\setminus (\beta D)\subset C\]
implies that $\mathbf x\in\beta D$.
On the other hand, $f(x)\in\mathcal A^D$ implies $f(x)\in(\mathcal A')^{\beta D}$.
Ultimately, we can deduce that $\lim{}^*f(\mathbf x)=f(x)\in(\mathcal A')^{\mathbf x}$.

If $\mathbf x\notin\beta A$ then $ \lim{}^*f|_{\beta A}(\mathbf x)=0\in(\mathcal A')^{\mathbf x}$
is clear. Hence, we see that $\lim{}^*f|_{\beta A}(\mathbf x)\in(\mathcal A')^{\mathbf x}$
except on the points of $C$.

Then, the weak negligibility of $C$ implies that we can a add a
small pointed envelope $\hat{\mathcal  C}$ to $\mathcal  A$ such that
\[\hat{\mathcal  A}=\mathcal  A'+\hat{\mathcal  C}\]
would approximate $\lim{}^*f|_{\beta A}$ and satisfy our statement.
\end{proof}
\end{lemma}

\begin{lemma}\label{lem:EQUIMAC}
\halfquote
Let $\mathfrak S$ be locally finite and
$X$ be a topological space compatible to $\mathfrak S$.

Assume that for each set $S,B\in\mathfrak S$

i.) The set
\[(B\cap\lim(\beta S)\setminus S)\cup(S\setminus\lim(\beta S) )\]
is weakly $L,\mu$-negligible.

ii.) The set
\[((\beta S)\setminus \lim{}^{-1}(S))\cup
(\beta B\cap\lim{}^{-1}(S)\setminus (\beta S))\]
is weakly $L,\beta_*\mu$-negligible.

Then, we claim, a function $f:X\rightarrow\mathcal V$ is extended Lebesgue-McShane
$L,\mu$-integrable if and only if it is McShane $L,\mu$-integrable.
In those cases the integrals are the same.
\begin{proof}
That immediately follows from Lemma \ref{lem:DEAT1} and characterization Lemma
\ref{lem:EQUIGEN}.$(\mathtt{\overline E_0})$.
\end{proof}
\end{lemma}
\begin{remark}
The awkward appearance of the sets $B$ in the statement above is due to the fact that we
have not defined extended (=local) negligibility.
\end{remark}
\begin{example}
Let $\mathfrak I$ be the set of possibly degenerate intervals on $\mathbb R$.
The statement above implies that classical Lebesgue integrability is equivalent to
classical McShane integrability.
\end{example}
~

\paragraph{\textbf{D. The Kurzweil-Henstock integral}}
~\\

The Kurzweil-Henstock integral can be defined in a  manner  analogous to the
the McShane integral. From practical viewpoint a major difference is that
in order to get a tagged Kurzweil-Henstock division the condition
\texttt{(DT2)} should be substituted by
\begin{itemize}

\item[\texttt{(DT2')}] $c:\mathfrak E\rightarrow \beta A$ is a function such that
\[c(E)\in\beta E\quad\text{and}\quad E\subset h(c(E))\]
for each $E\in\mathfrak E$.
\end{itemize}
\newpage\section{Locally convex integrals (in fibers)}\label{sec:MODIF}

The Lebesgue-McShane integral applies to the case of topological vector spaces,
even if the definition of the integral does not require a scalar multiplication structure.
\begin{conven}
When dealing with  topological vector spaces linearity is assumed unless told otherwise.
\end{conven}
~
\paragraph{\textbf{A. Locally convex envelopes}}~\\

We will show that if we integrate a locally convex space valued function the we can
make the assumption that the fibers of envelopes are convex.
\begin{remin}
Suppose that $\mathcal  V$ is a vector space, $C\subset \mathcal  V$. Let
\[\conv C\]
be the set of convex combinations of elements of $C$.

If $\mathcal V$ is a topological vector space then
\[\overline{\conv C}=\overline{\conv \overline C}=\conv\overline{\conv C}.\]
\end{remin}
\begin{lemma}\label{lem:APPCONV}
a.) Let $\mathcal V$ be a vector space, $\mathfrak S\subset\mathfrak P(X)$ be an interval
system. Suppose that $\mathcal A$ is a (pointed) envelope of   sets in $\mathcal  V$ with respect to $\mathfrak S$.

Then, we claim,
\[\conv \mathcal  A=\{(v,x)\in\mathcal V\,:\,v\in\conv (\mathcal A^x)\}\]
is a (pointed) envelope of sets in $\mathcal  V$ with respect to $\mathfrak S$.

b.) \triquartquote Suppose that $\mathcal  A$ is  an envelope  of pointed  sets in
$\mathcal  V$ with respect to $\mathfrak S$.

Then, we claim,
\[\intp L(\conv \mathcal  A,\mu)=\overline{\conv \intp L(\mathcal  A,\mu)}.\]
\begin{proof}
a.) Only the generator system property is nontrivial.
Through Lemma \ref{lem:FOREST}
we can find a generator system to $\widetilde{\mathfrak D}$ to $\mathcal A$ which is
a fully branching countable forest. If $(v,x)\in\conv \mathcal A$, $x\in\bigcup\widetilde{\mathfrak D} $ then
there exist finitely many elements $(v_k,x)\in\mathcal A$, $t_k\in\mathbb R$ $(k\in K)$ such that
\[\sum_{k\in K}t_kv_k=v.\]
By the generator system property there are elements $x\in D_k\in \widetilde{\mathfrak D}$ such that
$v_k\in \mathcal A^{D_k}$.
Let $D$ be the smallest of these $D_k$'s, ie. the intersection. Then $v_k\in \mathcal A^{D}$.
Consequently,
\[v=\sum_{k\in K}t_kv_k\in \conv \mathcal A^{\mathfrak D}\subset(\conv \mathcal A)^{\mathfrak D}.\]
On the other hand, if $x\in\bigcup\widetilde{\mathfrak D} $ then $(\conv \mathcal A)^x=\{0\}$
is obvious. That proves the generator system property for $\conv \mathcal A$.

b.) By linearity (Lemma \ref{lem:DECO}.b/b') it is enough to prove the statement for pointed
envelopes $\mathcal  A=\mathcal  C$.

The finite ``nondisjoint'' sums contributing to $\intp L(\conv \mathcal  C,\mu)$
are closed to convex combinations, so
\[\overline{\conv \intp L(\mathcal  C,\mu)}\subset\intp L(\conv \mathcal  C,\mu)\]
is clear.
Consider any finite ``disjoint'' sum
\[\sum_{j\in J} L(c_j,\mu(E_j))\]
contributing to $\intp L(\conv \mathcal  C,\mu)$. According to the standard methods
there is a refinement $\mathfrak D$ of $\{E_j\,:\,j\in J\}$ such that
for $E\supset D\in\mathfrak D$
\[c_j\in\conv \mathcal  C^D.\]
So, the finite ``disjoint'' refines to a countable ``disjoint'' sum.
Because of closure, it is enough to prove that finite disjoint sums
\[\sum_{j\in J} L(c_j,\mu(E_j))\]
such that
\[c_j\in\conv \mathcal  C^{E_j}\]
are in $\conv \intp L(\mathcal  C,\mu)$. Consider that case.
Here for all $j\in J$ there is a convex combination
\[c_j=\sum_{l\in L_j} t_{j,l}c_{j,l}\]
such that $c_{j,l}\in\mathcal  C^{E_j}$.
Taking common partitions of $1$ we can actually assume that that for all
$j\in J$ the family of $t_{j,l}$ is the same, say,
\[t_1,\ldots,t_s.\]
But then
\[\sum_{j\in J} L(c_j,\mu(E_j))\in \sum_{h=1}^s t_h \intp L(\mathcal  C,\mu)\subset
\conv \intp L(\mathcal  C,\mu).\]
\end{proof}
\end{lemma}
\begin{cor}
In the case when $\mathcal V$ is a vector space and is a
$\mathcal Z$ locally convex vector space then  the Lebesgue-McShane integral can be defined in terms of
envelopes of convex sets in
$\mathcal  V$ with respect to $\mathfrak S$.
\begin{proof}
The space $\mathcal  Z$ is locally convex and the convex sets form an
neighborhood base; hence nothing is lost if we pass to the convex closure of approximations.
\end{proof}
\end{cor}
\begin{remark}
If one is not interested in topological groups in general then one can
develop the integral in terms of envelopes of convex pointed sets from the start.
One convenient thing to do in this case is to pass to the notion of convex
semivariations
\begin{multline}
\csvar(L,c,\mu,A)=\overline{\biggl\{\sum_{j\in J}t_jL(c,\mu(A_j))\,:\,
J\text{ is finite},\,A_j\in\mathfrak S,\dots}\\\overline{\ldots \,
A_j\subset A,\text{ the $A_j$ are pairwise disjoint},\,0\leq t_j,\,\sum_{j\in J}t_j\leq1\biggr\}},\notag
\end{multline}
\begin{multline}
\csvar(L,f,\mu)=
\overline{\biggl\{\sum_{j\in J}\int t_jL(f|_{A_j},\mu)\,:
\,J\text{ is finite},\,A_j\in\mathfrak S,\dots}\\\overline{\ldots\,
A_j\subset A,\text{ the $A_j$ are pairwise disjoint},\,0\leq t_j,\,\sum_{j\in J}t_j\leq1  \biggr\}}.\notag
\end{multline}
It is easy to see that these are just the convex closures
of ordinary semivariations.

In general, using only envelopes of convex sets keeps us somewhat
closer to the spirit of the Riemann integral.
\end{remark}
~

\paragraph{\textbf{B. Integration of real-valued functions}}~\\

Here are some immediate benefits of the observations from the previous section:
\begin{conven}
Throughout this paragraph we will consider the case when
$\mathtt M:\mathbb R\times\mathcal  W\rightarrow \mathcal  W$
is the standard multiplication on a locally convex space $\mathcal W$.
(So, $\mathcal  V=\mathbb R$, $\mathcal  W=\mathbb Z$.)

In what follows here we always have the choice $L=\mathtt M$.
If it plays no specific role then we omit $\mathtt M$ form the notation
(for example, from semivariation expressions).
\end{conven}

\begin{lemma}[Beppo Levi's theorem]\label{lem:CONCL1}
\semiquote

Suppose that $f_\lambda:X\rightarrow\mathbb R$ ($\lambda\in\Lambda$, countable) are integrable,
and for each $0\in\mathcal  T\subset \mathcal  W$ neighborhood there is a finite set
$\Xi\subset\Lambda$ such that
\[\sum_{\lambda\in\Xi} \svar(f_\lambda,\mu)\subset \mathcal  T.\]
Then, we claim,
\[\sum_{\lambda\in\Lambda}^{\ovline} f_\lambda\]
(understood pointwise) is convergent except on a strongly negligible set, and for
\[f=\sum_{\lambda\in\Lambda}^{\ovline}f_\lambda\]
(defined arbitrarily on non-convergence points)
\[\int f\,\mu=\sum_{\lambda\in\Lambda}^{\ovline}\int f_\lambda\,\mu.\]
Furthermore, in this case
\[\svar(f,\mu)\subset\overline{ \sum_{\lambda\in\Lambda}\svar(f_\lambda,\mu)}.\]
\begin{proof}
Let $C$ be the points of non-convergence.
Let $0\in\mathcal  T'\subset \mathcal  W$ be an arbitrary neighborhood.
Let $0\in\mathcal  T\subset \mathcal  W$ be a neighborhood such that
\[\conv\mathcal  T-\conv\mathcal  T+\mathcal  T \subset\mathcal  T'.\]
Let us choose $\Xi$ as in our assumptions.
Then
\[\mathbb R\times C\subset\conv\biggl(\sum_{\lambda\in\Lambda\setminus\Xi}\mathcal  C_\lambda-
\sum_{\lambda\in\Lambda\setminus\Xi}\mathcal  C_\lambda\biggr).\]
On the other hand, we see that
\[\intp L\biggl(\conv\biggl(\sum_{\lambda\in\Lambda\setminus\Xi}\mathcal  C_\lambda-
\sum_{\lambda\in\Lambda\setminus\Xi}\mathcal  C_\lambda\biggr),\mu\biggr) \]
\[=\overline{\conv\biggl(\sum_{\lambda\in\Lambda\setminus\Xi}\intp L(\mathcal C_\lambda,\mu)
-\sum_{\lambda\in\Lambda\setminus\Xi}\intp L(\mathcal C_\lambda,\mu) \biggr) }\]
\[\subset\overline{\conv(\mathcal  T-\mathcal  T)}\subset{\conv \mathcal T}-
 {\conv \mathcal T}+\mathcal T\subset \mathcal T'.\]
Here $\mathcal  T'$ was arbitrary, hence strong negligibility is proven.

In the same setting, if $\mathcal A_\xi$ are envelopes approximating $f_\xi$ we find that
\[\sum_{\xi\in\Xi}\mathcal A_\xi+ \sum_{\lambda\in\Lambda\setminus\Xi}\mathcal  C_\lambda
-\sum_{\lambda\in\Lambda\setminus\Xi}\mathcal  C_\lambda+
\sum_{\lambda\in\Lambda\setminus\Xi}\mathcal  C_\lambda\]
is  an envelope approximating  $f$.
Estimating ( using a sufficiently small $\mathcal T$)
and applying the Cauchy criterium we obtain the sum statement.

Similarly, using the pointed envelopes
\[\sum_{\lambda\in\Lambda}\mathcal  C_\lambda
-\sum_{\lambda\in\Lambda\setminus\Xi}\mathcal  C_\lambda+
\sum_{\lambda\in\Lambda\setminus\Xi}\mathcal  C_\lambda\]
we obtain the semivariation statement.
\end{proof}
\end{lemma}
\begin{lemma}\label{lem:CONCL2}
\semiquote

If $f,g: X\rightarrow\mathbb R$ are integrable
then $f\vee g$, $f\wedge g$ are integrable. If $d\geq 0$ then
$f\wedge d$ is integrable, if $c\leq 0$ then $f\vee c$ is integrable. $|f|$ is integrable.
\begin{proof}
It is enough to prove the statement for $f\wedge d$, then the rest follows by the standard tricks.
Consider the map
\[(\wedge d)^*:\,r\mapsto r\wedge d.\]
Assume that
\[\mathcal  A=s+\mathcal  C\]
is  an envelope approximating $f$ such that $s$ is a simple step-function and $\mathcal  C$ is convex,
and $\intp L(\mathcal  A,\mu)$ is close to the integral.
(This can be achieved according Lemma \ref{lem:APPREFIN} and \ref{lem:APPCONV}.)

Then $(\wedge d)^*\mathcal  A$ (applied in the scalar component) is envelope approximating $f\wedge d$.
Indeed, it contains the simple step-function $s\wedge d$ and any generator system of $\mathcal  A$
will do as a generator system. From convexity,
\[(\wedge d)^*\mathcal  A-(\wedge d)^*\mathcal  A\subset \mathcal  C-\mathcal  C.\]
This and the Cauchy criterium implies the integrability of $f\wedge d$.
\end{proof}
\end{lemma}
\begin{lemma} \label{lem:CONCL3}
Let $\mathfrak S=\mathfrak P(X)$ be an interval system; $\mu:\mathfrak S\rightarrow\mathbb R$
be a measure.
Suppose that $\mu\geq0$.

Then, we claim,  $f\geq 0$ implies
\[\int f\mu\in\svar(f,\mu)\subset\biggl[0,\int f\mu\biggr].\]
More generally, if $f$ is integrable then $|f|^+$ and $|f|^-$ are integrable and
\[\int f\mu\in\svar(f,\mu)\subset \biggl[\int |f|^-\mu,\int |f|^+\mu\biggr].\]
The endpoints of the intervals are contained in the semivariation.
\begin{proof}
The first statement follows from the monotonicity statement in Lemma \ref{lem:SVMON}.
The second statement follows from Lemma \ref{lem:SVSUM}.b.
(In fact, Lemma \ref{lem:SVCHAR}.a is sufficient if we do not care about endpoints.)
\end{proof}
 \end{lemma}
\begin{lemma}\label{lem:CONSW}
\semiquote
Then, we claim,
weakly ($\mathtt M,\mu$-)negligible sets are strongly ($\mathtt M,\mu$-)negligible.
\begin{proof}
Suppose that $C$ is weakly negligible.
Let $0\in\mathcal  T'\subset Z$ be an arbitrary neighborhood.
Then there is a closed convex neighborhood $0\in\mathcal  T\subset\mathcal Z$ such that
$\mathcal  T\subset \mathcal  T'$. Let us divide $\mathcal  T$ by $\mathbb Z$.
For $n\in\mathbb Z$ let $\mathcal  C_n$ be  an envelope approximating $n|_C$ such that
\[\intp L(\mathcal  C_n,\mu)\subset \mathcal  T_n.\]
Let
\[\mathcal  C= \conv \sum_{n\in\mathbb Z}\mathcal  C_n.\]
Then
\[\intp L(\mathcal  C,\mu)=\overline{ \sum_{n\in\mathbb Z}\intp L(\mathcal  C_n,\mu) }\subset
\overline{\sum_{n\in\mathbb Z}\mathcal  T_n }=\overline{\mathcal  T}=\mathcal  T. \]
Then
\[\intp L(\conv \mathcal  C,\mu)=\overline{\conv \intp L( \mathcal  C,\mu)}\subset\mathcal T \subset\mathcal  T' \]
while $\conv \mathcal  C$ contains $\mathbb R\times C$.
\end{proof}
\end{lemma}
So, the simple term ``negligible'' can be used.
\begin{lemma}\label{lem:CONNEG}
\semiquote
Let $C$ be a negligible set.

Then, we claim,
for each neighborhood $0\in\mathcal  T\subset \mathcal  Z$
 there exist a countable family of sets $\{E_\lambda\}_{\lambda\in\Lambda}$
such that:

 For all $x\in X$ the  sum
\[\sum_{\lambda\in\Lambda} \chi_{E_\lambda}(x)\]
diverges (does not exist) on $C$ but
\[\overline{\sum_{\lambda\in\Lambda} \svar(E_\lambda,\mu)}\subset\mathcal  T.\]
Remark: In this case, $C$ is contained in the divergence set of
\[\sum_{\lambda\in\Lambda} \chi_{E_\lambda}(x)-\sum_{\lambda\in\Lambda} \chi_{E_\lambda}(x).\]
\begin{proof}
Let $\mathcal  C$ be  an envelope such that it contains $\mathbb R\times C $ but
\[\intp L(\mathcal C,\mu)\subset \mathcal  T.\]
According the generator property of envelopes for each $n\in\mathbb N$ the set
\[{}^{n}\mathcal  C=\{x\,:\,(n,x)\in\mathcal  C\}\]
is a countable union of elements of $\mathfrak S$.
In fact, according to Lemma \ref {lem:SETUN}.a
\[{}^{n}\mathcal  C=\bigcup_{\lambda\in\Lambda_n}^\updisjoint E_{n,\lambda},\]
where $E_{n,\lambda}\in\mathfrak S$.
Then let $\{E_\lambda\}_{\lambda\in\Lambda}$ be the disjoint union of the indexed
families $\{E_{n+1,\lambda}\}_{\lambda\in\Lambda_{n+1}}$
$(n\in\mathbb N)$.
Then, by Lemma \ref{lem:NEIGH},
\[\overline{\sum_{\lambda\in\Lambda}\svar(E_\lambda,\mu)}\subset
\intp L(\mathcal C,\mu)\subset \mathcal T.\]
The construction of the sets implies our statement.
\end{proof}
\end{lemma}
From Beppo Levi's theorem we can devise a ``approximate structure theorem'' for the Lebesgue integral:
\begin{lemma}\label{lem:STRUC}
Let $\mathfrak S\subset\mathfrak P(X)$ be an interval system and
$\mu:\mathfrak S\rightarrow\mathcal  W$ be a measure.
Assume that $\mathcal  W$ is a Fr\'echet $(\mathtt M_1)$ locally convex space.

Suppose that  $f:X\rightarrow\mathbb R$ is integrable.
Then, we claim, for arbitrary $0\in\mathcal  T\subset\mathcal  W$ there exist countable families
$c_\lambda\in\mathbb R$, $E_\lambda\in\mathfrak S$ ($\lambda\in\Lambda$) such that
the following hold:
\begin{itemize}
\item[i.)]
\[f(x)=\sum_{\lambda\in\Lambda}^{\ovline} c_\lambda\chi_{E_\lambda}(x)\]
except on  a negligible set.
\item[ii.)]
\[\overline{\sum_{\lambda\in\Lambda} \svar(c_\lambda,E_\lambda,\mu)}\subset
\svar(f,\mu)+\mathcal  T.\]
\item[iii.)]
and for each neighborhood $0\in\mathcal  T'\subset\mathcal  W$ there exists $\Xi\subset\Lambda$ such that
\[\overline{\sum_{\lambda\in\Lambda\setminus\Xi} \svar(c_\lambda,E_\lambda,\mu)}\subset
\mathcal  T'\]
\end{itemize}
\begin{proof}
Let us assume that $\{\mathcal  K_n\}_{n\in\mathbb N}$ is a countable neighborhood basis of
$\mathcal  Z$.

The by induction, we can define a series  of neighborhoods $\{\mathcal  B_n\}_{n\in\mathbb N}$
of $0$ such
that with the choice $\mathcal B_{-1}=\mathcal K_{-1}=\mathcal T$
\[\mathcal  B_{n}+\mathcal  B_{n}\subset\mathcal  B_{n-1},\qquad
\mathcal B_n\subset-\mathcal B_{n-1},\qquad\text{and}\qquad
\mathcal  B_n\subset \mathcal  K_n.\]
Then, by induction, we can define a series of simple step-functions
\[g_n=\sum_{j\in J_n}c_{j,n}\chi_{E_{j,n}}\]
$(n\in\mathbb N)$ such that
\[\svar(f-g_0-\ldots-g_n)\subset \mathcal  B_{n+2}.\]
Namely, from the integrability of $f-g_0-\ldots-g_{n-1}$ follows that we can find a very good
approximating simple step function $g_n$. Notice that
\[g_{n+1}=(f-g_0-\ldots-g_{n})-(f-g_0-\ldots-g_{n+1}),\]
which implies that
\[\svar(g_{n+1},\mu)\subset\overline{ \mathcal  B_{n+2}-\mathcal  B_{n+3}} \subset
\overline{\mathcal B_{n+2}+\mathcal B_{n+2}}\subset \overline {\mathcal B_{n+1}}\subset\mathcal B_{n}.\]

The simplicity of $g_n$ implies
\[\svar(g_n,\mu)=\overline{\sum_{j\in J_n}\svar(c_{j,n},E_{j,n},\mu)}.\]
Then, we can see that for $n_0\in\mathbb N$
\[\overline{\sum_{n> n_0,\, j\in J_n}\svar(c_{j,n},E_{j,n},\mu)}=
\overline{\sum_{n> n_0}\svar(g_n,\mu)}\subset \overline{\sum_{n> n_0}  \overline{\mathcal B_{n}}}
=\overline{\sum_{n> n_0} {\mathcal B_{n}}}\subset
\overline{\mathcal  B_{n_0}}\]
So, in particular,
\[\subset\mathcal  B_{n_0-1}\subset\mathcal  K_{n_0-1},\]
and this observation  implies iii.).

This latter semivariation statement implies that Beppo Levi's theorem can be applied and
\[f'= \sum_{n\in\mathbb N,\, j\in J_n}^{\ovline} c_{j,n}\chi_{E_{j,n}}\]
exists (converges) almost everywhere. Let us define $f'$ arbitrarily on the non-convergence places.
We claim that
\[\svar(f'-f,\mu)=0.\tag{\dag}\]
First of all,
\[f'-f= \biggl(\sum_{n>m,\, j\in J_n}^{\ovline} c_{j,n}\chi_{E_{j,n}}\biggr)-(f-g_0-\ldots-g_m)\]
holds  almost everywhere (from Beppo Levi's theorem).
Then the semivariation part of Beppo Levi's theorem implies that
\[\svar(f'-f)\subset \overline{
\overline{\sum_{n>m,\, j\in J_n}\svar(c_{j,n},E_{j,n},\mu)}
-\svar(f-g_0-\ldots-g_m,\mu)
}\subset \]
\[\subset\overline{\overline{\mathcal B_{m}} - \mathcal B_{m+2} }
\subset\overline{\mathcal B_{m}+{\mathcal B_{m+1}}}
\subset \mathcal B_{m}+\mathcal B_{m+1}+\mathcal B_{m+1}\subset\mathcal B_{m-1}  \subset \mathcal K_{m-1}.\]
That proves (\dag). Applying Beppo Levi's theorem to the infinite sum with equal terms $f'-f$
we obtain that $f'-f$ is almost everywhere $0$.
That proves part i.) of our statement.

From Lemma \ref{lem:SVCHAR}.a we know that
\[\svar(g_0,\mu)\subset\overline{\svar(f,\mu)-\svar(f-g_0,\mu)}. \]
Hence
\[\overline{\sum_{n\in\mathbb N,\, j\in J_n}\svar(c_{j,n},E_{j,n},\mu)}=
\overline{\sum_{n\in\mathbb N}\svar(g_n,\mu)}\]
\[\subset\overline{\svar(f,\mu)-\svar(f-g_0,\mu)+\sum_{n>1}\svar(g_n,\mu) }\]
\[\subset \overline{\overline{\svar(f,\mu)-\mathcal B_{2}}+\overline{\mathcal B_1}}
=\overline{{\svar(f,\mu)-\mathcal B_{2}}+{\mathcal B_1}}\]
\[\subset \overline{\svar(f,\mu)+\mathcal B_{1}+\mathcal B_{1}}\subset
\overline{\svar(f,\mu)+\mathcal B_{0}}\]
\[\subset\svar(f,\mu)+\mathcal B_{0}+\mathcal B_0\subset \svar(f,\mu)+\mathcal B_{-1} = \svar(f,\mu)+\mathcal T.\]
That proves part ii.)
\end{proof}
\end{lemma}
\begin{lemma}\label{lem:CONINT}
\semiquote

Consider the integrable extension
\[\hat\mu:\mathfrak S^\mu\rightarrow\mathbb R.\]

Then, we claim:

a.) If $\mathcal  A$ is  an envelope over $\mathfrak S^\mu$ and  then
for each neighborhood $0\in\mathcal  T\subset \mathcal  W$ there is  an envelope $\tilde{\mathcal  A}$
 over $\mathfrak S$
such that
\[\mathcal  A\subset\tilde{ \mathcal  A}\]
but
\[ \intp L(\tilde{\mathcal  A},\mu)\subset \conv\intp L(\mathcal  A,\hat\mu)+\mathcal  T.\]

b.)
A function $f:X\rightarrow\mathbb R$ is (locally) $\mu$-integrable if and only if
it is (locally) $\hat\mu$-integrable. In this case the integrals and the semivariation
sets are the same.
\begin{proof}
a.)
Let $0\in\mathcal  T'\subset \mathcal  W$ be a convex neighborhood such that
$\mathcal T'+\mathcal T'\subset\mathcal T$.
Suppose that $\mathfrak D$ is a generator set of $\mathcal A$ and
$s$ is (the graph of) a step-function contained in $\mathcal A$.
Then for each
$D\in \mathfrak D$ let
\[R_D=((\conv \mathcal A)^\mathfrak D\cap\mathbb Q)
\cup\{\max\, (\conv \mathcal A)^\mathfrak D\}\cup\{\min\, (\conv \mathcal A)^\mathfrak D\}\]
where the two last term are understood to be in the expression if only they exist.

Then, we can notice that
\[\conv \mathcal A^D=\conv R_D.\]
Let us define
\[\mathcal R=\bigcup_{D\in\mathfrak D}(R_D\times D)\cup s.\]
One can see that this is  an envelope over $\mathfrak S^\mu$ with countable many values in $\mathbb R$ and
\[\mathcal R\subset\conv \mathcal R=\conv\mathcal A .\]
According to Lemma \ref{lem:MECONTINT} there is  an envelope $\tilde {\mathcal A}$ over $\mathfrak S$
such that
\[\intp L(\tilde{\mathcal  A},\mu)
\subset
\intp L(\mathcal  R,\hat\mu)+\mathcal  T'.\]
But then we have
\[\subset \intp L(\conv \mathcal  A,\hat\mu)+\mathcal  T'\subset
\overline{\conv \intp L(\mathcal  A,\hat\mu)}+\mathcal  T'\subset
{\conv \intp L(\mathcal  A,\hat\mu)}+\mathcal  T.\]
b.) That is an immediate consequence of part a.
\end{proof}
\end{lemma}
\begin{lemma}\label{lem:CONNEVO}
 Suppose that $\mathfrak S\subset \mathfrak P(X)$ is an interval system,
$\mathcal V$ is a vector space.

Then, we claim,  for any pointed envelope $\mathcal C$,
there are countable families $c_\lambda\in\mathbb R$, $E_\lambda\in\mathfrak S$ $(\lambda\in\Lambda)$
such that
\[\conv\mathcal C=\sum_{\lambda\in\Lambda}\conv\var(c_\lambda,E_\lambda).\]
\begin{proof}
We can assume $\mathcal C=\conv\mathcal C$.
Then \[\mathcal C=\mathcal C^+-\mathcal C^-,\]
where
\[\mathcal C^+=\mathcal C\cap(\mathbb R^[0,+\infty)\times X)\quad\text{and}\quad
\mathcal C^-=(-\mathcal C)\cap([0,+\infty)\times X).\]
So, we can also assume that $\mathcal C=\mathcal C^+$.
Let us define the function
\[r:X\rightarrow\mathbb R\]
\[r(x)=\sup \mathcal C^x.\]
Along $\mathbb N$ let us list the elements of $\mathfrak D$ infinitely times.
Let $D_n\in\mathfrak D$ be the $n$th listed element.
By recursion, let us define a sequence $c_n\in[0,+\infty)$ such that
\[\sum_{n\in\mathbb N}^{\ovline}c_n\chi_{D_n}(x)=r(x).\tag{\ddag}\]
That will immediately yield our statement.

Indeed, such a sequence can be constructed by the prescription such that
\begin{itemize}
\item[i.)]
\[\max\biggl( 0, \inf_{x\in D_n} \biggl(r(x)-\sum_{0\leq j<n} c_j\chi_{D_j}(x) \biggr)-2^{-n}
 \biggr)\leq c_n\leq
\inf_{x\in D_n} \biggl(r(x)-\sum_{0\leq j<n} c_j\chi_{D_j}(x) \biggr);\]
\item[ii.)]
\[c_n=\min_{x\in D_n} \biggl(r(x)-\sum_{0\leq j<n} c_j\chi_{D_j}(x) \biggr)\]
if the right side exists.
\end{itemize}
(We can make such a sequence unique by an appropriate reference to dyadic numbers.)
One can prove that such a sequence $c_n$ has the property (\ddag).
\end{proof}
\end{lemma}
~

\paragraph{\textbf{C. Integration in strong variation}}~\\

A nonequivalent, but stronger fiber condition for locally convex spaces can be introduced as follows:
\begin{defin}
Suppose that $p$ is a seminorm on $\mathcal  V$.
We say that the function $f$ is $p$-strongly contained in the envelope $\mathcal  A$, if
\[\mathcal  A=s+\mathcal  C\]
such that
\[\{(x,y)\in X\times\mathcal  V\,:\,p(y)<p(f(x)-s(x))\}\subset \mathcal  C.\]
\end{defin}
\begin{defin}
Suppose that $\mathcal  V$ is a locally convex vector space.
We say that the integral exists in $p$-strong sense, if for each seminorm $p$ on $\mathcal  V$
for each $0\in\mathcal  T\subset Z$ there is  an envelope $\mathcal  A$ containing $f$ $p$-strongly
such that
\[\intp L(\mathcal  A,\mu)\subset a-\mathcal  T. \]
The integral exists in strong sense if it exists in $p$-strong sense for all $p$.
\end{defin}
\begin{remak}
The results of the previous section can be generalized for $p$-strong integrals, by claiming
convergence in $p$ almost everywhere.
If the integral exists in strong sense and $\mathcal  V$ space is Fr\'echet ($\mathtt M_1$) then
convergence can be claimed almost everywhere.
\end{remak}
~\\
\paragraph{\textbf{D. Infinite integrals}}~\\

If $\mathcal  Z=\mathbb R$ then we can define the integral with more general values:
\begin{defin}
Let $\mathcal  V$, $\mathcal  W$, be commutative
topological groups, $\mathfrak S\subset\mathfrak P(X)$ be an interval system,
$\mu:\mathfrak S\rightarrow\mathcal  W$ be a measure,
$L:\mathcal  V\times\mathcal  W\rightarrow\mathbb R$ be a biadditive pairing, which is
continuous in its second variable. We say
\[\int L(f,\mu)=+\infty\]
if for each $K\in\mathbb R$ there exist  an envelope $\mathcal  A$
 such that
\[\intp L(\mathcal  A,\mu)\subset [K,+\infty).\]
One can similarly define the integral with $-\infty$.
\end{defin}
\newpage\section{Finite type semivariations, sums }\label{sec:FINSVAR}
The title of this section can be dubbed as ``Open integration (in fibers)''.
~\\

\paragraph{\textbf{A. Fundamentals}}
\begin{defin}
If no confusion arises then  for $R\subset\mathcal  V$, $A\subset X$ we write
\[\var(R,A)=(R\times A)\cup(\{0\}\times X).\]
\end{defin}

\begin{defin}\label{def:LFSVAR}
\halfquote

a.) We say that the measure $\mu$ is of locally finite semivariation with respect to the pairing $L$ if
for each set $A\in\mathfrak S$ and neighborhood $0\in\mathcal  T\subset \mathcal  Z$ there exists
a neighborhood $0\in\mathcal  U\subset\mathcal   V$ such that
\[\intp L(\var(\mathcal  U,A),\mu)\subset \mathcal  T.\]

b.) We say that the measure $\mu$ is of locally $\sigma$-finite semivariation with respect to the pairing $L$
if for each set $A\in\mathfrak S$ and neighborhood $0\in\mathcal  T\subset \mathcal  Z$ there exists a
countable decomposition
\[A=\bigcup^\updisjoint_{\lambda\in\Lambda}A_{\lambda}\]
($A_\lambda\in\mathfrak S$) and neighborhoods $0\in\mathcal  U_\lambda\subset\mathcal  V$ such that
\[\intp L\biggl(\sum_{\lambda\in\Lambda}\var(\mathcal  U_\lambda,A_\lambda), \mu\biggr)\subset \mathcal  T.\]

c.) We say that the measure $\mu$ is of $\sigma$-locally finite semivariation with respect
to the pairing $L$
if for each set $A\in\mathfrak S$ there exists  an envelope $\mathcal  C$ of pointed  sets in
$\mathcal  V$ with respect to $\mathfrak S$ such that
\[\intp L(\mathcal  C, \mu)\subset \mathcal  T,\]
and for each $x\in A$
\[0\in (\mathcal  C^x)^\circ\]
(${}^{\circ}$ refers to interior.)
\end{defin}
\begin{lemma}\label{lem:FINCONS}
\halfquote

a.) If $\mu$ is of locally finite semivariation with respect to the pairing $L$ then
$\mu$ is of locally $\sigma$-finite semivariation with respect to the pairing $L$.

b.) If $\mu$ is of locally  $\sigma$-finite semivariation with respect to the pairing $L$ then
$\mu$ is of $\sigma$-locally finite semivariation with respect to the pairing $L$.

This latter implication can be reversed if $\mathcal  V$ is a Fr\'echet ($\,\mathtt M_1$) topological group.
\begin{proof}
Only the last statement is non-trivial. Let $0\in\mathcal  T'\subset\mathcal Z$ be an arbitrary neighborhood.
Let $0\in\mathcal  T\subset\mathcal Z$ be a neighborhood
 such that $\mathcal  T-\mathcal  T+\mathcal  T+\mathcal  T\subset \mathcal  T'$.
Consider a pointed envelope $\mathcal  C$ as in \ref{def:LFSVAR}.iii.
Let $\mathfrak D$ be a generator system for $\mathcal  C$.
For each $x\in A$
\[\bigcup_{D\in\mathfrak D,\,x\in D}\mathcal  C^D=\mathcal  C^x.\]
The set $\mathcal  C^x$ has a nonempty interior.
Any $\mathtt M_1$, complete uniform space has the Baire property, so there is a $D\in\mathfrak D$
such that $\overline{\mathcal  C^D}$ has nonempty interior. Hence
\[0\in (\overline{\mathcal  C^D}-\overline{\mathcal  C^D})^\circ\subset
\left({(\mathcal  C-\mathcal  C+\mathcal  C)^D }\right)^\circ.\]
We have thus proved that those $D$'s such that
$0\in\left({(\mathcal  C-\mathcal  C+\mathcal  C)^D }\right)^\circ$
cover $A$. Applying, say, Lemma \ref{lem:FOREST}, we see that there is a an exact decomposition
$\mathfrak D'$ of $A$ such that
\[0\in\left({(\mathcal  C-\mathcal  C+\mathcal  C)^D }\right)^\circ.\]
for all $D\in\mathfrak D$. But then
\[\sum_{D\in\mathfrak D'}\var((\mathcal  C-\mathcal  C+\mathcal  C)^D,D)\subset
\mathcal  C-\mathcal  C+\mathcal  C\]
and
\[\intp L(\mathcal  C-\mathcal  C+\mathcal  C,\mu)\subset
\intp L(\mathcal  C,\mu)-\intp L(\mathcal  C,\mu)+\intp L(\mathcal  C,\mu)+\mathcal  T
\subset \mathcal  T-\mathcal  T+\mathcal  T+\mathcal  T\subset\mathcal  T'\]
proves our statement.
\end{proof}
\end{lemma}
\begin{lemma}\label{lem:FINSIGMA}
\halfquote

a.) If $\mu$ is of  locally $\sigma$-finite semivariation with respect to the pairing $L$ then
for each set $A\in\boldsymbol\Sigma_0\mathfrak S$
and neighborhood $0\in\mathcal  T'\subset\mathcal Z$ there is an open set
$0\in\mathcal G\subset\mathcal V$
such that
\[\intp L(\var(\mathcal  G,A),\mu)\subset \mathcal  T'.\]

b.) If $\mu$ is of  locally $\sigma$-finite semivariation with respect to the pairing $L$ then
for each $A\in\boldsymbol\Sigma\mathfrak S$ and neighborhood $0\in\mathcal  T'\subset\mathcal Z$
there exists a
countable decomposition
\[A=\bigcup^\updisjoint_{\lambda\in\Lambda}A_{\lambda}\]
($A_\lambda\in\mathfrak S$) and open sets $0\in\mathcal  G_\lambda\subset\mathcal  V$ such that
for the pointed envelope
\[\mathcal C=\sum_{\lambda\in\Lambda}\var(\mathcal  G_\lambda,A_\lambda)\]
we have
\[\intp L(\mathcal C, \mu)\subset \mathcal  T'.\]

c.) If $\mu$ is of $\sigma$-locally finite semivariation with respect to the pairing $L$ then
for each $A\in\boldsymbol\Sigma\mathfrak S$ and neighborhood $0\in\mathcal  T'\subset\mathcal Z$
there exists an pointed envelope $\mathcal C$ such that
\[\mathcal C^x\]
is open for all $x\in A$ but
\[\intp L(\mathcal C,\mu)\subset\mathcal T' .\]
\begin{proof}
a.) Suppose that
\[A=\bigcup_{j\in J} A_j,\]
($A_j\in\mathfrak S$) where $J$ is finite. Let  $0\in\mathcal  T\subset\mathcal Z$ be a neighborhood such that
$\mathcal T+\mathcal T\subset\mathcal T'$. Let us divide $\mathcal T$ by $J$.

According to our assumption there a neighborhoods $0\in\mathcal U_j\subset\mathcal V$
such that
\[\intp L(\var(\mathcal U_j,A_j))\subset \mathcal T'_j.\]
Let
\[\mathcal G=\biggl(\bigcap_{j\in J}\mathcal U_j \bigg)^\circ.\]
Then $\mathcal G$ is an open $0$-neighborhood, yet
\[\intp L(\var(\mathcal G,A),\mu)\subset \intp L\biggl(\sum_{j\in J}\var(\mathcal U_j,A_j),\mu\biggr)
\subset\overline{\sum_{j\in J}\mathcal T_j}\subset\overline{\mathcal T}\subset\mathcal T'.\]

b.)
Suppose that
\[A=\bigcup_{\omega\in\Omega}^\updisjoint A_\omega,\]
($A_\omega\in\mathfrak S$)
where $\Omega$ is countable. Let  $0\in\mathcal  T\subset\mathcal Z$ be a neighborhood such that
$\mathcal T+\mathcal T\subset\mathcal T'$. Let us divide $\mathcal  T$ by $\Omega$.
Let us apply the requirement of Definition \ref{def:LFSVAR}.b. with respect to
$A_\omega$ and $\mathcal T_\omega$. That will yield countable decompositions
\[A_\omega=\bigcup^\updisjoint_{\lambda\in\Lambda_\omega}A_{\omega,\lambda},\]
and neighborhoods $\mathcal U_{\omega,\lambda}$.

Then there is a decomposition
\[A_\omega=\bigcup^\updisjoint_{\omega\in\Omega,\,\lambda\in\Lambda_\omega}A_{\omega,\lambda}.\]
Let us associate $\mathcal G_{\omega,\lambda}=(\mathcal U_{\omega,\lambda})^\circ$
to $A_{\omega,\lambda}$.

Then
\[\intp L\biggl(\sum_{\omega\in\Omega,\,\lambda\in\Lambda_\omega}
\var(\mathcal G_{\omega,\lambda},A_{\omega,\lambda} ), \mu\biggr)\subset
\overline{\sum_{\omega\in\Omega}\mathcal T_\omega  }
\subset\overline{\mathcal T}\subset\mathcal T'\]
proves our statement.

c.) Suppose that
\[A=\bigcup_{\omega\in\Omega}^\updisjoint A_\omega,\]
where $\Omega$ is countable. Let  $0\in\mathcal  T\subset\mathcal Z$ be a neighborhood such that
$\mathcal T+\mathcal T\subset\mathcal T'$. Let us divide $\mathcal  T$ by $\Omega\times\mathbb N$.
Let us apply the requirement of Definition \ref{def:LFSVAR}.b. with respect to
$A_\omega$ and $\mathcal T_{\omega,n}$. That will yield envelopes $\mathcal C_{\omega,n}$.
Then one can see that
\[\mathcal C=\sum_{\omega\in\Omega,\,n\in\mathbb N}\mathcal C_{\omega,n}\]
has the required openness property.
(The algebraic sum of infinitely many $0$-neighborhoods is always open, and so is any bigger algebraic sum.)
Yet
\[\intp L(\mathcal C, \mu)\subset
\overline{
\sum_{\omega\in\Omega,\,n\in\mathbb N}\mathcal T_{\omega,n} }
\subset\overline{\mathcal T}\subset\mathcal T'.\]
That proves our statement.\end{proof}
\end{lemma}
The finite type semivariation conditions are a quite restrictive:
\begin{lemma}\label{lem:FINCLO}
\halfquote
Suppose that $\mu$ is of $\sigma$-locally   finite semivariation with respect to the pairing $L$.

Then, we claim, for any envelope $\mathcal A$ of
sets in $\mathcal  V$ with respect to $\mathfrak S$
and any neighborhood $0\in\mathcal T\subset \mathcal Z$
 there exists a pointed envelope
$\mathcal C$ such that
\[\overline{\mathcal A^x}\subset(\mathcal A+\mathcal C)^x\]
for all $x\in X$ but
\[\intp L(\mathcal  A+\mathcal  C,\mu)\subset \intp L(\mathcal  A,\mu)+\mathcal  T.\]

In particular, if $f:X\rightarrow\mathcal V$ is a function such that
\[f(x)\in\overline{\mathcal  A^x},\]
then $f$ will be approximated by $\mathcal A+\mathcal C$.

\begin{proof}
Let $0\in\mathcal T'\subset \mathcal Z$ be a neighborhood such
that $\mathcal T'+\mathcal T'\subset\mathcal T$.
Let $\mathfrak D$ be a generator system of $\mathcal D$. Then let us apply
Lemma \ref{lem:FINSIGMA}.c to $A=\bigcup\mathfrak D$ in order to obtain the
pointed envelope $\mathcal C$. One can see that $\mathcal C$ has the required properties.
\end{proof}
\end{lemma}
\begin{cor}\label{cor:INTSIGMA}
\halfquote
Suppose that $\mu$ is of $\sigma$-locally finite semivariation with respect to the pairing $L$.

Then, integration can be developed using the weaker notion of approximations,
ie. requiring only the the closure property
\[f(x)\in\overline{\mathcal  A^x}.\]
\begin{proof}
According to the Lemma \ref{lem:FINCLO} we can always pass to an only slightly bigger proper approximation.
\end{proof}
\end{cor}
\begin{lemma}\label{lem:SIGCONT}
\halfquote
Suppose that $\mu$ is of $\sigma$-locally finite semivariation with respect to the pairing $L$.

Assume that $f:X\rightarrow\mathcal V$ is integrable, $\mathcal A$ is
 an envelope of sets in $\mathcal  V$ with respect to $\mathfrak S$, such that
\[f(x)\in\overline{\mathcal A^x}\]
for all $x\in X$. Then, we claim,
\[\int L(f,\mu)\in\intp L(\mathcal A,\mu).\]
\begin{proof}
That immediately follows either from Lemma \ref{lem:FINCLO} or Corollary \ref{cor:INTSIGMA}.
\end{proof}
\end{lemma}
As a consequence of Corollary \ref{cor:INTSIGMA} that we can formulate (an other version of) Beppo Levi's theorem:
\begin{defin}
Suppose that $\mathcal  V$ is commutative topological group, $\{a_\lambda\}_{\lambda\in\Lambda}$
is a countable family of elements of $\mathcal  V$. We say that $a$ is in the limit of the sum
\[\sum_{\lambda\in\Lambda} a_{\lambda}\]
if $a$ is an accumulation point of the filter base
\[\sum_{\lambda\in\Lambda}^{\mathrm{form}}a_\lambda.\]
\end{defin}
\begin{lemma}\label{lem:FINBP}
\halfquote
Suppose that $\mu$ is of $\sigma$-locally finite semivariation with respect to the pairing $L$.

Suppose that
$f_\lambda:X\rightarrow \mathcal  V$ (${\lambda\in\Lambda}$) is a family integrable functions and
$f:X\rightarrow \mathcal  V$
is a function such that $f(x)$ is in the limit of
\[\sum_{\lambda\in\Lambda}f_\lambda(x)\]
for all $x\in X$ except on a weakly negligible set $C$.

Also assume that for each neighborhood $0\in\mathcal  T\subset\mathcal  Z$ there exists
a finite set $\Xi\subset\Lambda$
such that
\[\sum_{\lambda\in\Lambda\setminus\Xi}\svar(L,f_\lambda,\mu)\subset\mathcal  T.\]

Then, we claim, $f$ is integrable and
\[\int L(f,\mu)=\sum_{\lambda\in\Lambda}^{\ovline}\int L(f_\lambda,\mu).\]
Furthermore
\[\svar(L,f,\mu)\subset\overline{\sum_{\lambda\in\Lambda}\svar(L,f_\lambda,\mu}).\]

\begin{proof}
First of all we can forget about the negligible set $C$, because we can set the
functions to be $0$ on $C$ without having any consequences with respect to integration or
semivariation.

One can translate the underlying condition according to Lemma \ref{lem:SVTRANS}.
Then one can use the proof of Lemma \ref{lem:ALGBP}, because weaker approximation sense is
enough.

The statements about semivariation can be proven in the following way:
Let $\mathcal C_\lambda$ pointed envelopes approximating $f_\lambda$. The $\mathcal C$ be their sum.
If $\intp L(\mathcal C_\lambda,\mu)$ was quite close to $\svar(L,f_\lambda,\mu,)$ then
\[\intp L(C,\mu)\subset \overline{\sum_{\lambda\in\Lambda}\svar(L,f_\lambda,\mu})+\mathcal T\]
can be assumed with an arbitrarily small $0$-neighborhood $\mathcal T$.
One the other hand, $\mathcal C$ approximates every restriction of $f$ in the weak sense,
hence
\[\svar(L,f,\mu)\subset \intp L(\mathcal C,\mu)\subset
\overline{\sum_{\lambda\in\Lambda}\svar(L,f_\lambda,\mu})+\mathcal T \]
proves our statement.
\end{proof}
\end{lemma}
The backward feature of  this lemma is that we had to assume that a the function $f$ exists.
Under more restrictive conditions the existence of  $f$ follows:
\begin{lemma}\label{lem:STRONGBP}
\halfquote
Suppose that $\mu$ is of $\sigma$-locally finite semivariation with respect to the pairing $L$.

Assume that $\mathcal  V$ is a Fr\'echet ($\mathtt M_1$) locally convex space.
Suppose that
$f_n:X\rightarrow \mathcal  V$ (${\lambda\in\Lambda}$) is a countable indexed
family of [strongly] integrable functions

Assume that for each $0\in\mathcal  T\subset Z$ and seminorm $p$
on $\mathcal  V$ there exists a finite subset  $\Xi\subset\Lambda$
and envelopes $\mathcal  C_\lambda$ ($\lambda\in\Xi$) such that
\[\{(v,x)\,:\, p(v)<p(f_\lambda(x)) \}\subset\mathcal  C_\lambda\] and
\[\sum_{\lambda\in\Lambda\setminus\Xi} \intp L(\mathcal  C_\lambda,\mu)\subset\mathcal  T.\]
Then
\[f(x)=\sum_{\lambda\in\Lambda}^{\ovline} f_n(x)\]
exists except on a strongly negligible set $C$, $f$ is [strongly] integrable and
\[\int L(f,\mu)=\sum_{\lambda\in\Lambda}^{\ovline}\int L( f_n,\mu).\]
\begin{proof}
Let $ C_p$ be the set of non convergence points in $p$.
Then, for a non-convergence point $x\in C$ and any finite set $\Xi\subset\Lambda$ the sum
\[\sum^{\ovline}_{\lambda\in\Lambda\setminus\Xi}p(f_\lambda(x))\]
diverges. Hence, form the properties of $\mathcal C_\lambda$ as above,
it follows that
\[\mathcal V\times C\subset \sum_{\lambda\in\Lambda\setminus\Xi} C_\lambda.\]
That proves the strong negligibility of $C_p$. Being $\mathcal  V$ a Fr\'echet space,
the set all non-convergence points $C$ is a countable union of some $C_p$'s, hence it is
still strongly negligible.

To prove the statements about [strong] integrability,
 it is sufficient to consider envelopes
\[\mathcal A=\sum_{\lambda\in\Xi} A_\xi+\sum_{\lambda\in\Lambda\setminus\Xi} C_\lambda,\]
where $\mathcal A_\xi$ [strongly] approximates $f_\xi$ and $C_\lambda$ is as above.
\end{proof}
\end{lemma}

\paragraph{\textbf{B. Consequences of separability}}
\begin{defin}\label{def:SEPAR}
a.) A topological space $X$ is separable if it contains a countable dense set.

b.) A commutative topological group $\mathcal V$ is weakly separable if for each neighborhood
$0\in\mathcal T\subset\mathcal V$ there is a countable set $R_{\mathcal T}\subset \mathcal V$
such that
\[R_{\mathcal T}+\mathcal T=\mathcal V. \]
\end{defin}

\begin{lemma}\label{lem:FINNEGL}
\halfquote

Suppose that at least one of the following two conditions holds:
\begin{itemize}

\item[x).]
The measure $\mu$ is of $\sigma$-locally finite semivariation with respect to the pairing $L$;
and $\mathcal  V$ is a separable  space.

\item[y).]
The measure $\mu$ is of locally $\sigma$-finite semivariation with respect to the pairing $L$;
and $\mathcal  V$ is a weakly separable  space.
\end{itemize}
Then, we claim, any weakly $L,\mu$-negligible set is strongly $L,\mu$-negligible.
\begin{proof}
Let $C$ be a weakly negligible set.
The restriction of a nontrivial constant function for $C$ is integrable; hence
we can assume that there is as set $A\in\boldsymbol\Sigma\mathfrak S$.

Let $0\in\mathcal T\subset \mathcal Z$ be an arbitrary neighborhood.
Let $0\in\mathcal T'\subset \mathcal Z$ be a neighborhood such that
$\mathcal T'+\mathcal T'+\mathcal T'\subset \mathcal T$.

In case (x): Let $\mathcal C$ be a pointed envelope as in Lemma \ref{lem:FINSIGMA}.c.
Let $R$ be a countable dense set in $\mathcal V$.

In case (y): Let $\mathcal C$ be a pointed envelope as in Lemma \ref{lem:FINSIGMA}.b,
and let
\[R=\bigcup_{\lambda\in\Lambda} R_{\mathcal G_\lambda}.\]
(cf. Definition \ref{def:SEPAR}).

Having $\mathcal C$ and $R$ defined let us continue as follows:

 Let us divide $\mathcal T'$ by $R$.
for each $r\in R$ we can consider the function $r|_C$. According to our assumptions there is a
pointed envelope $\mathcal C_r$ such that it approximates $r|_C$ but
\[\intp L(\mathcal C_r,\mu)\subset\mathcal T_r'.\]
Let \[\tilde{\mathcal C}=\sum_{r\in R}\mathcal C_r.\]
Then
\[\intp L(\tilde{\mathcal C},\mu)\subset \overline{\sum_{r\in R}\mathcal T'_r }\subset\overline{\mathcal T'}
\subset \mathcal T'+\mathcal T'.\]
One the other hand,
\[\var(\mathcal V,C)\subset \var(R,C)+\mathcal {\mathcal C} \subset
\tilde{\mathcal C}+\mathcal {\mathcal C}.\]
(only the the last term is  an envelope).
Then
\[\mathcal V\times C\subset\tilde{\mathcal C}+\mathcal {\mathcal C}\]
while
\[\intp L(\tilde{\mathcal C}+\mathcal {\mathcal C},\mu)=\overline{
\sum_{r\in R} \intp L(\mathcal C_r,\mu)+\intp L(\mathcal C,\mu)}\subset
\overline{\sum_{r\in R} \mathcal T'_r+\mathcal T'}\subset
\mathcal T'+\mathcal T'+\mathcal T'\subset\mathcal T.\]
That proves the strong negligibility of $C$.
\end{proof}
\end{lemma}
\begin{lemma}\label{lem:FINEXT}
\halfquote

Suppose that at least one of the following two conditions holds:
\begin{itemize}

\item[x).]
The measure $\mu$ is of $\sigma$-locally finite semivariation with respect to the pairing $L$;
and $\mathcal  V$ is a separable  space.

\item[y).]
The measure $\mu$ is of locally $\sigma$-finite semivariation with respect to the pairing $L$;
and $\mathcal  V$ is a weakly separable  space.
\end{itemize}

Consider the integrable extension
\[\hat\mu:\mathfrak S^\mu\rightarrow\mathbb R.\]

Then, we claim:

a.) If $\mathcal  A$ is  an envelope over $\mathfrak S^\mu$ and  then
for each neighborhood $0\in\mathcal  T\subset \mathcal  Z$ there is  an envelope $\tilde{\mathcal  A}$
 over $\mathfrak S$
such that
\[\mathcal  A\subset\tilde{ \mathcal  A}\]
but
\[ \intp L(\tilde{\mathcal  A},\mu)\subset \intp L(\mathcal  A,\hat\mu)+\mathcal  T.\]

b.)
A function $f:X\rightarrow\mathbb R$ is (locally) $\mu$-integrable if and only if
it is (locally) $\hat\mu$-integrable. In this case the integrals and the semivariation
sets are the same.
\begin{proof}
a.)
Let $0\in\mathcal  T'\subset \mathcal  Z$ be a neighborhood such that
$\mathcal T'+\mathcal T'+\mathcal T'+\mathcal T'-\mathcal T'\subset\mathcal T$.
Suppose that $\mathfrak D$ is a generator set of $\mathcal A$ and
$s$ is (the graph of) a step-function contained in $\mathcal A$.

Every integrable set from $D\in\mathfrak D\subset\mathfrak S^\mu$
is contained in a set $\boldsymbol\Sigma\mathfrak S$, hence we find that there is a set
$A\in \boldsymbol\Sigma\mathfrak S$ such that $\bigcup\mathfrak D\subset A$.

In case (x): Let $\mathcal C$ be a pointed envelope over $\mathfrak S$ as in Lemma \ref{lem:FINSIGMA}.c.
Let $R$ be a countable dense set in $\mathcal V$.

In case (y): Let $\mathcal C$ be a pointed envelope over $\mathfrak S$ as in Lemma \ref{lem:FINSIGMA}.b,
and let
\[R=\bigcup_{\lambda\in\Lambda} R_{\mathcal G_\lambda}.\]
(cf. Definition \ref{def:SEPAR}).

Having $\mathcal C$ and $R$ defined let us continue as follows:

According to Lemma \ref {lem:MECOS},
\[\intp L(\mathcal C,\hat\mu)=\intp L(\mathcal C,\mu)\subset\mathcal T'.\]
We define
\[\mathcal R=s\cup ((\mathcal A-\mathcal C)\cap(R\times X)).\]
Then one can see that $\mathcal R$ is  an envelope over $\mathfrak S^\mu$ and
\[\intp L(\mathcal R,\hat\mu)\subset\intp L(\mathcal A-\mathcal C,\hat\mu). \]
But then, according to Lemma \ref{lem:MECONTINT} there is  an envelope $\tilde{\mathcal R}\supset\mathcal R$
over $\mathfrak S$ such that
\[\intp L(\tilde{\mathcal R},\mu)\subset
\intp L(\mathcal R,\hat\mu)+\mathcal T'\subset\intp L(\mathcal A-\mathcal C,\hat\mu)+\mathcal T'.\]
Hence,
\[\intp L(\tilde{\mathcal R}+\mathcal C,\mu)\subset
\intp L(\mathcal A-\mathcal C,\hat\mu)+\mathcal T'+\intp L(\mathcal C,\mu)+\mathcal T'\]
\[\subset \intp L(\mathcal A,\hat\mu)-\mathcal T'+\mathcal T'+\mathcal T'+\mathcal T'+\mathcal T'
\subset \intp L(\mathcal A,\hat\mu)+\mathcal T.\]
One the other hand, for $x\in A$
\[\mathcal A^x\subset {((\mathcal A^x-\mathcal C^x)\cap R)+\mathcal C^x}
=\mathcal R^x+\mathcal C^x\subset \tilde{\mathcal R}^x+\mathcal C^x
=(\tilde{\mathcal R}+\mathcal C)^x.\]
This proves that $ \tilde{\mathcal A}=\tilde{\mathcal R}+\mathcal C$ is a good choice in our statement.

b.) That is an immediate consequence of part a.
\end{proof}
\end{lemma}
\begin{remark}
One very reasonable policy to adopt would be the following:
We should require that each envelope $\mathcal A$ is open in its fibers $\mathcal A^x$;
and we could demand that every pairing $L$ should be $\sigma$-locally final.
In general, that would restrict our attention to cases when $X\in\boldsymbol \Sigma\mathfrak S$
but extended integration circumvents all problems coming from here.

On the other hand, this convention would not invalidate any of our previous work: before this
section we essentially considered every $\mathcal V$ with the discrete topology.

However, some of our earlier statements would require extra work to adopt them to this
more general setting. Especially, formal integration would become less transparent.
\end{remark}
The following statement is a bit more geometric:
\begin{lemma}\label{lem:FINOPEN}
\halfquote

Suppose that $X\in\boldsymbol\Sigma\mathfrak S$, and $X$ is a topological space.
Also suppose that $\mu$ is of locally finite semivariation with respect to the pairing $L$;
and $\mathcal  V$ is a weakly separable  space.

Assume that the following condition holds:
\begin{itemize}
\item[\texttt{(Z1)}]  For each $E\in\mathfrak S$,
there is a countable subset $\mathfrak B_E\subset\mathfrak S$ such that
\[E\setminus E^\circ\subset\bigcup_{S\in\boldsymbol\Sigma_0(\mathfrak B_E)} S^\circ.\]
\item[\texttt{(Z2)}]
For each $E\in\mathfrak S$, $c\in\mathcal V$, and neighborhood $0\in\mathcal T\subset\mathcal Z$
there is a set $B\in\boldsymbol\Sigma\mathfrak S$ such that
\[E\setminus E^\circ\subset(E\cup B)^\circ\]
and
\[\svar(L,c,\bigcup \mathfrak B,\mu)\subset\mathcal T. \]
\end{itemize}

Then, we claim:

 If $\mathcal  A$ is  an envelope over $\mathfrak S$ and  then
for each neighborhood $0\in\mathcal  T\subset \mathcal  Z$ there is  an envelope $\tilde{\mathcal  A}$
 over $\mathfrak S$
such that
\[\mathcal  A\subset\tilde{ \mathcal  A};\]
and
\[\tilde{ \mathcal  A}\text{ is an open subset of $\mathcal V\times X$}\]
but
\[ \intp L(\tilde{\mathcal  A},\mu)\subset \intp L(\mathcal  A,\mu)+\mathcal  T.\]
\begin{proof}
Part 1. In this part we prove the following statement:
For each neighborhood $0\in\mathcal  U\subset \mathcal  Z$
there is a pointed envelope $\mathcal C$
such that
\[\{0\}\times X\subset (\mathcal C)^\circ \]
but
\[\intp L(\mathcal C,\mu)\subset\mathcal U.\]
Moreover, we can assume that there is a countable family of $0$-neighborhoods
\[\mathcal H_{\omega}\]
($\omega\in\Omega$) such that each point of $x\in X$ has neighborhood $\mathcal D$ in $X$
and an appropriate $\omega\in\omega$ such that
\[(0,x)\in\mathcal D\times\mathcal H_{\omega}\subset\mathcal C.\]

The proof goes as follows:

Let $0\in\mathcal T'\subset\mathcal V$ be a neighborhood such that
$\mathcal T'+\mathcal T'+\mathcal T'\subset\mathcal U$.
According to Lemma \ref{lem:FINSIGMA}.b there is a countable decomposition
\[X=\bigcup_{\lambda\in\Lambda}A_\lambda\]
and open 0-neighborhoods $\mathcal G_\lambda$ such that with the envelope
\[\tilde{\mathcal C}=\sum_{\lambda\in\Lambda}\var(\mathcal G_\lambda,A_\lambda)\]
we have
\[\intp L(\tilde{\mathcal C},\mu)\subset\mathcal T'.\]
Let $\mathfrak A=\{A_\lambda\,:\,\lambda\in\Lambda\} $.

We can see that $\tilde{\mathcal C}$ is already a neighborhood of each point $(0,x)$ such that
\[x\notin \bigcup_{A\in\mathfrak A} A\setminus A^\circ.\]
According to condition \texttt{(Z1)} let
\[\mathfrak B=\bigcup_{A\in\mathfrak A}\mathfrak B_A.\]
Let us divide $\mathcal T'$ by the countable set $\boldsymbol\Sigma_0\mathfrak B$.
According to Lemma \ref{lem:FINSIGMA}.a
for each element $B\in \boldsymbol\Sigma_0\mathfrak B$ there is an $0$-neighborhood
$\mathcal U_B$ such that
\[\intp L(\var(\mathcal U_B,B),\mu)\subset\mathcal T'_B.\]
Then let
\[\mathcal C=\tilde{\mathcal C}+\sum_{B\in \boldsymbol\Sigma_0\mathfrak B}\var(\mathcal U_B,B).\]
By construction this will be a neighborhood of any point of $\{0\}\times X$, yet
\[\intp L(\mathcal C,\mu)=\overline{\intp L(\tilde{\mathcal C},\mu)
+\sum_{B\in \boldsymbol\Sigma_0\mathfrak B}\intp L(\var(\mathcal U_B,B),\mu)}\]
\[\subset\overline{\mathcal T'+
\sum_{B\in \boldsymbol\Sigma_0\mathfrak B}\mathcal T'_B }\subset\mathcal T'+\mathcal T'+\mathcal T'
\subset\mathcal U.\]

As for the neighborhoods $\mathcal H_\omega$, we can see that the zusammensetzung of the $\mathcal G_\lambda$'s
and $\mathcal U_B$'s will be good.

Part 2. In this part we prove the following statement:
For every neighborhood $0\in\mathcal  U\subset \mathcal  Z$,
each countable subset $R\subset\mathcal V$,
and envelope
\[\mathcal R\subset R\times X\]
there is  an envelope $\tilde{\mathcal R}\supset \mathcal R$ such that
for all $r\in R$
\[{}^r\mathcal R\subset( {}^r\tilde{\mathcal R} )^\circ\]
but
\[\intp L(\tilde{\mathcal R},\mu)\subset\intp L({\mathcal R},\mu)+\mathcal U. \]
The proof goes as follows:

Let $0\in\mathcal T'\subset\mathcal V$ be a neighborhood such that
$\mathcal T'+\mathcal T'\subset\mathcal U$.
Let $\mathfrak D$ be a generator system of $\mathcal R$.
Let $s\subset\mathcal R$ be a step-function, with values from the finite set $K\subset\mathcal V$.
Let us divide $\mathcal T'$ by
$R\times K\times\mathfrak D$. According to assumption \texttt{(Z2)}, for $r\in R$, $k\in K$,$D\in\mathfrak D$
let $B_{r,k,D}\in\boldsymbol\Sigma\mathfrak S$ such that
\[D\setminus D^\circ\subset(D\cup B_{r,k,D})^\circ\]
and
\[\svar(L,r-k,B_{r,D},\mu)\subset \mathcal T'_{r,k,D}.\]
Then
\[\tilde{\mathcal R}=\mathcal R+\sum_{r\in R,\,D\in\mathfrak D} \var(r-k,B_{r,k,D})  \]
satisfies the local openness yet
\[\intp L(\tilde{\mathcal R},\mu)=\overline{\intp L({\mathcal R},\mu)+
\sum_{r\in R,\,k\in K,\,D\in\mathfrak D}\svar(L,r-k,B_{r,k,D},\mu)}\]
\[\subset\overline{\intp L({\mathcal R},\mu)+\sum_{r\in R,\,k\in K,\,D\in\mathfrak D}\mathcal T'_{r,k,D}}
\subset\intp L({\mathcal R},\mu)+\mathcal T'+\mathcal T'\subset \intp L({\mathcal R},\mu)+\mathcal U.\]

Part 3. Here we prove the following statement:

For every neighborhood $0\in\tilde{\mathcal  T}\subset \mathcal  Z$, and envelope
$\mathcal A$ there is  an envelope $\tilde{\mathcal A}$ such that
\[\mathcal A\subset(\tilde{\mathcal A})^\circ\]
but
\[ \intp L(\tilde{\mathcal  A},\mu)\subset \intp L(\mathcal  A,\mu)+\tilde{\mathcal  T}.\]

The proof is as follows:

Let $0\in\mathcal U\subset\mathcal V$ be a neighborhood such that
$-\mathcal U+\mathcal U+\mathcal U+\mathcal U\subset\tilde{\mathcal T}$.
Let $s$ be a step-function contained in $\mathcal A$, which takes values from a
finite set $K\subset\mathcal V$.
Let $\mathcal C$ as in Part 1.
Let
\[\mathcal R=\bigcup_{\omega\in\Omega} R_{\mathcal H_\omega}\]
(cf. Definition \ref{def:SEPAR}).
Let
\[\mathcal R=s\cup\Bigl((\mathcal A-\mathcal C)\cap(R\times X)\Bigr).\]
Let $\tilde{\mathcal R}$ as in Part 2.

Then, by construction,
\[\tilde{ \mathcal R}+ \mathcal C\]
is  an envelope such that for each $(v,x)\in \mathcal A$ there is  a neighborhood $x\in\mathcal D$
such that $\{v\}\times D\subset \tilde{ \mathcal R}+ \mathcal C$.

Hence, the envelope
\[ \tilde{\mathcal  A}= (\tilde{ \mathcal R}+ \mathcal C)+\mathcal C\]
will actually be a neighborhood of $\mathcal A$.

One the other hand,
\[\intp L(\tilde{\mathcal  A},\mu)\subset\overline{ \intp L(\tilde{\mathcal  R},\mu)+
\intp L(\tilde{\mathcal  C},\mu)+\intp L(\tilde{\mathcal  C},\mu)}\subset
\]
\[\subset \overline{ \intp L({\mathcal  A}-\mathcal C,\mu)+
\intp L(\tilde{\mathcal  C},\mu)+\intp L(\tilde{\mathcal  C},\mu)}\subset  \]
\[ \subset\intp L({\mathcal  A},\mu)-\mathcal U+\mathcal U+\mathcal U+\mathcal U
\subset\intp L({\mathcal  A},\mu)+\tilde{\mathcal T}.\]

Part 4. Here we prove the statement of the lemma.

Let $0\in\tilde{\mathcal T}\subset\mathcal V$ be a neighborhood such that
$\tilde{\mathcal T}+\tilde{\mathcal T}\subset\mathcal T$.
Let us divide $ \tilde{\mathcal T}$ by $\mathbb N$. Then, according to part 3, there is an  sequence
of envelopes
\[\tilde{\mathcal A}_n\]
such that
\[\tilde{\mathcal A}_0=\mathcal A;\qquad
\tilde{\mathcal A}_n\subset (\tilde{\mathcal A}_{n+1})^\circ ;\]
and
\[\intp L(\tilde{\mathcal A}_{n+1},\mu)\subset\intp L(\tilde{\mathcal A}_{n},\mu)+\tilde{\mathcal T}_n. \]

Then we see that for any $n\in\mathbb N$
\[ \intp L(\tilde{\mathcal A}_{n},\mu)\subset  \intp L(\mathcal A,\mu)+
\sum_{n\in\mathbb N } \tilde{\mathcal T}_n  \subset\intp L(\mathcal A,\mu)+\tilde{\mathcal T}. \]
In particular,
\[\overline{\bigcup_{n\in\mathbb N} \intp L(\tilde{\mathcal A}_{n},\mu) }
\subset\overline{\intp L(\mathcal A,\mu)+\tilde{\mathcal T}}\subset
\intp L(\mathcal A,\mu)+\tilde{\mathcal T}+\tilde{\mathcal T}
\subset\intp L(\mathcal A,\mu)+{\mathcal T}.\]
According to Lemma \ref{lem:UNIOLIM}, however, this implies that for
\[\tilde{\mathcal A}=\bigcup_{n\in\mathbb N}\tilde{\mathcal A}_{n}\]
we have
\[ \intp L(\tilde{\mathcal A},\mu)\subset\intp L(\mathcal A,\mu)+{\mathcal T}.\]
On the other hand, $\mathcal A\subset\tilde{\mathcal A}$ and
$\tilde{\mathcal A} $ is a neighborhood of its own points.
\end{proof}
\end{lemma}
\begin{remak}
If, in addition, the open sets of $X$ are from $\boldsymbol\Sigma\mathfrak S$ then
the term ``envelope'' can effectively be replaced by the term ``open subset of $\mathcal V\times X$
which contains the graph of a step-function''.
\end{remak}
\begin{lemma}\label{lem:LOCEX}
\halfquote

If $X=\mathbb R$ and $\mathfrak S=\mathfrak I$ is the system of the possibly degenerate finite
intervals in $\mathbb R$ then conditions \texttt{(Z1)} and \texttt{(Z2)} from Lemma \ref{lem:FINOPEN}
hold (even if not the other conditions).

\begin{proof} Condition \texttt{(Z1)} is trivial. Condition \texttt{(Z2)}  follows from
applying Lemma \ref{lem:MEAVAR} to set systems
\[\{[a-2^{-n},a-2^{-n-1})\,:\,n\in\mathbb N\}\cup\{\{a\}\}\cup\{(a+2^{-n-1},a+2^{-n}]\,; \,n\in\mathbb N\}\]
($a\in\mathbb R$).
\end{proof}
\end{lemma}

\begin{lemma}\label{lem:LOCOPEN}
\halfquote

Suppose that $X\in\boldsymbol\Sigma\mathfrak S$, and $X$ is a topological space.
Also suppose that $\mu$ is of locally finite semivariation with respect to the pairing $L$;
and $\mathcal  V$ is a weakly separable topological vector space,
$\mathcal Z$ is a locally convex space.

Assume that the following condition holds:
\begin{itemize}
\item[\texttt{(Z1)}]  For each $E\in\mathfrak S$,
there is a countable subset $\mathfrak B_E\subset\mathfrak S$ such that
\[E\setminus E^\circ\subset\bigcup_{S\in\boldsymbol\Sigma_0(\mathfrak B_E)} S^\circ.\]
\item[\texttt{(Z2)}]
For each $E\in\mathfrak S$, $c\in\mathcal V$, and neighborhood $0\in\mathcal T\subset\mathcal Z$
there is a set $B\in\boldsymbol\Sigma\mathfrak S$ such that
\[E\setminus E^\circ\subset(E\cup B)^\circ\]
and
\[\svar(L,c,\bigcup \mathfrak B,\mu)\subset\mathcal T. \]
\end{itemize}

Then, we claim:

 If $\mathcal  A$ is  an envelope over $\mathfrak S$ and  then
for each neighborhood $0\in\mathcal  T\subset \mathcal  Z$ there is  an envelope $\tilde{\mathcal  A}$
 over $\mathfrak S$
such that
\[\mathcal  A\subset\tilde{ \mathcal  A};\]
\[\tilde{ \mathcal  A}\text{ is an open subset of $\mathcal V\times X$};\]
\[\tilde{\mathcal  A}=\conv\tilde{\mathcal  A};\]
but
\[ \intp L(\tilde{\mathcal  A},\mu)\subset \conv\intp L(\mathcal  A,\mu)+\mathcal  T.\]
\begin{proof}
Let $0\in\mathcal U\subset\mathcal Z$ be a neighborhood such that
$\mathcal U+\mathcal U+\mathcal \subset \mathcal T.$
Let is divide $\mathcal U$ by $\mathbb N$. We can assume that the neighborhoods $\mathcal U_n$
are convex. The let us define the envelopes $\mathcal A_n$ by the following manner:
Let $\mathcal A_0=\conv \mathcal A$. Then let us define $\mathcal A_{n+1}$ such that
\[\conv \mathcal A_n\subset(\mathcal A_{n+1})\]
but
\[\intp L(\mathcal A_{n+1},\mu)\subset \intp (\conv \mathcal A_n,\mu)+\mathcal U_n.\]
Then, from Lemma \ref{lem:APPCONV} one can prove that for every $n\in\mathbb B$
\[\intp L(\mathcal A_{n},\mu)\subset
\overline{\conv\intp ( \mathcal A,\mu)+\sum_{n\in\mathbb N}\mathcal U_n }
\subset\overline{\conv\intp ( \mathcal A,\mu)+\mathcal U }.\]
Let
\[\tilde{\mathcal A}=\bigcup_{n\in\mathbb N}\mathcal A_n.\]
That envelope will will satisfy our requirements; even the last one:
According to Lemma \ref{lem:UNIOLIM}
\[\intp L(\mathcal A,\mu)
\subset\overline{\conv\intp ( \mathcal A,\mu)+\mathcal U }
\subset{\conv\intp ( \mathcal A,\mu)+\mathcal U+\mathcal U }
\subset \conv\intp ( \mathcal A,\mu)+\mathcal T.\]
\end{proof}
\end{lemma}
\begin{remak}
If the assumptions of Lemma \ref{lem:LOCOPEN} hold and we are in the situation when
$X=\mathbb R$, $\mathfrak S=\mathfrak I$ then one easily see that the following statement holds:

For such  an envelope $\tilde{\mathcal A}$ as in Lemma \ref{lem:LOCOPEN}, we claim,
\[\intp L(\tilde{\mathcal A},\mu)=\]\[
\overline{ \biggl\{\int L(f,\mu)\,:\,f\subset \tilde{\mathcal A};
\,f:\mathbb R\rightarrow\mathcal V \text{ is continuous, compactly supported}\biggr\}}.\]
That indicates that we can geometrize the situation to a great degree.

However, we will not dwell on that issue here; the closest thing in this direction will
be pursued in Section \ref{sec:EXTENS}, although in a rather special case.

\end{remak}
\newpage\section{Product measures}\label{sec:PRODMEAS}

Measures of $\sigma$-locally finite semivariation are also very pleasant with respect to products.

\begin{defin}\label{def:PRODEF}
Suppose that $\mathfrak R$ and  $\mathfrak S$ are interval systems.

i.) We define
\[\mathfrak R\boxtimes\mathfrak S=\{B\times A\,:\,B\in\mathfrak R,\,A\in\mathfrak S\}.\]

ii.) Moreover, suppose that $\mathcal V$, $\mathcal W$, $\mathcal Z$
are commutative topological groups,
$\nu:\mathfrak R\rightarrow\mathcal  V$ and $\mu:\mathfrak S\rightarrow\mathcal  W$
are measures,  $L:\mathcal  V\times\mathcal  W\rightarrow\mathcal  Z$ is a biadditive pairing.
We define
\[\mathfrak R\boxtimes_{L,\nu,\mu}\mathfrak S\]
as the set of all \[B\times A\] such that $B\in \mathfrak R$, $A\in\mathfrak S$, and
for each pair of decompositions
\[B=\bigcup_{\gamma\in\Gamma}^\updisjoint B_\gamma\qand
A=\bigcup_{\lambda\in\Lambda}^\updisjoint A_\lambda\]
in $\mathfrak R$ and $\mathfrak S$, respectively, the sum
\[\sum_{(\gamma,\lambda)\in\Gamma\times\Lambda}^{\ovline}L(\nu(B_\gamma),\mu(A_\lambda))\]
exists.

(If $B$ or $A$ is the empty set then the condition holds trivially; in what follows we
 will mostly ignore this trivial case.)
\end{defin}
\begin{lemma} \label{lem:PREDEF}
Let $\mathfrak R$ and  $\mathfrak S$ are interval systems,
$\mathcal V$, $\mathcal W$, $\mathcal Z$
be commutative topological groups,
$\nu:\mathfrak R\rightarrow\mathcal  V$ and $\mu:\mathfrak S\rightarrow\mathcal  W$
be measures, and let  $L:\mathcal  V\times\mathcal  W\rightarrow\mathcal  Z$ be a biadditive pairing
which is continuous in both variables (separately).

Assume that $B\in \mathfrak R$, $A\in\mathfrak S$,
$B\times A\in\mathfrak R\boxtimes_{L,\nu,\mu}\mathfrak S$, and
\[B=\bigcup_{\gamma\in\Gamma}^\updisjoint B_\gamma\qand
A=\bigcup_{\lambda\in\Lambda}^\updisjoint A_\lambda\]
in $\mathfrak R$ and $\mathfrak S$, respectively. Then, we claim,
\[L(\nu(B),\mu(A))=\sum_{(\gamma,\lambda)\in\Gamma\times\Lambda}^{\ovline}L(\nu(B_\gamma),\mu(A_\lambda)).\]
\begin{proof}
That follows from
\[\sum_{(\gamma,\lambda)\in\Gamma\times\Lambda}^{\ovline}L(\nu(B_\gamma),\mu(A_\lambda))
=\sum_{\gamma\in\Gamma}^{\ovline}\biggl(\sum_{\lambda\in\Lambda}^{\ovline}
L(\nu(B_\gamma),\mu(A_\lambda)) \biggr)=\]
\[=\sum_{\gamma\in\Gamma}^{\ovline}L(\nu(B_\gamma),\mu(A))=L(\nu(B),\mu(A))\]
(we used contraction of sums and continuity).
\end{proof}
\end{lemma}
\begin{lemma}\label{lem:PRODEF}
Suppose that $\mathfrak R$ and $\mathfrak S$ are interval systems. Then, we claim:

a.)  $\mathfrak R\boxtimes\mathfrak S$ is an interval system.

Moreover, if, $\mathcal V$, $\mathcal W$, $\mathcal Z$
are commutative topological groups,
$\nu:\mathfrak R\rightarrow\mathcal  V$ and $\mu:\mathfrak S\rightarrow\mathcal  W$
are measures, and   $L:\mathcal  V\times\mathcal  W\rightarrow\mathcal  Z$
is a biadditive pairing then the following holds:

b.) $\mathfrak R\boxtimes_{L,\nu,\mu}\mathfrak S$ is an interval system.

c.) If $S_1\in \mathfrak R\boxtimes\mathfrak S$, $S_2\in\mathfrak R\boxtimes_{L,\nu,\mu}\mathfrak S$
and $S_1\subset S_2$ then $S_1\in\mathfrak R\boxtimes_{L,\nu,\mu}\mathfrak S$.

d.) The equality
\[\mathfrak R\boxtimes\mathfrak S=\mathfrak R\boxtimes_{L,\nu,\mu}\mathfrak S,\]
holds if $L$ is continuous in both variables
and $\mu$ is of locally finite semivariation with respect to the pairing.

\begin{proof}
a.) Consider the identities
\[(B\times A)\cap(D\times C)=(B\cap D)\times(A\cap C)\]and
\begin{multline}
(B\times A)\setminus(D\times C)=\Bigl((B\setminus D)\times (A\setminus C)\Bigr)\,\dot\cup\,\\
\,\dot\cup\,\Bigl((B\cap D)\times (A\setminus C)\Bigr)\,\dot\cup\,
\Bigl((B\setminus D)\times (A\cap C)\Bigr).\notag
\end{multline}
As the sets $B\cap D$, $A\cap C$, $B\setminus D$, $A\setminus C$ all decompose in $\mathfrak R$
and $\mathfrak S$ respectively, it follows that $\mathfrak R\boxtimes\mathfrak S $ is an
interval system.

c.) Suppose that $ D\times C\in \mathfrak R\boxtimes\mathfrak S$,
$ B\times A\in\mathfrak R\boxtimes_{L,\nu,\mu}\mathfrak S$,
and $ D\times C \subset  B\times A$. We can assume the these sets are nonempty and
$D\subset B$, $C\subset A$.
 Now,  $B\setminus D$, $A\setminus C$ decompose in $\mathfrak R$
and $\mathfrak S$ respectively, so we extend the any exact decompositions of
$D$ and $C$ to $B$ and $A$, respectively.
Then, the big sum (as in Definition \ref{def:PRODEF}.ii) with respect to $B\times A$ will
converge, hence by the partial sum for $D\times C$ will certainly converge.

b.) The closedness in the case of $\mathfrak R\boxtimes_{L,\nu,\mu}\mathfrak S$ will follow from
point b.

d.) Let $0\in\mathcal  T\subset \mathcal  Z$ be an arbitrary neighborhood.
Let $0\in\mathcal  U\subset\mathcal  V$ be a neighborhood
such that
\[\intp L(\var(\mathcal U,A),\mu)\subset \mathcal  T.\]
Let $\Omega\subset\Gamma$ so that it is finite but
\[\svar\biggl(\bigcup_{\gamma\in\Gamma\setminus\Omega } B_\gamma,\nu\biggr)\subset\mathcal  U\]
(cf. Lemma \ref{lem:MEAVAR}).
Let us divide  $0\in\mathcal  T\subset \mathcal  Z$ by $\Omega$.

Now, for all $\omega\in\Omega$ there exists a $\Xi_\omega$ such that
\[\svar\biggr(L,\nu(B_\omega),\bigcup_{\lambda\in\Lambda\setminus\Xi_\omega }A_\lambda
,\mu\biggl)\subset\mathcal  T_\omega.\]
Let
\[\Xi=\bigcup_{\omega\in\Omega} \Xi_\omega.\]
This is a finite set.
Then $\Omega\times\Xi\subset \Gamma\times\Lambda$ is finite, yet
\[\sum_{(\gamma,\lambda)\in(\Gamma\times\Lambda)
\setminus(\Omega\times\Xi)}\{0,L(\nu(B_\gamma),\mu(A_\lambda))\}=\]
\[\sum_{\omega\in\Omega}\biggl(\sum_{\lambda\in\Lambda}\{0,L(\nu(B_\omega),\mu(A_\lambda))\}\biggr)
+\sum_{(\gamma,\lambda)\in(\Gamma\setminus\Omega )\times\Lambda}
\{0,L(\nu(B_\gamma),\mu(A_\lambda))\}\]
\[\subset
\sum_{\omega\in\Omega}\mathcal  T_\omega
+\intp L\biggl(\var\biggl(\sum_{\gamma\in\Gamma\setminus\Omega } \nu(B_\gamma),A\biggr),\mu\biggr)
\subset \mathcal T +\intp L(\var(\mathcal U,A),\mu)
\subset \mathcal  T+\mathcal  T.\]
Here $\mathcal  T$ was arbitrary, hence we proved the Cauchy property for the sum. This implies
that the sum is convergent.
\end{proof}
\end{lemma}

In order to treat the general case we need the following two technical statements:
\begin{lemma}\label{lem:ADDBAIRE}
Let $\mathcal  V$ be a commutative topological group.

Let $\{C_\lambda\}_{\lambda\in\Lambda}$ a countable family of closed pointed sets in $\mathcal  V$
such that for each neighborhood $0\in\mathcal  T\subset \mathcal  V$ there is a finite
subset $\Xi\subset\Lambda$ such that
\[\overline{\sum_{\lambda\in\Lambda\setminus\Xi}C_\lambda}\subset\mathcal  T.\]
Also suppose that $\{F_\omega\}_{\omega\in\Omega}$ a countable family of closed
sets in $\mathcal  V$ such that
\[\mathcal  U\subset\bigcup_{\omega\in\Omega}F_\omega\]
is a neighborhood of $0$ in $\mathcal  V$. Then, we claim, there exists a finite set
set $\Xi\subset\Lambda$, an element $\omega\in\Omega$, and
\[b\in\sum_{\lambda\in\Xi}C_\lambda\]
  such that
\[b+\overline{\sum_{\lambda\in\Lambda\setminus\Xi}C_\lambda}\subset F_\omega.\]
In particular,
\[\overline{\sum_{\lambda\in\Lambda\setminus\Xi}C_\lambda}\subset F_\omega-F_\omega.\]
\begin{proof}
Suppose that it is otherwise.
We can assume that $\Lambda=\Omega=\mathbb N$. By omitting a couple of terms we can assume that
\[\overline{\sum_{n\in\mathbb N}C_n}\subset\mathcal  U.\]
By induction we will define
$b_k\in\mathcal  V$ ($k\in\mathbb N$), neighborhoods a $\mathcal  T_k$ ($k\in\mathbb N$),
$\mathcal  T_{-1}=\mathcal  V$, a
sequence $s_k$ of natural numbers ($k\in\mathbb N$), $s_{-1}=-1$, such that
\[b_0+\ldots+b_k+\mathcal  T_k\subset\mathcal  U\setminus F_k\tag{c1}\]
and
\[s_{k-1}<s_k\tag{c2}\]
and
\[b_k\in\sum_{s_{k-1}<l \leq s_k}C_l\tag{c3}\]
and
\[\overline{\sum_{l>s_k}C_l}\subset\mathcal  T_k.\tag{c4}\]
Indeed, in a step to define  $b_k$, $\mathcal  T_k$, $s_k$, proceed as follows:
Then, by our indirect assumptions
\[b_0+\ldots+b_{k-1}+\overline{\sum_{l>s_{k-1}}C_l}\not\subset F_k.\]
Being the right side closed, there must be an element $b_k\in\sum_{l>s_{k-1}}C_l$
such that
\[b_0+\ldots+b_k\notin F_k,\]
or in other terms
\[b_0+\ldots+b_k\in\mathcal  U\setminus F_k.\]
Being the right side open we can certainly find an open subset $\mathcal  T_k$ satisfying
condition (c1).  Furthermore, if $s_k$ is large enough then the other conditions
are satisfied, too.

Having finished the induction process let us examine the sum
\[\sum_{k\in \mathbb N}^{\ovline} b_k.\]
First of all, the sum exists according to the Cauchy criterium.
(Cf. the assumptions of the Lemma and condition (c3).)

Furthermore, for arbitrary $k$ we find that
\[\sum_{k\in \mathbb N}^{\ovline} b_k=b_0+\ldots+b_{k}+\sum^{\ovline}_{l>k}b_k\in
b_0+\ldots+b_{k}+\mathcal  T_k\subset \mathcal  U\setminus F_k.\]
Consequently,
\[\sum^{\ovline}_{k\in \mathbb N} b_k\in\mathcal  U
\setminus\biggl(\bigcup_{k\in\mathbb N}F_k\biggr)=\emptyset. \]
That is a contradiction.
\end{proof}
\end{lemma}

\begin{lemma}\label{lem:PROTECH}
Let $\mathcal  V$, $\mathcal  W$, $\mathcal  Z$ be commutative
topological groups, $\mathfrak S\subset\mathfrak P(X)$,  $\mathfrak R\subset\mathfrak P(Y)$
be  interval systems,
$\mu:\mathfrak S\rightarrow\mathcal  W$, $\nu:\mathfrak R\rightarrow\mathcal  V$ be a measures,
$L:\mathcal  V\times\mathcal  W\rightarrow\mathcal  Z$ be a biadditive pairing, which is
continuous in both variables.
Assume that $\mu$ is of $\sigma$-locally finite semivariation with respect to the pairing.

Suppose that
\[B\times A=\bigcup_{\lambda\in\Lambda}^\updisjoint B_\lambda\times A_\lambda\]
is a countable decomposition in $\mathfrak R\boxtimes\mathfrak S$.
Also, let $\mathcal  C$ as in Definition \ref{def:LFSVAR}.c with respect to $\mathcal  T$.

Then, we claim, there exists a decomposition
\[A=\bigcup_{\omega\in\Omega}A'_\omega\]
such that for each $\omega\in\Omega$ there exists a finite subset
$\Xi_\omega\subset\Lambda$ such that
\[A_\omega'\subset A_\xi\]
for all $\xi\in \Xi_\omega$, and
\[\svar\biggl(B\setminus\bigcup_{\xi\in\Xi_\omega} B_{\xi},\nu\biggr)\subset
{(\mathcal  C-\mathcal  C+\mathcal  C)^{A_\omega'}}.\]
Moreover, in that case
\[\sum_{\omega\in\Omega}\var\biggr(\svar\biggl(B\setminus\bigcup_{\xi\in\Xi_\omega}
B_{\xi},\nu\biggr)\times A_\omega\biggr)\subset\mathcal  C-\mathcal  C+\mathcal  C\]
and
\[\intp L(\mathcal  C-\mathcal  C+\mathcal  C,\mu)\subset\mathcal  T-\mathcal  T+\mathcal  T+\mathcal  T.\]

\begin{proof} Let $\mathfrak D$ be a generator system for $\mathcal C$.
For $x\in A$ let
\[\Lambda_x=\{\lambda\in\Lambda\,:\,x\in A_\lambda\}.\]
Then \[B=\bigcup_{\lambda\in\Lambda_x}B_\lambda,\]
hence, by Lemma \ref{lem:MEAVAR},
\[C_\lambda=\svar(B_\lambda,\nu)\]
satisfies the conditions of Lemma \ref{lem:ADDBAIRE} (with $\Lambda_x$ in place of $\Lambda$).
Furthermore
\[\mathcal  C^x=\bigcup_{D\in\mathfrak D,\,x\in D}\mathcal  C^D\subset
\bigcup_{D\in\mathfrak D,\,x\in D}\overline{\mathcal  C^D}.\]
In particular, the sets $\overline{\mathcal C^D}$
 $(x\in D\in\mathfrak D)$ can be used in Lemma \ref{lem:ADDBAIRE}
in the place of $F_\omega$'s.

Hence, by Lemma \ref{lem:ADDBAIRE} there is a finite subset $\Xi_x\subset\Lambda_x$
 and $D\in\mathfrak D$, $x\in D$ such that
\[\svar\biggl(B\setminus \bigcup_{\lambda\in\Xi_x} B_\lambda,\nu\biggr)=
\svar\biggl(\bigcup_{\lambda\in\Lambda_x\setminus \Xi_x} B_\lambda,\nu\biggr)
\subset\overline{ \mathcal  C^D}-\overline{ \mathcal  C^D}
\subset \overline{ \mathcal  C^D- \mathcal  C^D}\subset\overline{ (\mathcal  C- \mathcal  C)^D}.\]
Being $\mathcal C^x$ is a neighborhood of $0$ for $x\in D$, we also have
\[\subset{(\mathcal  C-\mathcal  C+\mathcal  C)^D}.\]

Notice that in the statement above $\Lambda_x$ can be substituted by any $\Lambda_y$
with.
\[y\in D\cap\bigcap_{\xi\in\Xi}A_{\xi}.\]
Hence we have have proved that the sets
\[D\cap\bigcap_{\xi\in\Xi}A_{\xi}\]
($D\in\mathfrak D$, $\Xi\subset\Lambda$ is finite)
such that
\[\svar\biggl(B\setminus \bigcup_{\lambda\in\Xi} B_\lambda\biggr)\subset(
 \mathcal  C-\mathcal  C+\mathcal  C)^D.\]
cover $A$.

Applying Lemma \ref{lem:REFIN} to decompose the sets above,
and further, to get a  refinement which  decomposes $A$;
we find the desired decomposition.

The last two  statements follow immediately from the properties of envelopes.
\end{proof}
\end{lemma}

\begin{lemma}\label{lem:PROD}
Suppose that $\mathfrak R$ and  $\mathfrak S$ are interval systems;
 $\nu:\mathfrak R\rightarrow\mathcal  V$ and $\mu:\mathfrak S\rightarrow\mathcal  W$
are measures, the pairing $L:\mathcal  V\times\mathcal  W\rightarrow\mathcal  Z$ is continuous
in both variables, $\nu$ and $\mu$ are of $\sigma$-locally finite semivariation with respect to
the pairing L.

Then, we claim
\[L(\nu,\mu):\mathfrak R\boxtimes_{\nu,\mu}\mathfrak S\rightarrow\mathcal  Z\]
\[L(\nu,\mu)(B\times A)=L(\nu(B),\mu(A))\]
is a measure.
\begin{proof}
Suppose
\[B\times A=\bigcup_{\lambda\in\Lambda}^\updisjoint B_\lambda\times A_\lambda\]
is a countable decomposition.

Part 0. We prove the $\sigma$-additivity statement in the case when both the sets
$\mathfrak A=\{A_\lambda\,:\,\lambda\in\Lambda\}$ and $\mathfrak B=\{B_\lambda\,:\,\lambda\in\Lambda\}$
form a decomposition of $A$ and $B$ respectively. In fact, this is just
Lemma \ref{lem:PREDEF}.

Part 1. We prove the $\sigma$-additivity statement in the special case when
the sets $\{A_\lambda\,:\,\lambda\in\Lambda\}$ or $\{B_\lambda\,:\,\lambda\in\Lambda\}$
form a decomposition of $A$ or $B$ respectively.
By symmetry reasons it is enough to prove the statement only when the $B_\lambda$'s are so nice.

From the previous lemma it is clear that for neighborhood $0\in\mathcal  T'\subset \mathcal  Z$
there is a pointed envelope $\mathcal  C'$, an exact decomposition $\{A'_\omega\}_{\omega\in\Omega}$
of $A$, and finite subsets $\Xi_\omega\subset\Lambda$ such that
\[\sum_{\omega\in\Omega}\var\biggl(
\svar\biggl(B\setminus\bigcup_{\xi\in\Xi_\omega} B_{\xi},\nu\biggr), A_\omega\biggr)\subset\mathcal  C'
\qand \intp L(\mathcal  C',\mu)\subset\mathcal  T'.\]

Let us refine all $A_\lambda$'s with $A_\omega'$'s as in Lemma \ref{lem:SETUN}.c.
We can make the important assumption that if $A_\omega'\subset A_\lambda $ then the new decomposition of
$A_\lambda$ actually contains $A_\omega'$. That, in particular applies to all $\lambda\in\Xi_\omega$.

Through the decompositions of the $A_\lambda$'s
we get a decomposition
\[B\times A=\bigcup_{\theta\in\Theta}^\updisjoint B''_\theta\times A''_\theta.\]
Notice that in this decomposition all elements $B_\xi\times A_\omega $ such $\xi\in\Xi_\omega$
are contained.

Let
\[\mathfrak J=\{B_\xi\times A_\omega  \,:\,\xi\in\Xi_\omega \},\]
\[\mathfrak E_0=\{B_\lambda\times A_\lambda \,:\,\lambda\in\Lambda\}, \]
\[\mathfrak E_1=\{B''_\theta\times A''_\theta \,:\,\theta\in\Theta\}, \]
\[\mathfrak E_2=\{B_\lambda\times A'_\omega\,:\,\lambda\in\Lambda,\,\omega\in\Omega\}. \]
Then $\mathfrak J\subset\mathfrak E_1, \,\mathfrak E_2$.

By part 0,
\[\tag{a1}\label{FFa1}
 L(\nu(B),\mu(A))=\sum^{\ovline}_{D\times C\in \mathfrak E_2}L(\nu(D),\mu(C))=\]\[
=\sum^{\ovline}_{D\times C\in \mathfrak J}L(\nu(D),\mu(C))+
\sum^{\ovline}_{D\times C\in \mathfrak E_2\setminus\mathfrak J}L(\nu(D),\mu(C)).\]
We can see that
\[\tag{a2}\label{FFa2}
\sum_{D\times C\in \mathfrak E_2\setminus\mathfrak J}\{0,L(\nu(D),\mu(C))\} \subset\]
\[\intp L\biggl(\sum_{\omega\in\Omega}
\svar\biggl(B\setminus\bigcup_{\xi\in\Xi_\omega} B_{\xi},\nu\biggr)\times A_\omega
,\mu\biggr)\subset
\intp L(\mathcal  C',\mu)\subset\mathcal T'.\]
Hence, by closure,
\[\sum^{\ovline}_{D\times C\in \mathfrak E_2\setminus\mathfrak J}L(\nu(D),\mu(C))\in
\intp L(\mathcal  C',\mu)\subset\mathcal  T'.\]
Consequently, from (\ref{FFa1}),
\[\sum_{D\times C\in \mathfrak J}^{\ovline}L(\nu(D),\mu(C))
\subset L(\nu(B),\mu(A))-\mathcal  T'.\]

On the other hand, because of the convergence of the sum right up,
 we find that there is a finite set $J\subset\mathfrak J$  such that
\[\tag{a3}\label{FFa4}
\sum_{D\times C\in J}L(\nu(D),\mu(C))+
\sum_{D\times C\in \mathfrak J\setminus J}\{0,L(\nu(D),\mu(C))\}\]
\[\subset\sum^{\ovline}_{D\times C\in \mathfrak J}L(\nu(D),\mu(C))
+\mathcal  T'
\subset L(\nu(B),\mu(A))-\mathcal  T'+\mathcal  T'.\]

Now, similarly to (\ref{FFa2}), one can see that
\[\tag{a4}\label{FFa6}
{\sum_{D\times C\in \mathfrak E_1\setminus
\mathfrak J}\{0,L(\nu(D),\mu(C))\}}\subset \intp L\biggl(\sum_{\omega\in\Omega}
\svar\biggl(B\setminus\bigcup_{\xi\in\Xi_\omega} B_{\xi},\nu\biggr)\times A_\omega
,\mu\biggr)\]\[\subset
\intp L(\mathcal  C',\mu)\subset\mathcal  T'.\]
Hence, from (\ref{FFa4}), (\ref{FFa6}), we find
\[\sum_{D\times C\in J}L(\nu(D),\mu(C))+
\sum_{D\times C\in \mathfrak E_1\setminus J}\{0,L(\nu(D),\mu(C))\}=\]
\[\biggl(\sum_{D\times C\in J}L(\nu(D),\mu(C))+
\sum_{D\times C\in \mathfrak J\setminus J}\{0,L(\nu(D),\mu(C))\}\biggr)+
\sum_{D\times C\in \mathfrak E_1\setminus\mathfrak J}\{0,L(\nu(D),\mu(C))\}\]
\[\subset(L(\nu(B),\mu(A))-\mathcal  T'+\mathcal  T')+(\mathcal  T').\]
Taking closure,
\[\tag{a5}\label{FFa7}
\overline{\sum_{D\times C\in J}L(\nu(D),\mu(C))+
\sum_{D\times C\in \mathfrak E_1\setminus J}\{0,L(\nu(D),\mu(C))\}}\subset\]\[\subset
L(\nu(B),\mu(A))-\mathcal  T'+\mathcal  T'+\mathcal  T'+\mathcal T'.\]

However, every finite partial sum of
\[\sum_{B_\lambda\times A_\lambda\in\mathfrak E_0} L(\nu(B_\lambda),\mu(A_\lambda))\]
which contains the at least the set $J'\subset\mathfrak E_0$ of the majorizing sets of $J\subset \mathfrak J$
 is contained in
the closure set of (\ref{FFa7}) above, because the sets $B_\theta''\times A_\theta''$  were obtained by
decomposing set $B_\lambda\times A_\lambda$ in the second component.
That proves
\[\sum_{B_\lambda\times A_\lambda\in J'}L(\nu(B_\lambda),\mu(A_\lambda))+
\sum_{B_\lambda\times A_\lambda\in \mathfrak E_0\setminus J'}
\{0,L(\nu(B_\lambda),\mu(A_\lambda))\}\]
\[\subset L(\nu(B),\mu(A))-\mathcal  T'+\mathcal  T'+\mathcal  T'+\mathcal T'.\]
Here $\mathcal  T'$ was arbitrary. That implies
\[\sum_{B_\lambda\times A_\lambda\in\mathfrak E_0}^{\ovline}
 L(\nu(B_\lambda),\mu(A_\lambda))= L(\nu(B),\mu(A)),\]
which was our statement.

Part 2. The general case. Let us repeat the refinement construction of Part 1, but without the
assumption that the $B_\lambda$'s decompose $B$.

For each $\omega\in\Omega$ consider a decomposition
\[B\setminus\bigcup_{\xi\in\Xi_\omega} B_\xi=\bigcup_{\psi\in\Psi_\omega}^\updisjoint B_{\omega,\psi}.\]
Make the modification that instead of $\mathfrak E_2$ we take
\[\tilde{\mathfrak E_2}=\mathfrak E_0\cup \{B_{\omega,\psi}\times A'_\omega\} \]
Here, to the analogy of the first line of (\ref{FFa1}) in Part 1 we have
\[L(\nu(B),\mu(A))=\sum_{D\times C\in \bar{\mathfrak E_2}}L(\nu(D),\mu(C)), \]
 because the $A_\omega$'s nicely decompose $A$.

After that one can repeat the proof of Part 1.
\end{proof}
\end{lemma}

\newpage\section{Classical Lebesgue theory}\label{sec:CLASSLEB}
\begin{conven}
Here we suppose that  $\mathcal  V=\mathcal  W=\mathcal  Z=\mathbb R$
and $\mathtt M$ is ordinary multiplication.
\end{conven}
~

\paragraph{\textbf{A. Main statements}}~\\

The following statements are the cornerstone of classical Lebesgue theory:
\begin{lemma}[Beppo Levi's theorem]\label{lem:CL1}
Let $\mathfrak S\subset\mathfrak P(X)$ be an interval system and
$\mu:\mathfrak S\rightarrow\mathbb R$ be a measure.

Suppose that $f_n:X\rightarrow\mathbb R$ ($n\in\mathbb N$) are integrable,
\[\sum_{n\in\mathbb N} \svar(f_n,\mu)\]
is bounded. Then
\[\sum_{n\in\mathbb N}^{\ovline} f_n(x)\]
is almost everywhere convergent, and for
\[f=\sum_{n\in\mathbb N}^{\ovline}f_n\]
(defined arbitrarily on non-convergence points)
\[\int f\,\mu=\sum_{n\in\mathbb N}^{\ovline}\int f_n\,\mu.\]
Furthermore, in this case
\[\svar(f,\mu)\subset\overline{ \sum_{n\in\mathbb N}\svar(f_n,\mu)}.\]
\begin{proof} This is Lemma \ref{lem:CONCL1} reformulated.
\end{proof}
\end{lemma}
\begin{lemma}[Constructivity]\label{lem:CL2}
Let $\mathfrak S\subset\mathfrak P(X)$ be an interval system and
$\mu:\mathfrak S\rightarrow\mathbb R$ be a measure.

If $f,g:X\rightarrow\mathbb R$ are integrable
then $f\vee g$, $f\wedge g$ are integrable. If $c\geq 0$ then
$f\wedge c$ is integrable, if $c\leq 0$ then $f\vee c$ is integrable. $|f|$ is integrable.
\begin{proof} This is Lemma \ref{lem:CONCL2} reformulated.
\end{proof}
\end{lemma}
\begin{lemma}[Monotonicity] \label{lem:CL3}
Let $\mathfrak S\subset\mathfrak P(X)$ be an interval system and
$\mu:\mathfrak S\rightarrow\mathbb R$ be a measure.
Suppose that $\mu\geq 0$.

Then, we claim, if $f:X\rightarrow\mathbb R$ is integrable then $|f|^+$ and $|f|^-$ are integrable and
\[\int f\mu\in\svar(f,\mu)\subset\biggl[\int |f|^-\mu,\int |f|^+\mu\biggr].\]
Endpoints are contained in the semivariation
\begin{proof} This is Lemma \ref{lem:CONCL3} reformulated.
\end{proof}
 \end{lemma}
\begin{lemma}\label{lem:CL4}
Let $\mathfrak S\subset\mathfrak P(X)$ be an interval system and
$\mu:\mathfrak S\rightarrow\mathbb R$ be a measure.

For a set $C\subset X$ the following are equivalent:

i.) $C$ is strongly negligible

ii.) $C$ is weakly negligible

iii.)
For each $\varepsilon>0$
there is a countable indexed family of sets $\{E_\lambda\}_{\lambda\in\Lambda}$ from $\mathfrak S$ such that
\[\sum_{\lambda\in\Lambda}^{\ovline}\chi_E(x)\]
is divergent for all $x\in C$ but
\[\sum_{\lambda\in\Lambda}\svar(E_\lambda,\mu)\subset(-\varepsilon,\varepsilon).\]

iv.) There is an  $\varepsilon>0$, such that \dots (continued as in  iii.) )

v.) For each $\varepsilon>0$
 there is a countable indexed family of sets $\mathfrak E$ from $\mathfrak S$ covering $C$ such that
\[\sum_{E\in\mathfrak E }\svar(E_\lambda,\mu)\subset(-\varepsilon,\varepsilon). \]
\begin{proof}
The equivalence of i. and ii. follows from Lemma \ref{lem:CONSW},
while equivalence with iii. follows from Lemma \ref{lem:CONNEG}.

Statements iv. and v. are weaker statements the iii., but summing of the appropriate pointed envelopes
of semivariation contained in $(-\epsilon 2^{-n-1},\epsilon 2^{-n-1} )$ ($n\in\mathbb N$)
we obtain the reverse direction.
\end{proof}
\end{lemma}
 It is also notable that
\begin{lemma}[``Structure theorem'']\label{lem:CL5}
Let $\mathfrak S\subset\mathfrak P(X)$ be an interval system and
$\mu:\mathfrak S\rightarrow\mathbb R$ be a measure.

Suppose that  $f:X\rightarrow\mathbb R$ is integrable, and $\svar(f,\mu)$ is bounded.
Then, we claim, for arbitrary $\varepsilon>0$ there exist countable families
$c_\lambda\in\mathbb R$, $E_\lambda\in\mathfrak S$ ($\lambda\in\Lambda$) such that
for each $x\in X$
\[f(x)=\sum_{\lambda\in\Lambda}^{\ovline} c_\lambda\chi_{E_\lambda}(x)\]
or the right side is completely divergent [meaning that both the positive and negative terms
are unbounded]; and
\[\sum_{\lambda\in\Lambda} \svar(c_\lambda,E_\lambda,\mu)\subset
\svar(f,\mu)+(-\varepsilon,\varepsilon).\]
\begin{proof}
This is Lemma \ref{lem:STRUC} except that we have to provide full divergence on the negligible set
where the right side sum and the left side are not equal.

To do that, we can take
functions $\chi_{E_\lambda}$ as in Lemma \ref{lem:CL4}.iii plus the functions $-\chi_{E_\lambda}$
into the sum.
Adding these terms make the big sum divergent on a negligible set including $C$,
but leaves everything else as it was over other points $x$. If
\[\sum_{\lambda\in\Lambda}\svar(E_\lambda,\mu)\]
was small then the sum of all semivariations will not increase much.
\end{proof}
\end{lemma}
~

\paragraph{\textbf{B. Positive measures.}}~
\\
\begin{point}
The more classical part of Lebesgue integration deals with the case $\mu\geq0$.

Using Beppo Levi's theorem, constructivity and monotonicity, one can now prove
the usual theorems in the development of classical Lebesgue integration.
That includes Fatou's Lemma and Lebesgue's dominated convergence theorem.

The only major point missing then is to show that being measurable and
being dominated implies being integrable. That follows from Lebesgue's dominated
convergence theorem once measurability is properly defined.
\end{point}
\begin{point}
It is very easy to see that if $\mathfrak R$ and $\mathfrak S$ are interval systems then
\[\mathfrak R\boxtimes\mathfrak S=\{B\times A\,:\,B\in\mathfrak R,\,A\in\mathfrak S\}\]
is an interval system. Moreover, if $\mu,\nu\geq0$ are measures then the function
\[\nu\times\mu:\mathfrak R\boxtimes\mathfrak S\rightarrow\mathbb R\]
\[\nu\times\mu(B\times A)=\nu(B)\mu(A)\]
will define a measure.
Opposed to the lengthy discussions of Section\ref{sec:PRODMEAS},
 that is an immediate consequence of Beppo Levi's theorem.
After that, Fubini's theorem immediately follows from Lemma \ref{lem:CL5} and Beppo Levi's theorem.
\end{point}

\paragraph{\textbf{C. Signed measures.}}~
\\

For the sake of completeness we include the discussion about measure decomposition
without specifically applying to the Hahn decomposition theorem in terms of the underlying set $X$.
We restrict our attention to the locally finite case

\begin{lemma}\label{lem:SIGNPRE}
Let $\mathfrak S\subset\mathfrak P(X)$ be an interval system and
$\mu:\mathfrak S\rightarrow\mathbb R$ be a measure.

Then, we claim, $\mu$ is of locally finite semivariation if and only if $\mu$ is locally bounded,
ie. for each $A\in\mathfrak S$
\[\svar(A,\mu)\]
is bounded.
\begin{proof}
($\Rightarrow$) Let $\mathcal U$ be as in Definition \ref{def:LFSVAR}.a with any bounded neighborhood.
 Then we can rescale
$\mathcal U$ such that $1\in\mathcal U$ but
\[\intp \mathtt M(\svar(\mathcal U,A),\mu)\]
is still bounded. The set above contains $\svar(A,\mu)$, hence that is bounded.

($\Leftarrow$) According to Lemma \ref{lem:NEIGH} and Lemma \ref{lem:APPCONV} we know that
\[\intp \mathtt M(\svar([-1,1],A),\mu)=\overline{\conv \svar(A,\mu)- \conv\svar(A,\mu)}.\]
Then we can rescale.
\end{proof}
\end{lemma}
\begin{remak}
For the purposes of classical integration this boundedness
statement is more practical as definition.
\end{remak}
\begin{defin}\label{def:SIGNDECOMP}
Let $\mathfrak S\subset\mathfrak P(X)$ be an interval system and
$\mu:\mathfrak S\rightarrow\mathbb R$ be a measure of locally finite semivariation.
Then we define
\[|\mu|^+(A)=\sup\biggl\{\sum_{\lambda\in\Lambda} |\mu(A_\lambda)|^+\,:\,
\bigcup_{\lambda\in\Lambda}^\updisjoint A_\lambda=A\text{ in }\mathfrak S\biggr\}\]
\[|\mu|^-(A)=\sup\biggl\{\sum_{\lambda\in\Lambda} |\mu(A_\lambda)|^-\,:\,
\bigcup_{\lambda\in\Lambda}^\updisjoint A_\lambda=A\text{ in }\mathfrak S\biggr\}\]
\[|\mu|(A)=\sup\biggl\{\sum_{\lambda\in\Lambda} |\mu(A_\lambda)|\,:\,
\bigcup_{\lambda\in\Lambda}^\updisjoint A_\lambda=A\text{ in }\mathfrak S\biggr\}.\]
\end{defin}

\begin{lemma}\label{lem:SIGNDECOMP}
Let $\mathfrak S\subset\mathfrak P(X)$ be an interval system and
$\mu:\mathfrak S\rightarrow\mathbb R$ be a measure of locally finite semivariation.

The functions $|\mu|^+,|\mu|^-,|\mu|:\mathfrak S\rightarrow \mathbb R$ are non-negative
measures and
\[\mu(A)=|\mu|^+(A) -|\mu|^-(A) \qquad\qquad |\mu|(A)=|\mu|^+(A)+|\mu|^-(A).\]
Furthermore,
\[|\mu|^+(A)=\sup\svar(A,\mu)\qquad\qquad -|\mu|^-(A)=\inf\svar(A,\mu).\]

\begin{proof}
The fact that $|\mu|^+,|\mu|^-,|\mu|$ are measure follows from the fact that if
\[A=\bigcup_{\omega\in\Omega }^\updisjoint B_\omega\]
then computing $|\mu|^+(A), |\mu|^-(A), |\mu|(A)$ we can restrict our attention to
decompositions $A$ which refine the decomposition above. Indeed, we can always take
a common refinement of decompositions and $\sup$ will not suffer.

The decomposition equalities follow from the equalities
\[\sum_{\lambda\in\Lambda}\mu(A_\lambda)
=\sum_{\lambda\in\Lambda} |\mu(A_\lambda)|^+-\sum_{\lambda\in\Lambda} |\mu(A_\lambda)|^-\]
and
\[
\sum_{\lambda\in\Lambda}|\mu(A_\lambda)|
=\sum_{\lambda\in\Lambda} |\mu(A_\lambda)|^++\sum_{\lambda\in\Lambda} |\mu(A_\lambda)|^-.\]
Indeed, we can take successively refined decompositions such that  sums
\[\sum_{\lambda\in\Lambda} |\mu(A_\lambda)|^+,\quad \sum_{\lambda\in\Lambda} |\mu(A_\lambda)|^-,
\quad \sum_{\lambda\in\Lambda} |\mu(A_\lambda)|\]
will limit to $|\mu|^+(A), |\mu|^-(A), |\mu|(A)$, while the equalities stand.

The statements for the supremum and infininum of the semivariation follows from the definition.
\end{proof}
\end{lemma}
\begin{lemma}\label{lem:SIGN}
Let $\mathfrak S\subset\mathfrak P(X)$ be an interval system and
$\mu:\mathfrak S\rightarrow\mathbb R$ be a measure of locally finite semivariation.
Suppose that $f:X\rightarrow \mathbb R$ is an integrable function.

Then, we claim,
\[\int f|\mu|^+,\quad \int f|\mu|^-,\quad \int f|\mu|,\quad\]
exist and
\[\int f\mu=\int f|\mu|^+-\int f|\mu|^-,\qquad \int f|\mu|=\int f|\mu|^++\int f|\mu|^-.\]
\begin{proof}
From the definition of $|\mu|^+$ and $|\mu|^-$ one can see that
\[\svar(c\chi_{E},|\mu|^+)\subset \svar(c\chi_{E},\mu);\qquad
\svar(c\chi_{E},|\mu|^-)\subset -\svar(c\chi_{E},\mu). \]
Applying this to a sum as in Lemma \ref{lem:CL5} and using Beppo Levi's theorem
one immediately obtains the statements.
\end{proof}
\end{lemma}

\begin{lemma}
\label{lem:SIGNINT}
Let $\mathfrak S\subset\mathfrak P(X)$ be an interval system and
$\mu:\mathfrak S\rightarrow\mathbb R$ be a measure of locally finite semivariation.
Suppose that $f:X\rightarrow \mathbb R$ is an integrable function.

Then, we claim,
\[\int |f||\mu|=\int|f|^-|\mu|^- +\int|f|^+|\mu|^++\int|f|^-|\mu|^++\int|f|^+|\mu|^-\]
and
\[\int f\mu=\int|f|^-|\mu|^- +\int|f|^+|\mu|^+-\int|f|^-|\mu|^+-\int|f|^+|\mu|^-\]
(existence stated).
 Furthermore,
\[\svar(f,\mu)\subset
\biggl[-\int|f|^-|\mu|^+-\int|f|^+|\mu|^-,\int|f|^-|\mu|^- +\int|f|^+|\mu|^+\biggr].\]
The endpoints are in the semivariation.
In particular,
\[\subset\biggl[-\int |f||\mu| ,\int |f||\mu|\biggr]\]
follows.

\begin{proof}
According to the previous lemma $\int f\,|\mu|$ exists.
From this one obtains that $\int |f|\,|\mu|$ exists.
The formulae for the integrals follow from the earlier statements.
The formula for the semivariation immediately follows from
\[\svar(f,\mu)=\overline{\svar(|f|^+,\mu)+\svar(-|f|^-,\mu)}=
\overline{\svar(|f|^+,\mu)-\svar(|f|^-,\mu)}.\]
(The ``$\subset$'' part is easy while  ``$=$'' follows from Lemma \ref{lem:SVSUM}.b.)
\end{proof}
\end{lemma}
\begin{remak}
However, at this point it is probably easier to pass to a measure extension
and apply Hahn decomposition directly.
\end{remak}

\begin{remark}
Despite the great simplifying power of Hahn decomposition, we can treat measure product
just like in the non-negative case. The main difference is that when we apply for Beppo Levi's
theorem in order to prove $\sigma$-additivity then for each countable decomposition
\[B\times A=\bigcup^\updisjoint_{\lambda\in\Lambda} B_\lambda\times A_\lambda\]
we have to prove that
\[\sum_{\lambda\in\Lambda} \svar(\nu(B_\lambda)\chi_{A_\lambda})\]
is bounded, and we cannot use monotonicity as in the non-negative case.
Nevertheless, boundedness follows from
\[\overline{\sum_{\lambda\in\Lambda} \svar(\nu(B_\lambda)\chi_{A_\lambda},\mu)}=
\overline{\sum_{\lambda\in\Lambda}\intp \mathtt M(\var(\nu(B_\lambda),{A_\lambda}),\mu)}=\]
\[=\intp\mathtt M\biggl(\sum_{\lambda\in\Lambda}\var(\nu(B_\lambda),{A_\lambda}),\mu\biggr)\subset
 \intp\mathtt M(\var(\svar(B,\nu), A) ,\mu)\subset\]
 \[\subset\intp\mathtt M(\var(\conv\svar(B,\nu), A) ,\mu)
 \subset [-|\nu|(B)|\mu|(A),|\nu|(B)|\mu|(A) ]. \]
\end{remark}
~

\paragraph{\textbf{D. Comments.}}
~\\

The cases when
\[\int f\mu=+\infty\quad\text{or}\quad \int f\mu=-\infty\]
can be incorporated as it was indicated in Section \ref{sec:MODIF}.D.
See Section \ref{sec:LR} for statements covering the more general case,
and in fact, a more unified treatment.
\newpage
\section{Infinite formal sums}
\label{sec:INFORMAL}
\begin{conven}
In what follows the sum signs
\[\sum\]
will be used in the extended sum sense (as follows).
\end{conven}
\begin{remin}\label{rem:IFSREMIN}
In what follows we will use the extended set of the real numbers
\[\real^*=\real\cup\{+\infty,-\infty,\pm\infty\}.\]
The set of positive and nonnegative numbers are the intervals
$(0,+\infty]$ and $[0,+\infty]$, respectively.
To any extended real $a$ we assign a ``positive'' (nonnegative)
and ``negative'' (nonpositive) part $|a|^+$ and $|a|^-$ from $[0,+\infty]$ as usual.
We have $|\pm\infty|^+=|\pm\infty|^-=+\infty$.

The relation $a\leq b$ means that $a,b\in[-\infty,+\infty]$ and the equality holds
in the usual sense or it refers to $-\infty\leq\pm\infty$, $\pm\infty\leq\pm\infty$,
or $\pm\infty\leq+\infty$.

In general, the right-modified relation $a\prec' b$ means
that $a\prec b$ or $b=\pm\infty$. We will similarly use ${}'$ on the left.

We use two kinds of infinum and supremum; they are as usual except that in
the ordinary case
\[\inf\{\pm\infty\}=-\infty\qquad\text{and}\quad\sup\{\pm\infty\}=+\infty; \]
but in the extended case
\[\inf'\{\pm\infty\}=+\infty\qquad\text{and}\quad\sup'\{\pm\infty\}=-\infty . \]

Any countable family $\{a_\lambda\}_{\lambda\in\Lambda}$ of the extended real
numbers will yield one extended real number
\[\sum_{\lambda\in\Lambda}a_\lambda=\biggl(\sum_{\lambda\in\Lambda}|a_\lambda|^+
\biggr)-\biggl(\sum_{\lambda\in\Lambda}|a_\lambda|^-\biggr)\]
as a sum (generalizing the nonnegative sums). The case $\pm\infty$
corresponds to $(+\infty)-(+\infty)$, the other cases are as usual.
\end{remin}
A definition of lesser importance is:
\begin{defin}\label{def:IFS}
Suppose that $\{a_\lambda\}_{\lambda\in\Lambda}$ is an indexed family of numbers from $\mathbb R^*$.
Assume that $n\in\mathbb N$. Then we define
\[\sup^{(\prime)}{}_{\!-n}\,\{a_\lambda\}_{\lambda\in\Lambda}=\inf_{\substack{\Xi\subset\Lambda\\|\Xi|=n}}
\,\,\sup^{(\prime)}_{\lambda\in\Lambda\setminus\Xi}a_\lambda\]
and
\[\inf^{(\prime)}{}_{\!-n}\,\{a_\lambda\}_{\lambda\in\Lambda}=\sup_{\substack{\Xi\subset\Lambda\\|\Xi|=n}}
\,\,\inf^{(\prime)}_{\lambda\in\Lambda\setminus\Xi}a_\lambda\]
(Ie. supremum and infinum except $n$ elements.)
\end{defin}
The following definitions will be essential.
\begin{defin}\label{def:IFSCM}
Suppose that $\mathfrak S$ is a family of sets. A formal sum $\mathsf C$ over $\mathfrak S$
is a countable indexed family  $\{(c_\lambda,C_\lambda)\}_{\lambda\in\Lambda}$
from $\real\times\mathfrak S$. (Ie. multiplicities are allowed.)
A set $C_\lambda$ is called a coefficient set, a number $c_\lambda$ is called
a coefficient number, the pair $(c_\lambda,C_\lambda)$ is called a coefficient
pair. The coefficient set family of $\mathsf C$ is
\[\mathfrak C=\{C_\lambda\,:\,\lambda\in\Lambda\}.\]

A countable disjoint union of formal sums  $\mathsf C_\xi$  $(\xi\in\Xi)$
is the formal sum
\[\bigcup_{\xi\in\Xi}^\updisjoint\mathsf C_\xi=\{(c_{\xi,\lambda},C_{\xi,\lambda})
\}_{\xi\in\Xi,\lambda\in\Lambda_\xi}\]
which contains all the coefficient pairs from the original
formal sums (with appropriate reindexing).
The opposite formal sum $-\mathsf C$ is the formal sum where
the coefficient numbers $c_\lambda$ are substituted by $-c_\lambda$.
A formal sum $\mathsf C$ is positive (or nonnegative) if all the coefficients are from
$(0,+\infty]$ (or from $[0,+\infty]$).
\end{defin}
\begin{conven}
In what follows when we write $\mathsf C=-\mathsf C^+\,\dot\cup\,\mathsf C^+$ then we assume that
$\mathsf C^+$ and $\mathsf C^-$ are nonnegative formal sums. Elements of $0$ coefficient can belong
to any side (but they do not matter in general).
\end{conven}
\begin{defin}\label{def:IFSEVAL}
The value of the formal sum $\mathsf C$ at a point $x\in X$ is defined as
\[\mathsf C^\Sigma(x) =\sum_{\lambda\in\Lambda,\,x\in C_\lambda} c_\lambda .\]
\end{defin}
\begin{defin}\label{def:IFSDECO}
a.) A formal sum $\mathsf C$ number-decomposes into a formal sum
$\mathsf C'$ if each coefficient pair $(c_\lambda,C_\lambda)$ is substituted by
countably many other one $(c_{\lambda,\gamma},C_\lambda)$ where $\gamma\in\Gamma_\lambda$
such that
\[|c_\lambda|^+=\sum_{\gamma\in\Gamma_\lambda}|c_{\lambda,\gamma}|^+\qand
|c_\lambda|^-=\sum_{\gamma\in\Gamma_\lambda}|c_{\lambda,\gamma}|^-\,.\]
(Ie. nonnegative coefficients decompose to nonnegative coefficients, and similarly for
nonpositive ones. This allow the decomposition of a $0$ coefficients into $0$ many terms.)

We say ``finite-number'' instead of ``number'' if element of $\mathsf C$ actually number-decomposes
to finitely many parts.

b.) A formal sum $\mathsf C$ set-decomposes into a formal sum
$\mathsf C'$ if each coefficient pair $(c_\lambda,C_\lambda)$ is substituted by
countably many other one $(c_{\lambda},C_{\lambda,\gamma})$ where $\gamma\in\Gamma_\lambda$
and $\{C_{\lambda,\gamma}\,:\,\gamma\in\Gamma_\lambda\}$ exactly decomposes the set $C_\lambda$.
In particular, if a family of sets $\mathfrak F\subset\mathfrak S$ decomposes
$\mathfrak C=\{C_\lambda\,:\,\lambda\in\Lambda\}$ then one can naturally
set-decompose $\mathsf C$ along $\mathfrak F$.

c.) A formal sum $\mathsf C$ decomposes into a formal sum $\mathsf C'$
\[\mathsf C\rightarrow\mathsf C'\]
if it occurs after finitely many number- or set-decomposition steps.

d.) A formal sum $\mathsf C$ subdecomposes into a formal sum $\mathsf C'$
\[\mathsf C\xrightarrow{\sub}\mathsf C'\]
if it occurs after the following procedure:
\[\xymatrix@=5pt{
\mathsf C_{}^{}&\equiv&\mathsf C_0
\ar[dd]\ar[rr]&&\mathsf C_1\ar[dd]
\ar[rr]&&\mathsf C_2\ar[dd]\ar[rr]&&\ldots\\\\
&&\mathsf C_1'&\dot\cup&\mathsf C_2'&\dot\cup&\mathsf C_3'&\dot\cup&\ldots
&\equiv& \mathsf C'_{}.}\]
We start with $\mathsf C=\mathsf C_0$ and then we successively (finitely) decompose
$\mathsf C_{n-1}$ into $\mathsf C_n\,\dot\cup\,\mathsf C_n'.$ Then we take the union
of all $\mathsf C_n'$ and this will be $\mathsf C'$.

If $\mathsf C$ (sub)decomposes into $\mathsf C'$  and the coefficient pair
$(c'_\mu,C'_\mu)\in\mathsf C'$ is obtained from $(c_\lambda,C_\lambda)\in\mathsf C$
after finitely many steps then we say that $(c_\lambda,C_\lambda)\in\mathsf C$
is the progenitor of $(c'_\mu,C'_\mu)\in\mathsf C'$.
\end{defin}
\begin{lemma}\label{lem:IFSDECO}
Two (finite-)number-decomposition steps following each other
are equivalent to one (finite-)number-decomposition step.
Two set-decomposition steps following each other are equivalent to one set-decomposition steps.
\begin{proof}
That immediately follows from the definitions.
\end{proof}
\end{lemma}
\begin{conven}
In what follows we assume that every formal sum is over an interval system
$\mathfrak S\subset\mathfrak P(X)$; and every decomposition is supposed to be over $\mathfrak S$.
\end{conven}
The following procedure is a very intuitive one, but is not so easy to describe.
\begin{defin}[Description of Procedure 1] \label{def:IFSPROCONE}
Assume that $\mathfrak S\subset\mathfrak P(X)$ is an interval system.

Suppose that we have countably many formal sums
\[\{\mathsf C_n\}_{n\in N}=\{\{(a_{n,m}, A_{n,m})\}_{m\in N_n}\}_{n\in N}\]
over $\mathfrak S$, with  positive coefficients, where $N$ and $N_n$ are initial
segments of $\mathbb N$. Also list all possible pairs
$\langle{n,m}\rangle$ as $\langle n(j),m(j)\rangle$ $(j\in J)$ along an initial segment $J$ of $\mathbb N$.

In what follows we will describe a procedure which set-decomposes and then finite-number-decomposes each
formal sum $\mathsf C_n$, and
to each piece $(z,Z)$ resulted assigns a ``color'' from $N$ and a ``lift''
$[u,v)\times Z\subset[0,+\infty)\times Z$,
such that $|v-u|=z$.

First, consider the sequence,
\[A_{n(0),m(0)},\,A_{n(1),m(1)},\,A_{n(2),m(2)},\,A_{n(3),m(3)},\,\ldots\,.\]
Apply Lemma \ref{lem:REFIN}.a and b. in order to find a decomposition $\mathfrak A_{-1}$ of the union
of all these sets, and
successive refinements $\mathfrak A_j$ $(j\in J)$ such that $A_{n(j),m(j)}$ is already decomposed
by $\mathfrak A_j$.

Then, the procedure continues in steps along $J$. Suppose that $j\in J$ comes. First, we set-decompose
\[(a_{n(j),m(j)}, A_{n(j),m(j)})\]
along $\mathfrak A_j$. Suppose that $(a_{n(j),m(j)}, E)$ is one piece resulted.
Notice that for each element $A'$ of $\mathfrak A_0\cup\,\ldots\,\cup\mathfrak  A_j$
the set $E$ is either contained in $A'$  or $E$ is disjoint from $A'$.
Let
\[b_j^E=\sum_{\substack{1\leq k\leq j,\, n(k)=n(j)\\ E\subset A_{n(k),m(k)}}}a_{n(k),m(k)}
\qand
c_j^E=\sum_{\substack{1\leq k< j,\, n(k)=n(j)\\ E\subset A_{n(k),m(k)}}}a_{n(k),m(k)}.\]
Notice that $b_j-c_j=a_{n(j),m(j)}$.
More generally, let $S^E_j$ be the set of all positive numbers which occur as
\[\sum_{\substack{1\leq k\leq i,\, n(k)=n(i)\\ E\subset A_{n(k),m(k)}}}a_{n(k),m(k)}\]
with an $i<j$ chosen. We can consider all those values which fall into $(c_j^E,b_j^E)$.
Suppose that they are like
\[c_j^E<u_1<\ldots<u_h<b_j^E.\]
Now, what we do is that we finite-number-decompose $(a_{n(j),m(j)}, E)$ with coefficient numbers
\[u_1-c_j^E,\,u_2-u_1,\,\ldots,\, b_j^E-u_h.\]
Here, we assign the ``lifts''
\[[c_j^E,u_1)\times E,\,[u_1,u_2)\times E,\,\ldots\,[u_h,b_j^E)\times E,\]
respectively.
Moreover, to each piece $(u_{l+1}-u_l,E)$ we can assign the ``color'' $s$ if there are
$s$ many numbers $p$ such that
\[u_{l}<d_p=\sum_{\substack{1\leq k\leq j,\, n(k)=p\\ E\subset A_{n(k),m(k)}}}a_{n(k),m(k)}.\]
(We can take the conventions $u_0=c_j^E$ and $u_{h+1}=b_j^E$ here.)

From the viewpoint of each element $(a_{n,m},A_{n,m})$ a set-decomposition and a finite-number-decomposition
happened, so we can imagine the whole procedure so that we set- and finite-number-decompose each $\mathsf C_n$.
We can gather things at the end into coefficient systems
\[\mathsf C'_0,\,\mathsf C'_1,\,\mathsf C'_2,\,\ldots\]
as decompositions of $\mathsf C_n$'s, respectively, or into
\[\mathsf F_0,\,\mathsf F_1,\,\mathsf F_2,\,\ldots\,,\]
according to color.
\end{defin}
The next lemma describes some properties of Procedure 1.
\begin{lemma}\label{lem:PX}
Consider Procedure 1.

a.) Consider the elements $\mathsf F_n$ of color $n$. We claim that for each $x\in X$ there is a
possibly infinite sequence of non-negative real numbers
\[0=v_0<v_1<v_2<v_3<\ldots\]
and sets
\[E_0\supset E_2\supset E_1\supset E_2\supset E_3\supset \ldots\]
from $\mathsf S$, all containing $x$ such that
\[[v_j,v_{j+1})\times E_j\]
are lifts of elements of $\mathsf F_n$. Furthermore, the lifts account for all those lifts
\[[u,v)\times E\]
such that $x\in E$.

b.) We claim that
\[\mathsf F_n^\Sigma(x)=\sup{}_{-n}\{\mathsf C_0^\Sigma(x),\,\mathsf C_1^\Sigma(x),\
,\mathsf C_2^\Sigma(x),\,\ldots\}.\]

c.) Suppose that
\[[u,v)\times E\]
is a lift of an element of color $n$. Suppose that $0\leq m\leq n$. Then, there is exactly one element of
color $m$ whose lift
\[[u',v')\times E'\]
contains it, or even intersects.

d.) Suppose that
\[[u,v)\times E\]
is a lift of an element of color $n$. Suppose that $m>n$.
Let $\mathfrak E$ be the sets of all $E'$'s such that
\[[u',v)\times E'\]
is a lift of color $m$ (with an appropriate $u'$).

Then, we claim, if $\mathfrak E$ covers $E$, then it forms an exact decomposition of $E$.
Furthermore, $\mathfrak E$ they decomposes every $E''$ which occur from a lift
\[[u'',v'')\times E''\subset [u,v)\times E.\]

e.) The situation in the previous point occurs if
\begin{itemize}
\item[i.] For all $x\in E$ the inequality $\mathsf F_m^\Sigma(x)>v$ holds; or if
\item[ii.] $[u,v)\times E$ comes from $\mathsf C_s$ but $\mathsf F_{m-1}^\Sigma(x)>\mathsf C_s^\Sigma(x)$.
\end{itemize}
\begin{proof}
a.) From the construction, all $E$'s which contain $x$ form such a sequence as above,
and they occur in the order of steps of Procedure 1. In such a step if lifts
 of color $n$
 \[[w_0,w_1)\times E,\,\ldots,[w_{g-1},w_g)\times E \]
were created then it is just to say that the $(n+1)$th greatest value
of the partial sums $d_p$ coming from the various $\mathsf C_n$'s up to that point
increased from $w_0$ to $w_g$.

b.) As it was explained in the previous point, a cofinal sequence of $v_i$'s occur as
the $(n+1)$th greatest value of the $d_p$'s. Taking the limit will yield the statement.

c.) When a lift $[u,v)\times E$ of color $n>0$ is created then it is contained in a lift of
$[u',v')\times E'$ of color $n-1$. Simply, we know that $u$ was exceeded $n$ times before, and the
last time it was exceeded must have yielded a lift $[u',v')\times E'$ of color $n-1$,
which by construction contains $[u,v)\times E$. Also, that because of point a. all the lifts
of the same color are disjoint; that proves the rest of the statement.

d.) That statement immediately follows from the fact that all $E$'s occuring all either disjoint or
one which occurs later is contained in the earlier one.

e.) By construction, when the $(m+1)$the greatest $d_p$ exceeds $v$ at a point $x\in E$,
then such a nice lift occurs by definition.
In case i.) this is the direct requirement.

As for ii.), suppose that we are in the step  when $[u,v)\times E$ comes up. Assume that
the color $m$ never reaches $v$. Suppose, that later we are in a step when colors $0,\ldots, m-1$ have
exceeded not only $v$ but $\mathsf C_s^\Sigma(x)$ (according to ii.) ).
Consider $d_s$ in that step. It cannot be of color
$0,\ldots,m-1$, because it is too small. So it must be at least of color $m$. But this is a contradiction.
\end{proof}
\end{lemma}

The next lemma is the crucial point in order to prove the consistency of the Lebesgue integral.
Its statement can be proven ``by hand'', but the proofs simply follow from the
properties of Procedure 1.
\begin{lemma}\label{lem:FORMALDECOMP}
Consider an interval system $\mathfrak S\subset\mathfrak P(X)$.

a.) Suppose that $\mathsf C$ and $\mathsf D$ are positive formal sums, $\mathsf C$ is finite.
Also suppose and for all $x\in X$
\[\mathsf C^\Sigma(x)<\mathsf D^\Sigma(x)\quad\text{or}\quad \mathsf C^\Sigma(x)=0.\]
Then, we claim, there exist positive formal
sums $\mathsf C'$ and $\mathsf E$ such that $\mathsf C$ decomposes into $\mathsf C'$
and $\mathsf D$  decomposes into $\mathsf C'\,\dot\cup\,\mathsf E$:
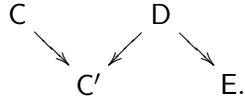
\[\xymatrix@=12pt{ \mathsf C\ar[dr] &  &\mathsf D\ar[dr]\ar[dl]
& \\& \mathsf C' & & \mathsf E. }\]
In fact, these decompositions above can be achieved by a set-decomposition step
followed by a finite-number-decomposition step for $\mathsf C$; and by a set/finite-number/set
decomposition for $\mathsf D$.

b.) Suppose that $\mathsf C$ and
$\mathsf D$ are positive formal sums and
 \[\mathsf C^\Sigma(x)\leq\mathsf D^\Sigma(x)\]
holds for all $x\in X$. Then, we claim, there exists a positive formal sum $\mathsf C'$
such that $\mathsf C$  decomposes into $\mathsf C'$ and
$\mathsf D$ subdecomposes into $\mathsf C'$:
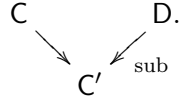
\[\xymatrix@=12pt{\mathsf C\ar[dr] &  & \mathsf D.\ar[dl]^{\sub}  \\ & \mathsf C' &} \]
In fact, the decomposition of $\mathsf C$ above can be achieved by a number/set/finite-number-decomposition.

c.) Suppose that  $\mathsf C_1$ and $\mathsf
C_2$ are positive formal sums such that
\[\mathsf C_1^\Sigma(x)= \mathsf C_2^\Sigma(x)\]
for all $x\in X$. Then there exists a positive formal sum
$\mathsf  C$ such that
 
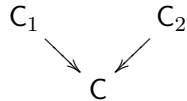
\[\xymatrix@=12pt{ \mathsf C_1\ar[dr] &  & \mathsf C_2\ar[dl] \\ & \mathsf C &} \]
both $\mathsf C_1$ and $\mathsf C_2$ can be decomposed into
$\mathsf C$. In fact, the decompositions above can be achieved by
a number/set/finite-number/set/finite-number-decomposition.

d.) Suppose that $\mathsf C$ and $\mathsf D$ are positive
formal sums. Then, we claim, there exists a decomposition
\[\mathsf C\,\dot\cup\,\mathsf D\rightarrow\mathsf  E_{\min}\,\dot\cup\, \mathsf E_{\max}\]
such that
\[ \mathsf E_{\min} {}^{\Sigma}(x)=\min( \mathsf C^\Sigma(x), \mathsf D^\Sigma(x) ) \]
and
\[ \mathsf E_{\max} {}^{\Sigma}(x)=\max( \mathsf C^\Sigma(x), \mathsf D^\Sigma(x) ). \]
The decomposition above can be achieved by a set-decomposition
step followed by a finite-number-decomposition step.

The statements of this lemma also valid for nonnegative sums instead of positive ones.
\begin{proof}
a.) Apply Procedure 1. Being $\mathsf C$ finite we can assume that in the listing along $J$
we list the elements of $\mathsf C$ and then the elements of $\mathsf D$.
Then, when
we do Procedure 1,
everything which comes from $\mathsf C$ gets color $0$, because only $d_0$ plays.
Suppose that the lift $[u,v)\times E$ comes from $\mathsf C$. According Lemma \ref{lem:PX}.e.ii,
$E$ is decomposed by an appropriate system $\mathfrak E$ as in Lemma \ref{lem:PX}.d .

Use a set-decomposition along $\mathfrak E$ with respect to all elements, whose lift is
is contained in $[u,v)\times E$. From the viewpoint of $\mathsf D$ we had a set-decomposition and a finite-number
decomposition followed by a set-decomposition. Let $\mathsf C'$ be the set of all
decomposition products from $\mathsf D'$ subordinated to such a $[u,v)\times E$ coming from $\mathsf C$.
And let the rest be $\mathsf E$.
 From the viewpoint of elements of $\mathsf C$, they  get set-decomposed one, then set-decomposed again,
 and we know that they finite-number-decompose into $\mathsf C'$. That amounts to a set-decomposition followed by
 a finite-number-decomposition.

b.) First, number-decompose each coefficient number of $\mathsf C$ into infinitely many
positive finite numbers. That way we obtain $\mathsf C''$.
Suppose that $\mathsf C''$ is indexed by $\natu$.
Now we start with $\mathsf D=\mathsf D_0$. Then, by recursion, in the $n$th step,
 using part a.), we obtain a coefficient sets $\mathsf C'_n$ and $\mathsf D_n$
 such that $\{(c''_{n-1},C''_{n-1})\}$ decomposes into $\mathsf C'_n$ and $\mathsf D_{n-1}$
 decomposes into $\mathsf C'_n\,\dot\cup\,\mathsf D_n$:
\[\xymatrix@=5pt{
\mathsf C\ar[rr]&&\mathsf C''&\equiv&\{(c''_0,C''_0)\}\ar[dd] &\dot\cup& \{(c''_1,C''_1)\}\ar[dd]
&\dot\cup&\{(c''_2,C''_2)\}\ar[dd]&\dot\cup& \dots \\\\
&&&&\mathsf C'_1&\dot\cup&\mathsf C'_2&\dot\cup&\mathsf C'_3&\dot\cup&\ldots&\equiv&\mathsf C'\\\\
&&&&\mathsf D=\mathsf D_0\ar[rr]\ar[uu]&& \mathsf D_1\ar[rr]\ar[uu]&&\mathsf D_2\ar[rr]\ar[uu]&&\ldots}\]
The reason that part a.) remains applicable is that we used that special number-decomposition first.
Then let
\[\mathsf C'=\bigcup^\updisjoint_{n\geq 1}\mathsf C_n'.\]
Part a.) also implies the statement about the possible number of number- and set-decompositions.

c.) First number-decompose each coefficient into an infinite sum of positive numbers, hence we obtain
$\mathsf C_1'$ and $\mathsf C_2'$. Now apply Procedure 1. We get element in color $0$ and $1$.
In fact, if $[u,v)\times E$ is a lift of color $0$, then according to Lemma \ref{lem:PX}.e.i
the situation of Lemma \ref{lem:PX}.d applies. If we apply the decompositions along those $\mathfrak E$
for each $0$-color element then we
we find that the resulted system $\mathsf C$ is such that both
$\mathsf C_1'$ and $\mathsf C_2'$
set/finite-number decomposes into it, either through as an element of color $0$ or color $1$.

d.) This is Procedure 1 applied directly; and separation by colors.

The closing remark follows from the fact the we can always finite-number-decompose
terms with $0$ coefficients into nothing.
\end{proof}
\end{lemma}
\begin{remak}
In Lemma \ref{lem:FORMALDECOMP}.a it is easy to improve both decompositions into set/finite-number
decompositions. Making detailed, optimal statements here, however, is not particularly worthwhile.
\end{remak}
\begin{cor}\label{cor:FORMALDECOMP}
Consider an interval system $\mathfrak S\subset\mathfrak P(X)$.

Suppose that
$\mathsf C=-\mathsf C^-\,\dot\cup \,\mathsf C^+ $
and $\mathsf D=-\mathsf D^-\,\dot\cup \,\mathsf D^+$ where
$\mathsf C^-,\,\mathsf C^+,$ $\mathsf D^-,\,\mathsf D^+$ are nonnegative formal sums.
Suppose that
\[\mathsf C^{\Sigma}(x)\,\,{}'\!\!='\mathsf D^{\Sigma}(x)\]
for all $x\in X$. (Ie. they are equal or we have $\pm\infty$ on one side.)

Then, we claim, there exists positive formal sums
$\mathsf E_{\mathsf C},\,\mathsf E^{-},\,\mathsf E^{+},\,\mathsf E_{\mathsf D}$
such that
$\mathsf C^+$ decomposes into $\mathsf E_{\mathsf C}\,\dot\cup \,\mathsf E^+$,
$\mathsf C^-$ decomposes into $\mathsf E_{\mathsf C}\,\dot\cup \,\mathsf E^-$,
$\mathsf D^+$ decomposes into $\mathsf E_{\mathsf D}\,\dot\cup \,\mathsf E^+$,
$\mathsf D^-$ decomposes into $\mathsf E_{\mathsf D}\,\dot\cup \,\mathsf E^-$,
ie.
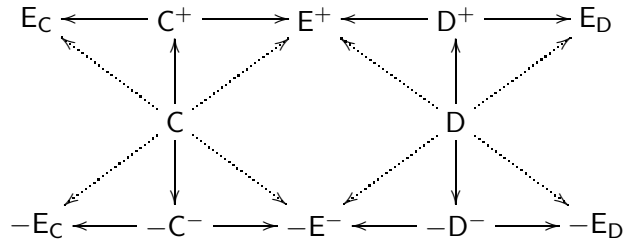
\[
\xymatrix{ \mathsf E_{\mathsf C}& \mathsf C^+ \ar[l]\ar[r]&
\mathsf E^+ & \mathsf D^+  \ar[l]\ar[r]&\mathsf E_{\mathsf D}\\
&\mathsf C\ar[u]\ar[d]\ar@{.>}[ul]\ar@{.>}[ur]\ar@{.>}[dl]\ar@{.>}[dr]
&&\mathsf D\ar[u]\ar[d]\ar@{.>}[ul]\ar@{.>}[ur]\ar@{.>}[dl]\ar@{.>}[dr]& \\
-\mathsf E_{\mathsf C}& -\mathsf C^- \ar[l]\ar[r]& -\mathsf E^- &
 -\mathsf D^-  \ar[l]\ar[r]& -\mathsf E_{\mathsf D} }
\]
\begin{proof}
We see that $\mathsf C^-\,\dot\cup \,\mathsf D^+ $ and $\mathsf C^+\,\dot\cup \,\mathsf D^- $
are nonnegative formal sums such that
\[(\mathsf C^-\,\dot\cup \,\mathsf D^+)^\Sigma(x)=(\mathsf C^+\,\dot\cup \,\mathsf D^-)^\Sigma(x) \]
for all $x\in X$. Then by Lemma \ref{lem:FORMALDECOMP}.c there is a common decomposition
$\mathsf E$ for $\mathsf C^-\,\dot\cup \,\mathsf D^+ $ and $\mathsf C^+\,\dot\cup \,\mathsf D^- $.
From $\mathsf E$,
let $\mathsf E_{\mathsf C}$ be the formal sum of those elements which come form $\mathsf C^+$ and $\mathsf C^-$;
let $\mathsf E^+$ be the formal sum of those elements which come form $\mathsf C^+$ and $\mathsf D^+$;
let $\mathsf E^-$ be the formal sum of those elements which come form $\mathsf C^-$ and $\mathsf D^-$;
let $\mathsf E_{\mathsf D}$ be the formal sum of those elements which come form $\mathsf D^+$ and $\mathsf D^-$.
\end{proof}
\end{cor}
\begin{remark}
The corollary above shows that if two formal sums $\,\,{}'\!\!='$ locally, then they are
sort of equal in global sense, too. In what follows we will not use this kind of
global arguments but this is the main idea behind the content of the next section.
\end{remark}

\newpage\section{The Lebesgue-Riesz integral}\label{sec:LR}
In this section we use interval systems, Section \ref{sec:IS}; the definition of measure,
Definition \ref{def:MEASURE}; and the decompositions from
Definition \ref{def:SIGNDECOMP} and Lemma  \ref{lem:SIGNDECOMP}.
\\

\paragraph{\textbf{A. Definition and consistence}}
\begin{defin}\label{def:FORMALSUM}
\lrquote

Suppose that $\mathsf C=\{(C_\lambda,c_\lambda)\}_{\lambda\in\Lambda}$ is a formal sum with respect to
$\mathfrak S$.

Then we define the positive variation as
\[\int^+\mu\, \mathsf C=\sum_\lambda |c_\lambda|^+|\mu(C_\lambda)|^++
|c_\lambda|^-|\mu(C_\lambda)|^-\in [0,+\infty];\]
the negative variation as
\[\int^-\mu\, \mathsf C=\sum_\lambda |c_\lambda|^-|\mu(C_\lambda)|^++
|c_\lambda|^+|\mu(C_\lambda)|^-\in [0,+\infty];\]
the total variation as
\[\int^\pm\mu\, \mathsf C= \int^+\mu\,\mathsf C+\int^-\mu\, \mathsf C \in [0,+\infty];\]
the overall sum as
\[\int\mu\,\mathsf C=\int^+\mu\,\mathsf C-\int^-\mu\, \mathsf C\in \real^*.\]
\end{defin}
\begin{lemma}\label{lem:FORMALSUM}
\lrquote

a.) If the formal sum $\mathsf C$ decomposes into the formal sum $\mathsf C'$ then
the positive variation, negative variation, total variation and the overall sum stays the same.

b.) If $\mu\geq0$ and the nonnegative formal sum $\mathsf D$ subdecomposes into a
nonnegative formal sum $\mathsf C$ then
\[\int\mu\,\mathsf C\leq\int\mu\,\mathsf D.\]
\begin{proof}
a.) Both number-decomposition and set-decomposition leaves the things above intact.
b.) Any finite sum from $\mathsf C$ is majorized.
\end{proof}
\end{lemma}

\begin{lemma}[Comparison]\label{lem:FORMALCOMPAR}
\lrquote
 Also assume that $\mathsf C$ and $\mathsf D$
are formal sums.

a.) If for all $x\in X$
\[\mathsf C^\Sigma(x)\,\,{}'\!\!='\mathsf D^\Sigma(x) \]
holds then, we claim,
\[\int \mu\,\mathsf C\,\,{}'\!\!='\int\mu\,\mathsf D. \]

b.) If $\mu\geq0$ and for all $x\in X$
\[\mathsf C^\Sigma(x)\,\,{}'\!\!\leq'\mathsf D^\Sigma(x) \]
holds then, we claim,
\[\int \mu\,\mathsf C\,\,{}'\!\!\leq'\int\mu\,\mathsf D. \]
\begin{proof}
We can assume that
$\mathsf C=-\mathsf C^-\,\dot\cup \,\mathsf C^+ $
and $\mathsf D=-\mathsf D^-\,\dot\cup \,\mathsf D^+$ where
$\mathsf C^-,\,\mathsf C^+,\,\mathsf D^-,\,\mathsf D^+$ are nonnegative formal sums.

a.) Our assumption
\[\mathsf C^+{}^\Sigma(x)-\mathsf C^-{}^\Sigma(x)=\mathsf C^\Sigma(x)\,\,{}'\!\!='\mathsf D^\Sigma(x)=
\mathsf D^+{}^\Sigma(x)-\mathsf D^-{}^\Sigma(x)\]
implies
\[(\mathsf C^+\,\dot\cup\,\mathsf D^-)^\Sigma(x)=\mathsf C^+{}^\Sigma(x)+\mathsf D^-{}^\Sigma(x)=
\mathsf C^-{}^\Sigma(x)+\mathsf D^+{}^\Sigma(x)= (\mathsf C^-\,\dot\cup\,\mathsf D^+)^\Sigma(x).\]
Hence, by Lemma \ref{lem:FORMALDECOMP}.c the formal sums
$\mathsf C^+\,\dot\cup\,\mathsf D^-$ and $\mathsf C^-\,\dot\cup\,\mathsf D^+$ have a common
decomposition and by then, by Lemma \ref{lem:FORMALSUM},
\[\int\mu\mathsf C^++\int\mu\mathsf D^-=
\int\mu(\mathsf C^+\,\dot\cup\,\mathsf D^-)=\int\mu(\mathsf C^-\,\dot\cup\,\mathsf D^+)=
\int\mu\mathsf C^-+\int\mu\mathsf D^+. \]
This, in turn, implies
\[\int\mu\mathsf C=\int\mu\mathsf C^+-\int\mu\mathsf C^-\,\,{}'\!\!='
\int\mu\mathsf D^+-\int\mu\mathsf D^-=\int\mu\mathsf D.  \]

b.) Our assumption
\[\mathsf C^+{}^\Sigma(x)-\mathsf C^-{}^\Sigma(x)=\mathsf C^\Sigma(x)\,\,{}'\!\!\leq'\mathsf D^\Sigma(x)=
\mathsf D^+{}^\Sigma(x)-\mathsf D^-{}^\Sigma(x)\]
implies
\[(\mathsf C^+\,\dot\cup\,\mathsf D^-)^\Sigma(x)=\mathsf C^+{}^\Sigma(x)+\mathsf D^-{}^\Sigma(x)\leq
\mathsf C^-{}^\Sigma(x)+\mathsf D^+{}^\Sigma(x)= (\mathsf C^-\,\dot\cup\,\mathsf D^+)^\Sigma(x).\]
Hence, by Lemma \ref{lem:FORMALDECOMP}.b,
there exists a  nonnegative
formal sum $\mathsf E$ such that
$\mathsf C^+\,\dot\cup\,\mathsf D^-$ decomposes into $\mathsf E$
and, in turn, $\mathsf C^-\,\dot\cup\,\mathsf D^+$ subdecomposes into $\mathsf E$.
Now, by Lemma \ref{lem:FORMALSUM},
\[\int\mu\mathsf C^++\int\mu\mathsf D^-=
\int\mu(\mathsf C^+\,\dot\cup\,\mathsf D^-)=\int\mu\mathsf E\leq
\int\mu(\mathsf C^-\,\dot\cup\,\mathsf D^+)=
\int\mu\mathsf C^-+\int\mu\mathsf D^+. \]
This, in turn, implies
\[\int\mu\mathsf C=\int\mu\mathsf C^+-\int\mu\mathsf C^-\,\,{}'\!\!\leq'
\int\mu\mathsf D^+-\int\mu\mathsf D^-=\int\mu\mathsf D.  \]
\end{proof}
\end{lemma}
\begin{defin}\label{def:RDEF}
\lrquote

Consider a function
\[f:X\rightarrow\real^*.\]
a.) Suppose that there is a formal sum $\mathsf C$ such that for each $x\in X$
\[f(x)='\mathsf D^\Sigma(x)\]
and
\[\int\mu\mathsf D\neq\pm\infty.\]
We say that this value gives the value of  Lebesgue-Riesz integral
\[\int^{(LR)}\mu\,f.\]
b.) Suppose that $\mu\geq0$. We say that the formal sum $\mathsf D$ upper bounds $f$ if  for each $x\in X$
\[f(x)\leq'\mathsf D^\Sigma(x).\]
We define the upper Lebesgue-Darboux integral of $f$ as
\[ \int^{(LD+)} \mu\,f=
 \inf'\left\{ \int\mathsf D\mu:\, {\text{$\mathsf D$ upper bounds $f$ }} \right\},\]
which is an element of $[-\infty,+\infty]$.

The lower Lebesgue-Darboux integral of $f$ can be defined in similarly.

c.) Suppose that $\mu\geq0$. If the lower and upper Lebesgue-Darboux integrals
are the same then we say that it is the Lebesgue-Darboux integral
\[\int^{(LD)} \mu\,f.\]
\end{defin}
\begin{lemma}\label{lem:RWELLDEF}
\lrquote
Consider a function
\[f:X\rightarrow\real^*.\]

a.) The Lebesgue-Riesz integral, if exists, is unique.

b.) Suppose that $\mu\geq0$. For the Lebesgue-Darboux integrals, which exist in all circumstances,
\[\int^{(LD-)} \mu \,f\leq\int^{(LD+)} \mu\,f.\]

c.) Suppose that $\mu\geq0$. If $f$ is Lebesgue-Riesz integrable then it is Lebesgue-Darboux
integrable and
\[\int^{(LD)} \mu\,f=\int^{(LR)} \mu\,f.\]
\begin{proof}
a.) If
\[f(x)='\mathsf C^\Sigma(x)\quad\qand\quad f(x)='\mathsf D^\Sigma(x)\]
then
\[\mathsf C^\Sigma(x)\,\,{}'\!\!='\mathsf D^\Sigma(x).\]
Hence, by Lemma \ref{lem:FORMALCOMPAR}
\[\int\mu\mathsf C\,\,{}'\!\!='\int\mu\mathsf D.\]
If those values are not $\pm\infty$ then they must be equal.

b.) Consider an arbitrary formal sum  $\mathsf C$ which lower bounds $f$ and an arbitrary
formal sum $\mathsf D$ which  upper bounds $f$. Then
\[\mathsf C^\Sigma(x)\,\,{}'\!\!\leq f(x)\qand f(x)\leq'\mathsf D^\Sigma(x)\]
implies
\[\mathsf C^\Sigma(x)\,\,{}'\!\!\leq'\mathsf D^\Sigma(x).\]
Then, by Lemma \ref{lem:FORMALCOMPAR} it yields
\[\int\mu\mathsf C\,\,{}'\!\!\leq'\int\mu\mathsf D.\]
Now, taking the extended supremum and infinum yields our statement.

c.) If
\[f(x)='\mathsf D^\Sigma(x)\]
and
\[\int\mu\mathsf D\neq\pm\infty\]
then $D$ both lower bound and upper bounds $f$ and sets the upper and lower Lebesgue-Darboux
integrals to be the same.
\end{proof}
\end{lemma}
We will prove the converse of point c. later for finite Lebesgue-Darboux integrals. The converse
is not true in the infinite case.
~\\

\paragraph{\textbf{B. Properties}}
\begin{lemma}\label{lem:ADDITMON}
\lrquote

a.) The Lebesgue-Riesz integral is additive, the upper Lebesgue-Darboux integral is subadditive,
the lower Lebesgue-Darboux integral is superadditive -- if the sums are not $\pm\infty$.

b.) There is the usual behavior with respect to scalar multiplication. Multiplication by negative numbers
interchanges the upper and lower Lebesgue-Darboux integrals.

c.) If $\mu\geq0$ then the integral is monotone.
\begin{proof}
We take always take the union of the formal sums.
The behavior with respect to scalar multiplication is straightforward.
The monotonicity follows from Lemma \ref{lem:FORMALCOMPAR}.b.
\end{proof}
\end{lemma}

\begin{defin}\label{def:RNEGL}
\lrquote

We say that a set $C\subset X$ is negligible with respect to $\mu$
there exists a formal sum  $\mathsf C$
such that
\[\int\mu\mathsf C=0,\]
but $\mathsf C^\Sigma(x)=\pm\infty$ at points of $C$.
\end{defin}
\begin{lemma}\label{lem:RNEGL}
a.) If $\mathsf C$ is a formal sum but
$\displaystyle\int\mu\mathsf C$
is finite then the set of points where $\mathsf C^\Sigma(x)$ is not finite is negligible.

a'.) If $\mu\geq 0$ and $\mathsf C$ is a formal sum and
$\displaystyle\int\mu\mathsf C>-\infty$
 then set of points where $\mathsf C^\Sigma(x)$ is $-\infty$ or $\pm\infty$ is negligible.

b.) If we change the values of Lebesgue-Riesz/etc. integrable function in a negligible set
then the resulted functions is still Lebesgue-Riesz/etc. integrable

c.) Subsets and countable unions of negligible sets are still negligible.
\begin{proof}
a.) Take $-\mathsf C\,\dot\cup\,\mathsf C$. That shows the statement.

a'.) Apply part a.) to the nonnegative part $\mathsf C^+$ of $\mathsf C$.

b.) If $\mathsf C$ is like in Definition \ref{lem:RNEGL} then we can add
it to any approximation/ upper bound, etc.

c.) Only the case of countable unions is nontrivial. Suppose that $C_0,\,C_2,\,C_2,\ldots$
are negligible sets with the respective $\mathsf C_0,\,\mathsf C_2,\,\mathsf C_2,\ldots$.
Multiplying by small nonzero numbers we can assume that the total variation of
$\mathsf C_n$ is less than $2^{-n}$. Then the disjoint union $\mathsf C$ of the $\mathsf C_n$ is of still
finite total variation. The places where $\mathsf C^\Sigma(x)$ is $\pm\infty$ contains all $C_n$ and
it is negligible by point a.
\end{proof}
\end{lemma}
Hence we have the liberty of considering our functions up to changes on negligible sets.

For the purpose of the next Lemma we use the notation
$f\vee g=\max(f,g)$, $f\wedge g=\min(f,g)$, except at places where
$f$ or $g$ are $\pm\infty$. At those places we allow any value.
 In the case of Lebesgue-Riesz-integrable functions those places are negligible,
 so they do not matter for the purpose of integration.

\begin{lemma}\label{lem:RMINMAX}
\lrquote

If $f,g:X\rightarrow\mathbb R$ are Lebesgue-Riesz-integrable,
such that the integrals are finite or $<+\infty$ or $>-\infty$. Also suppose that $c<0<d$.

Then, we claim, $f\vee g$, $f\wedge g$, $f\vee d$, $f\wedge d$ are Lebesgue-Riesz-integrable,
with integrals which are finite or $<+\infty$ or $>-\infty$, respectively.
\begin{proof}
First we prove the special cases $f\vee 0$ and $f\wedge0$.
Suppose that
$\mathsf C=-\mathsf C^-\,\dot\cup \,\mathsf C^+$ is a formal sum to $f$.
We can consider $\mathsf E$, formal sum to
$\min(\mathsf C^-{}^{\Sigma}(x),\mathsf C^+{}^{\Sigma}(x))$.
Then
\[-\mathsf E\,\dot\cup \,\mathsf C^+\qand-\mathsf C^- \,\dot\cup \,\mathsf E\]
are formal sums to $f\vee 0$ and $f\wedge 0$. Notice that the variations decreased.
Hence, we proved these special cases.
Now,
\[f\vee g=  (f\vee 0)\vee (g\vee0)+ (f\wedge 0)\vee(g\wedge 0)\]
\[f\wedge g=  (f\wedge 0)\wedge (g\wedge0)+ (f\vee 0)\wedge(g\vee 0)\]
\[f\vee c=(f\vee 0)+(f\wedge0)\vee c \]
\[f\wedge d=(f\wedge 0)+(f\vee0)\wedge d \]
show that it is enough to prove those cases when $f,g\geq 0$ or $f,g\leq 0$.

By symmetry we consider only the case when $f,g\geq 0$
and we prove that $f\vee g,\,f\wedge g,\,f\wedge d$ have those nice properties.

Suppose that
$\mathsf C=-\mathsf C^-\,\dot\cup \,\mathsf C^+$ and
$\mathsf D=-\mathsf D^-\,\dot\cup \,\mathsf D^+$
are formal sums to $f$ and $g$.
Notice that $\mathsf C^-,\,\mathsf D^-$ are of finite variation
($\mathsf C^+,\,\mathsf D^+$ subdecomposes into them, so infinite variation would make the
overall sum $\pm\infty$).
Now, we can find nonnegative formal sums, such that
$\mathsf E_{\max}$ is a formal sum to
\[\max((\mathsf C^+\,\dot\cup \,\mathsf D^-)^{\Sigma}(x),
(\mathsf C^-\,\dot\cup \,\mathsf D^+)^{\Sigma}(x));\]
and $\mathsf E_{\min}$ is a formal sum to
\[\min((\mathsf C^+\,\dot\cup \,\mathsf D^-)^{\Sigma}(x),
(\mathsf C^-\,\dot\cup \,\mathsf D^+)^{\Sigma}(x)).\]
Then,
\[-\mathsf D^-\,\dot\cup \,-\mathsf C^-\,\dot\cup \,\mathsf E_{\max}
\qand
-\mathsf D^-\,\dot\cup \,-\mathsf C^-\,\dot\cup \,\mathsf E_{\min}\]
are formal sums
\[f\vee g\qand f\wedge g.\]
Notice that the variations decreased. Remains the case $f\wedge d$.

Using Lemma \ref{lem:REFIN}.a we can find countable family $\mathfrak D$ of disjoint sets which cover
all the set from $\mathsf C=-\mathsf C^-\,\dot\cup \,\mathsf C^+$.
Let $\mathsf d$ be the formal sum which assign the coefficient everywhere.
Let $\mathsf E$ be a formal sum to
\[\min((\mathsf C^+)^{\Sigma}(x),
(\mathsf C^-\,\dot\cup \,\mathsf d)^{\Sigma}(x)).\]
Then
\[-\mathsf C^-\,\dot\cup \,\mathsf E\]
will be a formal sum to $f\wedge d$.

[Remark: Being the the integral possibly infinite, we cannot use the standard tricks for $f\vee g$.]
\end{proof}
\end{lemma}
\begin{defin}\label{def:PNINT}
\lrquote

We define
\[\int^-\mu f,\quad \int^+\mu f,\quad\int^\pm\mu f,\]
as the infinum of all
\[\int^-\mu\mathsf C,\quad\int^+\mu\mathsf C,\quad\int^\pm\mu\mathsf C,\]
respectively, such that $\mathsf C$ is a formal sum to $f$.
\end{defin}
\begin{lemma}\label{lem:PNINT}
\lrquote

If $ f$ is Lebesgue-Riesz integrable then
\[\int^{(LR)}\mu f=\int^+\mu f-\int^-\mu f\]
\[\int\mu^\pm f=\int^+\mu f+\int^-\mu f .\]
\begin{proof}
It immediately follows from the definitions and the consistency statement for
$\displaystyle\int\mu\mathsf C$ in Lemma \ref{lem:FORMALCOMPAR}.a.
\end{proof}
\end{lemma}
\begin{lemma}\label{lem:VARCONT}
Suppose that $f\geq0$ has a Lebesgue-Riesz integral.
Then we claim, that for any $\varepsilon>0$ there is a formal sum
$\mathsf C=-\mathsf C^-\,\dot\cup \,\mathsf C^+$ to $f$ such that
\[\int^\pm\mu \mathsf C^-<\varepsilon.\]
In particular, if $f\geq 0$ is Lebesgue-Riesz integrable and also $\mu\geq0$ then
\[\int^-\mu f=0. \]
\begin{proof}
Let $\mathsf C$ be a formal sum to $f$ such that $\int\mu\mathsf C\neq\pm\infty$.
First of all $\mathsf C^-$ has finite negative and positive variation; otherwise
$\mathsf C^-{}^\Sigma(x)\leq \mathsf C^+{}^\Sigma(x)$ would imply that
that $\mathsf C^+$ has at least the same negative and positive variation,
and that would make
$\int\mu\mathsf C=\pm\infty$.

So, the  variation in $\mathsf C^-$ it is finite.
We can number-decompose each positive coefficient in $\mathsf C^-$ into infinitely
many positive coefficients. That way we obtain $\mathsf C^{--}$.
The variation of $\mathsf C^{--}$ is still finite, and
it mainly comes only from finitely many terms, so
$\mathsf C^{--}=\mathsf D\,\dot\cup \,\mathsf C^-{}'$;  $\mathsf D$ is finite
and almost all variation comes from  $\mathsf D$ except $\varepsilon$.
Now $\mathsf D{}^{\Sigma}(x)<\mathsf C^-{}^{\Sigma}(x)\leq\mathsf C^+{}^{\Sigma}(x)$
or $\mathsf D{}^{\Sigma}(x)=0$. We can apply Lemma \ref{lem:FORMALDECOMP}.a, and
so we can substitute $\mathsf D$ by $\mathsf D'$ and
$\mathsf C^+$ by $\mathsf D'\,\dot\cup \,\mathsf C^+{}' $,
both decompositions.

Being $\mathsf D$ finite, for $\mathsf C'=-\mathsf C^-{}'\,\dot\cup \,\mathsf C^+{}'$
we see that
\[f(x)='\mathsf C'{}^\Sigma(x),\]
 and it has the required properties.
\end{proof}
\end{lemma}
The lemma above demonstrates the  applicability of
\begin{lemma}[Beppo Levi's theorem]\label{lem:RBLevi}
\lrquote

Suppose that $f_n$ are Lebesgue-Riesz integrable.
\[\sum_{n\in\mathbb N} \int^-\mu f_n<+\infty.\]
Then
\[\sum_{n\in\mathbb N}f_n\]
(which is, of course, is defined everywhere)
converges in classical sense to a finite number or to $+\infty$ or to $-\infty$ almost everywhere; the resulted function is Lebesgue-Riesz integrable and
\[\sum_{n\in\mathbb N} \int^{(LR)}\mu f_n = \int^{(LR)}\mu\sum_{n\in\mathbb N} f_n.\]
Moreover,
\[ \int^-\mu\sum_{n\in\mathbb N} f_n\leq \sum_{n\in\mathbb N} \int^-\mu f_n. \]
\begin{proof}
We can assume that $\mathsf C_n$ is a formal sum to $f_n$ such that
\[\int^-\mu \mathsf C_n< \int^-\mu f_n+2^{-n}\varepsilon.\]
Then take all $\mathsf C_n$ together into $\mathsf C$. For that
\[\int^-\mu\mathsf C<\sum_{n\in\mathbb N} \int^-\mu f_n+2\varepsilon.\]
Then $f=\mathsf C^{\Sigma}(x)$ is Lebesgue-Riesz integrable.
In particular, $f(x)=\mathsf C^{\Sigma}(x)\neq\pm\infty$ except on a negligible set $N$.
Also, $f_n(x)=\mathsf C_n^{\Sigma}(x)$ almost everywhere, but certainly outside of $N$.
Then, one can see $f(x)\,\,{}'\!\!\!\!=\sum_{n\in\mathbb N}f_n(x)$; which implies that
$f(x)=\sum_{n\in\mathbb N}f_n(x)$;
except on $N$.
In particular, that yields the classical convergence statement and the integrability of
$\sum_{n\in\mathbb N}f_n(x)$.

Our method immediately yields the estimate for $\int^-\mu\sum_{n\in\mathbb N} f_n$, because
taking formal sums together, their negative variation is taken together.
\end{proof}
\end{lemma}
We have seen that in the case $\mu\geq 0$ Lebesgue-Riesz integrability implies Lebesgue-Darboux integrability.
The reverse is also true if the Lebesgue-Darboux integral is finite:
\begin{lemma}\label{lem:LDREQUIV}
\lrquote

Suppose that $\mu\geq0$.
Then, we claim, if $f$ has a finite Lebesgue-Darboux integral then it is also  Lebesgue-Riesz integrable.
\begin{proof}
From the lower approximating sums there are functions $g_n'\leq f$ such that
\[\lim g_n'\]
is monotone increasing, the Lebesgue-Riesz integrals (ie. overall sums) limit to $c$.
We can actually take
\[g_n=f_0\vee f_1\vee\ldots\vee f_n.\]
Similarly with $h_n$ and the upper sums.
Then
\[f_-=g_0+\sum_{n\in\mathbb N}(g_{n+1}-g_n)\leq f\leq f_+=g_0+\sum_{n\in\mathbb N}(h_{n+1}-h_n)\]
almost everywhere. (Existence follows from  Beppo Levi's theorem.)
Furthermore,
\[\int\mu f^-=c=\int\mu f^+.\]
Considering
\[\sum_{n\in\mathbb N }(f^+-f^-)\]
and Beppo Levi's theorem we immediately see that not only $f^-\leq f^+$ but
\[f^-= f^+\]
almost everywhere. Hence
\[f^-= f=f^+\]
almost everywhere, and Lebesgue-Riesz integrability is clear.
\end{proof}
\end{lemma}
\begin{lemma}[Lebesgue's monotone convergence theorem]\label{lem:RLEBMON}
\lrquote

Suppose that $f_n,f\geq0$ are integrable; and
\[\sum_{n\in\mathbb N}f_n\leq f.\]
Then
\[\sum_{n\in\mathbb N}\int^-\mu f_n\leq \int^-\mu f \qand \sum_{n\in\mathbb N}\int^+\mu f_n\leq \int^+\mu f.\]
In particular,
\[\sum_{n\in\mathbb N}f_n\]
is integrable and
\[\sum_{n\in\mathbb N}\int^- \mu f_n= \int^-\mu\sum_{n\in\mathbb N}f_n
\qand \sum_{n\in\mathbb N}\int^+\mu f_n= \int^+\mu\sum_{n\in\mathbb N}f_n.\]
That also implies
\[\sum_{n\in\mathbb N}\int^\pm\mu f_n= \int^\pm\mu\sum_{n\in\mathbb N}f_n
\qand\sum_{n\in\mathbb N}\int^{(LR)} \mu f_n= \int^{(LR)}\mu\sum_{n\in\mathbb N}f_n.\]
\begin{proof}
i.) It is enough to prove the case of $\int^-$, otherwise we just take $-\mu$ instead of $\mu$.
It is enough to prove for a finite sum. In that case we can assume that
\[f_0+\ldots+f_s=f.\]
We can consider formal sums $\mathsf C_n=-\mathsf C_n^-\dot\cup\mathsf C_n$, and
$\mathsf C=-\mathsf C^-\dot\cup\mathsf C$ to $f_n$ and $f$. we can assume that
\[\int \mu^{\pm}\mathsf C^-_n<\varepsilon,\quad\int \mu^{\pm}\mathsf C^-<\varepsilon.\]
Then
\[\sum_{0\leq n\leq s }\int^- \mu f_n\leq\sum_{0\leq n\leq s }\int^- \mu\mathsf C^+_n+s\varepsilon\leq
\sum_{0\leq n\leq s }\int^- \mu\mathsf C^+_n+\int^-\mu\mathsf C^-+s\varepsilon=\]
\[\sum_{0\leq n\leq s }\int^-\mu\mathsf C^-_n+\int^-\mu\mathsf C^++s\varepsilon\leq
\int^-\mu\mathsf C^++2s\varepsilon\leq\int^-\mu\mathsf C^++\int^-\mu\mathsf C^-+2s\varepsilon \]
\[\leq\int^-\mu f+2s\varepsilon.\]
ii.) According to the integrability of $f$ either
\[\sum_{n\in\mathbb N}\int^-\mu f_n\quad\text{ or }\quad\sum_{n\in\mathbb N}\int^+\mu f_n\]
is bounded, hence Beppo Levi's theorem can be applied.
So the integrability of the sum  follows.

Hence, we can assume
\[f=\sum_{n\in\mathbb N}f_n.\]
Then, one has to prove the inequalities for the positive and negative integrals in the other direction.
These are nontrivial only in the unbounded case when they follow from Beppo Levi's theorem.
\end{proof}
\end{lemma}
From the these statements above one prove Fatou's lemma, Lebesgue's dominated (sequential) convergence theorems.
In general, having these basic statements one can continue with the Lebesgue-integral as usual.
Notice that Beppo Levi's theorem immediately implies Fubini's theorem in the present setting.

A fundamental link to the classical viewpoint on Lebesgue integration is:
\begin{lemma}\label{lem:RMEAS}
\lrquote

If $f$ is integrable then $c\chi_{\{f>d\}}$ is Lebesgue-Riesz integrable.
\begin{proof}
We can see that
\[c\chi_{\{f>d\}}=\sum_{n\in\mathbb N}\Bigr(((n+1)(f-f\wedge d))\wedge c-(n(f-f\wedge d))\wedge c\Bigl).\]
Notice that the terms of the sum are nonnegative. Then, the integrability follows from
Lemma \ref{lem:RLEBMON} with the majorizing function $m(f\vee 0)$, where $m$ is a sufficiently large integer.
\end{proof}
\end{lemma}
Then one can define integrable and measurable sets as usual.
\begin{remark}
If one is primary interested in the case of nonnegative measure
then it is much better to leave Procedure 1 out from the discussion and prove Lemma
\ref{lem:FORMALDECOMP}.a, b., d. ``by hand'', which is quite easy.
Then one case base integration on the Lebesgue-Darboux version.
\end{remark}
\begin{remark}
If one wants to extend the Lebesgue-Riesz integral to vector measures then one can define
\[\int\mu\mathsf C\]
as the common value of all
\[\sum_{\lambda\in\Lambda} c_\lambda'\mu(C_\lambda')\]
for all decomposition $\mathsf C'$ of $\mathsf C$ (if such a common value exists).
This definition is convenient to use if we assume a bounded variation property like we did above.
\end{remark}
\newpage\section{Extension theory}\label{sec:EXTENS}
The advantage of the picture presented in the previous sections is that
there is a full analogue in terms of Stonean vector lattices.
In fact, things are even somewhat nicer.
\\

\paragraph{\textbf{A. Positive vector lattices}}
\begin{defin}\label{def:EPOSLAT}
A positive (in fact, nonnegative) vector lattice $\mathcal  S$ is subset of $[0,+\infty)^X$ such that
\begin{itemize}
\item[\texttt{(Y1)}] $f,g\in\mathcal  S$ implies that $f+g$, $f\vee g$, $f\wedge g$, $(f-g)\vee 0$
are also elements of $\mathcal  S$.
\end{itemize}
The positive vector lattice $\mathcal  S$ is Stonean if
\begin{itemize}
\item[\texttt{(Y2)}] $f\in\mathcal  S$,  $d\geq 0$ implies that $f\wedge d$
is also an element of $\mathcal  S$.
\end{itemize}
\end{defin}
\begin{defin}\label{def:EDANF}
A Daniell functional on a positive vector lattice $\mathcal  S$ is a function
\[\tilde\mu:\mathcal  S\rightarrow\mathbb R\]
such that if $f_n$ ($n\in\mathbb N$), $f$ are in $\mathcal  S$ then
\[f=\sum_{n\in\mathbb N}f_n \]
(pointwise) implies
\[\tilde\mu\biggl(\sum_{n\in\mathbb N}f_n\biggr)=\sum_{n\in\mathbb N}\tilde\mu(f_n).   \]
\end{defin}
\begin{defin}\label{def:EBOUND}
A Daniell functional $\tilde\mu$ on a positive vector lattice $\mathcal  S$
is of locally bounded variation if for each $f\in\mathcal  S$
\[\{\tilde\mu(g)\,:\,0\leq g\leq f\}\]
is bounded.
\end{defin}
The fundamental statement is this matter is:
\begin{lemma}[Representability]\label{lem:EREP}
Let $\mathcal  S$ be a positive Stonean lattice.
Assume that $\tilde\mu$ is a Daniell functional of locally bounded variation on $\mathcal  S$.

Then, we claim, there exists a unique measure $\mu$ on the interval system
\[\mathfrak I(\mathcal S)=\biggl\{\{f>1\}\setminus\{g>1\}\,:\,f,g\in\mathcal S\biggr\}\]
such that
\[\tilde\mu(f)=\int \mu f .\]
\begin{proof} (Sketch, from \cite{kind}.)
1. Let $\mathfrak S'$ be the family of sets
\[A_{f,g}=\{(v,x)\in\mathbb R\times X\,:\,f(x)> v\geq g(x)\},\]
where $f,g\in\mathfrak S$.
One can check that these sets form an interval system.
Let $l$ be the interval length function. Then, we can notice
\[l(A_{f,g}{}^x)=((f-g)\wedge0)(x).\]
Let us define the measure
\[\mu'(A)=\tilde\mu(x\mapsto l(A^x));\]
Then the $\sigma$-additive property of $\mu'$ the is just an obvious
reformulation of the $\sigma$-additive property of $\tilde\mu$.

2. Now, let $m\geq b$, where $m\in\mathbb N$, $0\leq b\in\mathbb R$. Then
\[[0,b)\times\{f>1\}=\bigcup_{n\in\mathbb N}^{\updisjoint}
A_{((n+1)(f-f\wedge1))\wedge b,n(f-f\wedge1)}\subset A_{mf,0}\]
implies that the set above $\mu'$-integrable.
Taking differences it follows that the sets
\[[a,b]\times( \{f>1\}\setminus\{g>1\})\]
are all integrable.

3. One can simply check that $\mathfrak I(\mathcal S)$ forms an interval system.
Let
\[\mu(B)=\int\chi_{[0,1]\times B}\mu'.\]
This is certainly a measure, and it is of bounded variation.

One can simply prove that
\[2\int\chi_{[0,a]\times B}\mu'=\int\chi_{[0,2a]\times B}\mu'.\]
(we simply substitute all functions $f\in\mathbb S$ by $2f\in\mathbb S$),
from which it follows that
\[(b-a)\mu(B)=\int\chi_{[a,b]\times B}\mu'\]
if $a\leq b$ are nonnegative dyadic numbers.

4.)
Now, consider the decomposition
\[A_{f,0}=\bigcup_{n\in\mathbb Z,\,b\in\mathbb N} [2b2^n,(2b+1)2^n)\times
\Bigl(\{f>(2b+1)2^n\}\setminus\{f>2(b+1)2^n\}\Bigr),\]
apply $\mu'$ and Beppo Levi's theorem for $\mu$. Then we see that
\[\tilde\mu(f)=\int f\mu.\]

5. The unicity of $\mu$ is obvious from the construction.
\end{proof}
\end{lemma}
\begin{remark}
This method is a basic to obtain measures.

Having such a measure, it will yield an extension of $\mathcal  S$ to the
positive lattice of non-negative integrable functions.
It is called the Daniell extension.
\end{remark}

In the more traditional proofs, the associated integration theory (extension) is constructed directly.
Beside its philosophical importance, this latter method may be useful because the Stonean property can be
avoided for a good part, at the expense of representability.

In what follows our objective is to construct this extension  by imitating the
construction of the Lebesgue-Riesz integral.

\paragraph{\textbf{B. Formal sums in positive vector lattices }}
\begin{defin}\label{def:EESFORMAL}
a.) A formal sum $\mathsf  C$ in a positive Stonean vector lattice $\mathcal  S$ is just an indexed
family of elements $f$ and $-f$ where $f$ is from the lattice. The formal sum is nonnegative if
only elements of $\mathcal  S$ occur.

b.) A formal sum $\mathsf  C$ is decomposes to a formal sum $\mathsf  C'$ each term
$f\in\mathcal  S$ from $\mathsf  C$ is substituted by countably many terms
$\{f_\lambda\}_{\lambda\in\Lambda}$ ($f_\lambda\in\mathcal  S$) such that
\[f=\sum_{\lambda\in\Lambda }f_\lambda,\]
and terms like $-f$ are substituted by similar terms $\{-f_\lambda\}_{\lambda\in\Lambda}$
($f_\lambda\in\mathcal  S$).

c.) The subdecomposition of formal sums is defined similarly as before.

d.) The evaluations $\mathsf C^\Sigma(x)$ can be defined similarly.
\end{defin}
\begin{remin}
If we have finitely many real-valued function functions $f_j$ ($j\in J$) then the
$(n+1)$th greatest value function is
\[\max{}_{-n}\{f_j\,:\,j\in J\}=\bigwedge_{\substack{\Xi\subset J\\|\Xi|=n}}\bigvee_{j\in J\setminus\Xi}f_j.\]
The set notation is ambiguous  because we understand elements with multiplicities, but it will
not cause problems.
\end{remin}
One can describe the analogue of Procedure 1:
\begin{defin}[Description of Procedure 2]\label{def:EESPROCTWO}
Suppose that we have countably many formal sums
\[\{\mathsf C_n\}_{n\in N}=\{\{f_{n,m}\}_{m\in N_n}\}_{n\in N}\]
with  positive coefficients, where $N$ and $N_n$ are initial
segments of $\mathbb N$. Also list all possible pairs
$\langle{n,m}\rangle$ as $\langle n(j),m(j)$ $(j\in J)\rangle$ along an initial segment $J$ of $\mathbb N$.

In what follows we will describe a procedure which  finitely decomposes each
formal sum $\mathsf C_n$, and
lattice element $z$ resulted assigns a ``color'' from $N$ and a ``lift'' $u\leq v\in\mathcal  S$
such that $v-u=z$.

Then, the procedure is in steps along $J$. Suppose that $j\in J$ comes.
Let
\[b_j=\sum_{\substack{1\leq k\leq j,\, n(k)=n(j)}}f_{n(k),m(k)}
\qand
c_j=\sum_{\substack{1\leq k< j,\, n(k)=n(j)}}f_{n(k),m(k)}.\]
Notice that $b_j-c_j=f_{n(j),m(j)}$.

Let
\[S_j=\{n(k)\,;\,0\leq k\leq j\}.\]
We can consider
\[h_i=\sum_{\substack{1\leq k\leq j,\, n(k)=i}}f_{n(k),m(k)}\]
($i\in S_j$).
In fact, we can consider
\[h_n'=\max{}_{-n}\{h_i\,:i\in S_j\}\]
($0\leq n< |S_j|$), the $(n+1)$th greatest value function.
Here the set-notation is ambiguous, the $h_i$'s are understood with multiplicities.

More generally, we can consider all functions
\[b_i=\sum_{\substack{1\leq k\leq i,\, n(k)=n(i)}}f_{n(k),m(k)}\]
with an $i<j$ chosen.
Let
\[b_n'=\max{}_{-n}\{b_0,\ldots, b_j\},\]
the $(n+1)$th greatest's function.

Then defining $r_n=c_j\vee h_n'\wedge b_j $  ($0\leq n< |S_j|$) we find that
\[c_j=r_{|S_j|}\leq\ldots\leq r_1\leq r_0=b_j.\]
Moreover, defining $r_{n,k}=r_{n+1}\vee b_k'\wedge r_n$ ($0\leq k\leq j$) we find that
\[r_{n+1}= r_{n,j}\leq\ldots\leq r_{n,1}\leq r_{n,0}=r_n.\]

Now,
\[f_{n(j),m(j)}=b_j-c_j=\sum_{0\leq k<j,\, 0\leq n<|S_j|} (r_{n,k}-r_{n,k+1}).\]
Hence, we decompose $f_{n(j),m(j)}$ into these $(r_{n,k}-r_{n,k+1})$'s, with lifts
$r_{n,k}\leq r_{n,k+1}$.
The color given to such a piece is $n$.

From the viewpoint of each element $f_{n,m}$ a finite decomposition
happened, so we can imagine that we finitely decompose each $\mathsf C_n$.
We can gather things at the end into coefficient systems
\[\mathsf C'_0,\,\mathsf C'_1,\,\mathsf C'_2,\,\ldots\]
as decompositions of $\mathsf C_n$'s respectively, or into
\[\mathsf F_0,\,\mathsf F_1,\,\mathsf F_2,\,\ldots\]
according to color.
\end{defin}
This Procedure 2 has the nice properties similar to Procedure 1, except
that
\[[u,v)\times E\]
is supposed to be substituted by
\[\{(c,x)\,:\,u(x)\leq c<v(x)\}.\]
We will not spell these properties out here, the main matter is that in each step the sum of
elements of color $n$ will give $h_n'$, the $(n+1)$th greatest value.
The variant of Lemma \ref{lem:FORMALDECOMP} holds:
\begin{lemma}\label{lem:EFORMALDECOMP}
Consider a positive  vector lattice $\mathcal  S$.

a.) Suppose that $\mathsf C$ and $\mathsf D$ are nonnegative formal sums, $\mathsf C$ is finite.
Also suppose that for all $x\in X$
\[\mathsf C^\Sigma(x)\leq\mathsf D^\Sigma(x).\]
Then, we claim, there exist nonnegative formal sums $\mathsf C'$ and $\mathsf E$ such that
$C$ decomposes into $\mathsf C'$ and $D$  decomposes into $\mathsf C'\,\dot\cup\,\mathsf E$:
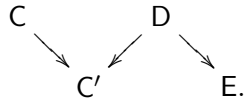
\[\xymatrix@=12pt{ \mathsf C\ar[dr] &  &\mathsf D\ar[dr]\ar[dl]& \\& \mathsf C' & & \mathsf E. }\]

b.) Suppose that $\mathsf C$ and $\mathsf D$ are nonnegative formal sums and
\[\mathsf C^\Sigma(x)\leq\mathsf D^\Sigma(x)\]
holds for all $x\in X$. Then, we claim, there exists a nonnegative formal sum $\mathsf C'$
such that $\mathsf C$  decomposes into $\mathsf C'$ and $\mathsf D$ subdecomposes into $\mathsf C'$:
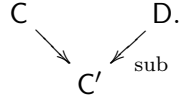
\[\xymatrix@=12pt{\mathsf C\ar[dr] &  & \mathsf D.\ar[dl]^{\sub}  \\ & \mathsf C' &} \]

c.) Suppose that  $\mathsf C_1$ and $\mathsf C_2$ are nonnegative formal sums such that
\[\mathsf C_1^\Sigma(x)= \mathsf C_2^\Sigma(x)\]
for all $x\in X$. Then there exists a nonnegative formal sum $\mathsf  C$ such that
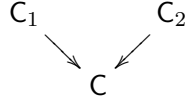
\[\xymatrix@=12pt{ \mathsf C_1\ar[dr] &  & \mathsf C_2\ar[dl] \\ & \mathsf C &} \]
both $\mathsf C_1$ and $\mathsf C_2$ can be decomposed into $\mathsf C$.

d.) Suppose that $\mathsf C$ and $\mathsf D$ are nonnegative formal sums.
Then, we claim, there exists a decomposition
\[\mathsf C\,\dot\cup\,\mathsf D\rightarrow\mathsf  E_{\min}\,\dot\cup\, \mathsf E_{\max}\]
such that
\[ \mathsf E_{\min} {}^{\Sigma}(x)=\min( \mathsf C^\Sigma(x), \mathsf D^\Sigma(x) ) \]
and
\[ \mathsf E_{\max} {}^{\Sigma}(x)=\max( \mathsf C^\Sigma(x), \mathsf D^\Sigma(x) ). \]
\end{lemma}

\paragraph{\textbf{C. Integration}}
\begin{defin}\label{def:ESMDECOMP}
Let $\mathcal S$ be a positive vector lattice.
Suppose that $\tilde\mu:\mathcal  S\rightarrow\mathbb R$ is Daniell functional. Then we define the
variations, then for $f\in\mathcal\mathcal  S$
\[|\tilde\mu|^+(f)=\sup\biggl\{\sum_{\lambda\in\Lambda} |\tilde\mu(f_\lambda)|^+\,:\,\{f_n\}_{\lambda\in\Lambda}
\text{ decomposes }f\biggr\} \]
\[|\tilde\mu|^-(f)=\sup\biggl\{\sum_{\lambda\in\Lambda} |\tilde\mu(f_\lambda)|^-\,:\,\{f_n\}_{\lambda\in\Lambda}
\text{ decomposes }f\biggr\} \]
\[|\tilde\mu|(f)=\sup\biggl\{\sum_{\lambda\in\Lambda} |\tilde\mu(f_\lambda)|\,:\,\{f_n\}_{\lambda\in\Lambda}
\text{ decomposes }f\biggr\} .\]
\end{defin}
One immediately notices that this definition is more difficult then it should be, however, in that
form invariance to decomposition due to Lemma \ref{lem:EFORMALDECOMP}.c is obvious.
In particular, it follows that $|\tilde\mu|^+$, etc. are Daniell-functionals themselves.
Also, from the definition it is straightforward that
\begin{lemma}\label{lem:ESMDECOMP}
Let $\mathcal S$ be a positive vector lattice.
Suppose that $\tilde\mu:\mathcal  S\rightarrow\mathbb R$ is a Daniell functional.

Then, we claim,
\[|\tilde\mu|(f)=|\tilde\mu|^+(f)+|\tilde\mu|^-(f),\qquad
\tilde\mu(f)=|\tilde\mu|^+(f)-|\tilde\mu|^-(f),\]
\[|\tilde\mu|^-(f)=-\inf\{\tilde\mu(g)\,:\,0\leq g\leq f\},\qquad
|\tilde\mu|^+(f)=\sup\{\tilde\mu(g)\,:\,0\leq g\leq f\}.\]
In particular, locally bounded variation means that $|\tilde\mu|(f)<+\infty$ for all $f\in\mathcal  S$.
\qed
\end{lemma}
\begin{defin}\label{def:ESVARDECOMP}
Let $\mathcal S$ be a positive vector lattice.
Suppose that $\tilde\mu:\mathcal  S\rightarrow\mathbb R$ is a Daniell functional.

If $\mathsf C=\{f_\lambda\}_{\lambda\in\lambda}$ is a formal sum then we define
\[\int^- \tilde\mu\mathsf C=\sum_{\lambda\in\Lambda,\,f_\lambda\in\mathcal  S}|\tilde\mu|^-(f) +
\sum_{\lambda\in\Lambda,\,-f_\lambda\in\mathcal  S}|\tilde\mu|^+(-f)\]
\[\int^+ \tilde\mu\mathsf C=\sum_{\lambda\in\Lambda,\,f_\lambda\in\mathcal  S}|\tilde\mu|^+(f) +
\sum_{\lambda\in\Lambda,\,-f_\lambda\in\mathcal  S}|\tilde\mu|^-(-f)\]
etc.
\end{defin}
\begin{point}
From here the theory of integration is essentially the same as before, proofs can be adapted
from earlier. One can define the Daniell integral
\[\int^{(D)}\tilde\mu f\]
from formal sums; and then one can prove its properties analogous the Lebesgue-Riesz integral.
\end{point}
\begin{point}
 There is one (in fact, the only) place when we have to apply for the Stonean property.
That is the proof of the analogue Lemma \ref{lem:RMINMAX}, in order to proof the integrability of
$f\wedge d$, ($d>0$ constant). For that purpose one can use the following lemma:
\end{point}
\begin{lemma}\label{lem:LEVELDECOMP}
Suppose that $f_n\in\mathcal  S$ ($n\in\mathbb N$), and $d>0$ is a constant.
Let
\[W=\bigcup_{n\in\mathbb N} \{x\,:\, f_n(x)>0\}.\]
Then, we claim, $d\chi_W$ decomposes in $\mathcal S$.
\begin{proof}\
Let
\[s_n=\left(\sum_{0\leq k< n} (n-k)f_k \right)\wedge d.\]
Notice that $s_0=0$.
Then one can see that the functions
\[(s_{n+1}-s_n)\]
realize the decomposition.
\end{proof}
\end{lemma}

\paragraph{\textbf{D. Equivalence of the extensions}}
~\\

In this paragraph we will not give full proofs.

\begin{lemma}\label{lem:DEAEQUIL1}
Consider the situation of Lemma \ref{lem:EREP}. Then, we claim:

a.) For all $S\in\mathfrak I(\mathcal S)$
\[|\mu|^+(S)=\int^+\tilde\mu\chi_S.\]
Similar statement apply to $|\mu|^-$ and $\int^-$; $|\mu|$ and $\int^\pm$; $\mu$ and $\int^{(D)}$.
All these values are finite.

b.) For all $f\in\mathcal S$
\[ |\tilde\mu|^+(f)=\int^+\mu f.\]
Similar statements apply to $|\tilde\mu|^-$ and $\int^-$; $|\tilde\mu|$ and $\int^\pm$;
$\tilde\mu$ and $\int^{(LR)}$.
All these values are finite.

\begin{proof}{(Indication)}
One has to use Lemma \ref{lem:VARCONT} and Lemma \ref {lem:RLEBMON} and their Daniell-analogues
several times.
\end{proof}
\end{lemma}

\begin{lemma}\label{lem:DEAEQUIL2}
Consider the situation of Lemma \ref{lem:EREP}.

Then, we claim, a function $f:X\rightarrow\mathbb R^*$ is Lebesgue-Riesz integrable
if and only if it is Daniell integrable. In that case the integrals are the same.

\begin{proof}{(Indication)}
That follows from Lemma \ref{lem:DEAEQUIL1}; and Lemma \ref {lem:RLEBMON}
and its Daniell-analogue because they control each other's variations
for formal sums.
\end{proof}
\end{lemma}

\end{document}